\documentclass[10pt]{amsart}
\usepackage{
	geometry
} 
\usepackage{graphicx}
\usepackage{dsfont}
\usepackage{amssymb}
\usepackage{amsmath}
\usepackage{epstopdf}
\usepackage{fullpage}
\usepackage{enumerate}
\usepackage{color, subcaption}
\usepackage{amsthm}
\usepackage{placeins}
\usepackage{gensymb} 
\usepackage{array}
\usepackage{blkarray}
\usepackage{multirow}
\usepackage{float}
\usepackage{morefloats}
\usepackage{pgfplots}
\usepackage[abs]{overpic}
\usepackage{xcolor}
\usepackage{amsopn}
\usepackage{diagbox}

\usepackage[colorlinks=true]{hyperref}
\usepackage[nameinlink, capitalise]{cleveref}
\crefname{equation}{}{}

\makeatletter
\def
\subsection{\@startsection{subsection}{3}%
\z@{.5\linespacing\@plus.7\linespacing}{.5\linespacing}%
{\bf}}
\makeatother

\makeatletter
\def
\subsubsection{\@startsection{subsubsection}{3}%
\z@{.5\linespacing\@plus.7\linespacing}{.5\linespacing}%
{\it}}
\makeatother

\let\oldtocsection=\tocsection

\let\oldtocsubsection=\tocsubsection

\let\oldtocsubsubsection=\tocsubsubsection

\renewcommand{\tocsection}[2]{\hspace{0em}\oldtocsection{#1}{#2}\textbf}
\renewcommand{\tocsubsection}[2]{\hspace{1em}\oldtocsubsection{#1}{#2}}
\renewcommand{\tocsubsubsection}[2]{\hspace{2em}\oldtocsubsubsection{#1}{#2}}

\newcommand{\dv}{\text{\rm div}}
\renewcommand{\o}{\text{\rm o}}
\newcommand{\hhat}{\widehat{h}}
\renewcommand{\d}{\text{\rm d}}

\newcommand{\Cov}{\text{\rm Cov}}

\newcommand{\tr}{\text{\rm tr}}

\newcommand{\e}{\varepsilon}
\newcommand{\I}{\text{\rm I}}
\newcommand{\Id}{\text{\rm Id}}

\newcommand{\calA}{{\mathcal{A}}}
\newcommand{\calP}{{\mathcal{P}}}

\newcommand{\calC}{{\mathcal{C}}}
\newcommand{\calF}{{\mathcal{F}}}
\newcommand{\calH}{{\mathcal{H}}}
\newcommand{\calD}{{\mathcal{D}}}

\newcommand{\calO}{{\mathcal{O}}}

\newcommand{\calT}{{\mathcal{T}}}

\newcommand{\Vol}{\text{\rm Vol}}

\newcommand{\lcor}{l_{\text{\rm cor}}}

\newcommand{\Winfty}{W^{1,\infty}(\mathbb{R}^d,\mathbb{R}^d)}

\newcommand{\R}{{\mathbb{R}}}
\renewcommand{\P}{{\mathbb{P}}}
\newcommand{\Ptrue}{{\mathbb{P}}_{\text{\rm true}}}
\renewcommand{\S}{{\mathbb{S}}}
\newcommand{\Q}{{\mathbb{Q}}}

\newcommand{\Uad}{{\mathcal{U}}_{\text{\rm ad}}}
\newcommand{\Jmean}{J_{\text{\rm mean}}}
\newcommand{\Jdr}{J_{\text{\rm dr}}}

\newcommand{\JW}{J_{\text{\rm W}}}
\newcommand{\JM}{J_{\text{\rm M}}}
\newcommand{\VaR}{\textrm{\texttt{VaR}}}
\newcommand{\CVaR}{\textrm{\texttt{CVaR}}}

\usepackage{soul, xcolor}
\setstcolor{violet}

\begin{document}
\newtheorem{theorem}{Theorem}[section]
\newtheorem{problem}{Problem}[section]
\newtheorem{remark}{Remark}[section]
\newtheorem{example}{Example}[section]
\newtheorem{definition}{Definition}[section]
\newtheorem{lemma}{Lemma}[section]
\newtheorem{corollary}{Corollary}[section]
\newtheorem{proposition}{Proposition}[section]
\numberwithin{equation}{section}

\title{Distributionally robust shape and topology optimization}
\author{
C. Dapogny\textsuperscript{1}, J. Prando\textsuperscript{2} and B. Thibert\textsuperscript{2}
}
\maketitle

\begin{center}
\emph{\textsuperscript{1} Sorbonne Universit\'e, Universit\'e Paris Cit\'e, CNRS, Inria, Laboratoire Jacques-Louis
Lions, LJLL, F-75005 Paris, France}.\\
\emph{\textsuperscript{2} Univ. Grenoble Alpes, CNRS, Grenoble INP, LJK, 38000 Grenoble, France}.\\
\end{center}

\begin{abstract}
This article aims to introduce the paradigm of distributional robustness from the field of convex optimization to tackle optimal design problems under uncertainty. 
We consider realistic situations where the physical model, and thereby the cost function of the design to be minimized depend on uncertain parameters. The probability distribution of the latter is itself known imperfectly, through a nominal law, reconstructed from a few observed samples. 
The distributionally robust optimal design problem is an intricate bilevel program which consists in minimizing the worst value of a statistical quantity of the cost function (typically, its expectation) when the law of the uncertain parameters belongs to a certain ``ambiguity set''. We address three classes of such problems: firstly, this ambiguity set is made of the probability laws whose Wasserstein distance to the nominal law is less than a given threshold; secondly, the ambiguity set is based on the first- and second-order moments of the actual and nominal probability laws. Eventually, a statistical quantity of the cost other than its expectation is made robust with respect to the law of the parameters, namely its conditional value at risk.
Using techniques from convex duality, we derive tractable, single-level reformulations of these problems, framed over augmented sets of variables. Our methods are essentially agnostic of the optimal design framework; they are described  in a unifying abstract framework, before being applied to multiple situations in density-based topology optimization and in geometric shape optimization. Several numerical examples are discussed in two and three space dimensions to appraise the features of the proposed techniques.
\end{abstract}

\bigskip
\bigskip
\hrule
\tableofcontents
\vspace{-0.5cm}
\hrule
\bigskip
\bigskip

\section{\textbf{Introduction and Related Works}}\label{sec.intro}

\noindent Driven by the need for raw material and energy savings and enabled by the constant progress of high-performance computing, the recent advances in shape and topology optimization now make it possible to optimize devices as diverse as mechanical structures, airfoils and magnetic circuits from the early stages of the design process. However, the systematic deployment of optimal design techniques in industry faces critical challenges, related to their ability to accommodate the requirements of fabrication processes such as milling, casting and additive manufacturing \cite{langelaar2019topology,liu2018current,michailidis2014manufacturing,wang2020topology,zhou2002progress} and, closer to the purpose of this article, the inevitable uncertainties over the knowledge of real-life conditions.

\medskip\noindent\textit{Uncertainties in optimal design.} Optimal design is about minimizing a cost function depending on the design via its ``physical behavior'', under constraints of the same nature.
The mathematical models governing this behavior are inherently built from parameters: for instance, the displacement of a mechanical structure depends on the applied loads and on the properties of the constituent material (its Young's modulus and Poisson's ratio);
likewise, the velocity and pressure within a fluid are influenced by its viscosity and density, the electric and magnetic fields of a wave passing through a photonic device are functions of the operating frequency and the refraction index of the medium, etc.
In practice, these physical parameters are often subject to uncertainty, 
either because they are inaccurately measured (which is typical of the loads applied to a mechanical structure), or because they experience uncontrolled fluctuations over time: changes in the ambient medium, wear, etc. Designing systems that are robust to such uncertainties is crucial, as even minor deviations from the assumed values of the physical parameters during the optimization process can severely degrade the performance of the design.
This sharp sensitivity is strikingly illustrated by a situation exhibited in \cite{cherkaev1999optimal}, where the best elastic microstructure to withstand a specific loading scenario becomes the worst possible configuration when even infinitesimal perturbations of these loads occur.
Understandably enough, the ambition to predict designs with stable performance under uncertain conditions has inspired a thriving literature lately, see for instance \cite{maute2014topology} and \cite{acar2021modeling} for a recent overview.

\medskip\noindent\textit{Existing approaches for optimal design under uncertainties.}  Various paradigms can be envisioned to incorporate robustness in an optimal design problem, depending on the assumptions about the uncertain physical parameters of the model.
\begin{itemize}
    \item Worst-case approaches are usually adopted when no information is available, aside from a maximum bound on their amplitude: one minimizes the worst (i.e. maximum) value of the cost functional over all the possible perturbations. Such problems generally incur an intense computational effort due to their min-max bilevel structure \cite{guo2009confidence}.
They are tractable in a few very particular situations; for instance, the quadratic dependence of the compliance of a structure on the applied loads allows to reformulate its minimization under load uncertainty as a semi-definite program \cite{amstutz2016notion,de2008shape,holmberg2015worst}.  
In general, however, one has to resort to heuristic, approximate linearization techniques based on the ``smallness'' of the maximum amplitude of the uncertainties \cite{allaire2014linearized,guo2013robust}.
 In addition to their high computational cost, worst-case approaches are often overly pessimistic, as they result in designs showing poor nominal performances, for the sake of anticipating an unlikely worst-case scenario.    

\item Probabilistic approaches are naturally preferred in the presence of information about the uncertain parameters, e.g. about their law or some of its moments. 
Typically, the expected value of the cost functional, its standard deviation, or a weighted sum of both are minimized. Alternatively, so-called reliability-based formulations include constraints on the probability of failure of the system, i.e. on the probability that the cost exceed a certain safety threshold, see for instance \cite{hasofer1974exact,zhao1999general} about the popular First- and Second-Order Reliability Methods (FORM and SORM) in engineering.
    Both types of probabilistic frameworks raise the need to calculate integrals over the (high-dimensional) parameter space. This task can be carried out by expansive Monte Carlo sampling methods \cite{cardoso2019robust,chen2010level,hamdia2022multilevel}, accelerated by efficient stochastic collocation techniques and parallelized GPU implementations \cite{martinez2016large,martinez2016robust}.
 A different strategy consists in using the stochastic gradient algorithm to mitigate the number of samples required by the Monte-Carlo strategy, to the expanse of an increase in the number of optimization iterations \cite{de2019topologyoptimizationuncertaintyusing,de2021reliabilitybasedtopologyoptimizationusing,de2021topologyoptimizationmicroscaleuncertainty,jofre2021rapidaerodynamicshapeoptimization,pflug2024stochastic}, see also \cite{grieshammer2024continuous,grieshammer2024continuous2,uihlein2025140} for related strategies. Yet another possibility is to construct a reduced-order model for the solution of the boundary-value problem at play: assuming a finite-dimensional structure for the uncertain parameters (which typically stems from a truncated Karhunen-Lo\`eve expansion), this solution is approximated by a polynomial chaos expansion  \cite{keshavarzzadeh2017topology}, or a stochastic Galerkin method is used for its computation \cite{tootkaboni2012topology}. Eventually, in this setting also, linearization techniques based on the ``small'' amplitude of the uncertainties, are applicable to build deterministic approximations for the considered statistical quantities of the cost function of interest, see e.g. our previous work \cite{allaire2015deterministic} or \cite{guest2008structural}.
\end{itemize}

\medskip\noindent\textit{Distributionally robust optimization.} Classical probabilistic optimal design formulations themselves suffer from a conceptual shortcoming:
in realistic applications, the ``true'' probability distribution $\Ptrue$ of the
uncertain parameters is unknown; at best, it can be roughly approximated
based on (often scarce) observations. Distributionally robust optimization has
recently emerged as a sensible paradigm to enforce robustness with
respect to uncertainties in a way overcoming this objection. In
a nutshell, the worst-case of the expectation (or another statistical quantity, such as the standard deviation) of the cost function is minimized over all the probability laws $\Q$ within a so-called ambiguity set $\calA$: $\calA$ is made of the laws that are ``close'' to a given nominal scenario $\P$, which is in practice reconstructed from a few samples. Ambiguity sets of various natures can be thought of. The perhaps most intuitive practice
is to define $\calA$ as the set of laws $\Q$ lying within a certain distance from
$\P$; this distance can be measured in terms of
the Kullback-Leibler divergence \cite{hu2013kullback}, a maximum mean discrepancy
\cite{staib2019distributionally}, or the Wasserstein distance
\cite{esfahani2018data}. Other ambiguity sets
are made of laws whose first moments are close to those of $\P$ \cite{delage2010distributionally,nakao2021distributionally}.
We refer to
\cite{lin2022distributionally,rahimian2019distributionally,rahimian2022frameworks}
for overviews of distributionally robust optimization, to the book \cite{chen2020distributionally}
for a focus on machine learning, and to the recent surveys
\cite{blanchet2024distributionally,kuhn2024distributionallyrobustoptimization}.

This recent literature about distributional robustness arises in the
rich and structured framework of convex optimization. Unfortunately,
optimal design problems often lack such convenient features; moreover, they raise
specific difficulties such as the large cost entailed by the numerical evaluations
of the objective function and its derivative, as they depend on the solution to one or several physical boundary-value problems. 
Partly for these reasons, despite its natural character, distributionally robust optimization has
been seldom considered in engineering design: beyond our note
\cite{dapogny2023entropy}, we are only aware of very few contributions
in this direction, dealing with particular situations at that. In
\cite{kapteyn2019distributionally}, a model problem of parametric optimal
design is considered, where the uncertain parameter is allowed to take
a finite number of values. The experiments conducted in there clearly demonstrate
the ability of distributionally robust formulations to predict optimal designs
that perform well even in scenarios that are not present in the nominal
law. The article \cite{Kanno_2021} considers truss structures in a context
where the optimized thickness of the bars is subjected to uncertainty; a semi-definite program is derived to handle
the worst value of the failure probability of the system when the
moments of the law are ``close'' from reference values.

\medskip\noindent\textit{Summary of our contributions.}  
This article is the natural sequel of our preliminary note \cite{dapogny2023entropy}. Inspired by recent ideas in robust convex optimization and optimal transport theory, we develop tractable distributionally robust versions for a wide range of realistic optimal design problems, that are compatible with the particularities of this field, and notably its lack of mathematical structure and the large computational cost of evaluating the optimization criteria at play. 
More precisely, we consider three situations: 
\begin{itemize}
    \item The ambiguity set $\calA$ is the set of probability laws $\Q$ whose Wasserstein distance from the nominal law $\P$ is lower than a certain threshold;
    \item The ambiguity set $\calA$ is made of laws whose first- and second-order moments are close to reference values; 
    \item The statistical quantity which is made robust is no longer the mean value of the cost but its conditional value at risk -- a useful notion that allows to handle failure probability constraints.
\end{itemize}
In the spirit of \cite{azizian2022regularization,wang2021sinkhorn}, in all three cases, we rely on suitable entropy-based penalizations, that enable the use of convex duality techniques. Importantly, the proposed methodology is  agnostic of the shape and topology optimization framework: it is applicable to both density-based topology optimization and geometric shape optimization based on the method of Hadamard.
To emphasize this genericity, the discussion is framed within a formal, abstract optimal design setting insofar as possible. 

\medskip\noindent\textit{Organization of the article.}  
Using an abstract optimal design setting, \cref{sec.absframework} introduces the notion of distributional robustness, 
as well as the three distributionally robust problems analyzed in the sequel and their tractable reformulations. The next sections are devoted to the application of these concepts in two different optimal design settings:
\cref{sec.TO} considers the popular engineering framework of density-based topology
optimization. In there, we notably conduct several numerical experiments to appraise the
main features of our formulations. Then, \cref{sec.SO} deals with geometric
shape optimization. A conclusion and a few research perspectives are outlined in \cref{sec.concl}. The article ends with three appendices, devoted to
formal sketches of the proofs of the mathematical results underlying our developments.

\section{\textbf{Distributionally Robust Optimal Design Problems}}\label{sec.absframework}

\noindent
This section presents the three types of distributionally robust optimal design problems considered
in this article and details their tractable reformulation allowed by convex duality; it takes place in a formal and abstract setting, that we first introduce. 

\subsection{Presentation of the abstract distributionally robust optimal design framework}\label{sec.genpb}

\noindent
Let $\Uad$ be a subset of a vector space $\calH$, made of admissible designs $h$. 
In \cref{sec.TO}, $\Uad=L^{\infty}(D,[0,1])$ is made of the density functions $h$ over a given computational domain $D \subset \R^{d}$ ($d=2,3$). In \cref{sec.SO}, with a slight abuse of this framework, $\Uad$ is a set of ``true'' shapes $\Omega \subset \R^{d}$.

The mathematical model characterizing the physical behavior of a design $h \in \Uad$ involves parameters;
these are gathered in a vector $\xi$, lying in a compact subset
$\Xi$ of a finite-dimensional space $\R^{k}$. Assuming a perfect
knowledge of $\xi$, the optimal design problem of interest reads:
\begin{equation} \label{eq.absoptpb}
\min\limits_{h \in \Uad}\calC(h,\xi),
\end{equation}
where $\calC(h,\xi)$ denotes the cost of a design $h$ when the physical
parameters at play equal $\xi$. For simplicity of the exposition,
this formulation omits constraints; as demonstrated by the examples
in \cref{sec.TO,sec.SO}, their treatment does not entail any additional
difficulty from the conceptual viewpoint.

In practice, the parameters $\xi$ are only known imperfectly: they are
modeled as a $k$-dimensional random vector $\xi \equiv \xi(\omega)$, where
the event variable $\omega$ belongs to an abstract probability space $(\calO
,\calF,\mu)$. They are described by their law $\Ptrue$, that is an
element in the set $\calP(\Xi)$ of probability measures on $\Xi$: 
$$\text{For each subset } A \subset \Xi, \:\:\Ptrue(A) \in [0,1] =  \text{probability that }\xi \text{ belong to A.}$$
 As we have mentioned, a simple and popular method to make the nominal problem \cref{eq.absoptpb} aware of these uncertainties
is to optimize the mean value of the cost functional $\calC(h,\xi)$, i.e.:
\begin{equation}
	\min\limits_{h \in \Uad}\Jmean(h), \text{ where }\Jmean(h)= \int_{\Xi}
	\calC(h,\xi) \:\d \Ptrue(\xi).\label{eq.pbmean}
\end{equation}
Unfortunately, in realistic applications, the ``true'' law $\Ptrue$ of
$\xi$ is itself unknown; at best, an estimate $\P \in \calP(\Xi)$ is
available, which is for instance reconstructed from a few observed samples
$\xi^{1},\ldots,\xi^{N}\in \Xi$:
$$\P= \frac{1}{N}\sum\limits_{i=1}^{N}\delta_{\xi^{i}},$$
where $\delta_{\xi^i} \in \calP(\Xi)$ denotes the Dirac mass at $\xi^{i} \in \Xi$, $i=1,\ldots,N$.

Acknowledging this uncertainty about the law of $\xi$ leads to the distributionally robust version of \cref{eq.absoptpb}, minimizing
the worst possible value of the expected cost:
\begin{equation}\label{eq.pbdrogen}
\min\limits_{h \in \Uad}\Jdr(h), \text{ where }\Jdr
	(h) = \sup\limits_{\Q \in \calA}\int_{\Xi}\calC(h,\xi) \:\d \Q(\xi).
\end{equation}
This formulation features an ambiguity set $\calA \subset \calP(\Xi)$, made
of laws $\Q \in \calP(\Xi)$ which are ``close'' to $\P$. Depending on
the available information about $\P$, the set $\calA$ may take different forms. This article proposes three variations around this general concern:
\begin{itemize}
\item In \cref{sec.wdro}, we consider problems where the ambiguity set $\calA$
is made of the probability laws $\Q$ which are ``close'' to $\P$ in
terms of (a suitable version of) the Wasserstein distance. 
\item In \cref{sec.momdro},
$\calA$ is defined as the set of probability laws
whose first and second-order moments are ``close'' to given reference values up to a certain
threshold.
\item \cref{sec.cvardro} eventually deals with a slightly different situation,
where another measure of risk than the mean value of the cost function
is made robust, namely, its conditional value at risk.
\end{itemize}

\begin{remark}\label{rem.Adeph}
\noindent
\begin{itemize}
\item In certain applications, the ambiguity set $\calA \equiv \calA(h)$ in \cref{eq.pbdrogen} depends on the design $h$ itself,
see e.g. \cite{drusvyatskiy2023stochastic}. In the (scarce) literature devoted to such problems, \cref{eq.pbdrogen} is treated by ``freezing'' the ambiguity set at the current iteration of the optimization process. Thus, the solution of a version of \cref{eq.pbdrogen} featuring a design-dependent ambiguity set boils down to that of a series of such problems where this set does not depend on $h$. For this reason, throughout the sequel, we limit ourselves to considering design-independent ambiguity sets $\calA$. 
\item The assumption that the uncertain parameter $\xi$ should lie in a finite-dimensional vector space $\R^k$ may seem restrictive at first glance.
In practice, though, infinite-dimensional random vectors are often approximated by a finite-dimensional reduced-order model, obtained e.g. via a truncated Karhunen-Lo\`eve expansion, see \cref{subsec.matdro}.
\end{itemize}
\end{remark}
\subsection{Ambiguity sets based on the Wasserstein distance to a nominal law}\label{sec.wdro}

\noindent
This section considers the situation where the ambiguity set $\calA$ in
\cref{eq.pbdrogen} is defined as the set of probability laws which are ``close''
to the nominal law $\P$ in terms of a suitable distance between probability
measures, namely, the (entropy-regularized) Wasserstein distance from optimal
transport theory. We start with a brief reminder about this notion.

\subsubsection{Preliminaries about the Wasserstein distance}

\noindent
Recalling that $\Xi$ is a compact subset of the finite-dimensional
space $\R^{k}$, let $\P$, $\Q$ be two elements in the set $\calP(\Xi)$
of probability measures on $\Xi$. The Wasserstein distance $W(\P,\Q)$ between
$\P$ and $\Q$ is defined as the minimum cost of sending the mass of $\P$
onto that of $\Q$, precisely:
\begin{equation}\label{eq.defW}
W(\P,\Q) = \min \left\{ \int_{\Xi \times \Xi}c(\xi,\zeta
	) \: \d \pi(\xi,\zeta), \: \pi \in \calP(\Xi \times \Xi), \:\: \pi_{1}
	= \P, \: \pi_{2}= \Q \right\}.
\end{equation}
This definition involves couplings $\pi$, that is, probability measures
on the product set $\Xi \times \Xi$; the first and second marginals of
such a coupling are respectively denoted by $\pi_{1}$ and $\pi_{2}$. These are defined through their action on the set $\calC(\Xi)$ of continuous functions on $\Xi$:
$$\forall \varphi \in \calC(\Xi), \quad \int_{\Xi} \varphi(\xi) \:\d\pi_{1}(\xi)
	= \int_{\Xi \times \Xi}\varphi(\xi) \:\d\pi(\xi,\zeta), \text{ and }\int
	_{\Xi} \varphi(\zeta) \:\d\pi_{2}(\zeta)= \int_{\Xi \times \Xi}\varphi(\zeta)
	\:\d\pi(\xi,\zeta).$$
Intuitively, $\pi(\xi,\zeta)$ represents the amount of mass transferred from
$\xi$ to $\zeta$ according to the transport plan $\pi$, $\pi_{1}(\xi)$ is the total mass leaving $\xi$ and
$\pi_{2}(\zeta)$ is that arriving at $\zeta$, so that the marginal constraints in \cref{eq.defW} encode that $\pi$ transports the ``mass distribution'' $\P$ onto $\Q$. The ``ground cost'' $c(\xi,\zeta)$ evaluates the cost of transporting one unit of mass from $\xi$
to $\zeta$; throughout the article, we use the quadratic cost 
\begin{equation}\label{eq.quadcost}
\forall \xi, \zeta \in \Xi, \quad c(\xi,\zeta) = \lvert \xi - \zeta\lvert^{2}.
\end{equation}

The Wasserstein distance is a principled notion of distance between probability
measures. It has recently aroused much enthusiasm to compare shapes
under extraneous formats in imaging \cite{feydy2020analyse} and as a
loss function between discrete probability distributions in machine
learning \cite{montesuma2024recent}. Notably, it allows to smoothly
discriminate probability distributions with different supports, contrary to the Kullback-Leibler divergence. We refer for
instance to
\cite{merigot2021optimal,peyre2019computational,santambrogio2015optimal,villani2009optimal}
for mathematical background about optimal transport theory.

\begin{remark}
Strictly speaking, it is the square root of the quantity $W(\P,\Q)$ in \cref{eq.defW} (in the present situation where the ground cost $c(\xi,\zeta)$ is given by \cref{eq.quadcost}) that is the distance over $\calP(\Xi)$ called Wasserstein distance. For simplicity of the terminology, we abusively refer to $W(\P,\Q)$ as the Wasserstein distance.
\end{remark}

In our applications, we use a regularized version of
the Wasserstein distance. Let us define the entropy $H(\pi)$ of a
coupling $\pi \in \calP(\Xi \times \Xi)$ by:
\begin{equation}\label{eq.Hpi}
H(\pi) = \left\{
	\begin{array}{cl}
		\displaystyle\int_{\Xi \times \Xi} \log \frac{\d \pi}{\d \pi_{0}} \:\d\pi & \text{if } \pi \text{ is absolutely continuous w.r.t. } \pi_0, \\
		+ \infty                                                                  & \text{otherwise.}
	\end{array}
	\right.
\end{equation}
Following \cite{azizian2022regularization}, the above definition of $H(\pi)$ features a reference coupling $\pi_{0} \in \calP(\Xi \times \Xi)$, which is given by:
\begin{equation}\label{eq.refcouplingW}
\pi_{0}(\xi,\zeta) = \P(\xi) \otimes \nu_{\xi}
	(\zeta), \text{ where }\nu_{\xi}(\zeta) = \alpha_{\xi}e^{-\frac{\lvert
	\xi-\zeta\lvert^{2}}{2\sigma^2}}\mathds{1}_{\Xi}(\zeta),
\end{equation}
and the normalization constant $\alpha_{\xi}$ is tailored so that
$\nu_{\xi}\in \calP(\Xi)$. More explicitly, the action of $\pi_{0}$
on an arbitrary continuous function $\varphi \in \calC(\Xi \times \Xi)$
reads:
$$\int_{\Xi \times \Xi}\varphi(\xi,\zeta) \:\d\pi_{0}(\xi,\zeta) = \int_{\Xi}
	\left( \int_{\Xi}\varphi(\xi, \zeta ) \:\d\nu_{\xi}(\zeta)\right)\:\d\P
	(\xi).$$
In \cref{eq.Hpi}, a coupling $\pi$ is called absolutely continuous with respect to $\pi_0$ if there exists an integrable function $\alpha \in L^1(\Xi \times \Xi; \d\pi_0)$ such that $\pi = \alpha(\xi,\zeta) \pi_0$; the latter is often denoted by $\alpha = \frac{\d\pi}{\d\pi_0}$.
	
We eventually define the entropy-regularized Wasserstein distance $W_{\e}$ via a minimization problem similar to that in \cref{eq.defW}, up to an additional penalization by the entropy of couplings: 
\begin{equation}\label{eq.RegWass}
W_{\e}(\P,\Q) = \min \left\{ \int_{\Xi \times \Xi}c(
	\xi,\zeta) \:\d \pi(\xi ,\zeta) + \e H(\pi), \: \pi \in \calP(\Xi \times
	\Xi), \:\: \pi_{1}= \P, \: \pi_{2}= \Q \right\}.
\end{equation}
This addition of a small, strictly convex log-barrier term to the optimal transport problem was proposed in \cite{cuturi2013sinkhorn}, with the purpose to ease numerical computations. Since then, this ``blurred'' optimal transport has found an interest of its own; for instance, it is often preferred to the original optimal transport formulation in supervised machine learning, where it helps to avoid overfitting. We refer to the lecture notes \cite{nutz2021introduction} about entropy-regularized optimal transport.

\begin{remark}
The parameter $\sigma$ in the definition \cref{eq.refcouplingW} of $\nu_\xi(\zeta)$ reflects the degree of confidence of the user in the nominal law $\P$. Loosely speaking, when $\sigma$ is ``small'', $\pi_0$ reduces to the trivial coupling $\P(\xi) \otimes \delta_\xi(\zeta)$, sending $\P$ onto itself by leaving in place each mass unit $\xi \in \Xi$. Then, the entropy penalization $H(\pi)$ in the regularized Wasserstein distance \cref{eq.RegWass} imposes that only couplings $\pi$ very close to this ``trivial'' one should be considered.
When $\sigma$ gets larger, $\pi_0$ acts by ``spreading the mass'' of $\P$ at each point $\xi \in \Xi$ over a neighborhood with increasing size. The couplings $\pi$ with ``small'' entropy values will then have second marginals corresponding to more and more ``blurred'' versions of $\P$. 
\end{remark}

\begin{remark} \label{rem.Wenodist} 
Contrary to what is suggested by its misleading name, the quantity $W_{\e}(\P,\Q)$ is solely an approximation of (the square of) a distance. For instance, $W_{\e}(\P,\P)$ is often different from $0$. The function $W_{\e}$ can be modified into a proper (squared) distance by setting:
$$S_{\e}(\P,\Q) := W_{\e}(\P,\Q) - \frac{1}{2}W_{\e}(\P,\P) - \frac{1}{2} W_{\e}(\Q,\Q).$$
This so-called Sinkhorn distance is significantly more complicated to handle;
in particular, it does not allow for convenient reformulations of distributionally robust problems similar to those derived below using $W_{\e}$. We refer to \cite{feydy2020analyse,sejourne2019sinkhorn} for more details about the Sinkhorn distance.
\end{remark}

\subsubsection{Distributionally robust optimization using a Wasserstein ambiguity set} \label{sec.wdroformula}

\noindent
In this section, we study the instance of the distributionally robust
optimal design problem \cref{eq.pbdrogen} where the ambiguity set $\calA$
is that $\calA_{\text{W}}$ defined by:
\begin{equation}\label{eq.calAW}
\calA_{\text{W}}= \Big\{\Q \in \calP(\Xi), \:\: W_{\e}
	(\P,\Q) \leq m \Big\},
\end{equation}
i.e. $\calA_{\text{W}}$ is made of the probability laws $\Q \in \calP(\Xi)$ that are within a distance $m>0$ from the given nominal law $\P$.
The distributionally robust problem \cref{eq.pbdrogen} rewrites in this case: 
\begin{equation}\label{eq.pbwdro}
\min\limits_{h \in \Uad}\JW(h), \text{ where }\JW(h) = \sup\limits_{\Q \in \calA_{\text{W}}}\int_{\Xi}\calC(h,\xi) \:\d \Q(\xi).
\end{equation}
The idea of considering ambiguity sets based on the Wasserstein distance to a fixed nominal law was originally introduced in \cite{pflug2007ambiguity}, see \cite{esfahani2018data,kuhn2019wasserstein} for the subsequent developments of this idea and the review \cite{blanchet2021statistical}. The precise formulation involving the entropy-regularized ambiguity set in \cref{eq.RegWass,eq.calAW} was proposed in \cite{azizian2022regularization,wang2021sinkhorn} and then reused in our note \cite{dapogny2023entropy}. 
We also refer to \cite{azizian2023exact,fournier2013rateconvergencewassersteindistance} for an account of its generalization guarantees and to \cite{vincent2024texttt} for a recent, open-source implementation of this framework.

The formulation \cref{eq.pbwdro} is a priori intricate, as it is made of nested minimization and maximization problems.
The key point of our methodology is that the worst-case functional
$\JW(h)$ can be conveniently reformulated as an infimum, thanks to the
following result from convex duality. The latter is proved rigorously in
\cite{azizian2022regularization}. For completeness, and since its main argument
will be reused, a short and formal sketch of the proof is provided in
\cref{app.wdro}.

\begin{proposition}\label{prop.formulaJdir} 
Let $\Xi$ be a compact subset of $\R^k$, and let $f : \Xi \to \R$ be a continuous function. 
Recalling the definition \cref{eq.calAW} of the ambiguity set $\calA_{\text{W}}$,
the following formula holds true:
\begin{equation}\label{eq.reformspwdro}
\sup\limits_{\mathbb{Q}\in \mathcal{A}_{\text{\rm W}}}\int_{\Xi}f(\xi
		) \:\d \Q(\xi) = \inf\limits_{\lambda \geq 0}\left( \lambda m + \lambda
		\e \int_{\Xi}\log\left( \int_{\Xi}e^{\frac{f(\zeta) - \lambda c(\xi,\zeta)}{\lambda\e}}
		\:\d\nu_{\xi}(\zeta) \right)\:\d\P(\xi)\right).
\end{equation}
\end{proposition}

\begin{remark}\label{rem.lambda}
As highlighted by the argument in \cref{app.wdro}, the optimization variable $\lambda$ in the right-hand minimization problem of \cref{eq.reformspwdro} is the Lagrange multiplier for the inequality constraint $W_\e(\P,\Q) \leq m$ of the maximization problem at the left-hand side. Hence, the optimal point $\lambda^*$ accounts for the sensitivity of the optimal value of the quantity $\int_{\Xi} f(\xi) \:\d\Q(\xi)$ when the threshold $m$ of this constraint is perturbed. In particular, according to the complementarity relations of the first-order necessary conditions of optimality for \cref{eq.reformspwdro},
\begin{itemize}
\item If an optimal law $\Q^*$ in the maximization problem at the left-hand side of \cref{eq.reformspwdro} does not saturate the constraint, i.e. if $W_\e(\P,\Q^*) < m$, then $\lambda^* =0$: the optimal value of $\int_{\Xi} f(\xi) \:\d\Q(\xi)$ does not change if $m$ is slightly modified. 
\item On the contrary, if $W_\e(\P,\Q^*) = m$, then $\lambda^*$ may differ from $0$, and it often does in practice: the optimal value of $\int_{\Xi} f(\xi) \:\d\Q(\xi)$ may change significantly as soon as $m$ is perturbed. 
\end{itemize}
We refer to e.g. \cite{nocedal2006numerical} about the classical interpretation of the Lagrange multiplier in optimization theory. 
\end{remark}

\cref{prop.formulaJdir} allows to rewrite the distributionally robust problem \cref{eq.pbwdro}
in the following way:
\begin{equation}\label{eq.wdrogen}
	\min\limits_{(h,\lambda) \in \Uad \times \R_+}\calD_{\text{W}}(h,\lambda
	) \text{ where }\calD_{\text{W}}(h,\lambda) = \lambda m + \lambda\e \int
	_{\Xi}\log\left( \int_{\Xi}e^{\frac{\calC(h,\zeta) - \lambda c(\xi,\zeta)}{\lambda\e}}
	\:\d\nu_{\xi}(\zeta) \right)\:\d\P(\xi).
\end{equation}
This new expression is much more amenable to the numerical treatment than the original one \cref{eq.pbwdro}, as it features a single minimization, over the
product space $\Uad \times \R_{+}$. Problem \cref{eq.wdrogen} can be solved owing to standard numerical methods, such as an alternating descent algorithm, see \cref{sec.TO,sec.SO} for illustrations.

\begin{remark}
\noindent
\begin{itemize}
\item A result similar to \cref{prop.formulaJdir} holds when the ambiguity set $\calA_{\text{W}}$ is constructed from the classical Wasserstein distance, i.e. when
$W_\e$ is replaced by $W$ in \cref{eq.calAW}. The main difference lies in the replacement of the ``log-sum-exp" operation in \cref{eq.reformspwdro} by its formal limit as $\e \to 0$, which is a supremum.  
\item Different notions of distance between probability measures could be employed in the definition of the ambiguity set $\calA$, such as the Kullback-Leibler divergence, the total variation distance, etc. We refer to \cite{kuhn2024distributionallyrobustoptimization} for a more exhaustive list of the possibilities and in-depth discussions about the features of these formulations.
	\end{itemize}
\end{remark}

%
%
%
%
%

\subsection{Ambiguity sets based on moments}\label{sec.momdro}

\noindent
This section assumes that the only available information about the law of the uncertain parameter $\xi$ is about its mean
value $\mu_{0}\in \R^{k}$ and covariance matrix
$\Sigma_{0}\in \R^{k\times k}$, as was initially proposed in \cite{delage2010distributionally},
see also \cite{hanasusanto2015distributionally,liu2018note}. The ambiguity
set $\calA_{M}$ containing the possible probability laws for $\xi$ is defined by:
\begin{equation}\label{eq.calAmoments}
\calA_{\text{M}}= \left\{ \Q \in \calP(\Xi), \:\:
	 \left\lvert\int_{\Xi}\xi \:\d \Q(\xi) -\mu_{0}\right\lvert \leq m_{1}
	, \text{ and } \int_{\Xi}(\xi - \mu_{0}) \otimes (\xi - \mu_{0}) \:\d
	\Q(\xi) \leq m_{2}\Sigma_{0}\right \}.
\end{equation}
Here, $m_1>0$ and $m_2 >0$ are given bounds and the second inequality is understood in the sense of positive semi-definite
matrices.

In this situation, taking inspiration from the analysis conducted in \cref{sec.wdro}, we shall derive a tractable alternative form for the distributionally robust problem \cref{eq.pbdrogen} by penalizing the mean value $\int_\Xi \calC(h,\xi) \:\d\Q(\xi)$ of the cost function by a ``small'', strictly convex entropy term. The latter,  that we denote by $H(\Q)$ (not to be confused with the entropy $H(\pi)$ of couplings in \cref{eq.Hpi}), is defined by:
\begin{equation}\label{eq.entropymoments}
H(\Q) = \left\{
	\begin{array}{>{\displaystyle}cl}
		\int_\Xi \log\frac{\d\Q}{\d\Q_{0}} \:\d\Q & \text{if } \Q \text{ is absolutely continuous w.r.t. } \Q_0, \\
		+\infty                                   & \text{otherwise}.
	\end{array}
	\right.
\end{equation}
Here, the reference law $\Q_{0} \in \calP(\Xi)$ is (the restriction to $\Xi$ of) the Gaussian distribution with center $\mu_0$ and covariance matrix $\Sigma_0$: 
$$\Q_0 = \alpha e^{-\frac12 (\xi-\mu_0)\Sigma_0^{-1} (\xi-\mu_0)} \:\mathds{1}_\Xi (\xi) \:\d \xi, \text{ where } \alpha \text{ is a normalization constant.} $$
These ingredients now allow to introduce the following moment-based distributionally robust optimal design problem:
\begin{equation}\label{eq.minJM}
\min\limits_{h \in \Uad}\JM(h), \text{ where }\JM(h) :
	= \sup\limits_{\Q \in \calA_{\text{M}} }\left( \int_{\Xi}\calC(h,\xi)
	\:\d\Q(\xi) - \e H(\Q)\right).
\end{equation}
Again, this a priori intricate problem is made tractable thanks to the following proposition, which converts the supremum functional $J_{\text{M}}(h)$ into an infimum;
we refer to \cite{prando2025distributionally} for a complete statement and a rigorous proof, a sketch of which is provided in \cref{app.moment}.

\begin{proposition}\label{prop.formulamoment} 
Let $\Xi\subset \mathbb{R}^{k}$ be compact, and let $f : \Xi \to \R$ be a continuous function. 
The ambiguity set $\calA_{\text{\rm M}}$ being defined by \cref{eq.calAmoments} with $m_{1}>0, m_{2}>0$, it holds:
\begin{multline}\label{eq.dualmoments}
\sup\limits_{\Q \in \calA_{\text{\rm M}} }\left
		( \int_{\Xi}f(\xi) \:\d \Q(\xi) -\e H(\Q)\right) = \inf\limits
		_{\lambda \geq 0, \: \lvert \tau \lvert \leq 1 ,\atop
		{ S \in \mathbb{S}^k_+}}\Bigg\{ \lambda m_{1}- \lambda \tau \cdot \mu
		_{0}+ m_{2}S : \Sigma_{0}\\[-1em] + \e \log \left( \int_{\Xi}\left(e
		^{\frac{f(\xi) + \lambda \tau \cdot \xi - S : (\xi-\mu_{0}) \otimes
		(\xi-\mu_{0})}{\e}}\right) \:\d \Q_{0}(\xi)\right) \Bigg\},
\end{multline}
where $\mathbb{S}_{+}^{k}$ is the set of $k \times k$ symmetric and positive semi-definite matrices.
\end{proposition}
The above proposition paves the way to the following alternative formulation for \cref{eq.minJM}:
\begin{multline*}
	\min\limits_{h \in \Uad, \:\lvert  \tau  \lvert\leq 1, \atop
	{\lambda \geq 0 , \: {S \in \mathbb{S}_+^k} } }\calD_{\text{M}}(h,\lambda
	,\tau,S), \text{ where }\\[-1em] \calD_{\text{M}}(h,\lambda,\tau,S):=
	\lambda m_{1}- \lambda \tau \cdot \mu_{0}+ m_{2}S : \Sigma_{0}+ \e \log
	\left( \int_{\Xi}\left(e^{\frac{\calC(h,\xi) + \lambda \tau \cdot \xi
	- S : (\xi-\mu_{0}) \otimes (\xi-\mu_{0})}{\e}}\right) \:\d \Q_{0}(\xi
	)\right).
\end{multline*}
Again, this minimization problem, posed over an augmented set, can be solved by standard numerical algorithms for optimal design. 

\subsection{Distributional uncertainties using the conditional value at risk}\label{sec.cvardro}

\noindent
The previous \cref{sec.wdro,sec.momdro} have considered distributionally robust functionals accounting for the worst-case value of the expectation of the cost $\calC(h,\xi)$ of a design $h$ when the law of $\xi$ belongs to a certain ambiguity set $\calA$. Here, we derive robust
versions of probabilistic optimal design problems involving a statistical quantity of the cost which is different from its expectation,
namely, its conditional value at risk. In spite of its popularity, for instance in mathematical finance, the use
of this notion in optimal design seems to be relatively confidential, see
however \cite{eigel2018risk} and \cite{kouri2016risk} in the more general context of optimal control of partial differential equations. We exemplify this use in the perspective of reliability-based
optimal design, where the conditional value at risk plays the role of a
surrogate for a probability of failure.

\subsubsection{A few basic facts about value at risk and conditional value at risk}\label{sec.cvarbg}

\noindent
This section collects a few useful facts about the notions of value at risk
and conditional value at risk in the context of interest to us; we
refer to \cite{rockafellar2002conditional,rockafellar2000optimization}
for details. For the moment, let us assume that the probability law $\P \in \calP(\Xi)$ of the parameters $\xi$ of the physical model is known perfectly. 
For a fixed design $h \in \Uad$, we denote by $t\mapsto \Psi(h,t)$ the cumulative
distribution function of the cost $\calC(h,\cdot)$ when the latter is seen as a function of the
random parameter $\xi \in \Xi$, that is:
$$\forall t \in \R, \quad \Psi(h,t) = \P \Big\{ \xi \in \Xi, \:\: \calC(h,\xi) \leq t \Big\}.$$
We recall that the function $t \mapsto \Psi(h,t)$ is non decreasing, with limits:
$$\lim\limits_{t \to -\infty }\Psi(h,t) = 0, \text{ and }\lim\limits_{t \to +\infty}\Psi(h,t) = 1.$$
This function is continuous from the right, but it fails to be continuous from the left at those points $t \in \R$ such that $\calC(h,\cdot)$ equals $t$ on a subset of $\Xi$ with positive measure, if any \cite{billingsley2017probability}. For simplicity, we suppose throughout this presentation that $\Psi(h,\cdot)$ is continuous; although the definitions and intuitions provided below require technical adjustments, the conclusions of our analyses remain valid in the case where it is not. \par\medskip

Let us now define the main objects at stake.

\begin{definition} \label{def.varcvar} 
Let $\beta \in [0,1]$ be a given threshold; then,
\begin{itemize}
\item The $\beta$-value at risk (or $\beta$-quantile) {\normalfont $\VaR_{\beta}(\calC(h,\cdot))$} of the cost $\calC(h,\cdot)$ is the smallest value $t \in \R$ greater than $\calC(h,\cdot)$ with probability $\beta$, i.e. 
$${\normalfont \VaR_\beta(\calC(h,\cdot))}= \inf\limits \Big\{ t \in \R , \:\: \Psi(h,t) \geq \beta \Big\}.$$
\item If $\beta < 1$, the $\beta$-conditional value at risk \normalfont $\CVaR_{\beta}(\calC(h,\cdot))$ \it of $\calC(h,\cdot)$ is the average of $\calC(h,\cdot)$ over those events where it exceeds the threshold {\normalfont $\VaR_{\beta}(\calC(h,\cdot))$}:
\begin{equation}\label{eq.defcvar}
{\normalfont \CVaR_\beta(\calC(h,\cdot))}= \frac{1}{1-\beta}\int_{\left\{ \xi \in \Xi, \:\: \calC(h,\xi) \geq {\normalfont \VaR_\beta(\calC(h,\cdot))} \right\}}\calC(h,\xi) \:\d \P(\xi).
\end{equation}
\end{itemize}
\end{definition}

\begin{figure}[ht] \centering
\includegraphics[width=0.7\textwidth]{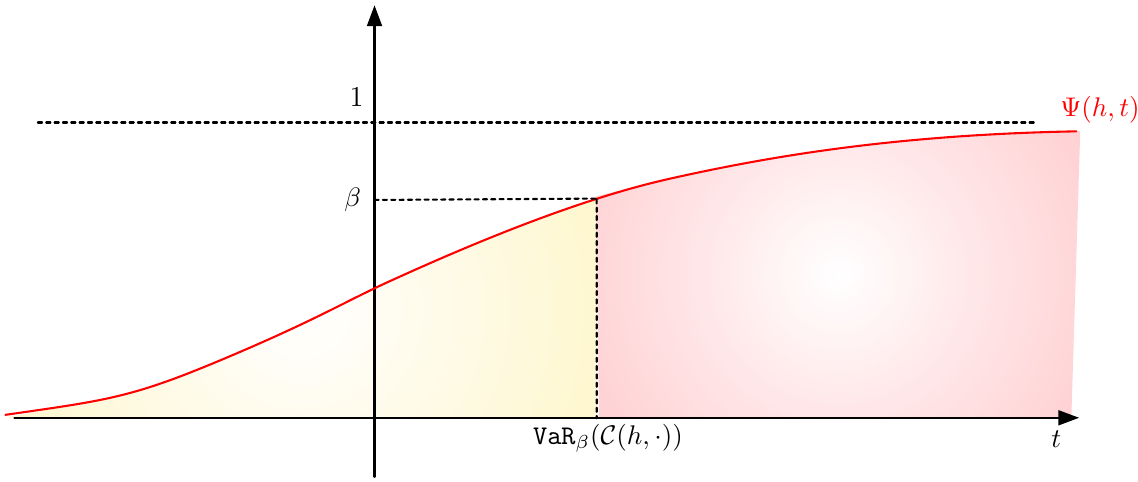}
\caption{\it The $\beta$-value at risk ${\normalfont \VaR_\beta(\calC(h,\cdot))}$ of the cost $\calC(h,\cdot)$ is the lowest level where the cumulative distribution function $\Psi(h,\cdot)$ takes values greater than $\beta$; the $\beta$-conditional value at risk ${\normalfont \CVaR_\beta(\calC(h,\cdot))}$ is the integral of $\calC(h,\cdot)$ over those events where it exceeds this value (red area).}
\label{fig.PhiVaR}
\end{figure}

Note that, by definition of {\normalfont $\VaR_{\beta}(\calC(h,\cdot))$}, it holds:
\begin{equation}\label{eq.probavar}
\P \Big\{ \xi \in \Xi, \:\: \calC(h,\xi) \geq{\normalfont \VaR_\beta(\calC(h,\cdot))}\Big\} = 1-\beta,
\end{equation}
and so \cref{eq.defcvar} is indeed an average. In particular, it holds:
\begin{equation}
	\label{eq.VaRCVaR}{\normalfont \VaR_\beta(\calC(h,\cdot))}\leq{\normalfont \CVaR_\beta(\calC(h,\cdot))}.
\end{equation}
Intuitively, the cost $\calC(h,\cdot)$ takes values lower than the $\beta$-value at risk 
${\normalfont \VaR_\beta(\calC(h,\cdot))}$ with probability $\beta$, see \cref{fig.PhiVaR}. 
This quantity is not fully satisfactory, as it does not supply any information about the tail of the distribution of $\calC(h,\cdot)$, i.e. those events $\xi$ where $\calC(h,\xi)$ exceeds {\normalfont $\VaR_{\beta}(\calC(h,\cdot))$}. 
On the contrary, the conditional value at risk ${\normalfont \CVaR_\beta(\calC(h,\cdot))}$ is precisely the average of $\calC(h,\cdot)$ over these scenarios.

The handling of conditional value at risk is eased by the following
representation as a minimum value. The latter is proved in \cite{acerbi2002coherence,rockafellar2002conditional}; for the sake
of completeness, a formal and intuitive argument is provided in \cref{app.cvar}.

\begin{proposition}\label{prop.cvarrep} 
The following representation formula holds true:
$${\normalfont \CVaR_\beta(\calC(h,\cdot))}= \inf \limits_{\alpha \in
		\R}\left\{ \alpha + \frac{1}{1-\beta}\int_{\Xi}\left[ \calC(h,\xi) -
		\alpha \right]_{+}\:\d \P(\xi)\right\}, \text{ where } [t]_+ := \max(t,0).$$
The minimum in the above program is uniquely attained at $\alpha = {\normalfont \VaR_\beta(\calC(h,\cdot))}$.
\end{proposition}
\subsubsection{Reliability-based optimal design using the conditional value at risk}\label{sec.relbasedOD}

\noindent The notion of conditional value at risk paves the way to a convenient reformulation of optimization
problems featuring failure probability constraints, of the form
$$\P \Big\{ \xi \in \Xi , \:\: \calC(h,\xi) \geq C_{T} \Big\} \leq 1 - \beta,$$
for a parameter $\beta \in (0,1)$ and a given safety threshold $C_T \in \R$, see e.g. \cite{chaudhuri2020risk}.
Indeed, it follows from the definition of $\Psi(h,t)$ that:
$$\P \Big\{ \xi \in \Xi , \:\: \calC(h,\xi) \geq C_{T} \Big\} \leq 1 - \beta \text{ if and only if }\Psi(h,C_{T}) \geq \beta.$$
Since $\VaR_{\beta}(\calC(h,\cdot))$ is by definition the smallest value $t \in \R$ such that
the cost is lower than $t$ with probability $\beta$, this is in turn equivalent to:
$$\VaR_{\beta}(\calC(h,\cdot)) \leq C_{T}.$$
Now invoking \cref{eq.VaRCVaR}, a conservative surrogate for this constraint is:
$$\CVaR_{\beta}(\calC(h,\cdot)) \leq C_{T}.$$
\par\medskip 

Let us apply this observation to a model optimal design problem featuring an objective function $J(h)$ and a constraint on the probability that the cost exceed the threshold
$C_{T}$:
$$\min\limits_{h \in \Uad}J(h) \text{ s.t. }\P\Big\{\xi \in \Xi \text{ s.t. } \calC(h,\xi) \geq C_{T} \Big\} \leq 1-\beta.$$
The foregoing discussion naturally leads to consider the following problem, featuring a more stringent constraint:
$$\min\limits_{h\in \Uad}\; J(h) \: \text{ s.t. }\texttt{CVaR}_{\beta}(\mathcal{C}(h,\cdot)) \leq C_{T}.$$
Using the representation formula for the Conditional Value at Risk in \cref{prop.cvarrep},
this problem rewrites as a more conventional optimization
problem over the pair $(h,\alpha)$:
\begin{equation}\label{eq.CVARRBSO}
\min\limits_{h\in \Uad, \atop \alpha \in \R}\; J(h)
	\text{ s.t. }\;\alpha + \frac{1}{1-\beta} \int_{\xi\in \Xi}\left[\mathcal{C}
	(h,\xi)-\alpha\right]_{+}\;\mathrm{d}\mathbb{P}(\xi) \leq C_{T}.
\end{equation}


\begin{remark}
In numerical practice, the non smooth function $\left[ \cdot \right]_{+}$ featured in \cref{eq.CVARRBSO} is replaced by the smooth approximation:
$$\left[t \right]_{+} \approx \dfrac{t}{1+e^{-\gamma t}}, \text{ where }\gamma := 20.$$
\end{remark}

\subsubsection{Distributionally robust formulation of problems involving the conditional value at risk}

\noindent
Let us now address the situation where the law of the parameter $\xi$ is itself uncertain. To set ideas, we assume that this law belongs to an
ambiguity set $\calA_{\text{W}}$ of Wasserstein type around a nominal law $\P$, see \cref{eq.calAW}.
The distributionally robust counterpart of \cref{eq.CVARRBSO} then reads:
\begin{equation}
    \label{eq.CVARDRO}
	\min\limits_{h\in \Uad}\; J(h) \text{ s.t. }\; \sup\limits_{\Q \in
	\calA_{\text{W}}}\min\limits_{\alpha \in \R}\left(\alpha + \dfrac{1}{1-\beta}
	\int_{\xi\in \Xi}\left[\mathcal{C}(h,\xi)-\alpha\right]_{+}\;\mathrm{d}
	\Q(\xi)\right) \leq C_{T}.
\end{equation}
We then use the following proposition, whose proof is similar to that
of \cref{prop.formulaJdir}.

\begin{proposition}
Let $\Xi$ be a compact subset of $\R^k$, and let $f : \Xi \to \R$ be a continuous function. Recall that $\calA_{\text{\rm W}}$  denotes the subset of $\calP(\Xi)$ defined by \cref{eq.calAW}. The following equality holds true:
	\begin{multline*}
		\sup\limits_{\Q \in \calA_{\text{\rm W}}}\min\limits_{\alpha \in \R}\left
		(\alpha + \dfrac{1}{1-\beta}\int_{\xi\in \Xi} f(\xi) \;\mathrm{d}\Q(\xi)\right) =\\ \inf\limits
		_{\lambda \geq 0, \atop \alpha \in \R}\left( \alpha + \frac{\lambda
		m}{1-\beta}+ \frac{ \lambda \e }{1-\beta}\int_{\Xi}\log\left( \int_{\Xi}
		e^{\frac{f(\zeta) - \lambda c(\xi,\zeta)}{\lambda\e}}
		\:\d\nu_{\xi}(\zeta) \right)\:\d\P(\xi) \right) .
	\end{multline*}
\end{proposition}

As a result, the distributionally robust optimal design problem \cref{eq.CVARDRO} rewrites:
\begin{multline*}
	\min\limits_{h\in \Uad, \atop \lambda \geq 0, \alpha \in \R}\; J(h) \:
	\text{ s.t. }\; \calD_{\text{C}}(h,\lambda,\alpha) \leq C_{T}, \text{
	where }\\ \calD_{\text{C}}(h,\lambda,\alpha) = \alpha + \frac{\lambda m}{1-\beta} + \frac{\lambda \e}{1-\beta}\int_{\Xi}\log\left( \int_{\Xi}
	e^{\frac{\left[\mathcal{C}(h,\zeta)-\alpha\right]_{+} - \lambda c(\xi,\zeta)}{\lambda\e}}
	\:\d\nu_{\xi}(\zeta) \right)\:\d{\P}(\xi).
\end{multline*}
This is a classical constrained optimal design problem, posed over the product space $\Uad \times \R_+ \times \R$, which can be solved by standard numerical algorithms.

\section{\textbf{Distributionally Robust Density-Based Topology Optimization}}\label{sec.TO}

\noindent
This section slips into the framework of density-based optimal design of mechanical
structures, about which we refer to the classical book \cite{bendsoe2013topology}, see also \cite{sigmund2013topology}.
We present various concrete applications of the general methods exposed in \cref{sec.absframework} in this setting. 
The first \cref{sec.elasSIMP,sec.simpnum} introduce the physical context and the numerical framework, respectively. Then, \cref{sec.simpwdroload} deals with topology optimization
problems where the uncertain parameters are the loads applied to the structure,
and the ambiguity set $\calA$ is based on the Wasserstein distance to a
nominal law. The next \cref{sec.simpmdro} illustrates the use of ambiguity sets based on the moments of the uncertain probability law. After that in \cref{subsec.matdro}, we look at another type of uncertainty, concerning the material coefficients. Finally, in \cref{sec.TOcvardro}, we handle reliability-based optimization problem using the notion of conditional value at risk.

\subsection{Topology optimization of linearly elastic structures}\label{sec.elasSIMP}

\noindent
Let $D$ be a fixed hold-all domain in $\R^d$ ($d=2,3$); the design variable under scrutiny is a density function
$h \in L^{\infty}(D, [0,1])$, representing a relaxed version of the characteristic function of a subset of $D$. Specifically,
\begin{itemize}
\item The regions where $h(x) = 1$ are made of an isotropic elastic material with Hooke's law $A$ given by:
\begin{equation}\label{eq.Hookelaw}
	\text{For any symmetric }d \times d \text{ matrix }e, \quad Ae = 2\mu
	e + \lambda \tr(e) \I.
\end{equation}
Here, $\lambda$ and $\mu$ are the Lam\'e parameters of the material; these are related to the more physical Young's modulus $E$ and Poisson ratio $\nu$ via the following relations:
\begin{equation}\label{eq.mulambdaE}
\mu = \frac{E}{2(1+\nu)}, \text{ and }\lambda = \left\{
	\begin{array}{cl}
		\frac{E\nu}{(1+\nu)(1-2\nu)} & \text{if } d =3,                        \\[0.5em]
		\frac{E}{1-\nu^{2}}          & \text{if } d =2 \text{ (plane stress).} \\
	\end{array}
	\right.
\end{equation}
\item The function $h$ equals $0$ at points $x \in D$ that are completely surrounded by a very soft material with Hooke's law $\eta A$, $\eta \ll 1$, mimicking void; 
\item The intermediate regions where $h(x) \in (0,1)$ are composed of a mixture of material and void: their material properties interpolate between those the bulk and void phases.
\end{itemize}
These features are incorporated into an inhomogeneous Hooke's tensor $A(h)$ associated to the density function $h$. 
According to the Solid Isotropic Material Penalization approach (SIMP), the latter is defined by:
\begin{equation}\label{eq.Ah}
\forall x \in D, \quad A(h)(x) = \Big(\eta + (1 - \eta) h(x)^{p}\Big) A,
\end{equation}
where the penalization exponent $p$ is typically set to $3$.

The considered structures are clamped on a region $\Gamma_{D}$ of the boundary $\partial D$, and external loads $g : \Gamma_{N} \to \R^{d}$ are applied on a disjoint region $\Gamma_{N}\subset \partial D$.
In this situation, the displacement of the structure is the unique solution $u \in H^{1}(D)^{d}$  to the system of linear elasticity:
\begin{equation}
	\label{eq.elassimp}\left\{
	\begin{array}{cl}
		-\text{div}(A(h)e(u)) = 0 & \text{in } D,                                                                 \\
		u = 0                     & \text{on } \Gamma_D,                                                          \\
		A(h)e(u)n = g             & \text{on } \Gamma_N,                                                          \\
		A(h)e(u)n = 0             & \text{on } \partial D \setminus \overline{\Gamma_N} \cup \overline{\Gamma_D}.
	\end{array}
	\right.
\end{equation}
Here, we have denoted by $e(u) := \frac{1}{2}(\nabla u + \nabla u^{T})$
the strain tensor associated to a displacement field $u : D \to \R$. 
We have also omitted body forces (gravity) for simplicity. 

In the sequel, we shall endow the elastic displacement $u$ and the associated physical quantities of interest with various indices or arguments to highlight their dependence on the relevant parameters of this system, notably, the design $h$, the loads $g$ and the material parameters $E$, $\nu$.

\begin{remark}
In practice, the definition of the Hooke's law $A(h)$ in \cref{eq.Ah} involves various filters, imparting effects such as smoothness, among others, which are not mentioned for brevity. Note however that the parameters of these filters are updated in the course of the optimization, which partly explains the large jumps that are sometimes visible in the convergence histories of our numerical results.
\end{remark}

A general topology optimization problem in the above context reads:
\begin{equation}\label{eq.optpbnomcplyvol}
\min_{h\in\Uad}J(h) \text{ s.t. }G(h) = 0,
\end{equation}
where the objective and constraint functions $J(h)$, $G(h)$ may depend on $h$ via the elastic displacement $u_h$. Examples of such functionals are:
\begin{itemize}
\item The compliance $C(h)$ of the structure, which equivalently accounts for the elastic energy stored in the structure, or the work done by the external loads:
\begin{equation}\label{eq.defcplyTO}
C(h) = \int_D A(h)e(u_h) : e(u_h)\:\d x = \int_{\Gamma_N} g \cdot u_h \: \d s ; 
\end{equation}
\item The difference between the volume $\Vol(h) = \int_D h \:\d x$ of $h$ and a target value $V_T$.
\end{itemize}

The differentiation of functions depending on the design via boundary-value problems is a difficult, albeit classical issue in optimal design, see e.g. \cite{bendsoe2013topology,cea1986conception,plessix2006review} about the celebrated adjoint method. 
For completeness, we recall the following derivative formulas, see e.g. \cite{allaire2020shape,bendsoe2013topology}.
\begin{proposition}
The compliance and volume functionals are differentiable at an arbitrary density function $h \in \Uad$, and their derivatives equal:
$$ C^\prime(h)(\widehat{h}) = -\int_D A^\prime(h)( \widehat{h})  e(u_h) : e(u_h) \:\d x, \text{ and } \Vol^\prime(h)(\widehat{h}) = \int_D \widehat{h} \:\d x. $$
\end{proposition}
 
\subsection{Numerical Strategy}\label{sec.simpnum}

\noindent Our numerical framework for solving problems of the form \cref{eq.optpbnomcplyvol} is based on the open-source implementation \cite{dapogny2022tuto}, and more exactly on a version of the latter dedicted to topology optimization \cite{totuto}. 
In a nutshell, it is a steepest-descent strategy based on the derivatives of the objective and constraint functionals $J(h)$ and $G(h)$. The physical boundary-value problems characterizing the elastic displacement $u_h$ in their expressions are solved by the Finite Element Method, using the open-source library  \texttt{FreeFem} \cite{hecht2012new}. 
A descent direction is then inferred from these data thanks to the Null-Space algorithm \cite{feppon2024tutoNS, feppon2020null}, a constrained optimization algorithm in infinite dimension which is suitable for shape and topology optimization applications, see also \cite{nullspace} for a \texttt{python} open-source implementation. These independent operations are orchestrated by a \texttt{python} script, and all our computations are performed on a standard \texttt{macbook pro} laptop with 16Gb of RAM.  

The treatment of distributionally robust problems such as \cref{eq.wdrogen} raises the need to evaluate the probabilistic integrals over the parameter set $\Xi$ involved in the expressions of the distributionally robust functionals and their derivatives. To achieve this task, we rely on Monte Carlo sampling: every integral of the form $\int_\Xi f(\xi) \:\d \P(\xi)$ is approximated as:
$$\int_\Xi f(\xi) \:\d \P(\xi) = \frac{1}{N}\sum\limits_{i=1}^N f(\xi^i),$$
where the $\xi^i$ are $N$ independent and identically distributed samples drawn from the probability distribution $\P$. In practice, in order to keep the computational effort moderate and since our illustrations do not target industrial complexity, we use a maximum number of $N=10$ vectors; these evaluations can be conducted in parallel thanks to the \texttt{numpy} vectorization solutions.
On a different note, according to \cite{blanchard2021accurately}, we use the so-called ``log-sum-exp trick'' to overcome the problems posed by the smallness of the integrand, and thus stabilize computations.
Admittedly, this stochastic gradient descent strategy for the resolution of a problem such as \cref{eq.wdrogen} is overly simplistic, and much more advanced techniques could be used to improve its efficiency, starting from improvements of the sampling strategy by variance reduction techniques, about which we refer to the review article \cite{mohamed2020monte}, or more advanced versions of the stochastic gradient method, such as that introduced in \cite{grieshammer2024continuous,grieshammer2024continuous2}. Although crucial in real-life applications, we believe that these topics are beyond the scope of the present study.

\subsection{Wasserstein distributionally robust topology optimization under load uncertainty}\label{sec.simpwdroload}

\noindent Our first set of numerical examples arises in situations where the loads applied to the structure are uncertain, i.e., $\xi = g$, and the considered cost function  is the compliance $C(h,\xi)$ defined in \cref{eq.defcplyTO}.

\subsubsection{Loads uncertainty with a nominal distribution defined from a single measurement}\label{sec.simpwdro1load}

\begin{figure}[H]
	\begin{tabular}{cc}
	\begin{minipage}{0.53\textwidth}
	 \begin{overpic}[width=1.0\textwidth]{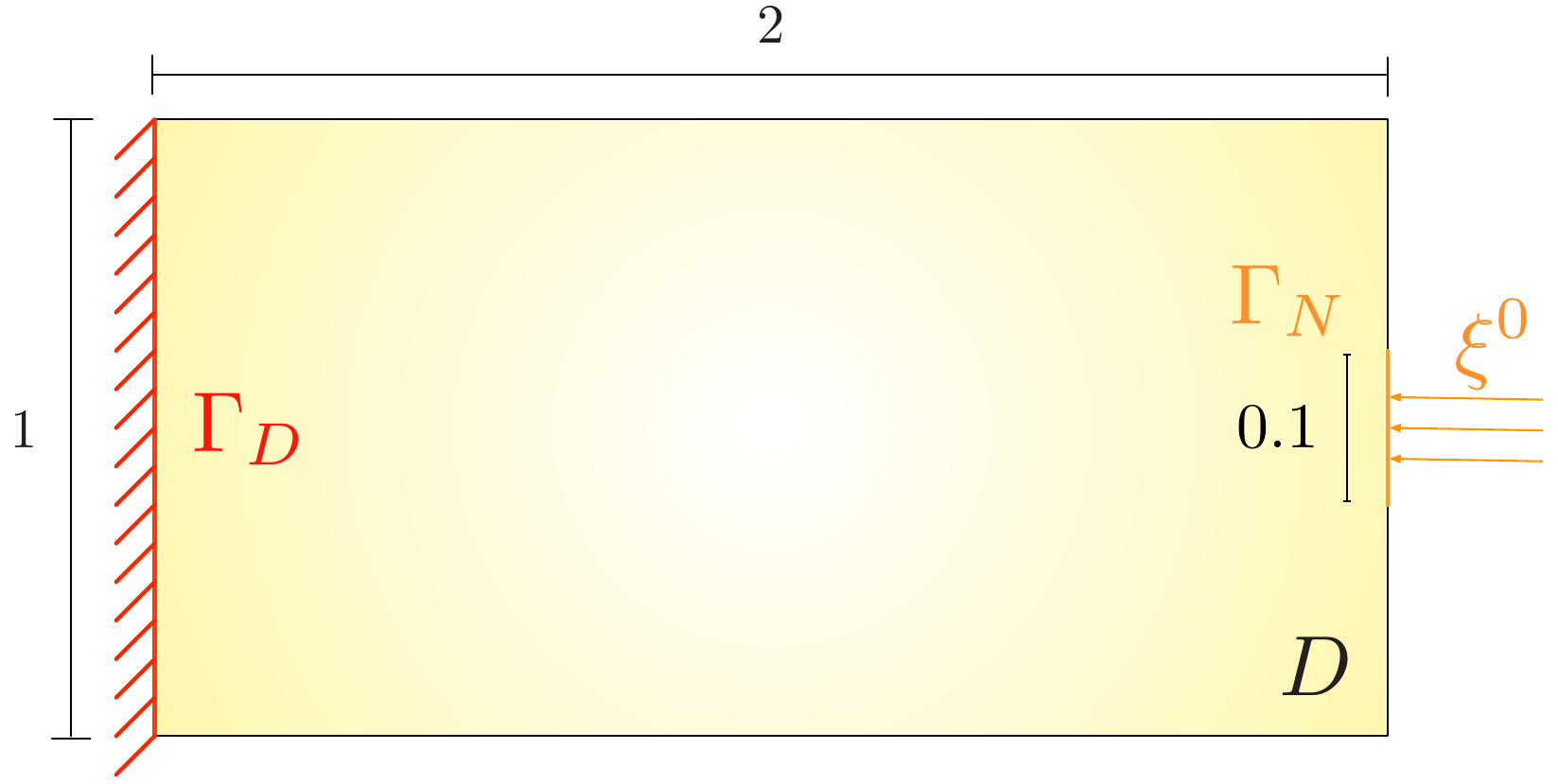} 
	 \put(2,5){\fcolorbox{black}{white}{a}}
	 \end{overpic}
	 \end{minipage} &
	\begin{minipage}{0.45\textwidth}
	\begin{overpic}[width=1.0\textwidth]{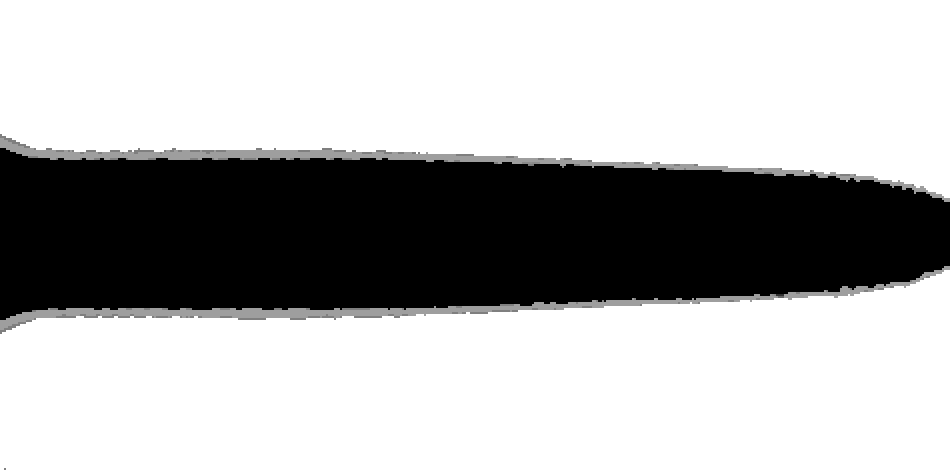}
	 \put(2,5){\fcolorbox{black}{white}{b}}
	 \end{overpic}
	 \end{minipage} \\
	  \begin{minipage}{0.45\textwidth}
	  \begin{overpic}[width=1.0\textwidth]{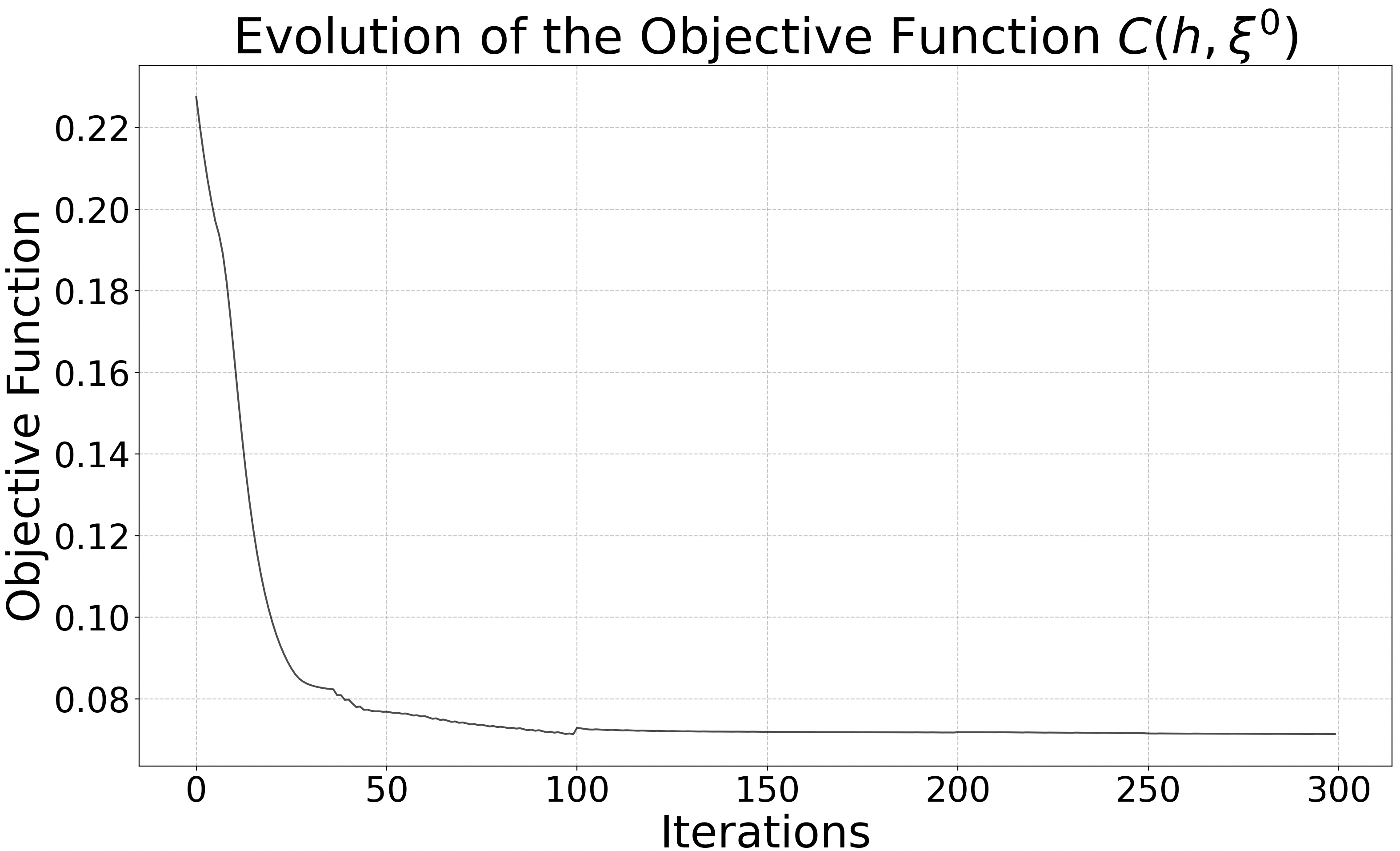} 
	  \put(2,5){\fcolorbox{black}{white}{c}}
	  \end{overpic}
	  \end{minipage}&
        \begin{minipage}{0.45\textwidth}
        \begin{overpic}[width=1.0\textwidth]
        {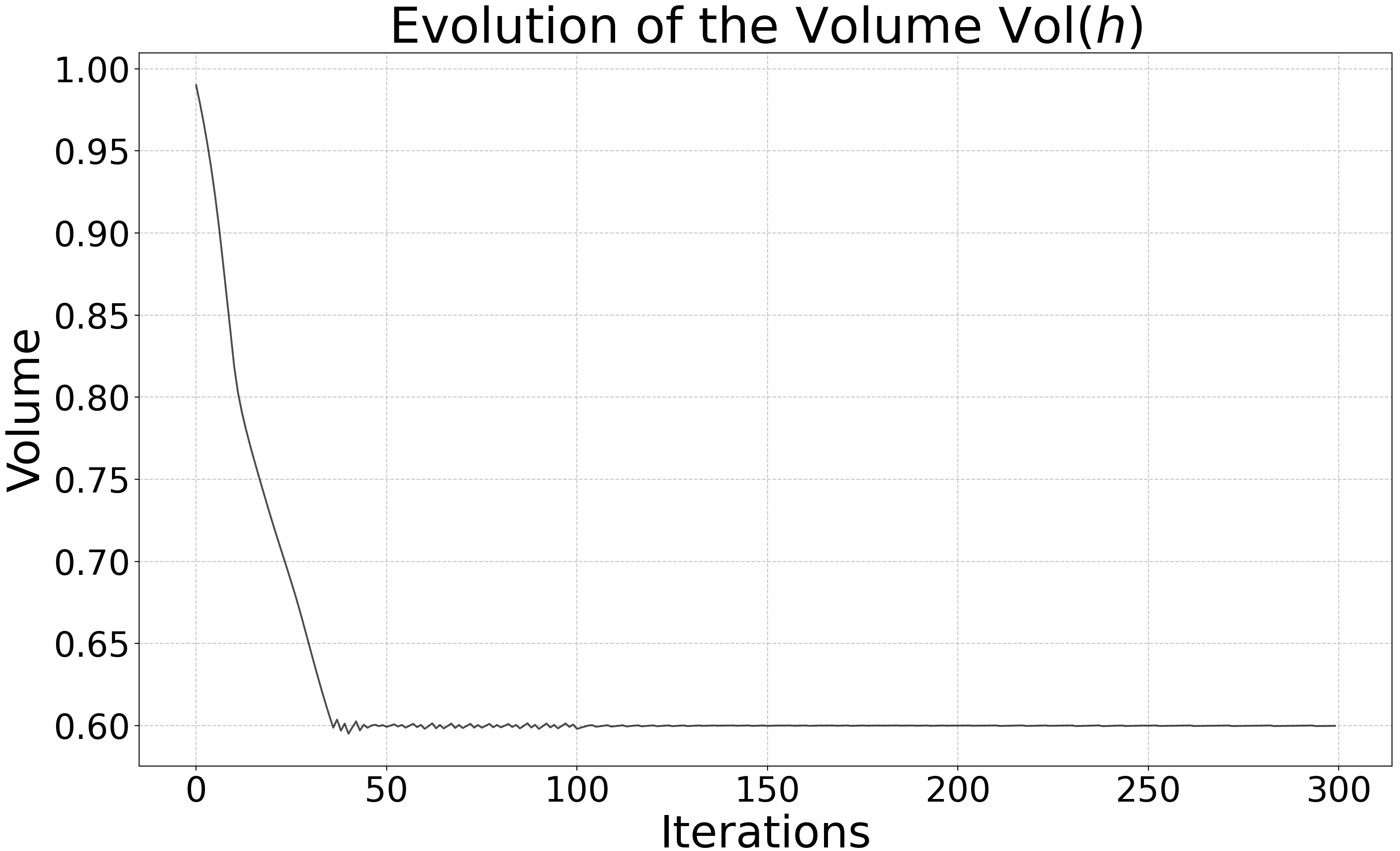}\put(2,5){\fcolorbox{black}{white}{d}}
        \end{overpic}            
        \end{minipage}
	\end{tabular}
	\caption{\it (a) Setting of the cantilever beam optimization problem considered in \cref{sec.simpwdro1load}; (b) Optimized structure for the ideal problem \cref{eq.optpbcantiideal}; (c) Convergence history of the objective $C(h,\xi^0)$; (d) Convergence history of the volume $\Vol(h)$.}
	\label{fig.cantinom}
\end{figure}

\noindent The structure under scrutiny is a 2d cantilever beam; it is contained in a fixed computational domain $D$ with size $2 \times 1$, which is in practice equipped with a triangular mesh made of about $9,\!000$ vertices ($\approx 17,\!000$ triangles). The structure is clamped on the left-hand side $\Gamma_D$ of $\partial D$, and surface loads are applied on a small region $\Gamma_N$ near the center of its right-hand side. The corresponding setting is shown in \cref{fig.cantinom} (a).

In an ideal situation, the structure is subjected to a horizontal load $\xi^0=(-1,0)$; we then consider the deterministic optimization problem:
\begin{equation}\label{eq.optpbcantiideal}
\min \limits_{h \in \Uad} C(h,\xi^0) \text{ s.t. } \Vol(h) = V_T,
\end{equation}
where the target volume $V_T$ equals $0.6$. The outcome of its resolution by the algorithm outlined in \cref{sec.simpnum} is depicted on \cref{fig.cantinom} (b,c,d). 
The computation requires 300 iterations, and it lasts about 17 mn.

We now assume that the load $\xi$ applied on $\Gamma_N$ is unknown: it belongs to a sufficiently large ball $\Xi$ in $\R^2$, and its probability law is itself unknown: 
the latter is estimated solely from the assumed ideal value $\xi^0$. Hence, introducing the nominal probability distribution ${\P}=\delta_{\xi^0}$ made from the single sample $\xi^0$, we consider the following Wasserstein distributionally robust problem:
\begin{equation}\label{eq.wdrocantiprimal}
\min \limits_{h \in \Uad}  \sup\limits_{\Q \in \calA_{\text{W}}} \int_\Xi C(h,\xi) \:\d\Q(\xi) \text{ s.t. } \Vol(h) = V_T. 
\end{equation}
Let us recall from the discussion in \cref{sec.wdroformula} that this problem has the following tractable equivalent form:
\begin{multline}\label{eq.wdrocanti}
	\min\limits_{h \in \Uad, \atop \lambda \geq 0} \mathcal{D}_{\text{W}}(h,\lambda)
	\quad \text{subject to} \quad \text{Vol}(h) = V_T,\\  \text{where} \quad 
	\mathcal{D}_{\text{W}}(h,\lambda) := \lambda m + \lambda \varepsilon \log\left(
	\int_{\Xi} \exp\left(\frac{C(h,\zeta) - \lambda c(\xi^0,\zeta)}{\lambda \varepsilon}\right)
	\,\mathrm{d}\nu_{\xi}(\zeta) \right), 
\end{multline}
and the probability measure $ \nu_{\xi}(\zeta) = \mathds{1}_{\Xi}(\zeta) \alpha_{\xi} e^{-\frac{|\xi-\zeta|^2}{2\sigma^2}} \d\zeta$ is that featured in the definition \cref{eq.refcouplingW} of the reference coupling $\pi_0$. We conduct several numerical investigations in this setting.

\bigskip\noindent\textit{Influence of the variance $\sigma^2$ and the entropy penalization parameter $\e$.} We solve \cref{eq.wdrocanti} for several values of the parameters $\sigma^2$ and $\e$, for a fixed value $m=0$ of the maximum bound on the radius of the Wasserstein ball featured in the definition of $\mathcal{A}_{\text{W}}$. 
Let us point out that the outcome of this analysis is already not trivial since, as noted in \cref{rem.Wenodist}, the quantity $W_{\varepsilon}(\cdot,\cdot)$ is not a true distance, so that the version of \cref{eq.wdrocanti} associated with $m=0$ may not coincide with the nominal problem \cref{eq.optpbcantiideal} unless $\e$ tends to $0$.

\begin{figure}[ht]
	\begin{tabular}{ccc}
		\begin{minipage}{0.3\textwidth}\begin{overpic}[width=1.0\textwidth]{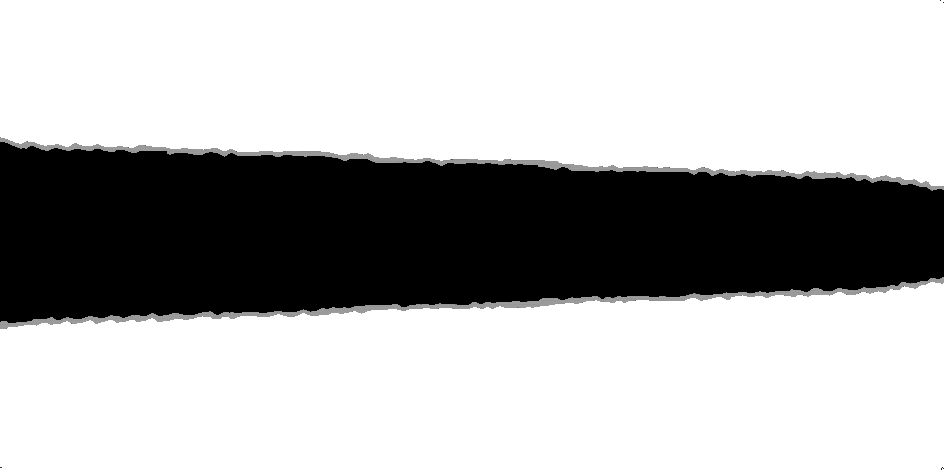} \put(2,5){\fcolorbox{black}{white}{$\sigma^2= 2e{-3}$}}\end{overpic}\end{minipage} & \begin{minipage}{0.3\textwidth}\begin{overpic}[width=1.0\textwidth]{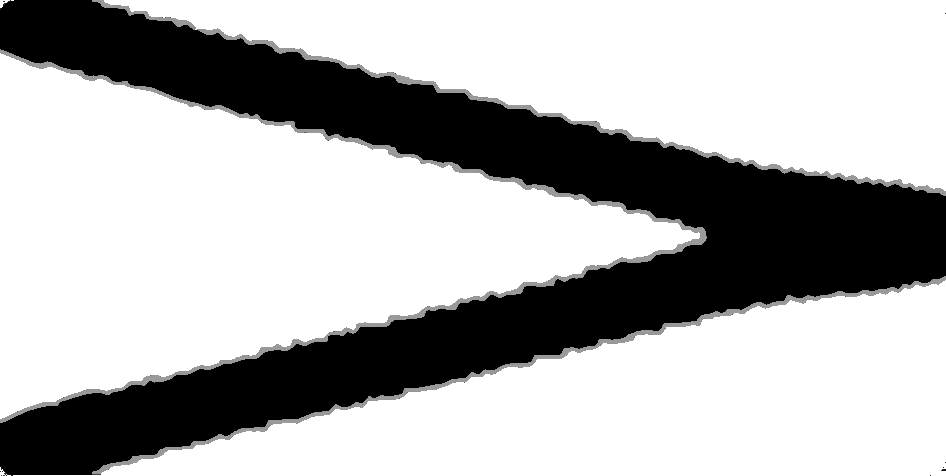} \put(2,5){\fcolorbox{black}{white}{$\sigma^2 = 2e{-2}$}}\end{overpic}\end{minipage} & \begin{minipage}{0.3\textwidth}\begin{overpic}[width=1.0\textwidth]{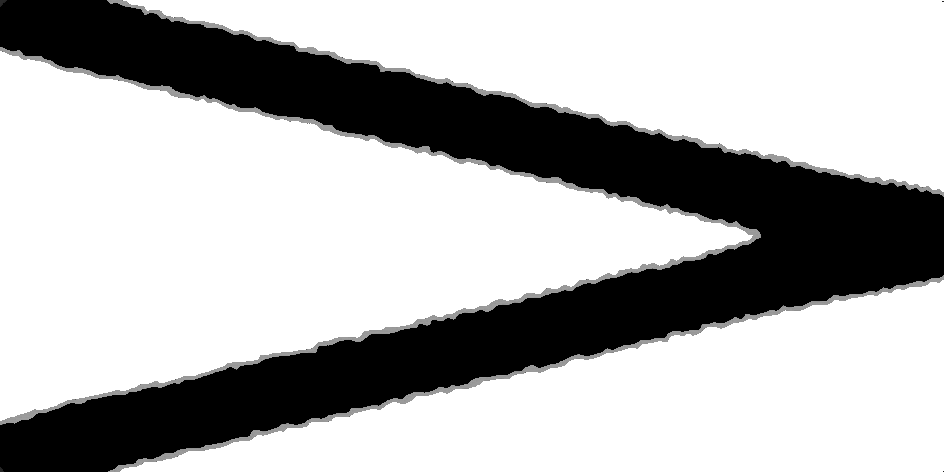} \put(2,5){\fcolorbox{black}{white}{$\sigma^2= 1e{-1}$}}\end{overpic}\end{minipage} \\
	\end{tabular}
    \par
    \medskip
    \begin{tabular}{ccc}
		\begin{minipage}{0.3\textwidth}\begin{overpic}[width=1.0\textwidth]{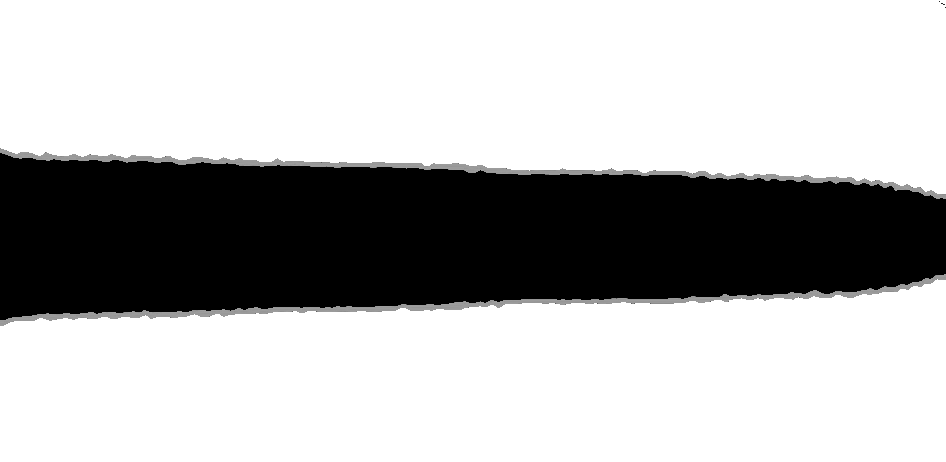} \put(2,5){\fcolorbox{black}{white}{$\sigma^2= 2e{-3}$}}\end{overpic}\end{minipage} & \begin{minipage}{0.3\textwidth}\begin{overpic}[width=1.0\textwidth]{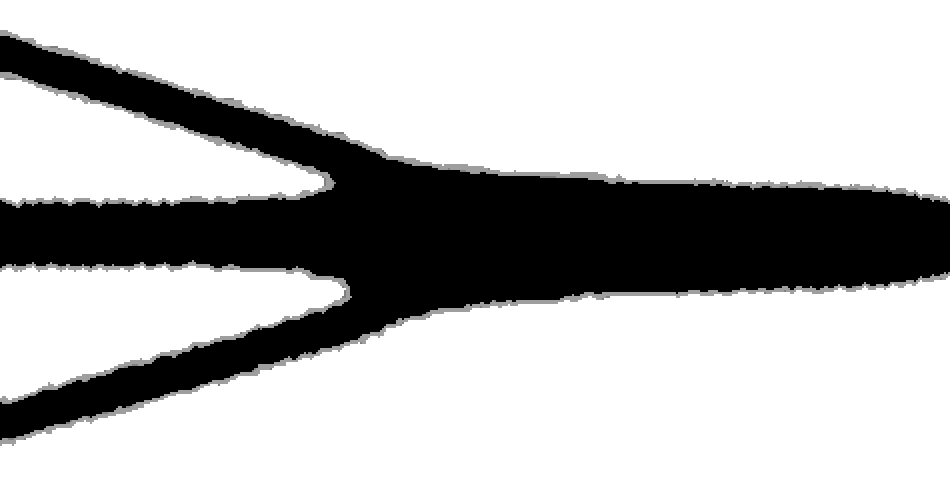} \put(2,5){\fcolorbox{black}{white}{$\sigma^2 = 2e{-2}$}}\end{overpic}\end{minipage} & \begin{minipage}{0.3\textwidth}\begin{overpic}[width=1.0\textwidth]{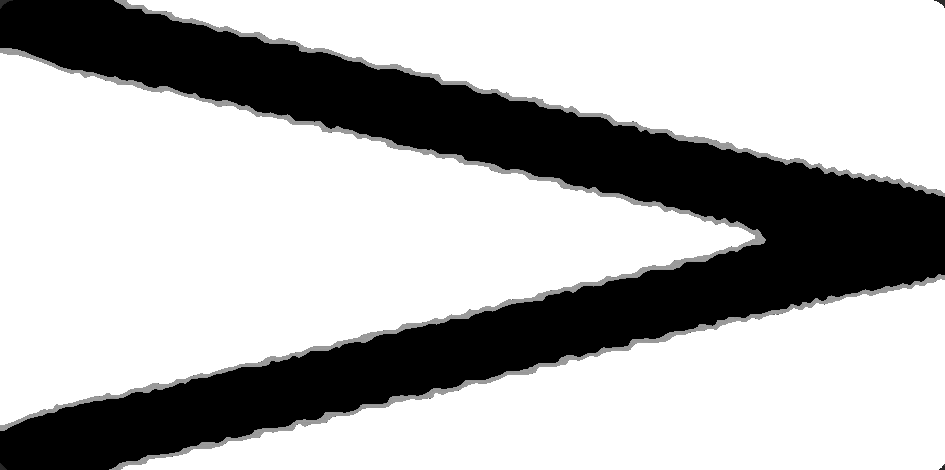} \put(2,5){\fcolorbox{black}{white}{$\sigma^2= 1e{-1}$}}\end{overpic}\end{minipage} \\
	\end{tabular}
	\caption{\it Optimized designs of the cantilever beam of \cref{sec.simpwdro1load} obtained by solving the distributionally robust problem \cref{eq.wdrocanti} for $m=0$ and different values of $\sigma^2$ when $m=0$; (upper row) $\varepsilon=1e{-2}$; (lower row) $\varepsilon=1e{-4}$.}
    \label{fig.cantiwdroVar}
\end{figure}

The results of the computation are presented on \cref{fig.cantiwdroVar}; each computation involves 600 iterations of our numerical algorithm, for a CPU time of about 2h 45 mn. As expected, when $\e$ and $\sigma^2$ are both ``very small'', the optimized design is similar to that obtained in the solution of the ideal problem \cref{eq.optpbcantiideal}, displayed on \cref{fig.cantinom} (b). However, as $\e$ or $\sigma^2$ increases, it tends to develop diagonal bars, that are more typical of structures withstanding vertical loads \cite{bendsoe2013topology}. These observations are in line with the intuitive meanings of these parameters, discussed in \cref{sec.wdro}. The variance $\sigma^2$ can be interpreted as (the inverse of) the degree of confidence of the user in the nominal law $\P$: for ``small'' values of $\sigma^2$, the only laws $\Q$ giving rise to a relatively small value of the Wasserstein distance $W_\e(\P,\Q)$ in \cref{eq.RegWass} are those very close to $\P$, whereas laws $\Q$ corresponding to increasingly smeared versions of $\P$ are captured in the situation where $\sigma^2$ is ``not too small''. Likewise, an increase in the value of the entropic regularization parameter $\e$ causes a more and more significant ``blurring'' of the true Wasserstein distance.

\begin{figure}[ht]
	\begin{tabular}{ccc}
		\begin{minipage}{0.33\textwidth}\begin{overpic}[width=1.0\textwidth]{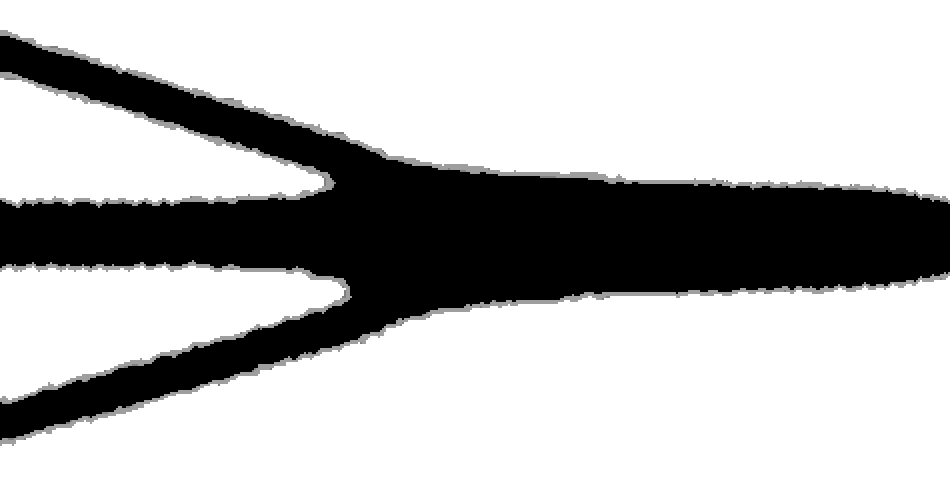} \put(2,5){\fcolorbox{black}{white}{$m=0$}}\end{overpic}\end{minipage} & \begin{minipage}{0.33\textwidth}\begin{overpic}[width=1.0\textwidth]{ 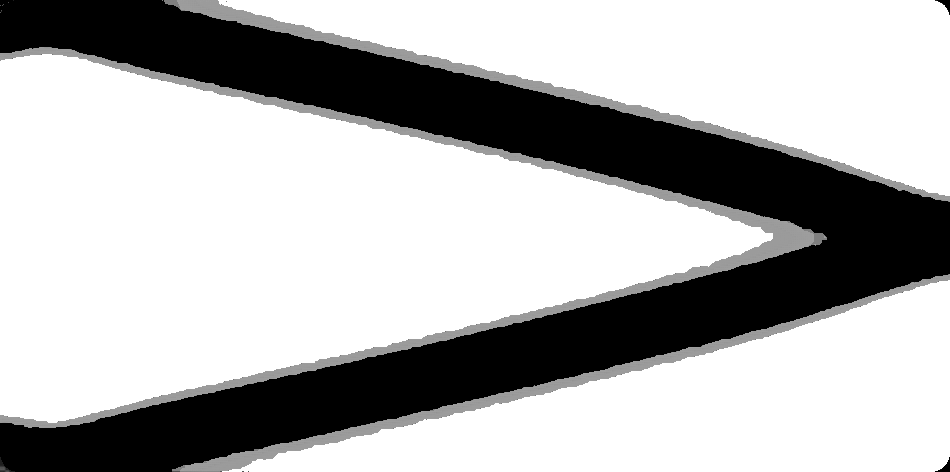 } \put(2,5){\fcolorbox{black}{white}{$m=1$}}\end{overpic}\end{minipage} & \begin{minipage}{0.33\textwidth}\begin{overpic}[width=1.0\textwidth]{ 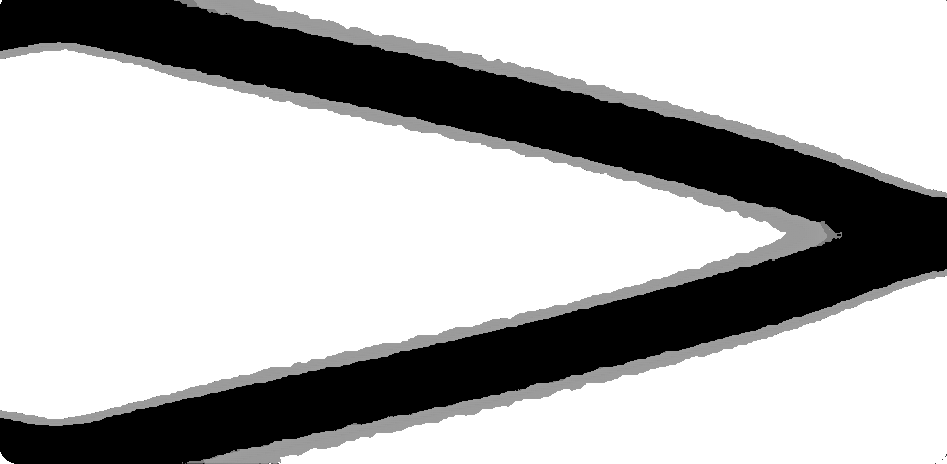 } \put(2,5){\fcolorbox{black}{white}{$m=5$}}\end{overpic}\end{minipage}
	\end{tabular}
    \par
    \medskip
    \begin{tabular}{ccc}
		\begin{minipage}{0.33\textwidth}\begin{overpic}[width=1.0\textwidth]{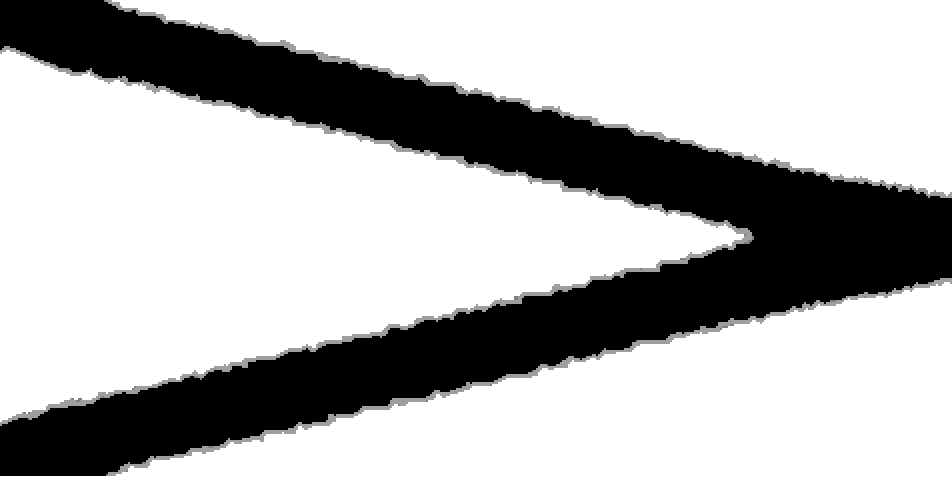} \put(2,5){\fcolorbox{black}{white}{$m=0$}}\end{overpic}\end{minipage} & \begin{minipage}{0.33\textwidth}\begin{overpic}[width=1.0\textwidth]{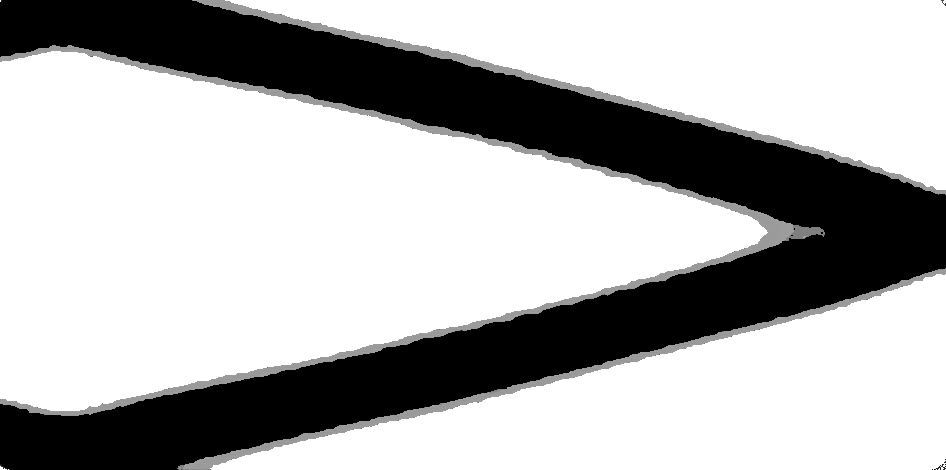} \put(2,5){\fcolorbox{black}{white}{$m=0.375$}}\end{overpic}\end{minipage} & \begin{minipage}{0.33\textwidth}\begin{overpic}[width=1.0\textwidth]{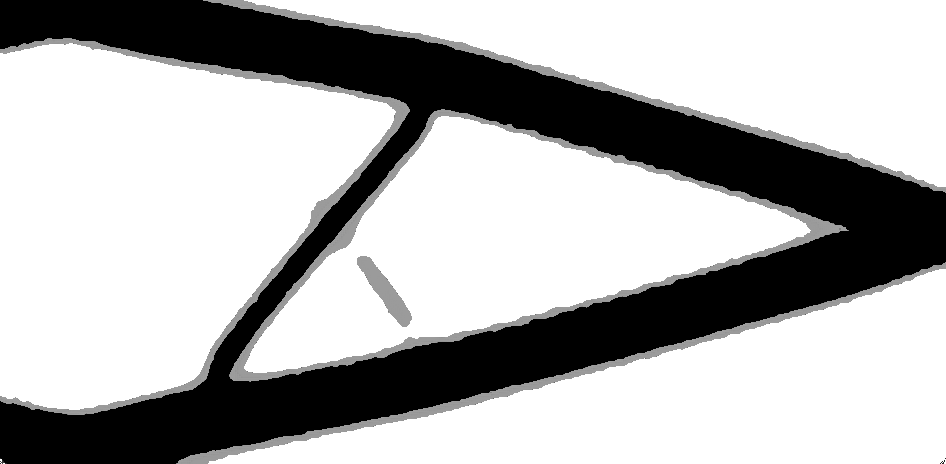} \put(2,5){\fcolorbox{black}{white}{$m=0.4$}}\end{overpic}\end{minipage}
	\end{tabular}
	\par
	\medskip
	\begin{tabular}{ccc}
		\begin{minipage}{0.33\textwidth}\begin{overpic}[width=1.0\textwidth]{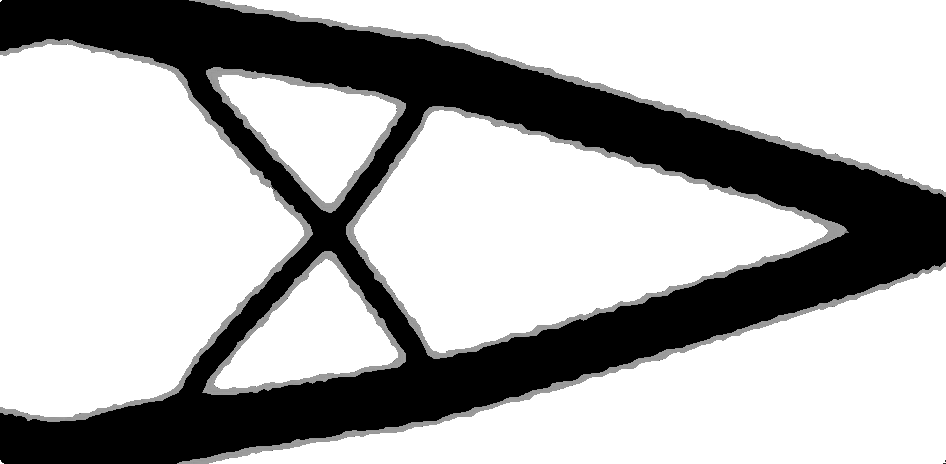} \put(2,5){\fcolorbox{black}{white}{$m=0.45$}}\end{overpic}\end{minipage} & \begin{minipage}{0.33\textwidth}\begin{overpic}[width=1.0\textwidth]{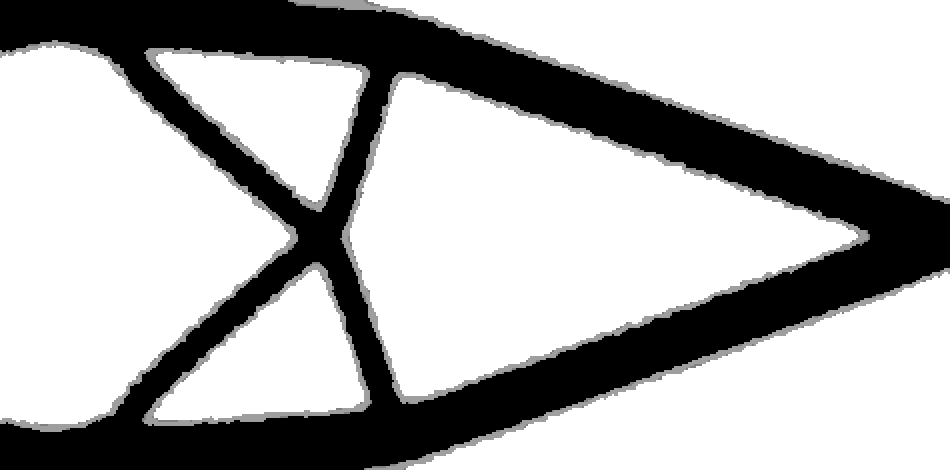} \put(2,5){\fcolorbox{black}{white}{$m=0.675$}}\end{overpic}\end{minipage} & \begin{minipage}{0.33\textwidth}\begin{overpic}[width=1.0\textwidth]{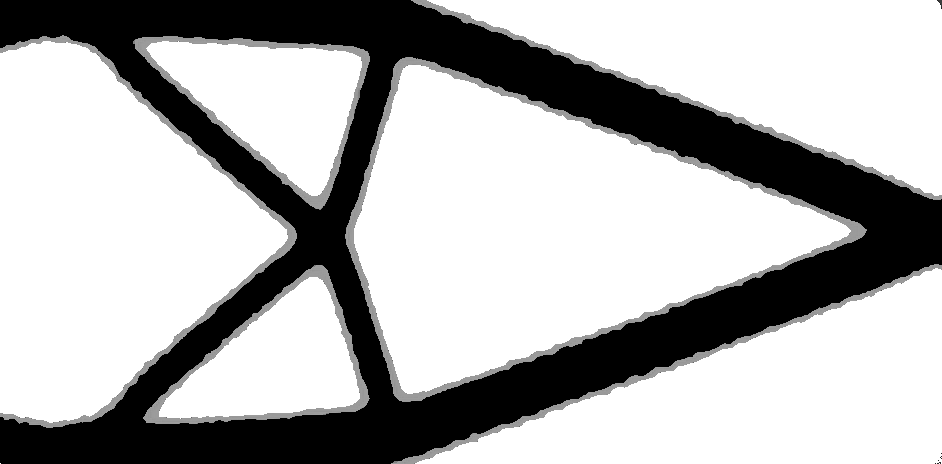} \put(2,5){\fcolorbox{black}{white}{$m=1$}}\end{overpic}\end{minipage} \\
	\end{tabular}
	\caption{\it Optimized designs of the cantilever beam of \cref{sec.simpwdro1load} obtained by solving the distributionally robust problem \cref{eq.wdrocanti} for several values of the Wasserstein radius $m$; (upper row) $\sigma^2 = 2e{-2}$, $\e =1e{-4}$;
   (middle and bottom rows) $\sigma^2 = 1e{-1}$, $\e=1e{-4}$.}
   \label{fig.cantiwdrom}
\end{figure}

\medskip\noindent\textit{Influence of the Wasserstein radius $m$.} We now solve \cref{eq.wdrocanti} for several values of the radius $m$ of the Wasserstein ball; the results are depicted on \cref{fig.cantiwdrom} for $\sigma^2 = 2e{-2}$, $\e =1e{-4}$, and then
$\sigma^2 = 1e{-1}$, $\e=1e{-4}$.

\begin{figure}[ht]
	\begin{tabular}{cc}
		\begin{minipage}{0.45\textwidth}
		\begin{overpic}[width=1.0\textwidth]{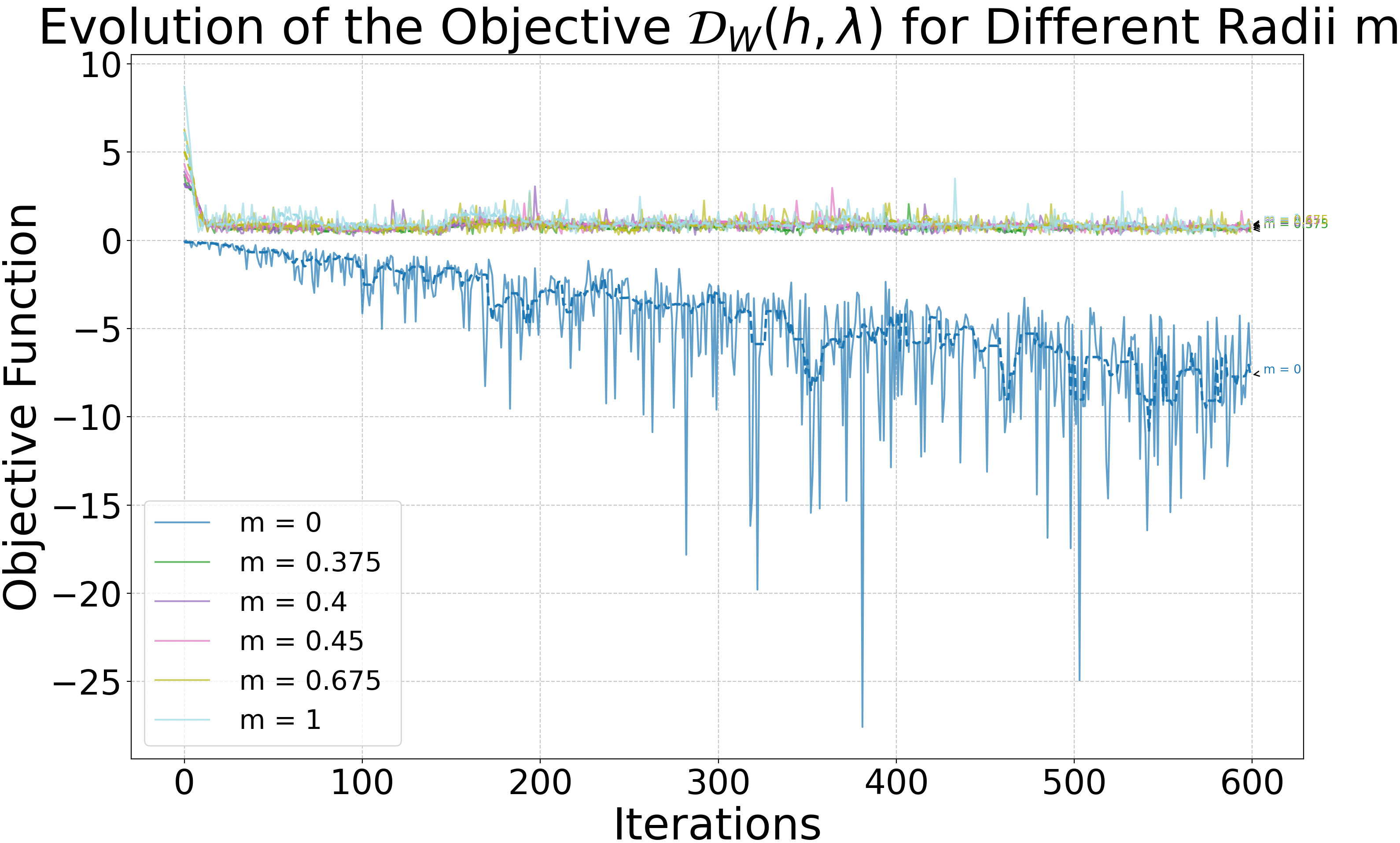} 
		\put(2,5){\fcolorbox{black}{white}{a}}
		\end{overpic} 
		\end{minipage} & 
		\begin{minipage}{0.45\textwidth}
		\begin{overpic}[width=1.0\textwidth]{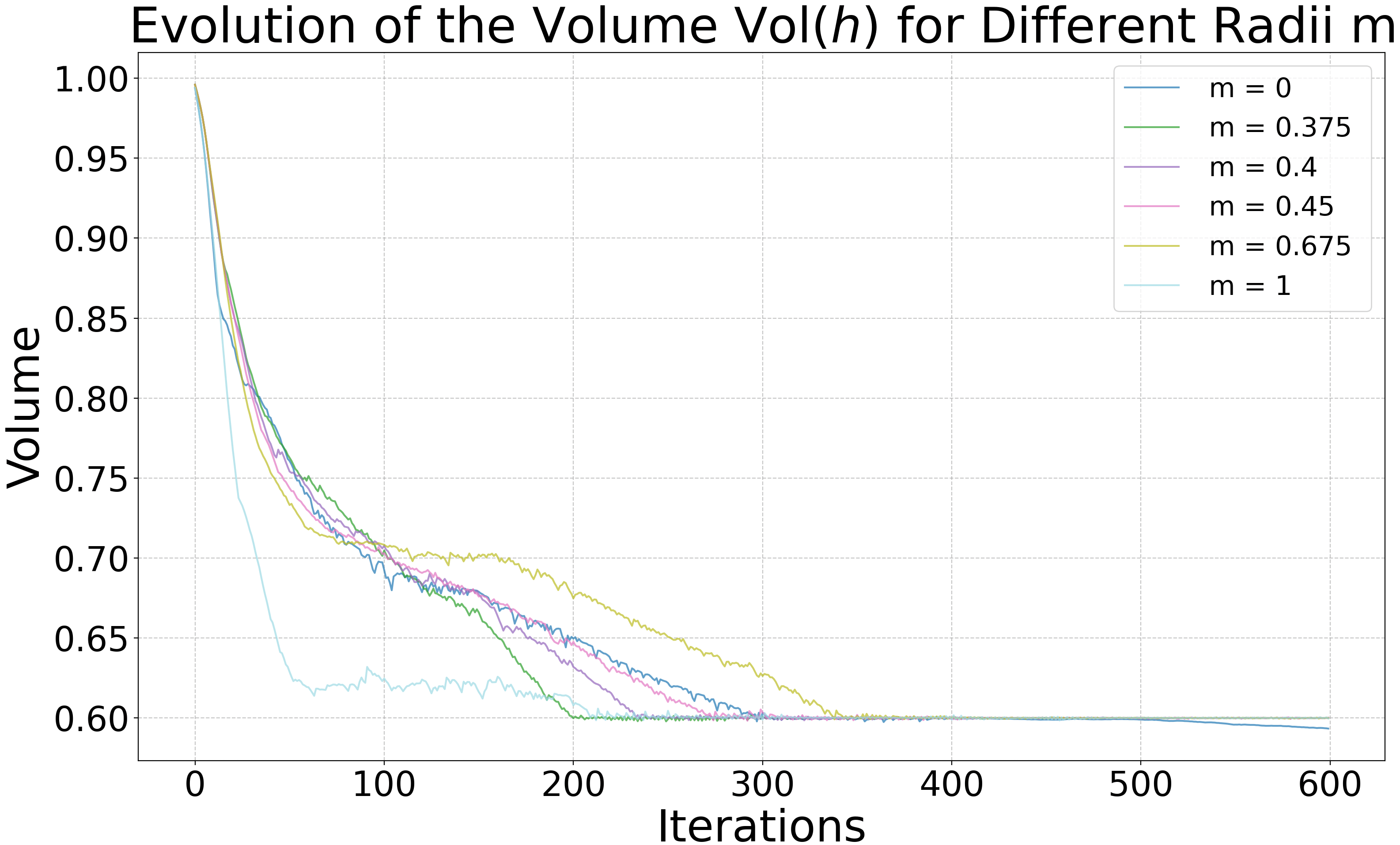} 
		\put(2,5){\fcolorbox{black}{white}{b}}
		\end{overpic} 
		\end{minipage} 	
		\end{tabular}
   \caption{\it Convergence histories in the distributionally robust optimization of the cantilever considered in \cref{sec.simpwdro1load} for the value $\sigma^2 = 1e{-1}$ and several values of $m$; (a) Objective function; (b) Volume.}
   \label{fig.cvcantidro}
\end{figure}

As the Wasserstein radius $m$ increases, the ``bad'' vertical perturbations of loads are anticipated, urging the formation of diagonal reinforcements that are reminiscent of the typical design of an optimized cantilever beam subjected to vertical loads, see again \cite{bendsoe2013topology}.
 The convergence histories of the objective function $\calD_{\text{W}}(h,\lambda)$ and of the volume $\Vol(h)$ are presented on \cref{fig.cvcantidro} in the case where $\sigma^2 = 1e{-1}$ and $\e = 1e{-4}$.
The optimization path is very noisy, but it eventually shows convergence to a local optimum of \cref{eq.wdrocanti}, as is typical of a stochastic descent strategy. 

\begin{figure}[ht]
    \centering
	\begin{tabular}{ccc}
		\begin{minipage}{0.33\textwidth}\begin{overpic}[width=1.0\textwidth]{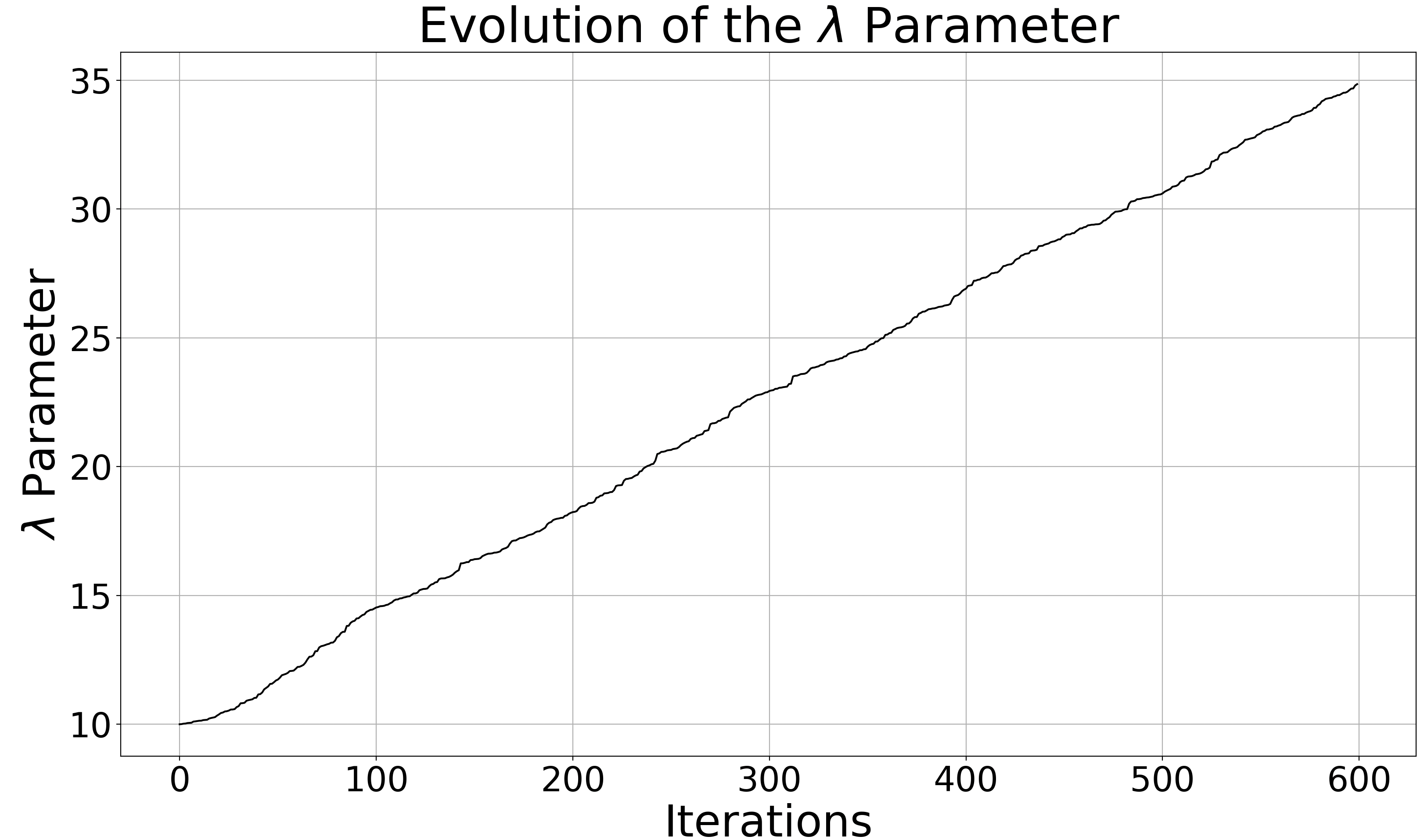} \put(2,5){\fcolorbox{black}{white}{$m=0$}}\end{overpic}\end{minipage} & \begin{minipage}{0.33\textwidth}\begin{overpic}[width=1.0\textwidth]{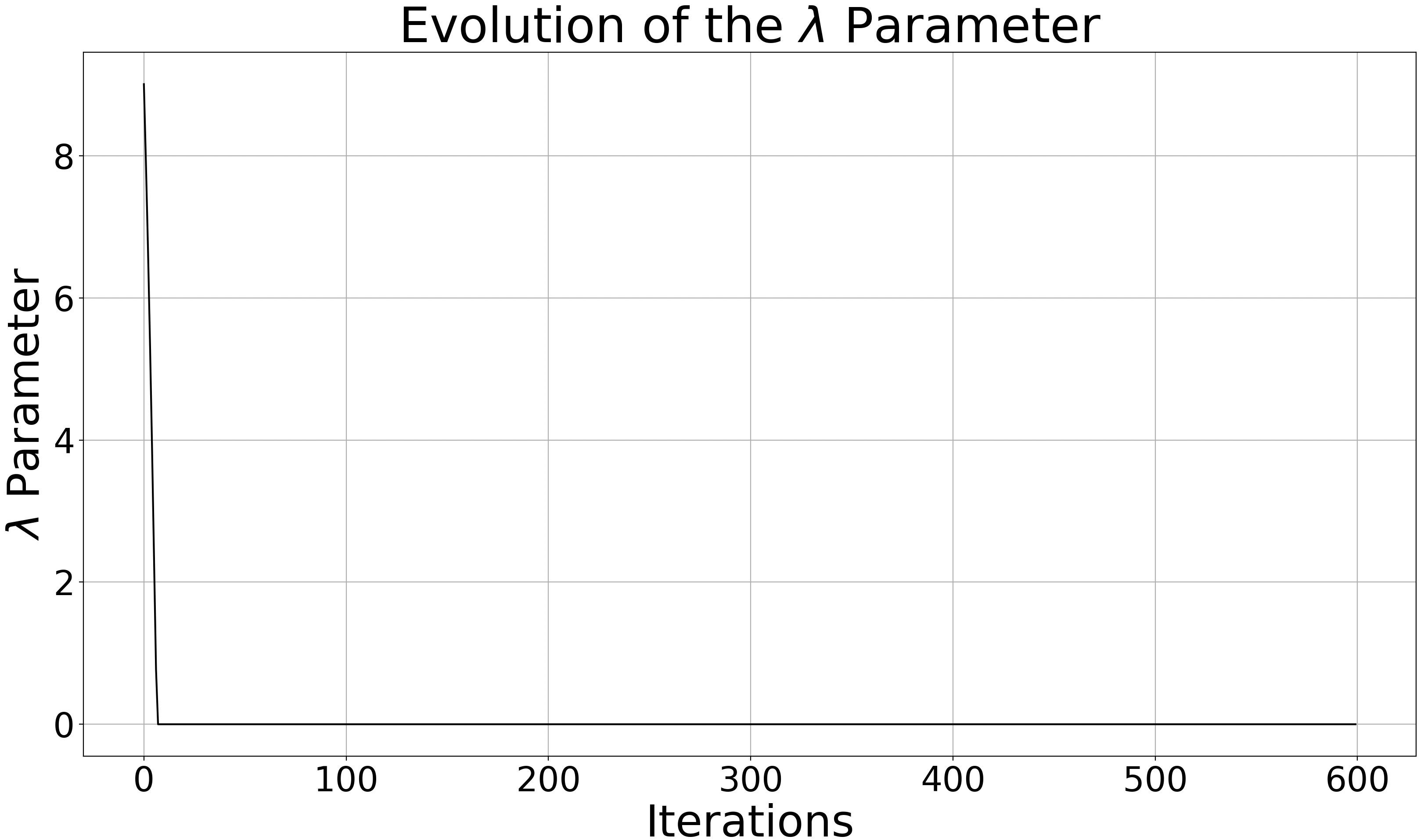} \put(2,5){\fcolorbox{black}{white}{$m=1$}}\end{overpic}\end{minipage} & \begin{minipage}{0.33\textwidth}\begin{overpic}[width=1.0\textwidth]{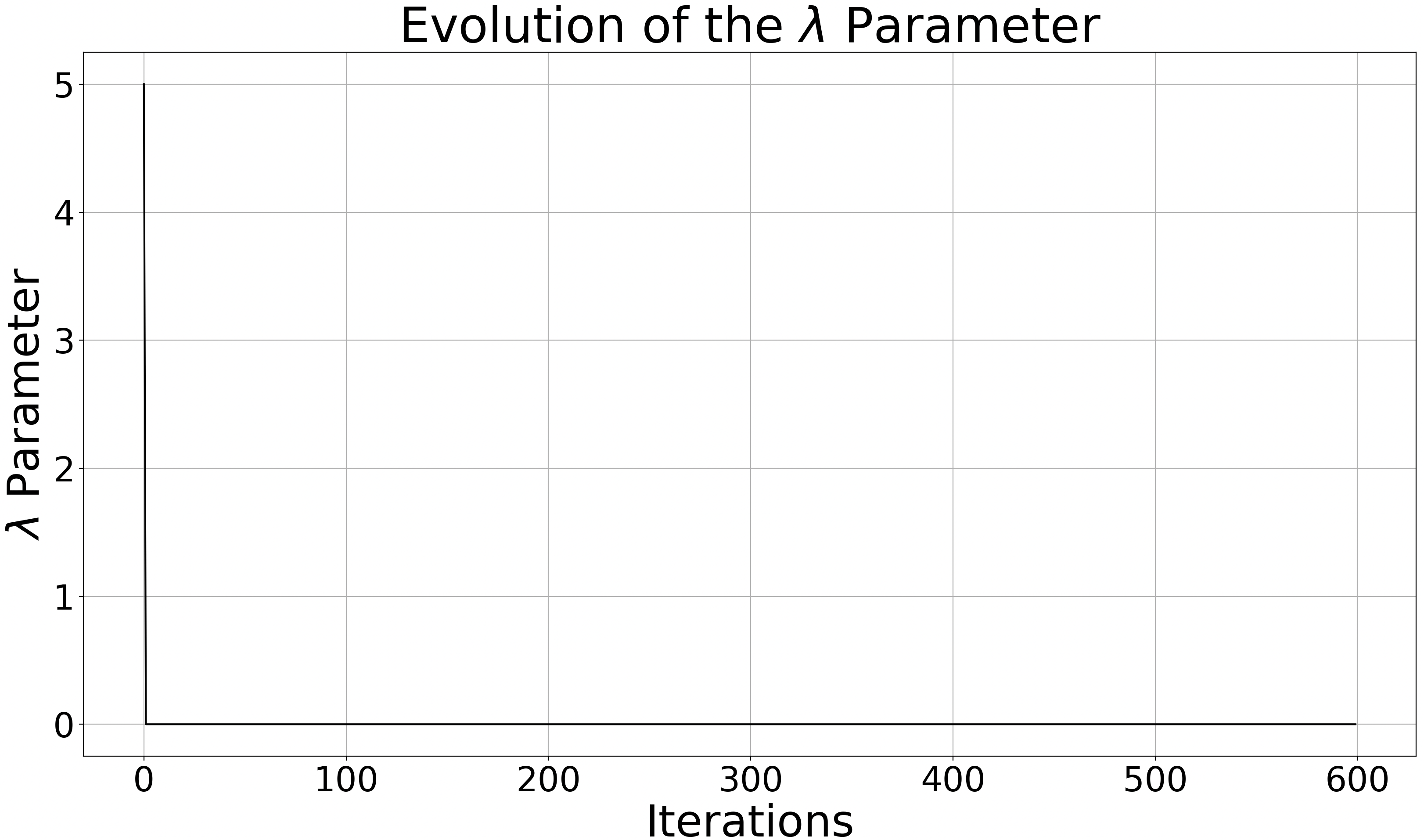} \put(2,5){\fcolorbox{black}{white}{$m=5$}}\end{overpic}\end{minipage}
	\end{tabular}
    \par
    \medskip
    \begin{tabular}{ccc}
		\begin{minipage}{0.33\textwidth}\begin{overpic}[width=1.0\textwidth]{ 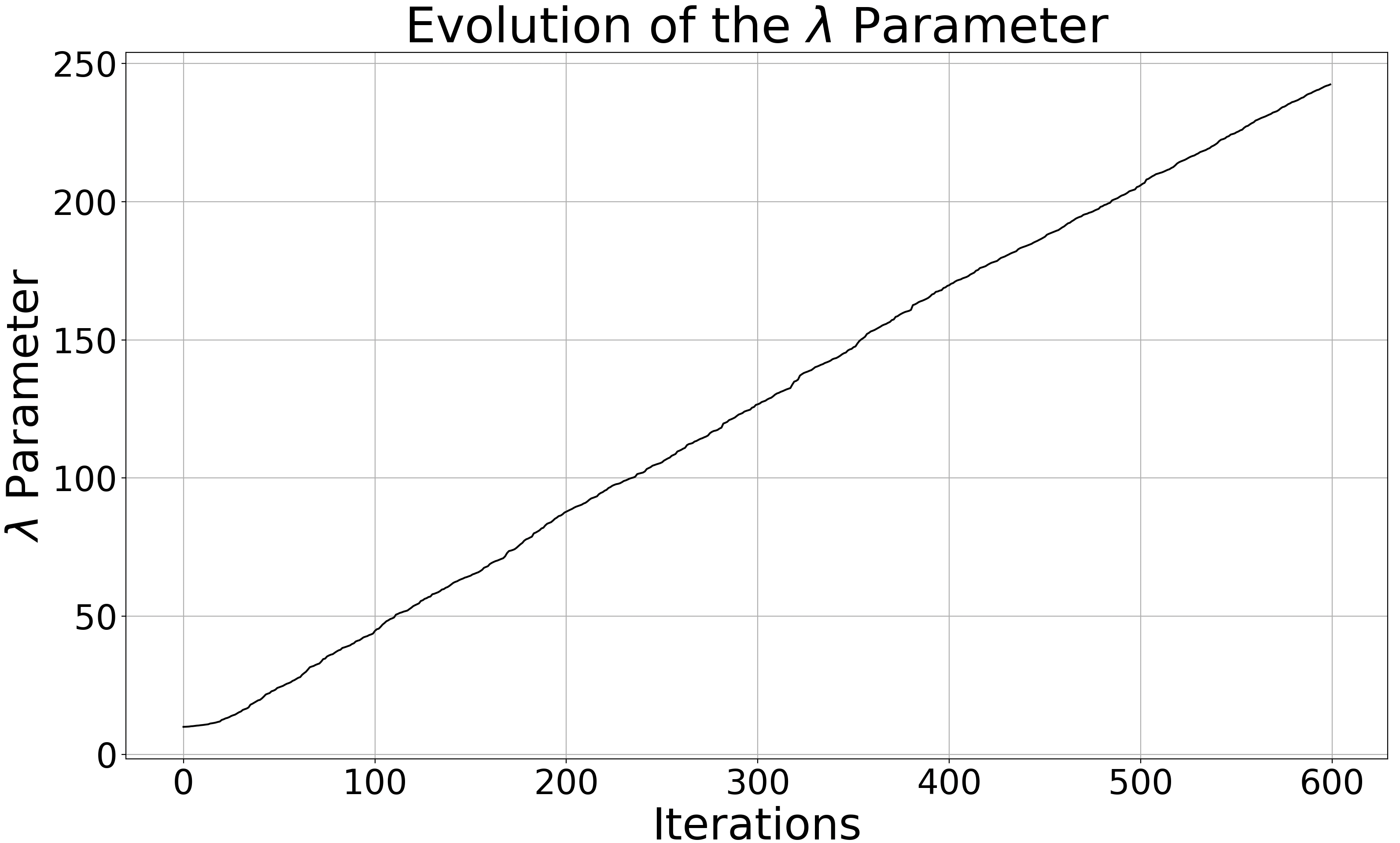 } \put(2,5){\fcolorbox{black}{white}{$m=0$}}\end{overpic}\end{minipage} & \begin{minipage}{0.33\textwidth}\begin{overpic}[width=1.0\textwidth]{ 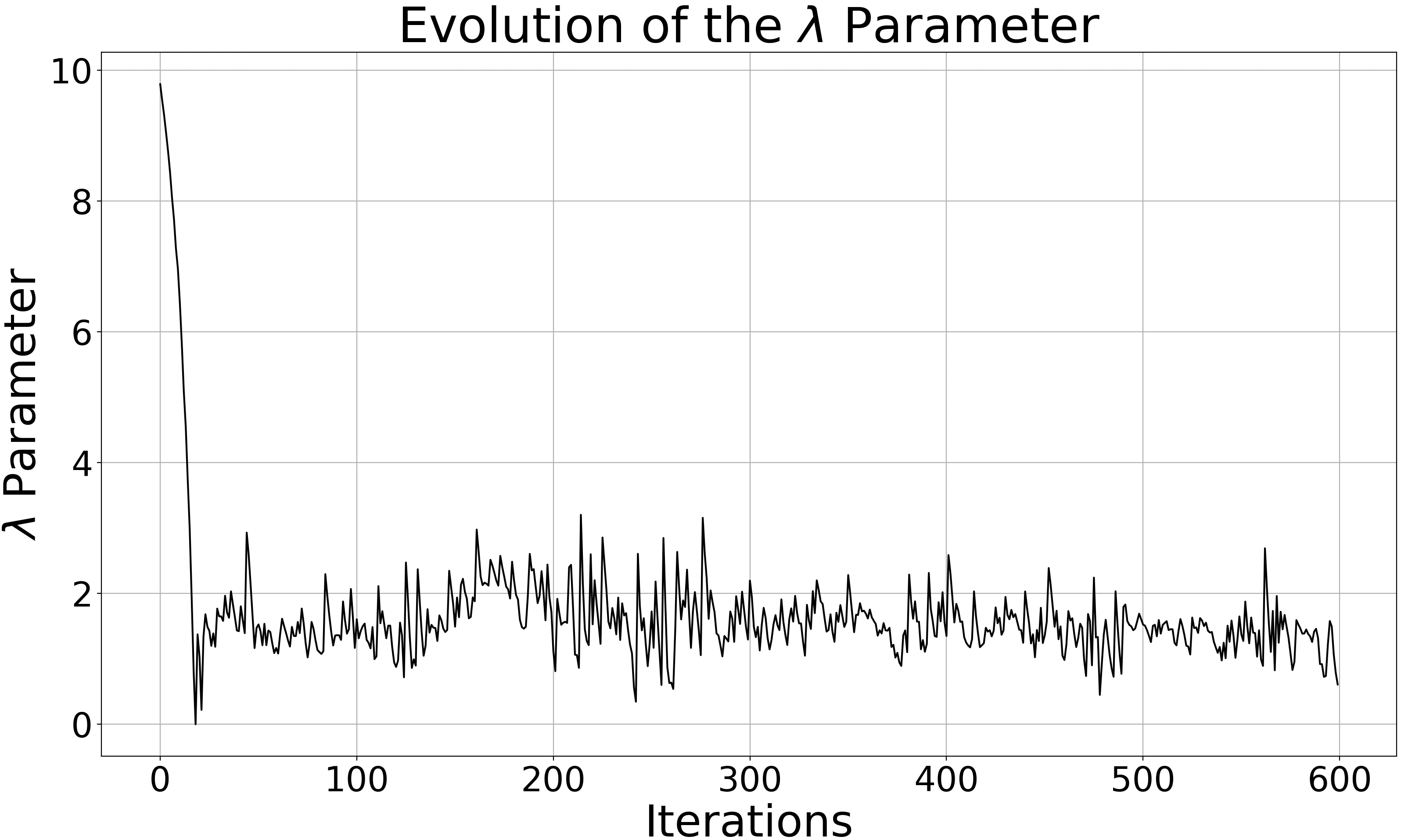 } \put(2,5){\fcolorbox{black}{white}{$m=0.25$}}\end{overpic}\end{minipage} & \begin{minipage}{0.33\textwidth}\begin{overpic}[width=1.0\textwidth]{ 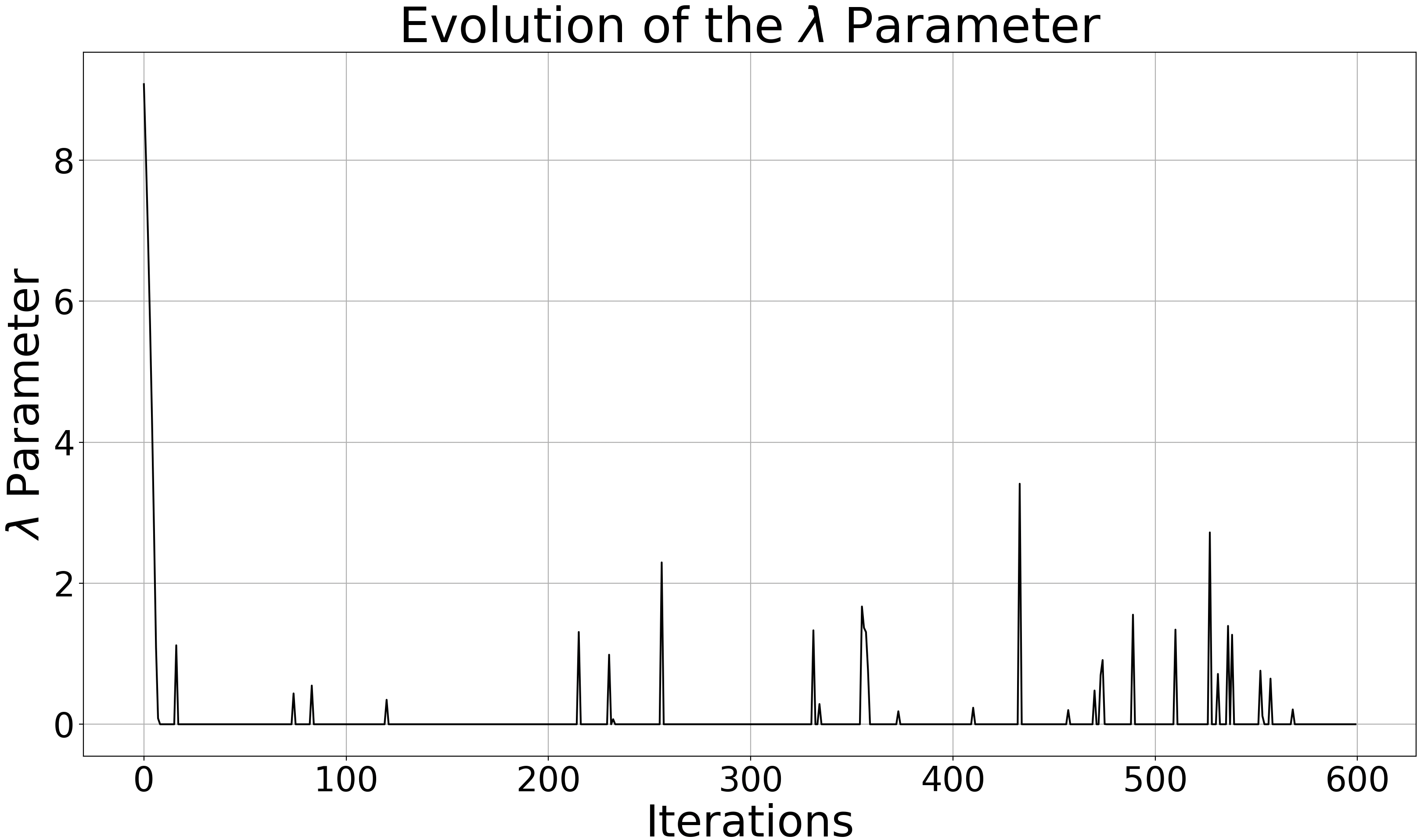 } \put(2,5){\fcolorbox{black}{white}{$m=1$}}\end{overpic}\end{minipage}
	\end{tabular}
	\caption{\it Evolution of the Lagrange multiplier $\lambda$ during the optimization of the cantilever conducted in \cref{sec.simpwdro1load} for various values of the Wasserstein radius $m$; (upper row) $\sigma^2=2e{-2}$, $\varepsilon=1e{-4}$; (lower row) $\sigma^2= 1e{-1}$, $\varepsilon =1e{-4}$.}
	\label{fig.histolambdacanti}
\end{figure}

\bigskip\noindent\textit{Evolution of the Lagrange multiplier $\lambda$.} Eventually, it is interesting to observe the behavior of the added optimization variable $\lambda$ in \cref{eq.wdrocanti} with respect to the original formulation \cref{eq.wdrocantiprimal}. As noted in \cref{rem.lambda}, $\lambda$ can be interpreted as the Lagrange multiplier for the Wasserstein distance constraint in \cref{eq.wdrocantiprimal}, or equivalently as the sensitivity of the maximum value of the mean $\int_\Xi \calC(h,\xi) \:\d \Q$ with respect to perturbations of the Wasserstein radius $m$. 
The evolution of $\lambda$ is reported on \cref{fig.histolambdacanti} for various values of the parameters $\sigma^2$, $\e$ and $m$. 
Looking first at the case where $\sigma^2= 2e{-2}$ and $\varepsilon =1e{-4}$ (upper row in \cref{fig.histolambdacanti}), $\lambda$ keeps increasing through the optimization process when $m=0$, whereas it rapidly converges to $0$ when $m=1$ or $m=5$. This is because the Wasserstein distance constraint is saturated in the former case; hence, a perturbation of $m$ is expected to have a large impact on the optimal value in \cref{eq.wdrocantiprimal}. On the contrary, this constraint is no longer active in the latter two cases, and so modifying the maximum bound on the Wasserstein distance $W_\e(\P,\Q)$ has little to no impact on the optimal value of the objective. Accordingly, we observe that the optimized designs are similar for $m=1$ and $m=5$. 
In the situation where $\sigma^2= 1e{-1}$, $\varepsilon =1e{-4}$ (lower row in \cref{fig.histolambdacanti}), $\lambda$ retains positive values when $m=0$ or $m=0.25$, which is consistent with the fact seen from \cref{fig.cantiwdrom} that the optimized shape is still sensitive to the level of the imposed constraint. 

\subsubsection{Loads uncertainty with a more complex nominal distribution} \label{sec.simpwdromload}

\begin{figure}[ht]
\begin{tabular}{cc}
\centering
		\begin{minipage}{0.4\textwidth}
		\centering
		\begin{overpic}[width=0.9\textwidth]{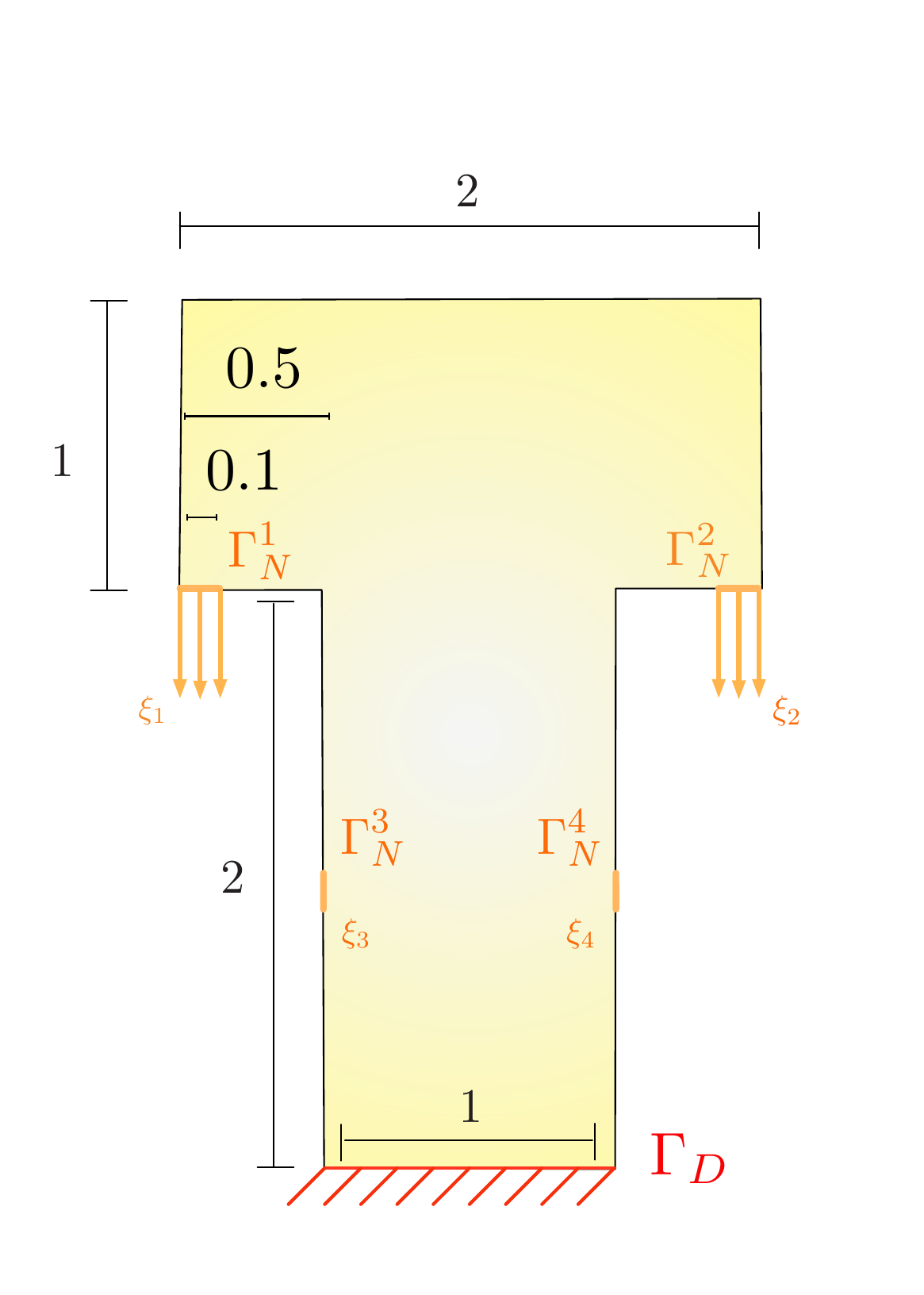} 
		\put(2,5){\fcolorbox{black}{white}{a}}
		\end{overpic}
		\end{minipage} &
	\begin{minipage}{0.4\textwidth}
	\centering
	\begin{overpic}[width=0.65\textwidth]{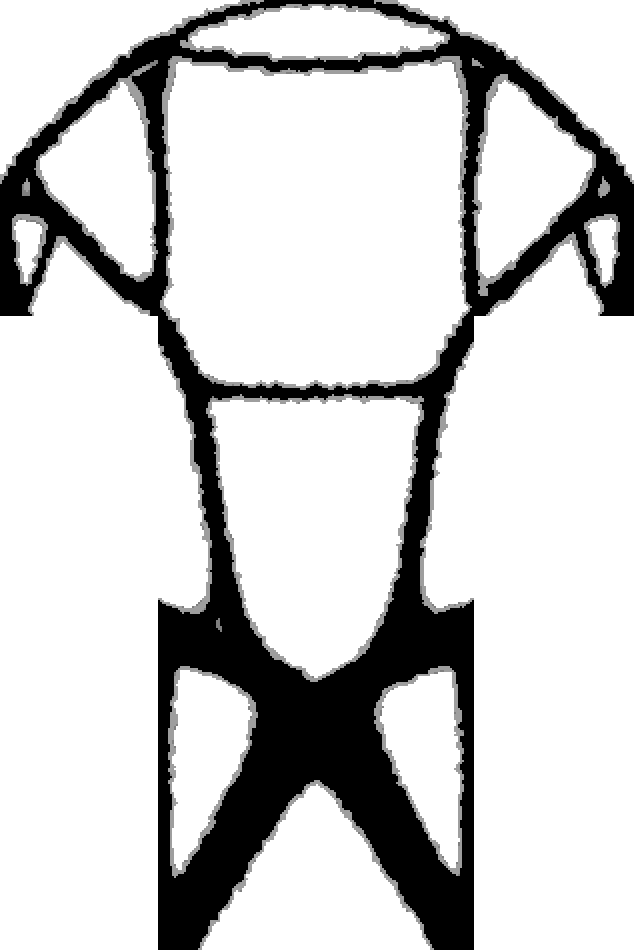} 
	\put(2,5){\fcolorbox{black}{white}{b}}
	\end{overpic}
	\end{minipage}
	\end{tabular}
	\caption{\it (a) Setting of the mast problem considered in \cref{sec.simpwdromload}; (b) Optimized design for the problem \cref{eq.mastJmean}, where the objective function is the mean value of the compliance under the nominal law.}
    \label{fig.mastNom}
\end{figure}

\noindent Our second example illustrates a situation where the nominal law of the uncertain parameter is built from several observations. It concerns the optimal design of a mast, contained in a T-shaped computational domain $D$ with size $2 \times 3$, see \cref{fig.mastNom} (a).
In practice, $D$ is equipped with a mesh composed of approximately $9,\!000$ vertices and $17,\!000$ triangles. The structure is clamped along its bottom boundary $\Gamma_D$, and different loads are applied on four disjoint regions collectively denoted by $\Gamma_N= \bigcup_{i=1}^4 \Gamma_{N,i}$.
These loads are accordingly described by a vector $\xi = \left\{ \xi_i \right\}_{i=1,\ldots,4}$, contained in a sufficiently large ball $\Xi \subset (\R^2)^4$; $\xi$ is made of four two-dimensional components $\xi_i = (\xi_{i,1},\xi_{i,2})$, representing the load applied on the part $\Gamma_{N,i}$, for $i = 1,\ldots,4$. More precisely,
\begin{itemize}
    \item The loads $\xi_1$ and $\xi_2$ applied on the arms of the mast represent the weight of the hanging cables.
    \item The loads $\xi_3$ and $\xi_4$ applied on the lower part of the structure model the effect of wind.
\end{itemize} 

In this context, the displacement of the structure is the solution $u_{h,\xi}$ to the following boundary-value problem:
\begin{equation}\label{eq.elasMast}
\left\{
	\begin{array}{cl}
		-\text{div}(A(h)e(u)) = 0 & \text{in } D,                                                                 \\
		u = 0                     & \text{on } \Gamma_D,                                                          \\
		A(h)e(u)n = \xi_i             & \text{on } \Gamma_{N,i} ,                                                          \\
		A(h)e(u)n = 0             & \text{on } \partial D \setminus \overline{\Gamma_N} \cup \overline{\Gamma_D}.
	\end{array}
\right.
\end{equation}
As far as the  loads are concerned, we observe three scenarios: 
\begin{itemize}
\item $\xi^{1} = ((0, -1), (0, -1), (0, 0), (0, 0))$, i.e. two vertical loads are applied on the arms $\Gamma_{N,1}$, $\Gamma_{N,2}$ of the mast and the regions $\Gamma_{N,3}$ and $\Gamma_{N,4}$ are free;
\item $\xi^{2} = ((0, -1),( 0, -1),( -1, 0),( -1, 0))$: the same vertical loads are applied on $\Gamma_{N,1}$, $\Gamma_{N,2}$, and $\Gamma_{N,3}$ and $\Gamma_{N,4}$ are subjected to a horizontal force directed to the left, mimicking the effect of wind;
 \item $\xi^{3} = ((0, -1),( 0, -1),( 1, 0), (1, 0))$: this scenario is similar to the previous one $\xi^2$, except that the force applied on $\Gamma_{N,3}$ and $\Gamma_{N,4}$ is directed to the right.
\end{itemize}\par\medskip

We first empirically attribute the probabilities $\frac12$, $\frac14$ and $\frac14$ to these respective scenarios, so that the law $\P \in \calP(\Xi)$ of the uncertain parameter $\xi$ is set to:
\begin{equation}\label{eq.Pmast}
\P = \frac{1}{2} \delta_{\xi^1} + \frac{1}{4} \delta_{\xi^2} + \frac{1}{4} \delta_{\xi^3}.
\end{equation}
We then optimize the design of the mast with respect to the mean value of its compliance $C(h,\xi)$ defined in \cref{eq.defcplyTO}, i.e. we solve:
\begin{multline}\label{eq.mastJmean}
\min\limits_{h \in \Uad} \Jmean(h), \text{ s.t. } \Vol(h) = V_T, \text{ where } V_T = 0.3\Vol(D) \text{ and }\\
 \Jmean(h) = \int_\Xi C(h,\xi) \:\d\P(\xi) = \frac12C(h,\xi^1) + \frac14C(h,\xi^2) + \frac14C(h,\xi^3). 
\end{multline}
We run $400$ iterations of the optimization algorithm outlined in \cref{sec.simpnum}, corresponding to approximately $30$ mn of computation; the resulting optimized density function  $h^*_{\text{\rm{mean}}}$ is depicted on \cref{fig.mastNom} (b). \par\medskip
%

\begin{figure}[ht]
\centering
\begin{tabular}{ccc}
		\begin{minipage}{0.3\textwidth}\begin{overpic}[width=1.0\textwidth]{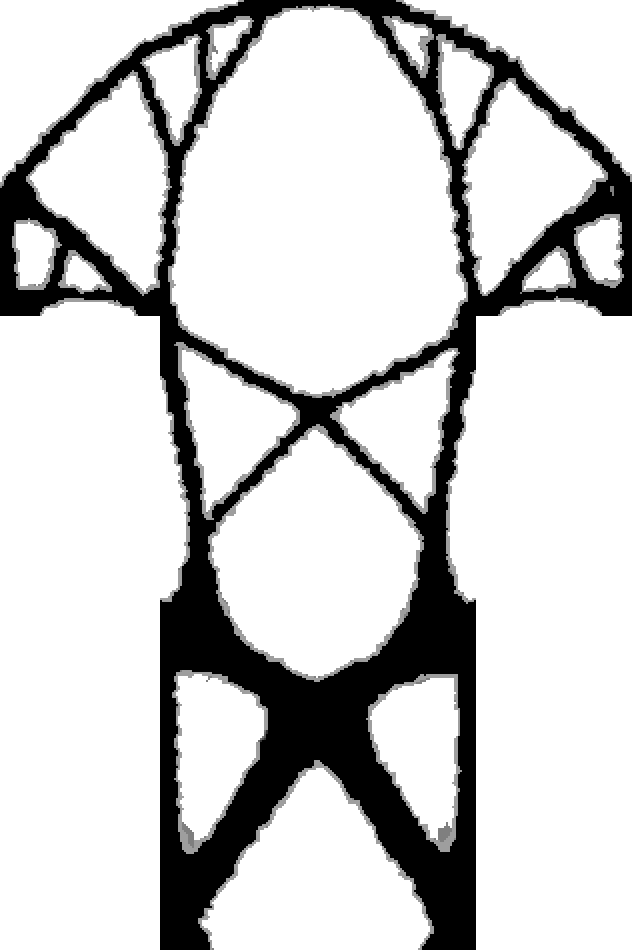} \put(2,5){\fcolorbox{black}{white}{$m=0$}}\end{overpic}\end{minipage} & \begin{minipage}{0.3\textwidth}\begin{overpic}[width=1.0\textwidth]{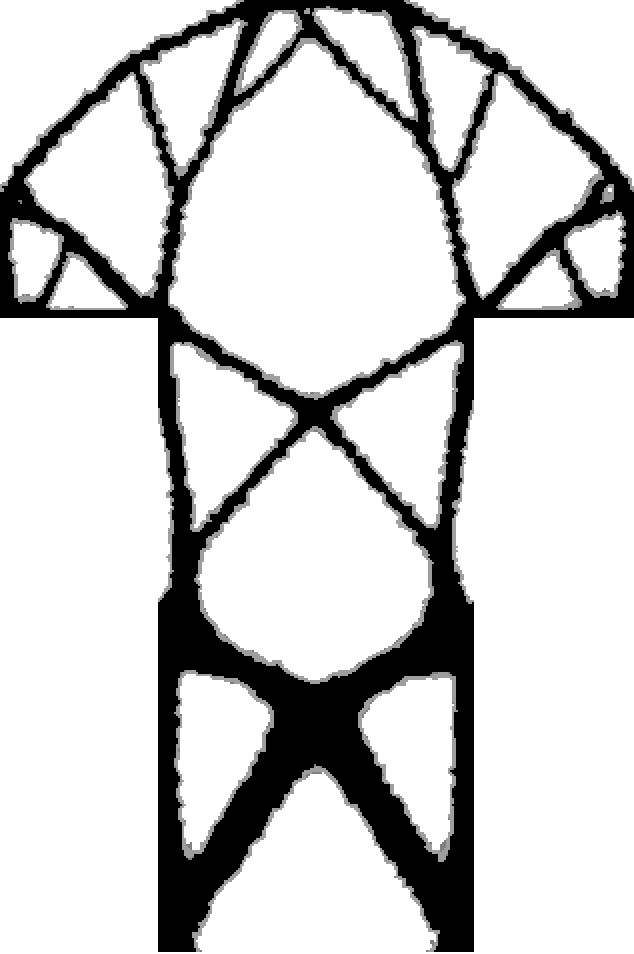} \put(2,5){\fcolorbox{black}{white}{$m=5$}}\end{overpic}\end{minipage} & \begin{minipage}{0.3\textwidth}\begin{overpic}[width=1.0\textwidth]{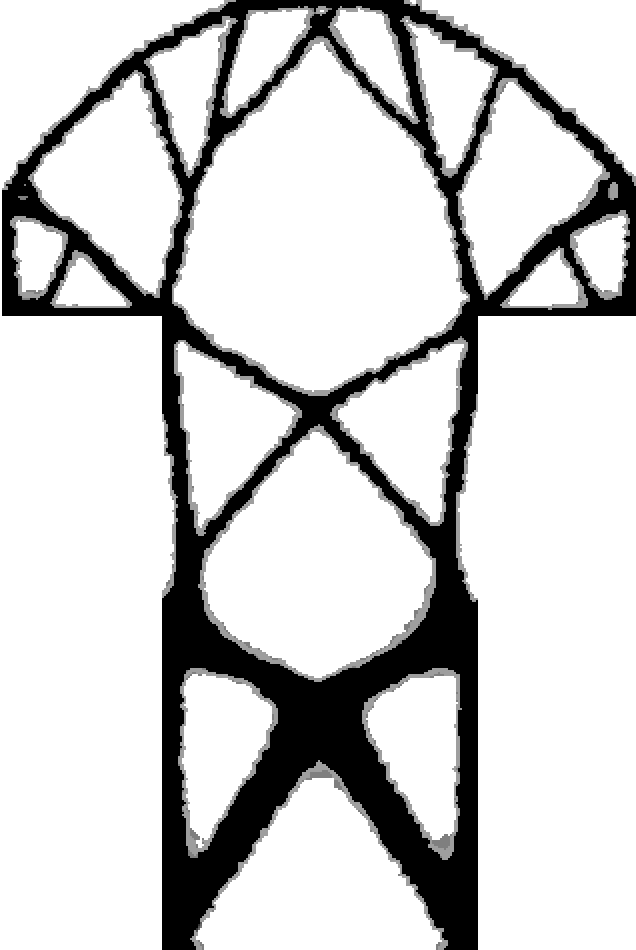} \put(2,5){\fcolorbox{black}{white}{$m=10$}}\end{overpic}\end{minipage}
	\end{tabular}
	\caption{\it Distributionally robust optimized designs $h^*_{m=0}$, $h^*_{m=5}$ and $h^*_{m=10}$ of the electric mast considered in \cref{sec.simpwdromload} associated to the values $m=0$, $5$ and $10$.}
    \label{fig.MultiMast}
\end{figure}

We now come to the situation of interest, where the probability law of the uncertain parameter $\xi$ is uncertain. Having little faith in the empirical reconstruction $\P$ in \cref{eq.Pmast}, we now minimize the worst value of the expectation of $C(h,\xi)$ when the actual law $\Q$ of $\xi$ is ``close'' to $\P$ in the sense that it belongs to the Wasserstein ambiguity set $\calA_{\text{W}}$ in \cref{eq.calAW}: 
\begin{equation}\label{eq.dromast}
\min\limits_{h \in \Uad} \sup\limits_{\Q \in \calA_{\text{W}}} \int_\Xi C(h,\xi) \:\d\Q(\xi), \text{ s.t. } \Vol(h) = V_T.
\end{equation}
Now invoking the material in \cref{sec.wdroformula}, this problem rewrites:
\begin{multline*}
	\min\limits_{h \in \Uad, \atop \lambda \geq 0} \calD_{\text{W}}(h,\lambda)  \text{ s.t. } \quad \text{Vol}(h) = V_T, \quad \text{ where }\\
	 \calD_{\text{W}}(h,\lambda) = \lambda m + \lambda \e \int_{\Xi} \log\left(
	\int_{\Xi} e^{\frac{C(h,\zeta) - \lambda c(\xi,\zeta)}{\lambda \e}}
	\,\d\nu_{\xi}(\zeta) \right)\,\d\P(\xi),
\end{multline*}
and $\nu_\xi(\zeta) \in \calP(\Xi)$ is the probability measure \cref{eq.refcouplingW}.  
This example is more involved than the previous one: since the nominal law $\P$ is made of several observations, the objective function $\mathcal{D}_{\text{W}}$ now involves two nested integrals over the parameter space $\Xi$. In practice, the innermost integral is approximated using $10$ samples, so that an evaluation of $\mathcal{D}_{\text{W}}$ requires the solution to 30 versions of \cref{eq.elasMast}, associated to as many different load cases; in practice, these operations can easily be conducted in parallel. 
Three designs, $h^*_{m=0}$, $h^*_{m=5}$ and $h^*_{m=10}$, optimized for the respective values $m=0,5$ and $10$ of the Wasserstein radius, are displayed in \cref{fig.MultiMast}. Each computation takes about 600 iterations of our numerical strategy, for a computational time of about 10h.

The three optimized designs resulting from this experiment present less obvious visual differences than in the example of \cref{sec.simpwdro1load}. This is because, in the present situation, contrary to the previous one, the three loads $\xi^1$, $\xi^2$, $\xi^3$ featured by the nominal law are already quite harmful for the structure compliance, and so the nominal distribution is already not that far from capturing the worst-case situation. 
However, we observe that as $m$ increases, the arches near the arms of the mast $h^*_{\text{\rm{mean}}}$ turn to flat connections, in response to the anticipation of possible horizontal perturbations of the loads in there. Moreover, the regions of the structure near $\Gamma_{N,3}$ and $\Gamma_{N,4}$ become thicker, in order to better withstand vertical perturbations of the loads.

\bigskip\noindent\textit{Out-of-sample analysis.}
In order to compare the design $h^*_{\text{\rm{mean}}}$ resulting from the minimization \cref{eq.mastJmean} of the expected compliance with the distributionally robust designs $h^*_{m=0}$, $h^*_{m=5}$ and $h^*_{m=10}$, we compute their compliance $C(h,\xi)$ for several load scenarios; the results are registered in \cref{tab.tablecompmulti}. The first three scenarios $\xi^1$, $\xi^2$, $\xi^3$ are those present in the nominal distribution; the four next vectors $(\xi^i)_{4 \leq i\leq 7}$ model situations where, in addition to the weight of the cables hanging from the arms, horizontal stretching or compression effects are added. The final ten vectors $(\xi^i)_{8 \leq i \leq 17}$ have random entries, drawn in the set $\{-1, 1, 0\}$. 
The design $h^*_{\text{\rm{mean}}}$ optimized for the expectation of the compliance performs best in the first 3 load situations, which is expected since they are present in the nominal law ${\P}$; on that contrary, $h^*_{\text{\rm{mean}}}$ shows poor performance in the other situations. 
In contrast, the distributionally robust designs $h^*_{m=0}$, $h^*_{m=5}$ and $h^*_{m=10}$ show lesser performance in the first 3 events, but their efficiency remains stable in the case of loads that are not represented in the nominal distribution.
This confirms the ability of the distributionally robust approach to anticipate events that are not present in the nominal distribution, but ``close'' to some of them, see again \cite{kapteyn2019distributionally} for an account of this capability. 
Note that, in some scenarios, the designs optimized at larger Wasserstein radius have better performance; this happens when the applied load corresponds to a perturbation which is ``far'' from those in the nominal law, and is thus only effectively captured when the radius of the Wasserstein constraint is large enough. On the contrary, when such perturbations are already captured at the initial radius, increasing the radius further may decrease the impact of this capture, ultimately resulting in lower performance.
 
\renewcommand{\arraystretch}{1.25}
\begin{table}[ht]
	\centering
	\centering
	\begin{tabular}{|c|c|c|c|c|}
		\hline
		\backslashbox[60mm]{$\xi$}{$h$}     & $h^*_{\text{\rm{mean}}}$ & $h^*_{m=0}$ &  $h^*_{m=5}$ &  $h^*_{m=10}$ \\
		\hline
		$\xi^{1} =((0, -1),( 0, -1), (0, 0), (0, 0))$                   & 1.45263              & 1.58169           & 1.80382           & 1.84381           \\
		\hline
		$\xi^{2}=((0, -1),( 0, -1), (1, 0), (1, 0))$                    & 3.96816              & 4.54952           & 4.95328           & 5.04692           \\
		\hline
		$\xi^{3}= ((0, -1), (0, -1), (-1, 0), (-1, 0))$                 & 3.97635              & 4.62913           & 4.93366           & 5.02842        \\
		\hline
		$\xi^{4} = ((0, -1), (0, -1),( -1, 0), (1, 0))$    & 2.30295              & 2.56223           & 2.99777            & 3.02278         \\
		\hline
		$\xi^{5} =$ ((0, -1), (0, -1), (1, 0), (-1, 0))    & 2.0726               & 1.99442           & 2.18521           & 2.21876         \\
		\hline
		$\xi^{6} =$ ((-1, -1),( -1, -1), (-1, 0), (-1, 0)) & 90.9951              & 18.2647           & 16.5193           & 16.5381  \\
		\hline
		$\xi^{7} =$ ((1, -1), (1, -1), (-1, 0), (-1, 0))   & 81.7526              & 8.31552           & 6.92599           & 6.76357           \\
		\hline
		$\xi^{8} =$ ((1, -1), (1, 1), (-1, -1), (1, 0))    & 53.2508              & 15.6051           & 10.7477           & 10.1949           \\
		\hline
		$\xi^{9} =$ ((1, 1),( 1, 1), (-1, 1),( 0, -1))     & 87.6837              & 12.8845           & 11.0046           & 10.9336           \\
		\hline
		$\xi^{10} =$ ((0, 0), (0, 0), (-1, -1), (0, 0))    & 2.55941              & 2.30058           & 2.22696           & 2.22671           \\
		\hline
		$\xi^{11} =$ ((-1, -1), (-1, 1), (1, 1), (0, 1))   & 196.139              & 28.1766           & 19.6703           & 19.3135           \\
		\hline
		$\xi^{12} =$ ((0, 0), (1, 0), (1, 1),( -1, 1))     & 29.1637              & 6.60941           & 5.6305           & 5.56953          \\
		\hline
		$\xi^{13} =$ ((-1, 1), (1, -1),( -1, 1),(0, -1))  & 59.6293              & 18.1618           & 12.9192           & 12.609          \\
		\hline
		$\xi^{14} =$ ((0, 0),( 0, -1), (0, -1),( 0, 1))    & 11.3029              & 3.96029           & 2.79272           & 2.72636
          \\
		\hline
		$\xi^{15} =$ ((1, -1), (1, 1),( 0, -1), (0, 0))    & 51.8462              & 14.3922           & 9.50305           & 8.9527           \\
		\hline
		$\xi^{16} =$ ((0, 0),( 0, 0), (0, -1),( 1, 0) )    & 1.51937              & 1.35681           & 1.40098           & 1.40919
         \\
		\hline
		$\xi^{17} =$ ((0, 0), (0, 0), (0, 1), (1, 0))      & 1.89735              & 1.81272           & 1.77383           & 1.77781           \\
		\hline
	\end{tabular}
	\caption{\it Values of the compliance $C(h,\xi)$ of the optimized masts in \cref{sec.simpwdromload}, associated to several values of the load parameter $\xi$.} 
	\label{tab.tablecompmulti}
\end{table}

\subsection{Moment-based distributionally robust optimization of the compliance of an L-shaped beam}\label{sec.simpmdro}

\noindent
This section discusses a situation where the ambiguity set $\calA \equiv \calA_{\text{M}}$ used in the formulation of the distributionally robust problem is defined from the moments of probability distributions, see \cref{eq.calAmoments} in \cref{sec.momdro}. 

\begin{figure}[ht]
\centering
	\begin{tabular}{ccc}
		\begin{minipage}{0.31\textwidth}
		\begin{overpic}[width=1.0\textwidth]{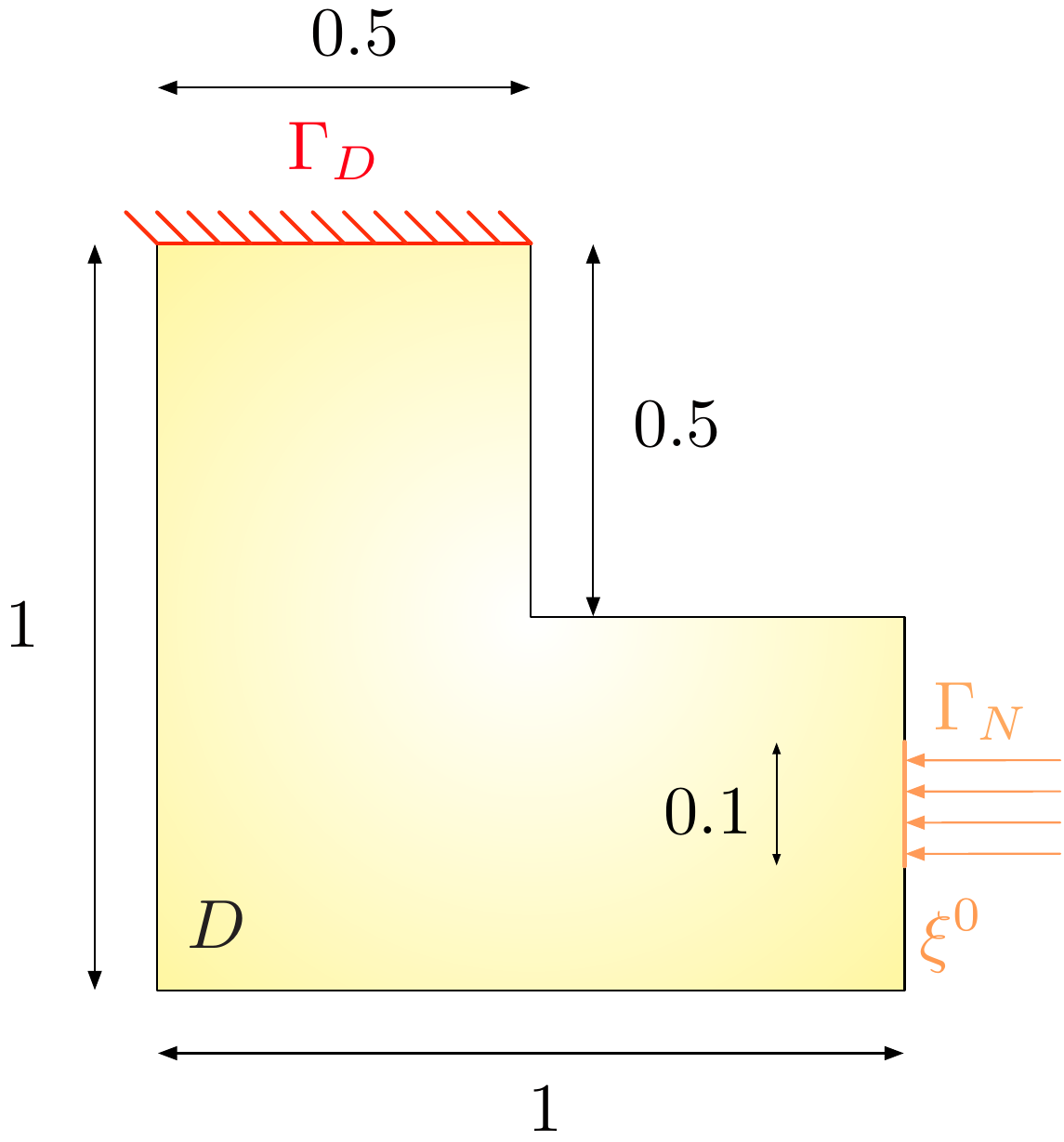} 
		\put(2,5){\fcolorbox{black}{white}{a}}
		\end{overpic}
		\end{minipage} & 
	\begin{minipage}{0.31\textwidth}
	\begin{overpic}[width=0.9\textwidth]{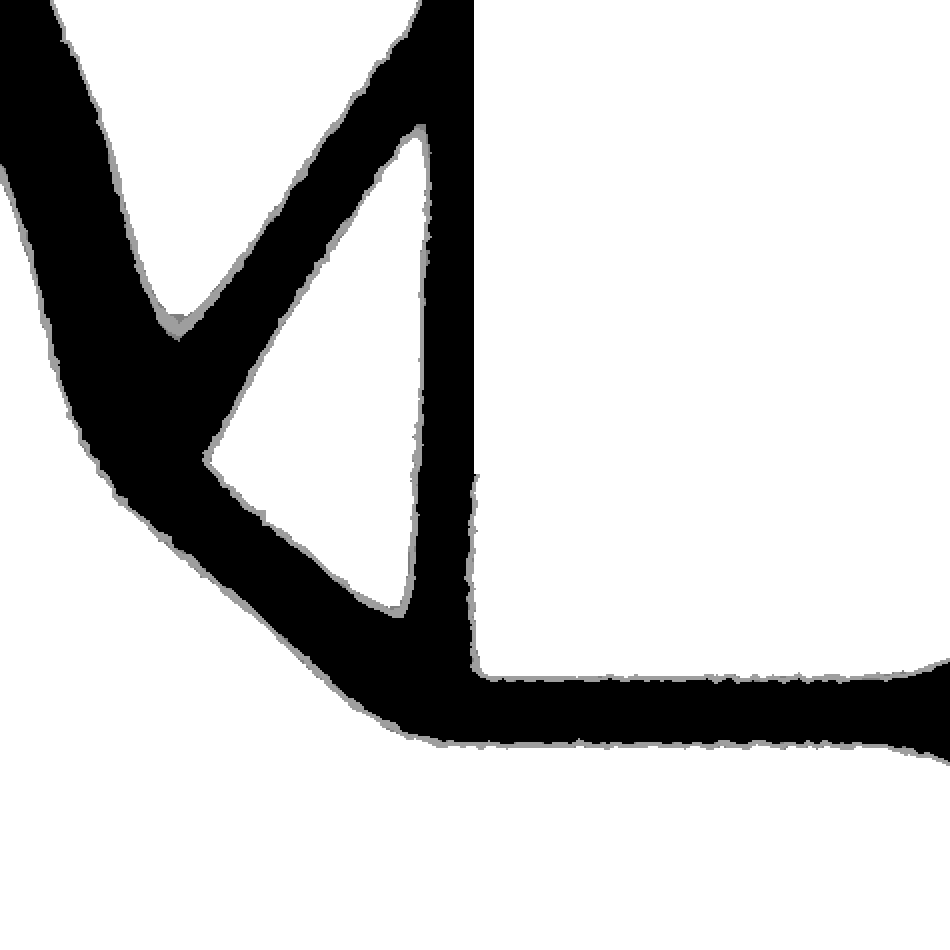}
	\put(2,5){\fcolorbox{black}{white}{b}}
	\end{overpic}
	\end{minipage} & 
		\begin{minipage}{0.31\textwidth}
	\begin{overpic}[width=0.9\textwidth]{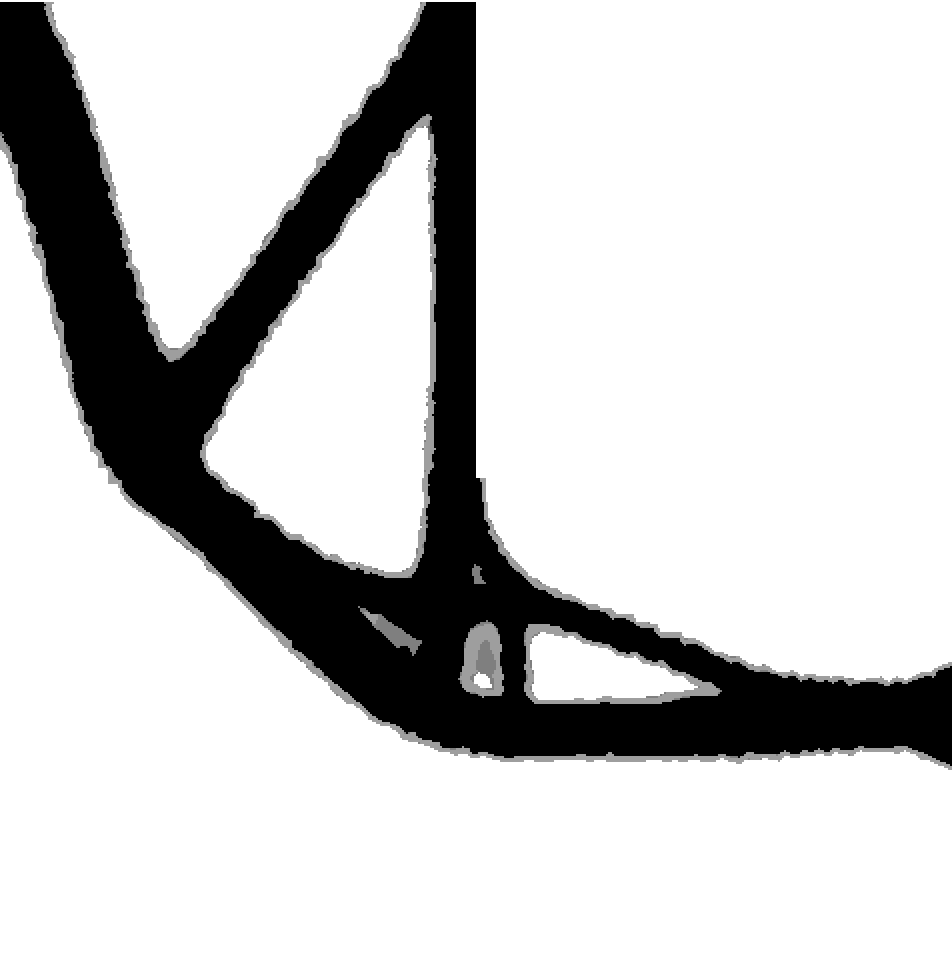}
	\put(2,5){\fcolorbox{black}{white}{c}}
	\end{overpic}
	\end{minipage}
	\end{tabular}
	\caption{\it Optimal design example of an L-shaped beam considered in \cref{sec.simpmdro}; (a) Setting of the test-case; (b) Optimized design $h^*_{\text{\rm det}}$ in the unperturbed situation; (c) Optimized design  $h^*_{\text{\rm mean}}$ for the mean value of the compliance $C(h,\xi)$ under the nominal law $\P$.}
    \label{fig.Lshapenom}
\end{figure}

As illustrated on \cref{fig.Lshapenom} (a), we now consider beams, contained in an L-shaped computational box $D$ with dimensions $1 \times 1$; the domain $D$ is equipped with a mesh comprising approximately $8,500$ vertices and $17,000$ triangles. The optimized structures are clamped on the upper side $\Gamma_D$ of $\partial D$, and a load $\xi$ is applied on a small region $\Gamma_N$ at its right-hand side.

In ideal circumstances, $\xi$ is known perfectly: it is the unit horizontal vector $\xi^0 = (-1, 0)$. We then solve the problem:
\begin{equation}\label{eq.minJdetLbeam} 
\min\limits_{h \in \Uad} C(h, \xi^0) \text{ s.t. } \Vol(h) = V_T.
\end{equation}
Here, $C(h,\xi)$ stands for the compliance \cref{eq.defcplyTO} of the beam when it is subjected to the load $\xi$, and the target volume equals $V_T = 0.2$. 
About $600$ iterations of the optimization algorithm sketched in \cref{sec.simpnum} are applied to this end; the resulting density $h^*_{\text{det}}$ is depicted on \cref{fig.Lshapenom} (b). 

\par\medskip

We now assume that the parameter $\xi$ is uncertain, but that its probability law $\P$ is perfectly known: $\xi$ is a Gaussian vector whose first- and second-order moments $\mu_0$ and $\Sigma_0$ are given by:
$$ \mu_0 = \xi^0, \text{ and } \Sigma_0 = \sigma^2 \I, \:\: \sigma^2 = 1e{-2}.$$ 
We then minimize the mean value of the cost functional with respect to $\P$, i.e. we solve:
\begin{equation}\label{eq.minJmeanLbeam}
    \underset{h\in \Uad}{\min}\; \int_{\Xi} C(h, \xi) \:\d\P(\xi), \text{ s.t. } \quad \Vol(h)=V_T.
\end{equation}
Again applying $600$ iterations of our algorithm ($\approx$ 1h of computation), the optimized density $h^*_{\text{\rm{mean}}}$ is that depicted in \cref{fig.Lshapenom} (c).\par\medskip

We now turn to a situation where the probability law of $\xi$ is uncertain, and only its first- and second-order moments are known.
Following \cref{sec.momdro}, we then consider the following distributionally robust counterpart of \cref{eq.minJdetLbeam,eq.minJmeanLbeam}:
$$ \min\limits_{h \in \Uad} \:\left(\sup\limits_{\Q \in \calA_{\text{\rm M}}} \int_\Xi C(h,\xi)\:\d\Q(\xi) - \e H(\Q) \right) \text{ s.t. } \Vol(h) = V_T,$$
featuring the moment-based ambiguity set $\calA_{\text{\rm M}}$ defined in \cref{eq.calAmoments}; the entropy $H(\Q)$ of a law $\Q \in \calP(\Xi)$ is defined by \cref{eq.entropymoments} and the reference law $\Q_0$ in this formula is $\Q_0 = \P$. The parameter $\e$ is set to $\e=0.01$.
According to \cref{sec.momdro}, this problem equivalently rewrites:
\begin{multline*}
	\min\limits_{h \in \Uad, \: \lvert \tau \lvert\leq 1, \atop
	{\lambda \geq 0, {S \in \mathbb{S}_+^m} } }\calD_{\text{M}}(h,\lambda
	,\tau,S), \text{ s.t. } \Vol(h)=V_T \text{ where }\\[-1em] \calD_{\text{M}}(h,\lambda,\tau,S):=
	\lambda m_{1}- \lambda \tau \cdot \mu_{0}+ m_{2}S : \Sigma_{0}+ \e \log
	\left( \int_{\Xi}\left(e^{\frac{C(h,\xi) + \lambda \tau \cdot \xi
	- S : (\xi-\mu_{0}) \otimes (\xi-\mu_{0})}{\e}}\right) \:\d \P(\xi
	)\right).
\end{multline*}

\begin{figure}[ht]
	\centering
	\begin{tabular}{ccc}
		\begin{minipage}{0.3\textwidth}\begin{overpic}[width=1.0\textwidth]{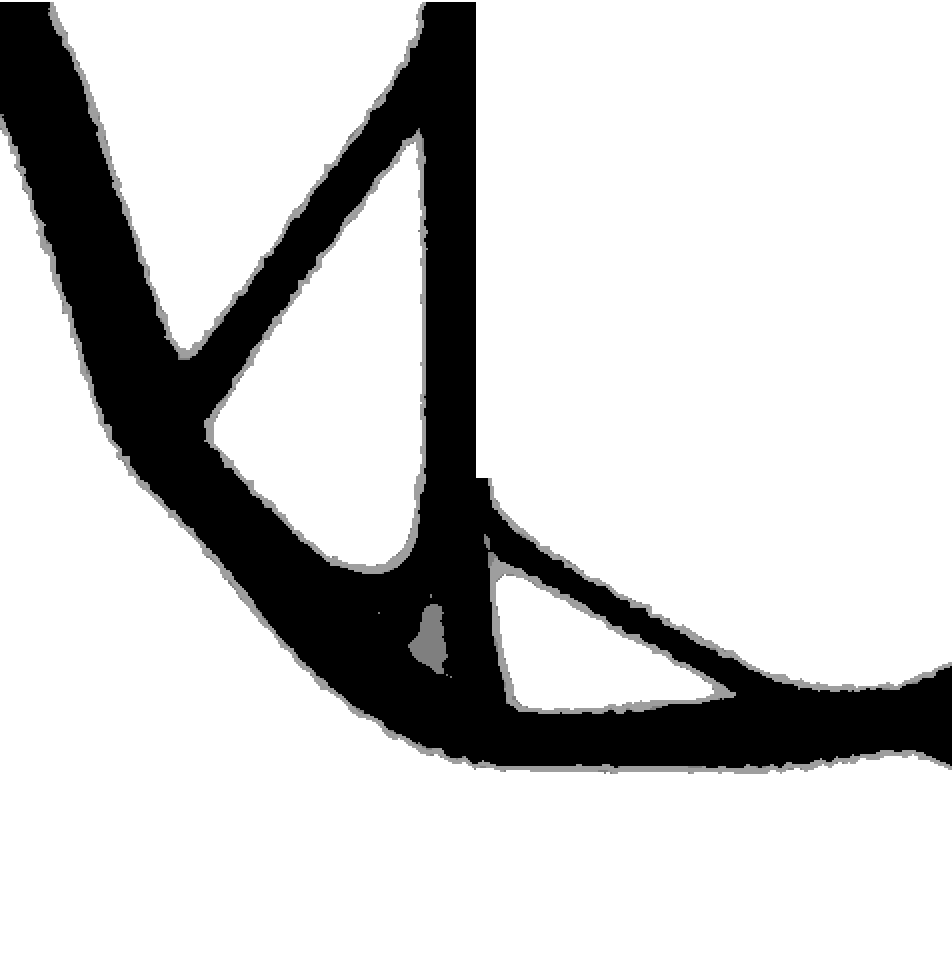} \put(2,5){\fcolorbox{black}{white}{$m_{1}=0$, $m_{2}=1$}}\end{overpic}\end{minipage} & \begin{minipage}{0.3\textwidth}\begin{overpic}[width=1.0\textwidth]{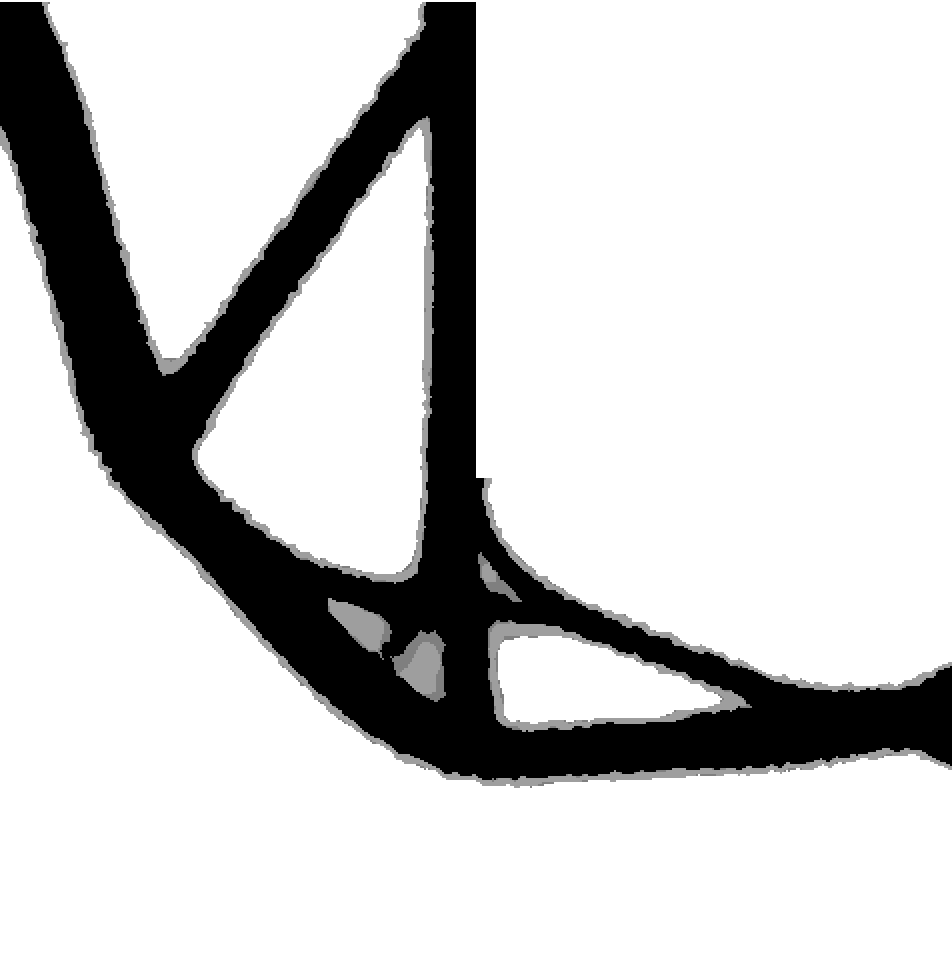} \put(2,5){\fcolorbox{black}{white}{$m_{1}=1$, $m_{2}=1$}}\end{overpic}\end{minipage} & \begin{minipage}{0.3\textwidth}\begin{overpic}[width=1.0\textwidth]{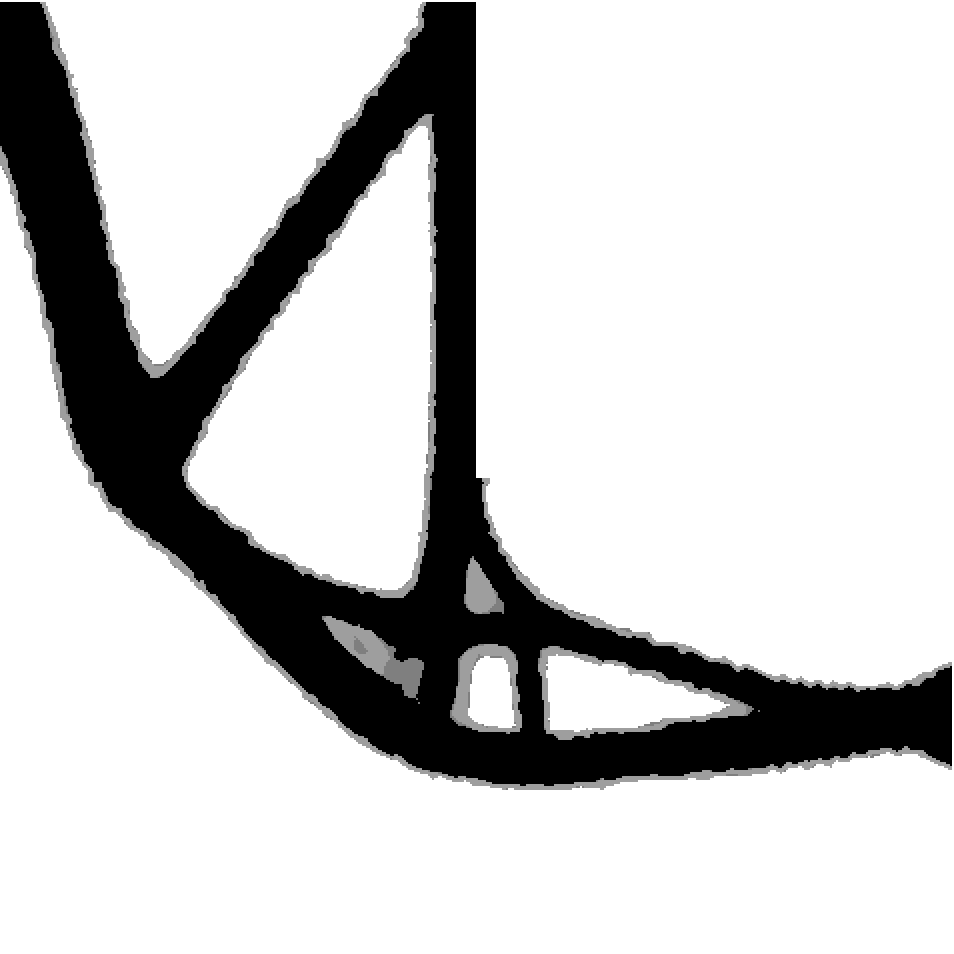} \put(2,5){\fcolorbox{black}{white}{$m_{1}=5$, $m_{2}=1$}}\end{overpic}\end{minipage}
	\end{tabular}
	\par
	\medskip
	\begin{tabular}{ccc}
		\begin{minipage}{0.3\textwidth}\begin{overpic}[width=1.0\textwidth]{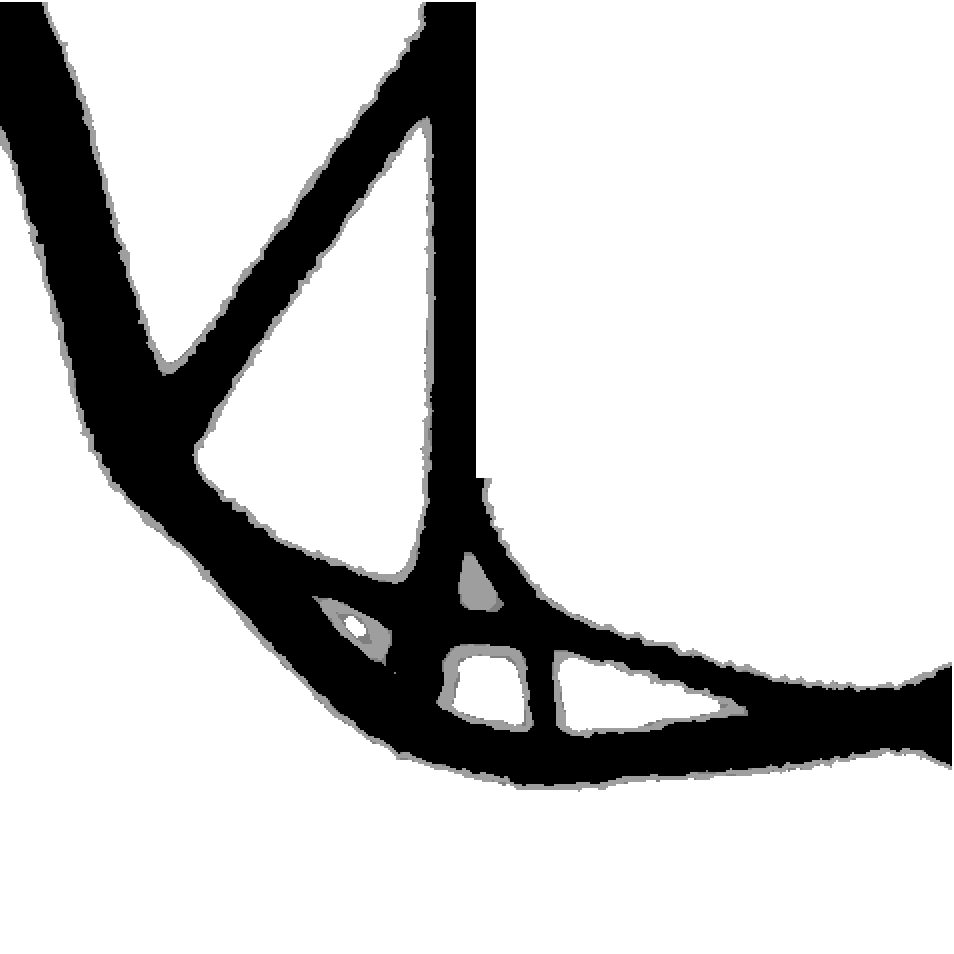} \put(2,5){\fcolorbox{black}{white}{$m_{1}=1$, $m_{2}=2$}}\end{overpic}\end{minipage} & \begin{minipage}{0.3\textwidth}\begin{overpic}[width=1.0\textwidth]{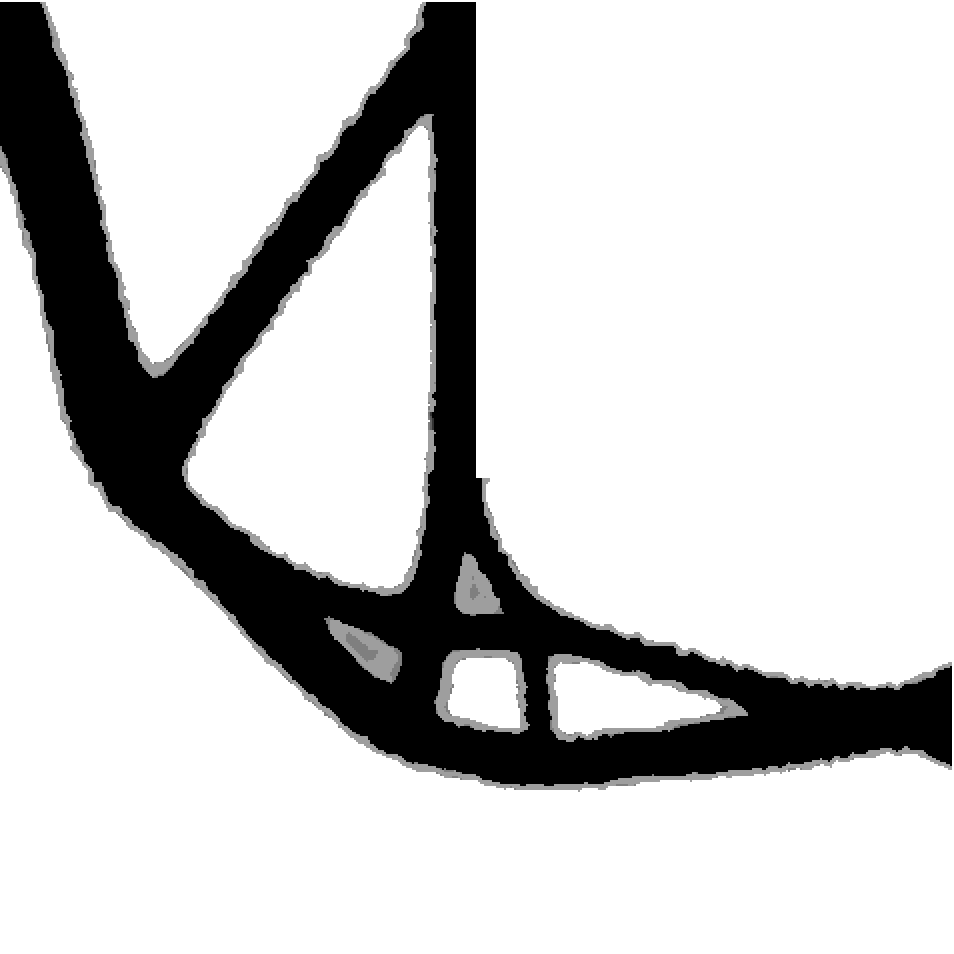} \put(2,5){\fcolorbox{black}{white}{$m_{1}=5$, $m_{2}=5$}}\end{overpic}\end{minipage} & \begin{minipage}{0.3\textwidth}\begin{overpic}[width=1.0\textwidth]{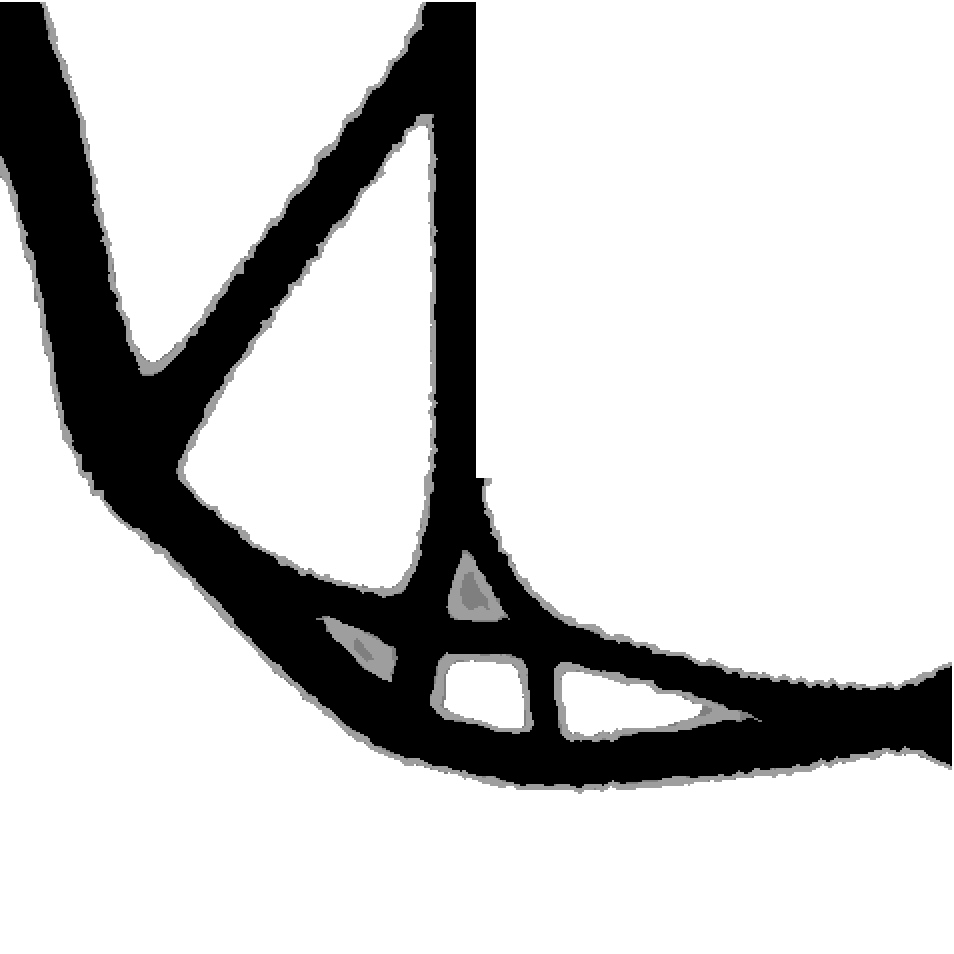} \put(2,5){\fcolorbox{black}{white}{$m_{1}=10$, $m_{2}=10$}}\end{overpic}\end{minipage}
	\end{tabular}
	\caption{\it Distributionally robust optimized designs obtained in the L-shaped beam example of \cref{sec.simpmdro}, associated to various values of the parameters $m_1$, $m_2$ of the moment-based ambiguity set ${\mathcal A}_{\text{M}}$.}
    	\label{fig.Ldro}
\end{figure}

We solve this problem for several values of the bounds $m_1$, $m_2$ defining the ambiguity set $\mathcal{A}_{\text{M}}$; the optimized designs $h^*_{m_1,m_2}$ are presented in \cref{fig.Ldro}. Each computation requires 600 iterations of our algorithmic strategy, for a CPU time of about 70 mn.

According to these results, as the ambiguity set gets wider, the optimized designs tend to develop a rounded re-entrant corner. 
This feature is well-known to better attenuate the stress concentrations induced by vertical loads in this region which is constrained by the geometry of the design domain $D$, see e.g. \cite{allaire2008minimum}. 
As expected, the compliance $C(h^*_{m_1,m_2},\xi^0)$ of these designs in the situation where the ideal load $\xi^0$ is applied deteriorates as the ambiguity set increases, as shown in \cref{tab.lshapeopti}. 
As far as the effect of loads different from $\xi^0$ is concerned, an out-of-sample analysis in the spirit of that conducted in the previous \cref{sec.simpwdromload}, yields similar conclusions. 
Omitting full details for brevity -- we refer to the thesis \cite{prando2025distributionally} about these --, let us solely mention that again, the designs $h^*_{\text{det}}$ and  $h^*_{\text{\rm{mean}}}$ are less efficient than those resulting from the distributionally robust formulation when they are submitted to a load which is different from $\xi^0$. This happens even for loads that are ``close'' to $\xi^0$, although the difference gets of course much more significant when they get farther.
This observation again underlines the ability of distributionally robust formulations to accommodate uncertain events that are not part of the observations used for reconstructing the nominal law $\P$.

\begin{table}[ht]
    \centering
    \resizebox{\textwidth}{!}{%
    \begin{tabular}{|c|c|c|c|c|c|c|c|c|}
        \hline
              & $h^*_{\text{det}}$ &  $h^*_{\text{\rm{mean}}}$ & $h^*_{m_{1}=0, m_{2}=1}$       & $h^*_{m_{1}=1,m_{2}=1}$       & $h^*_{m_{1}=5, m_{2}=1}$       & $h^*_{m_{1}=1, m_{2}=2}$      & $h^*_{m_{1}=5, m_{2}=5}$       & $h^*_{m_{1}=10, m_{2}=10}$   \\
        \hline
        $C(h,\xi^0)$ & 0.5409451 & 0.595818 & 0.63298 & 0.661997 & 0.666666 & 0.645945 & 0.659693 & 0.670644 \\
        \hline
    \end{tabular}
    }
    \caption{\it Values of the compliance $C(h,\xi^0)$ of the distributionally robust designs of the L-shaped beam considered in \cref{sec.simpmdro} when the nominal load $\xi^0$ is applied.}
    \label{tab.lshapeopti}
\end{table}


\subsection{Distributionally robust topology optimization under material uncertainties}\label{subsec.matdro}

\noindent The present section considers uncertainties about the properties of the constituent material of the design. 
This situation departs from \cref{sec.simpwdroload,sec.simpmdro} on at least two aspects: on the one hand, 
the uncertain parameter $\xi$ does not naturally pertain to a finite-dimensional space; rather, it stems from a truncated Karhunen-Lo\`eve expansion of an infinite-dimensional random field. On the other hand,
the considered optimal design problem features a cost function which is not self-adjoint, see \cref{prop.Dhprime} below.  \par\medskip

The designs under scrutiny are gripping mechanisms, see \cref{fig.gripnom} (a); they are enclosed in a fixed computational domain
$D$ with size $1\times 1$, equipped with a mesh composed of $9,000$ vertices ($\approx 17,000$ triangles). A load $g = (0.1,0)$ is applied on a region $\Gamma_N$ of the left-hand side of $\partial D$.

\begin{figure}[ht]
\centering
\begin{tabular}{cccccc}
\begin{minipage}{0.16\textwidth}
\begin{overpic}[width=\linewidth]{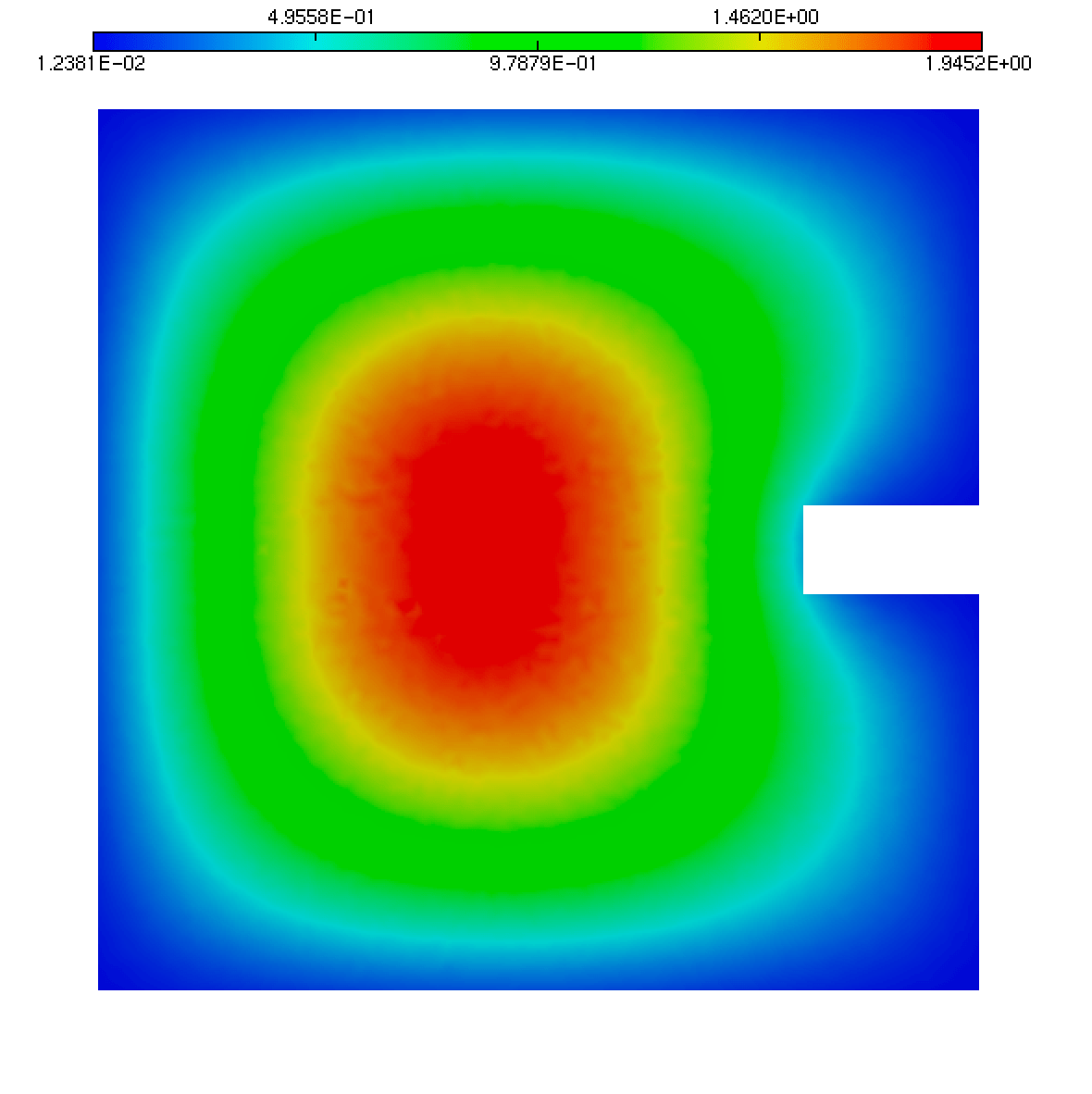} 
\put(2,5){\fcolorbox{black}{white}{$\lambda_{1}=0.7411$}}
\end{overpic}
\end{minipage} & 
\begin{minipage}{0.16\textwidth}
\begin{overpic}[width=\linewidth]{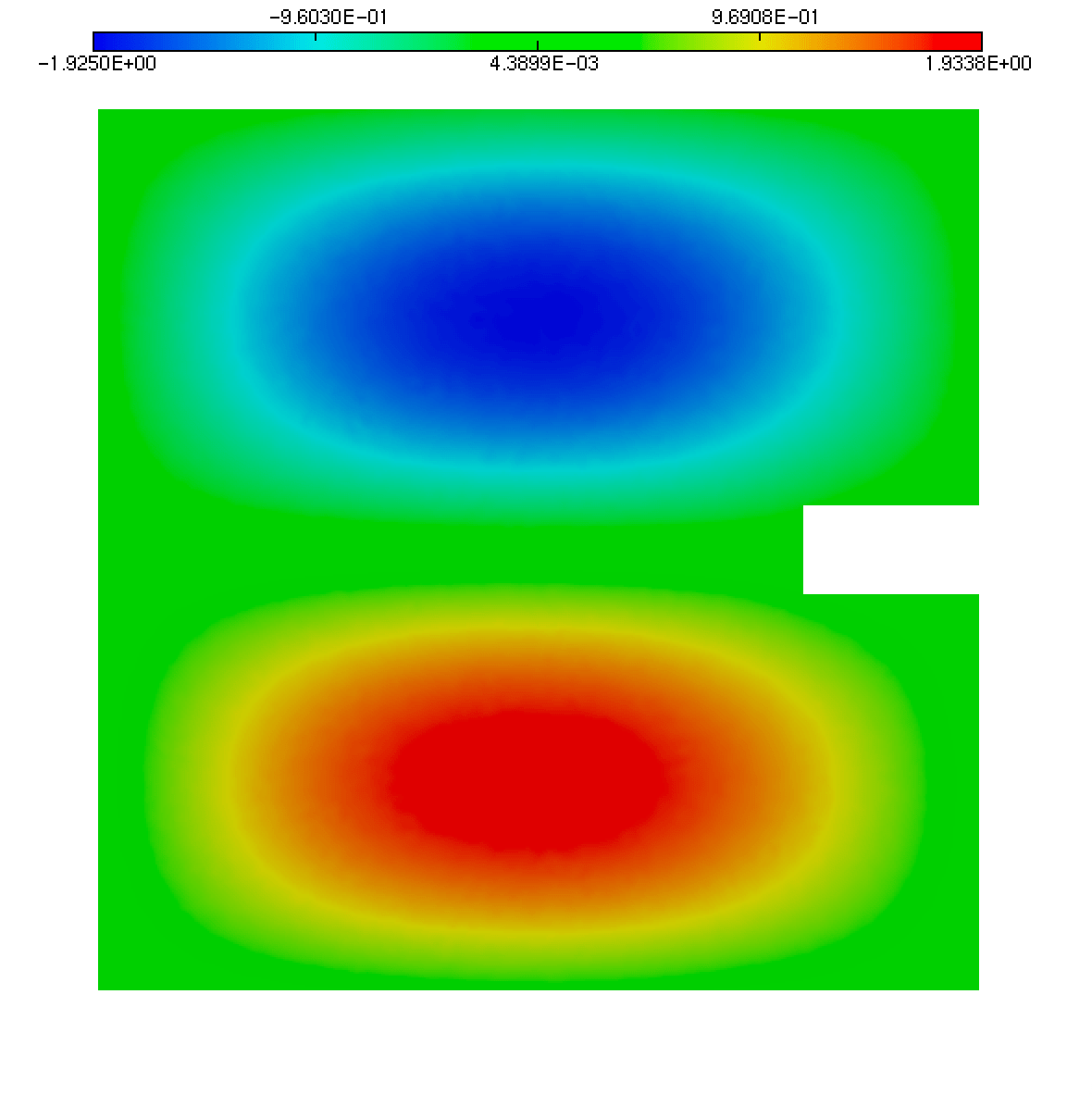} 
\put(2,5){\fcolorbox{black}{white}{$\lambda_{2}=0.7250$}}
\end{overpic}
\end{minipage} & 
\begin{minipage}{0.16\textwidth}
\begin{overpic}[width=\linewidth]{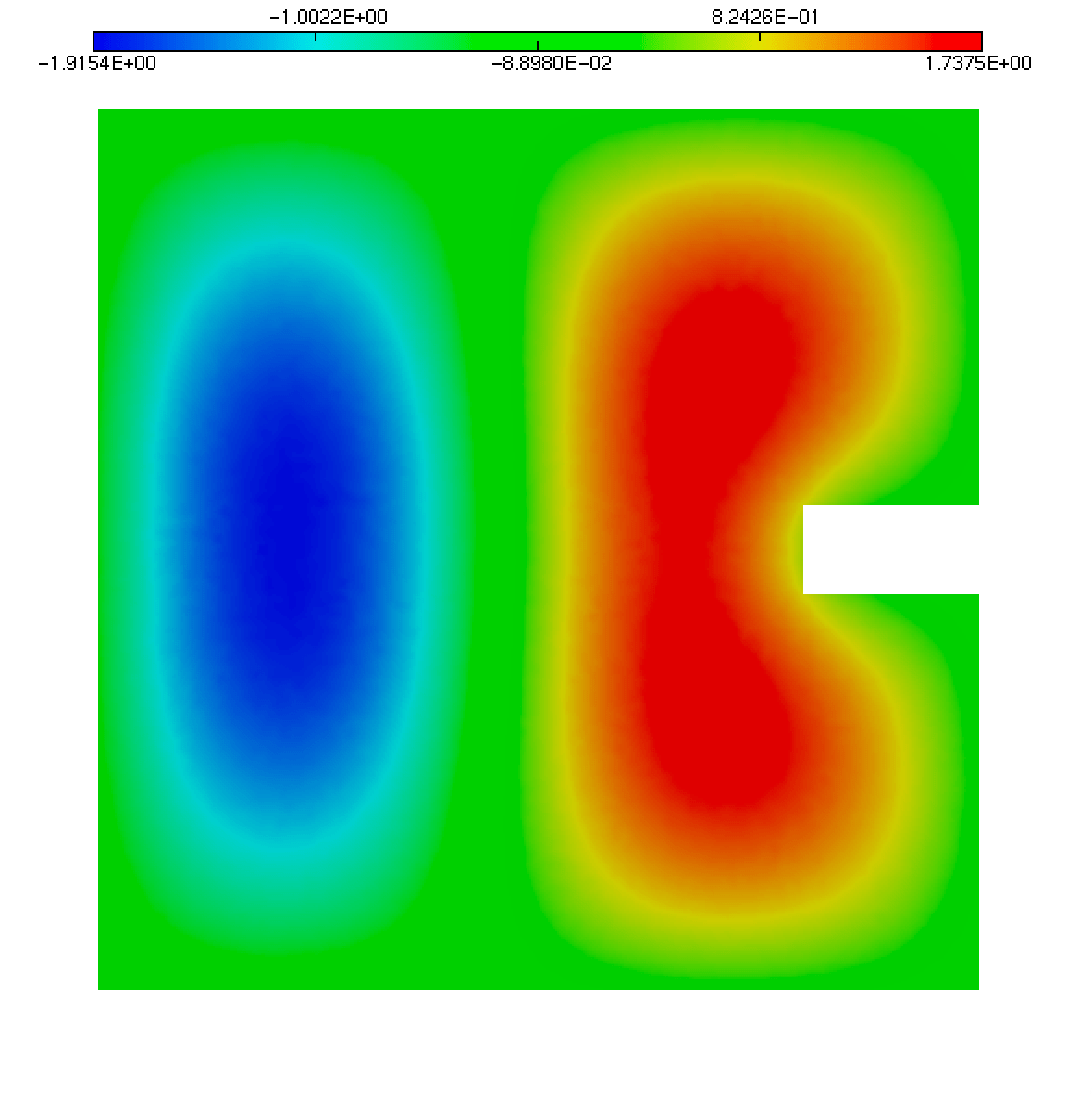} 
\put(2,5){\fcolorbox{black}{white}{$\lambda_{3}=0.7190$}}
\end{overpic}
\end{minipage} & 
\begin{minipage}{0.16\textwidth}
\begin{overpic}[width=\linewidth]{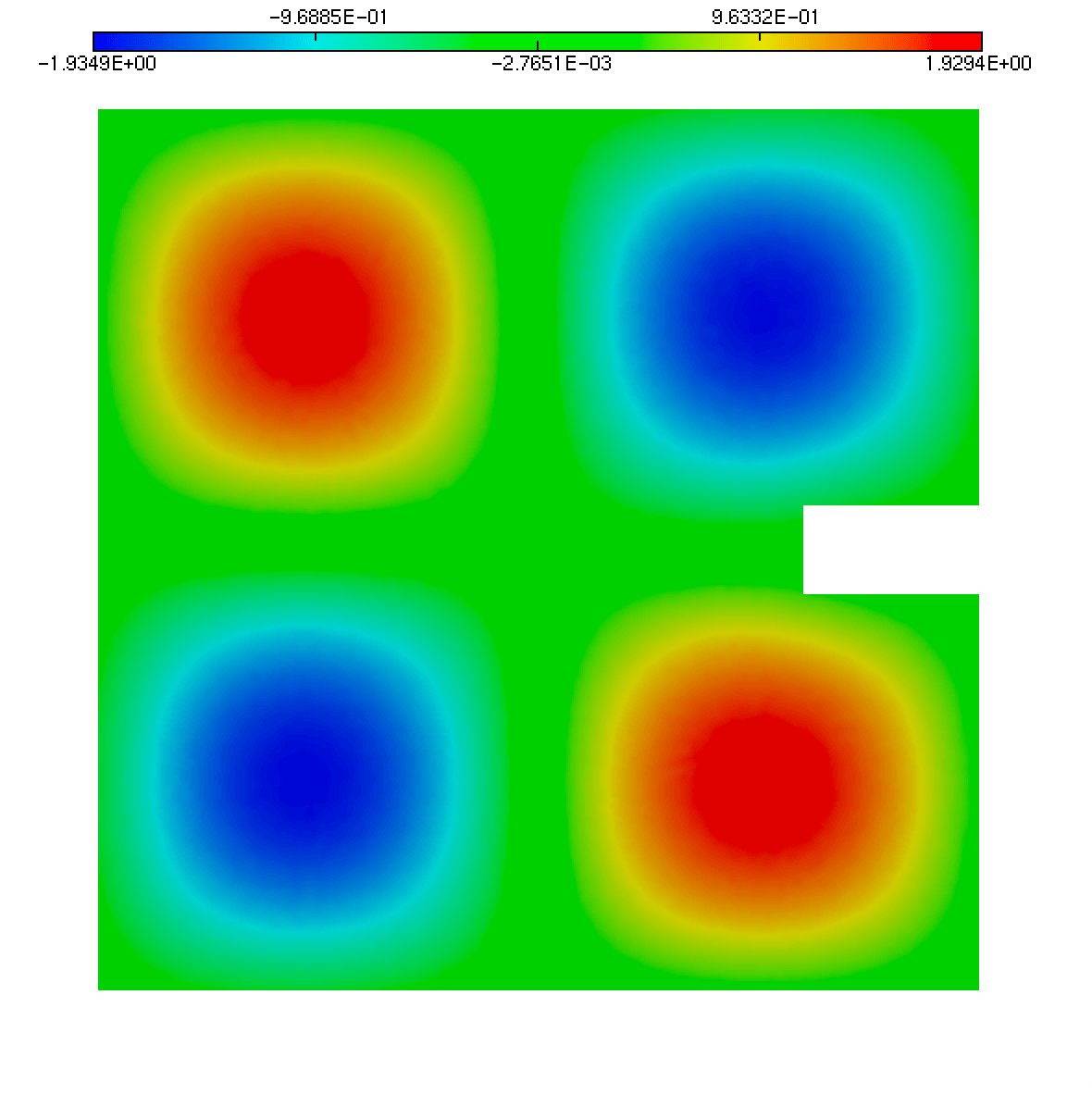} 
\put(2,5){\fcolorbox{black}{white}{$\lambda_{4}=0.7066$}}
\end{overpic}
\end{minipage} & 
\begin{minipage}{0.16\textwidth}
\begin{overpic}[width=\linewidth]{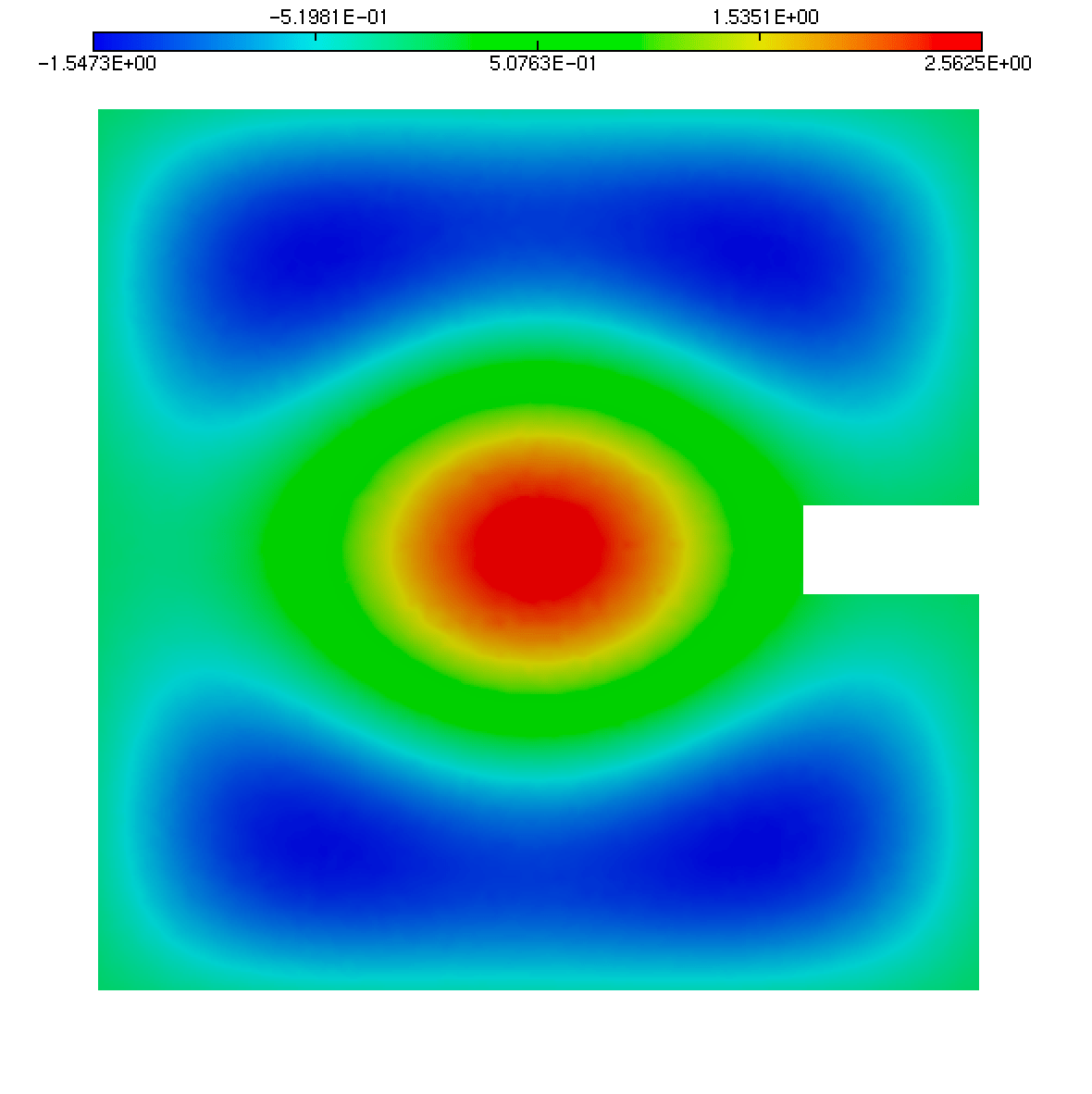} 
\put(2,5){\fcolorbox{black}{white}{$\lambda_{5}=0.6939$}}
\end{overpic}
\end{minipage} & 
\begin{minipage}{0.16\textwidth}
\begin{overpic}[width=\linewidth]{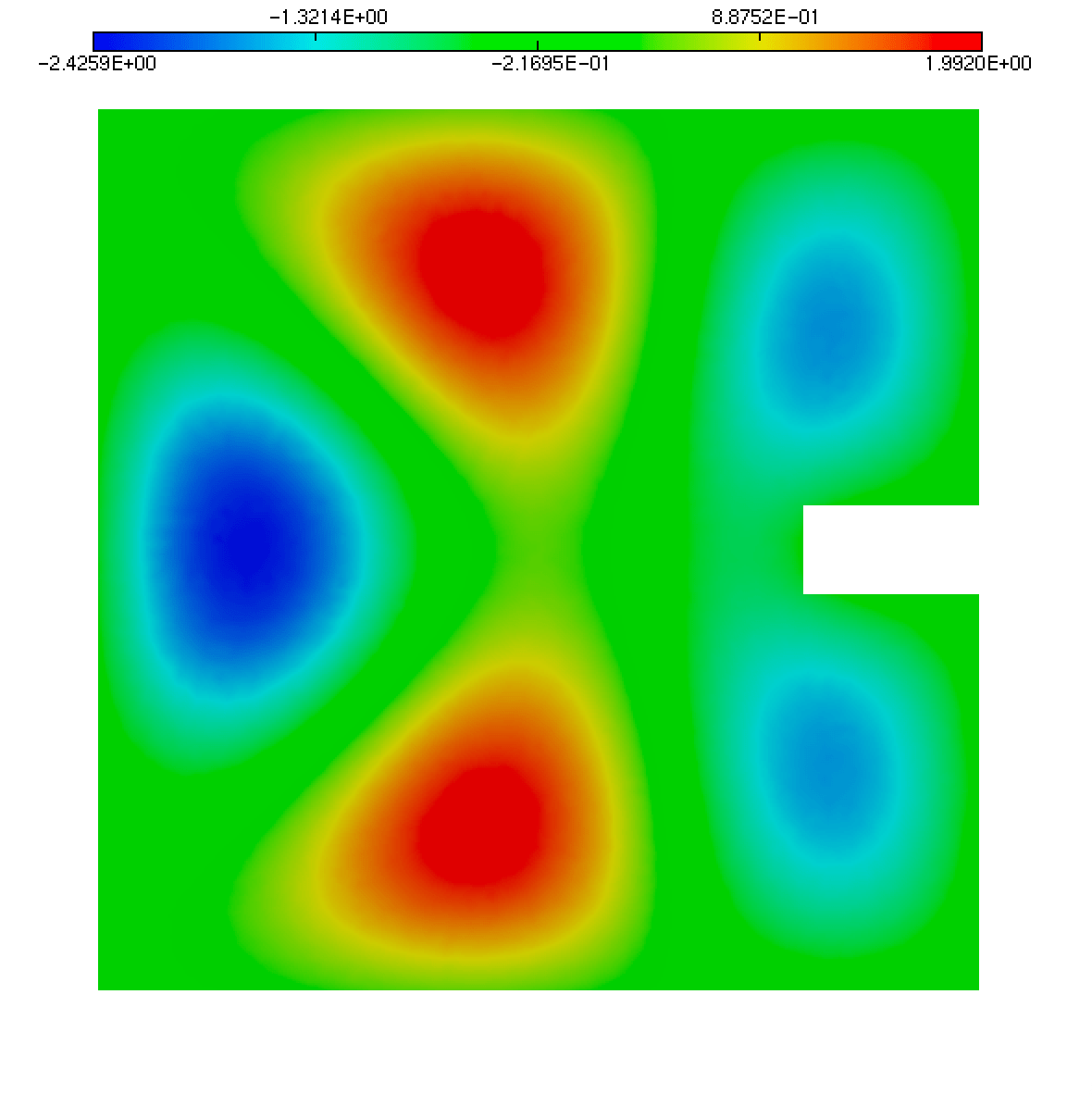} 
\put(2,5){\fcolorbox{black}{white}{$\lambda_{6}=0.6868$}}
\end{overpic}
\end{minipage}
\end{tabular} \par\medskip
\begin{tabular}{ccccc}
    \begin{minipage}{0.19\textwidth}
        \includegraphics[width=\linewidth]{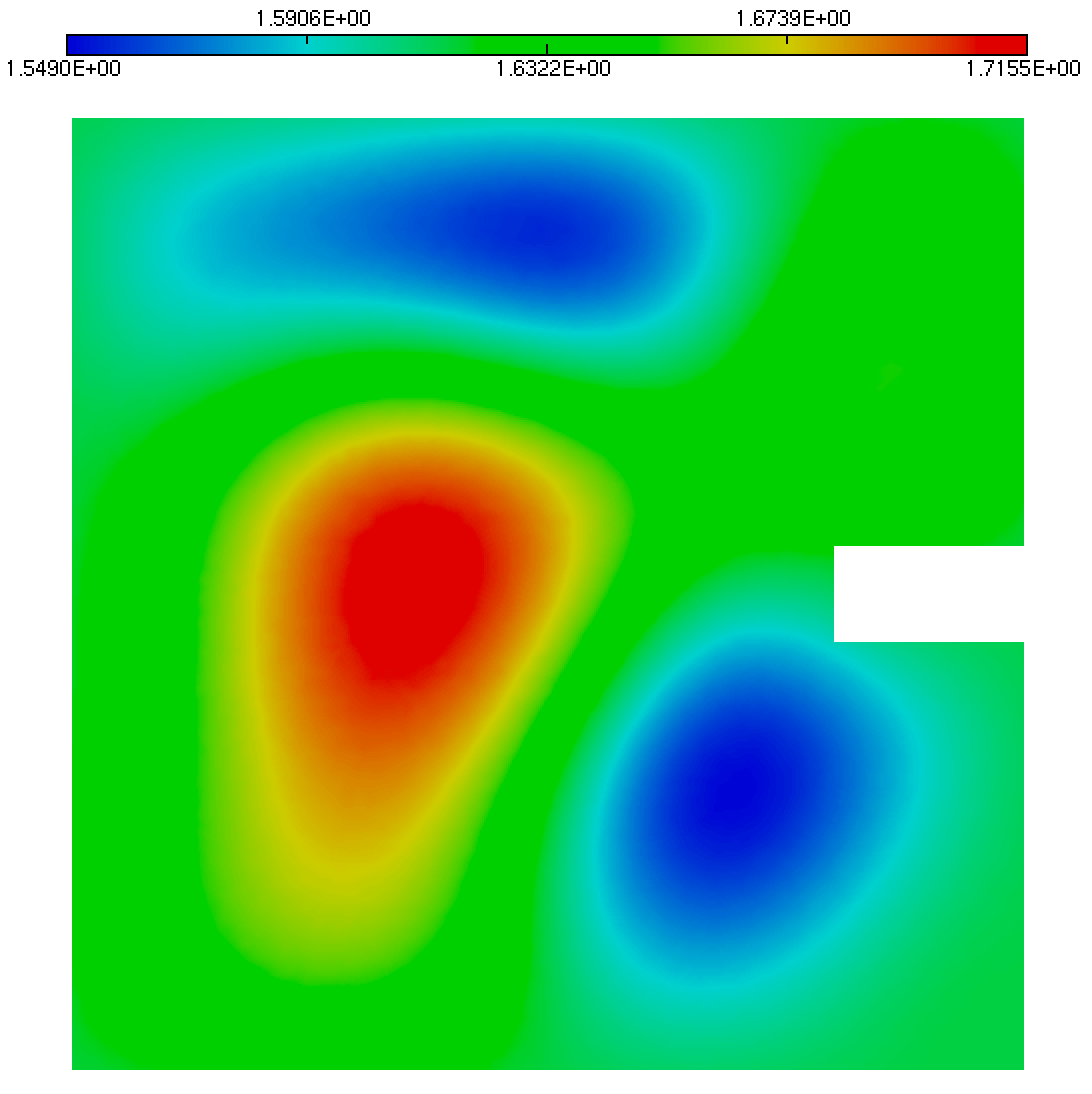}
    \end{minipage}&
    \begin{minipage}{0.19\textwidth}
        \includegraphics[width=\linewidth]{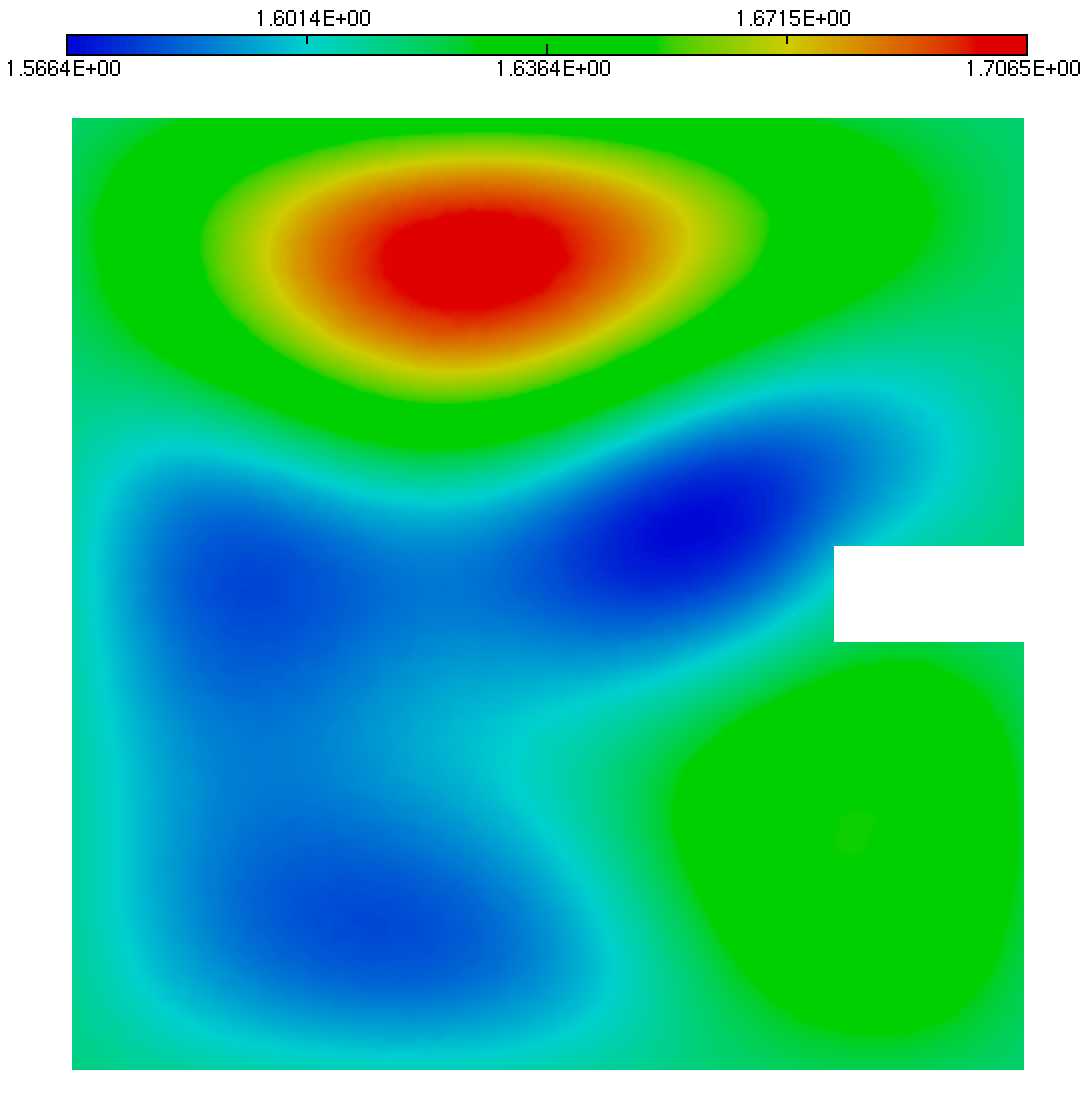}
    \end{minipage}&
    \begin{minipage}{0.19\textwidth}
        \includegraphics[width=\linewidth]{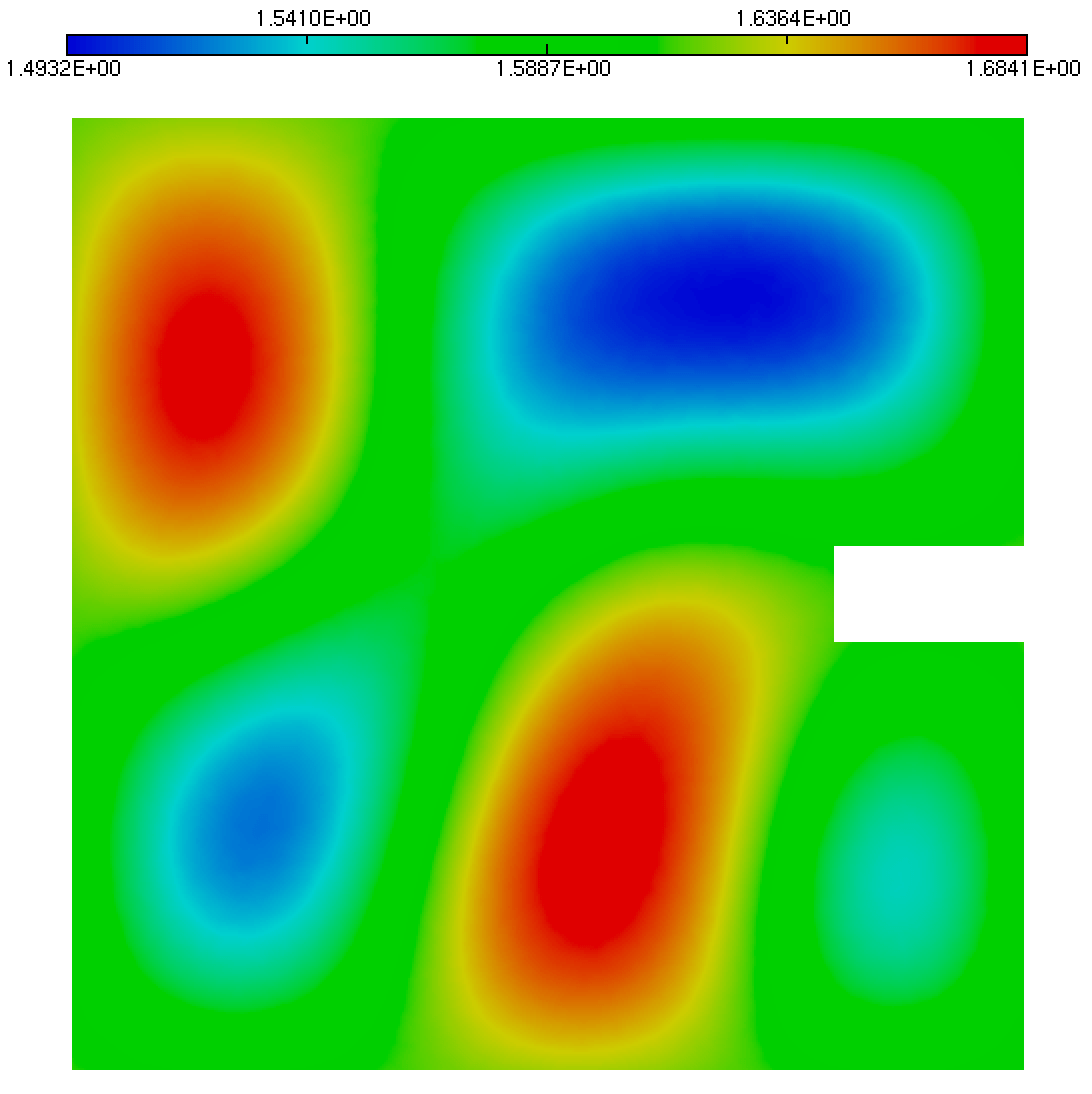}
    \end{minipage}&
    \begin{minipage}{0.19\textwidth}
        \includegraphics[width=\linewidth]{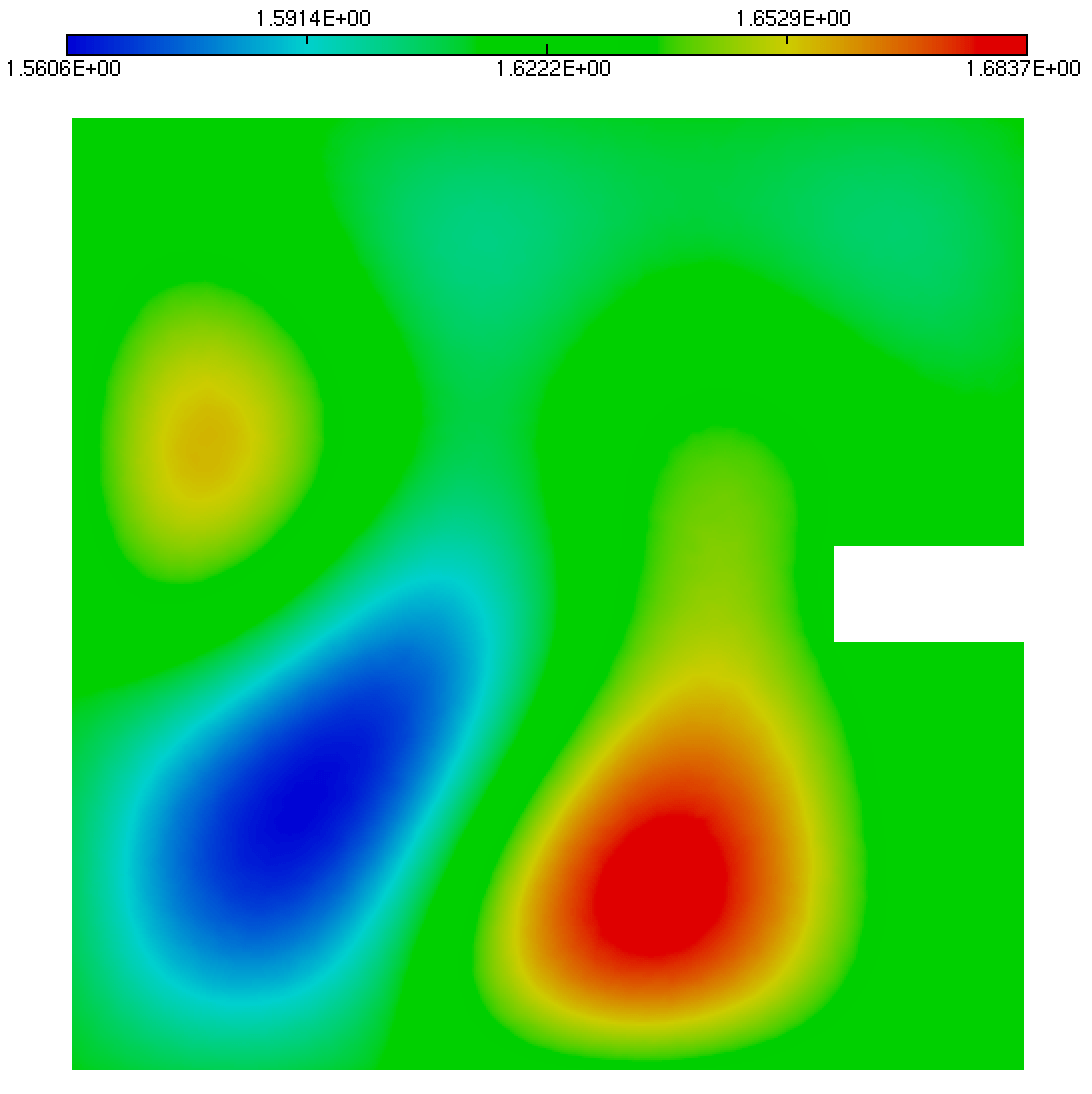}
    \end{minipage}&
    \begin{minipage}{0.19\textwidth}
        \includegraphics[width=\linewidth]{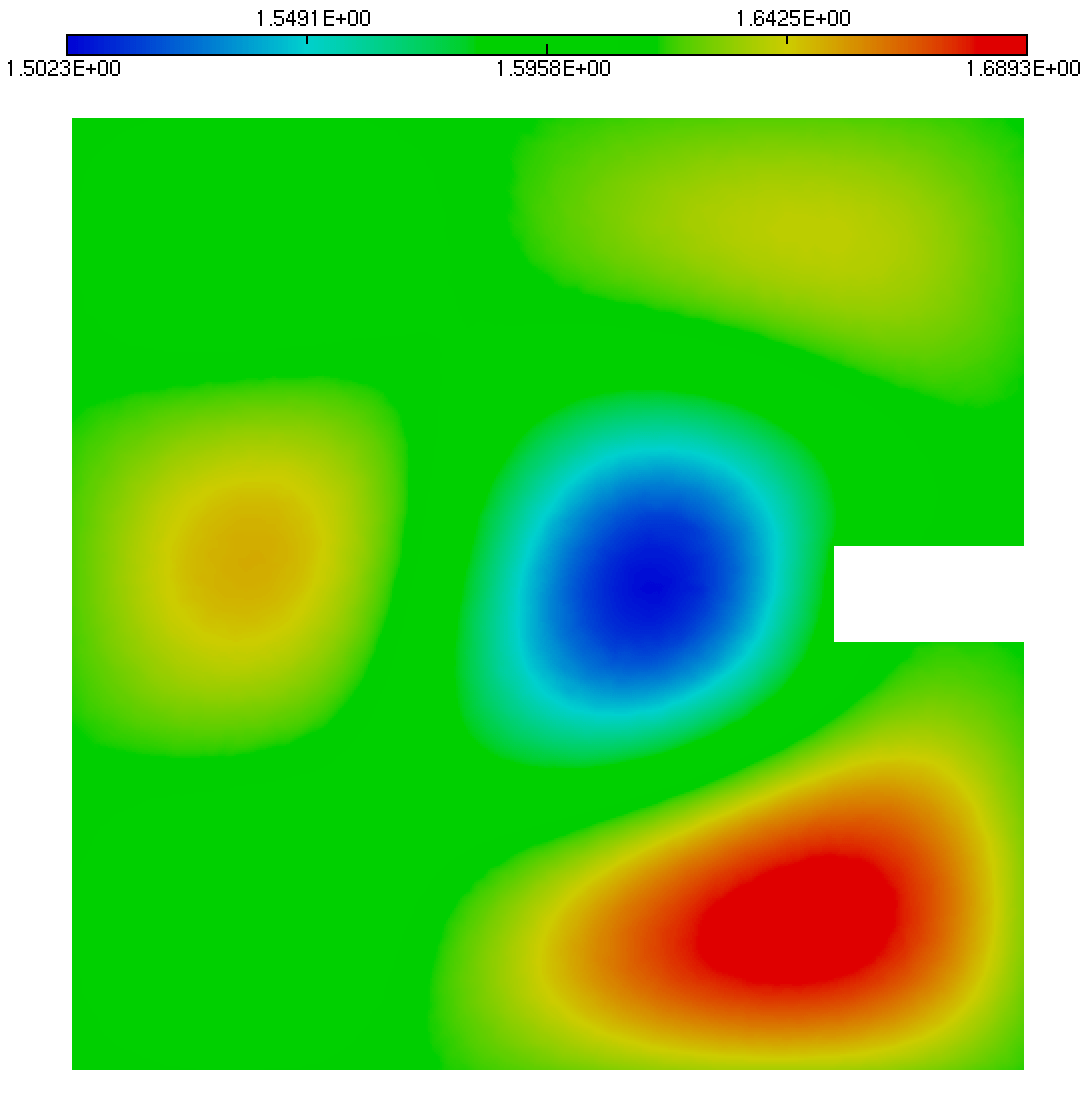}
    \end{minipage}
\end{tabular}
\caption{\it Distributionally robust optimization under material uncertainties: (Upper row) First six eigenpairs $(E_i,\lambda_i)$ of the covariance kernel \cref{eq.Cov}; (Lower row) Realizations of the random field $E(x,\xi)$ in \cref{eq.Exxi}, associated to five different instances of the parameter $\xi$.}
    \label{fig.KLev}
\end{figure}

\begin{figure}[ht]
\begin{tabular}{ccc}
		\begin{minipage}{0.33\textwidth}
		\begin{overpic}[width=1.0\textwidth]{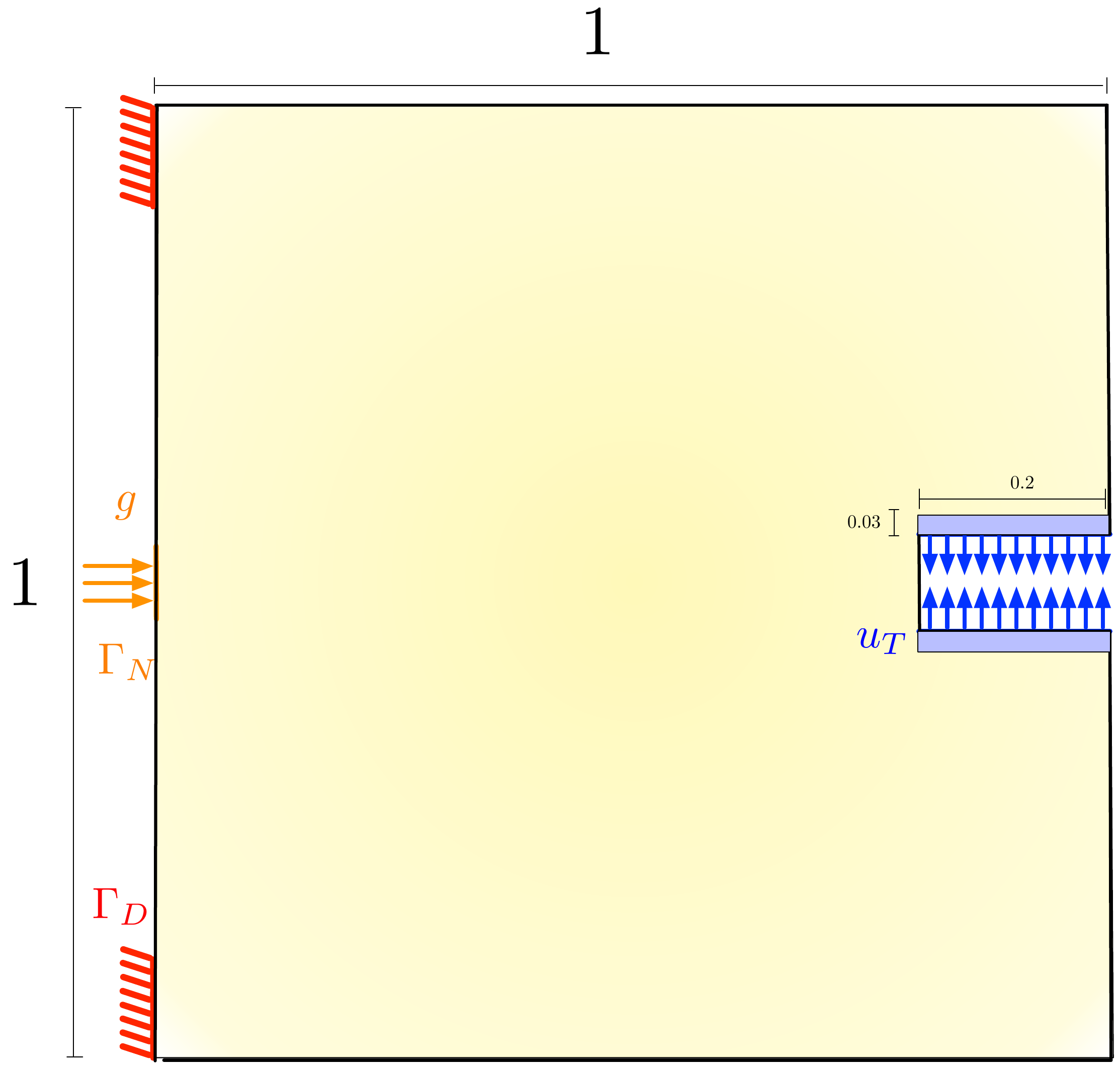} 
		\put(2,5){\fcolorbox{black}{white}{a}}
		\end{overpic}
		\end{minipage} &
		\begin{minipage}{0.33\textwidth}
		\centering
		\begin{overpic}[width=0.8\textwidth]{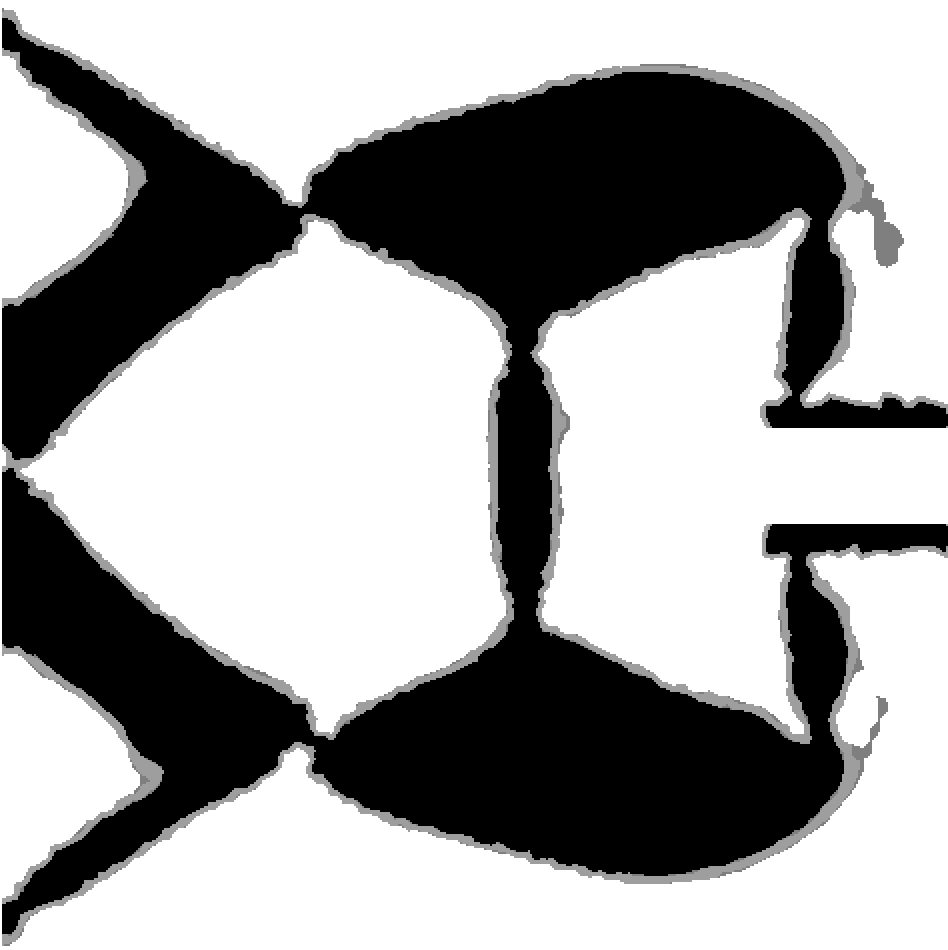} 
		\put(2,5){\fcolorbox{black}{white}{b}}
		\end{overpic}
		\end{minipage} & 
		\begin{minipage}{0.35\textwidth}
		\begin{overpic}[width=1.0\textwidth]{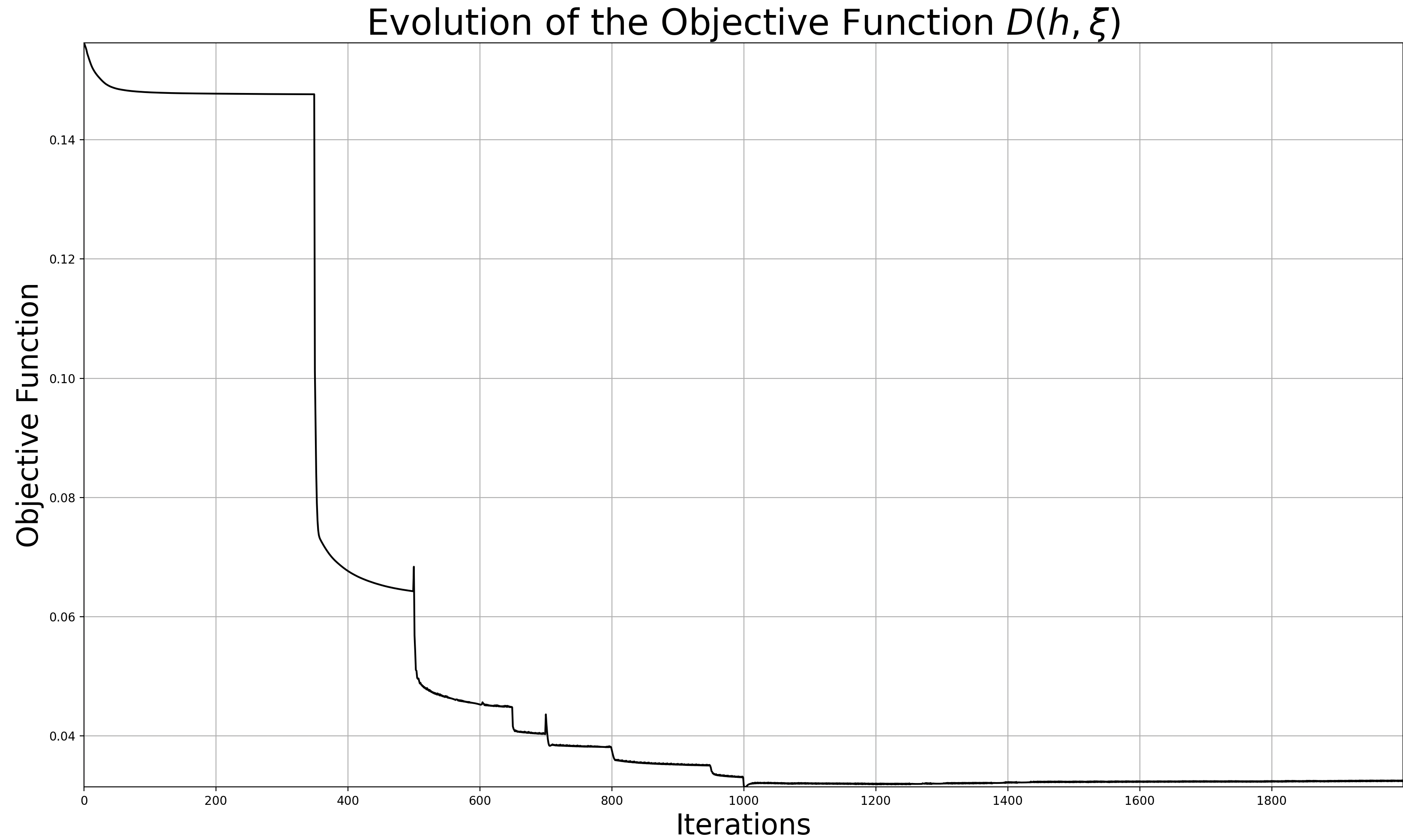} 
		\put(2,5){\fcolorbox{black}{white}{c}}
		\end{overpic}
		\end{minipage}
	\end{tabular}
    \caption{\it (a) Setting of the optimal design example of a gripping mechanism considered in \cref{subsec.matdro}; (b) Optimized design in the ideal situation; (c) Convergence history.}
	\label{fig.gripnom}
\end{figure}

The uncertainties over the physical model concern the (inhomogeneous) Young's modulus $E$ of the constituent elastic material. Following \cite{grigoriu1998simulation,lazarov2012topology}, this function is assumed to be of the form:
\begin{equation}\label{eq.Exxi}
E(x,\xi) = G(\widetilde{E}(x,\xi)), \text{ where }\widetilde
	{E}(x,\xi) = E_{0}(x) + \sum\limits_{i=1}^{k} \sqrt{\lambda_{i}}E_{i}(
	x) \xi_{i}, \quad x \in D, \:\: \xi \in \Xi. 
\end{equation}
Here,
\begin{itemize}
\item The $E_i(x)$ are given, deterministic elements in $L^2(D)$; 
\item The uncertain parameter $\xi$ lies in a compact subset $\Xi$ of $\R^k$. In the nominal situation, its law $\P$ is the (restriction to $\Xi$ of the) normalized $k$-variate Gaussian law. 
\item The function $G$ equals $G:=F^{-1}\circ \Phi$, where $\Phi$ is the cumulative distribution function of the $k$-variate Gaussian law and $F$ is the cumulative distribution function of each real-valued random variable $E(x,\cdot)$, $x \in D$. In practice, $F$ is the cumulative distribution function of the uniform law on $[0.1,1.9] \subset \R$; in particular, $E(x,\xi)$ is uniformly bounded away from $0$ and $\infty$.  
\end{itemize}
Let us briefly explain how the structure \cref{eq.Exxi} of $\widetilde{E}$ shows up in practice, see e.g. \cite{chen2010level,guest2008structural,lazarov2012topology}. The quantity $\widetilde E$ is modeled as a random field $\widetilde E(x,\omega)$, where the event variable $\omega$ belongs to the abstract probability space $(\calO,\calF,\mu)$: intuitively, a random variable $\widetilde E(x,\cdot)$ is attached at each spatial position $x\in D$. The available information about $\widetilde E(x,\omega)$ concerns its mean value $E_0(x)$ and covariance kernel $\Cov(x,y)$, defined by:
$$\forall x, y \in D, \quad E_{0}(x) = \int_{\calO} \widetilde E(x,\omega) \:\d\mu(\omega), \text{ and }  \Cov(x,y) = \int_{\mathcal{O}}(\widetilde E(x,\omega) - E_{0}(x))(\widetilde E(y,\omega)-E_{0}(y))\;\mathrm{d}\mu(\omega).$$
Here, we assume that $E_0(x) \equiv 1$ and that $\Cov(x,y)$ is given by:
\begin{equation}\label{eq.Cov}
\forall x,y\in D, \quad \Cov(x,y)=\alpha e^{-\frac{\lvert x-y \lvert}{\lcor}}, \text{ where }\lcor=2e{-2} \text{ and } \alpha=100.
\end{equation}
The finite-dimensional structure of $\widetilde E(x,\xi)$ in \cref{eq.Exxi} now follows from a Karhunen-Lo\`eve decomposition of the covariance kernel $\Cov(x,y)$, truncated after the $k^{\text{th}}$ term, see e.g. \cite{loeve2017probability,schwab2006karhunen}: the pairs $(\lambda_i,E_i)$ are the eigenelements of the operator
$$T: L^{2}(D) \to L^{2}(D), \quad T\varphi(x) = \int_{D}\Cov(x,y) \varphi(y) \:\d y.$$
The first six eigenfunctions of \cref{eq.Cov} are displayed on \cref{fig.KLev}, together with a few realizations of $E(x,\xi)$.
  
\begin{remark}
In practice, the covariance kernel $\Cov(x,y)$ of the random field $\widetilde{E}(x,\omega)$ is itself unknown and it has to be reconstructed from (a few) observed instances of $E(x,\cdot)$ (or $\widetilde{E}(x,\cdot)$). Hence, in principle, one should take into account the fact that the basis functions $E_i(x)$ in the decomposition \cref{eq.Exxi} are themselves uncertain. This topic of distributionally robust covariance estimators, which goes beyond the scope of this article, is discussed in \cite{dehghannasiri2018intrinsically,nguyen2022distributionally,yue2024geometric}. 
\end{remark}

We denote by $u_{h,\xi}$ the solution to the linear
elasticity system \cref{eq.elassimp} when the Young's modulus featured
in the Hooke's law $A \equiv A(h,\xi)$ in \cref{eq.Hookelaw} is that $E(x,\xi)$ given in
\cref{eq.Exxi}, see \cref{eq.mulambdaE}.
In this situation, we aim to optimize the design $h$ of the structure so that it adopts a prescribed displacement $u_T$ on a region $\Gamma_T$ located near the jaws of the mechanism. To achieve this, we consider the following problem: 
\begin{equation} \label{eq.target_material}
\min\limits_{h} D(h,\xi), \text{ where }
D(h,\xi) = \int_{D}\chi(x) \lvert u_{h,\xi}-u_{T} \lvert ^{2}\;\mathrm{d}x
\end{equation}
and $\chi(x)$ is a given weight factor, which equals $1$ in a thin region around the jaws (depicted in blue in \cref{fig.gripnom} (a)) and $0$ elsewhere.

The derivative of the least-square functional $D(h,\xi)$ is given by the following proposition, whose classical proof is omitted for brevity, see \cite{bendsoe2013topology}.

\begin{proposition}\label{prop.Dhprime}
The function $D(h,\xi)$ is differentiable with respect to the design variable $h$, and its derivative reads: 
$$D^{\prime}(h,\xi)(\hhat) = \int_{D}A^\prime(h)(\hhat)e(u_{h,\xi}) : e(p_{h,\xi}) \:\d x,$$ 
where the adjoint state $p_{h,\xi}$ is the unique $H^1(D)^2$ solution to the following boundary-value problem:
\begin{equation} \label{eq.adjTarg}\left\{
	\begin{array}{cl}
		-\dv(A(h,\xi)e(p_{h,\xi})) = -2\chi(x) (u_{h,\xi} - u_T) & \text{in } D,                                        \\
		p_{h,\xi} = 0                                               & \text{on } \Gamma_D,                                 \\
		A(h,\xi)e(p_{h,\xi})n = 0                                   & \text{on } \partial D \setminus \overline{\Gamma_D}.
	\end{array}
	\right.
\end{equation}
\end{proposition}

We first consider the optimization of the mechanism in the ideal situation where the Young's modulus is known perfectly, i.e. we solve \cref{eq.target_material} with $\xi=0$, and thus $E(x,\xi)= E_0(x)$. The computation requires 2,000 iterations of our optimization algorithm, for a CPU time of about 3h 15 mn. The resulting design and the convergence history are reported on \cref{fig.gripnom} (b) and (c).

\begin{remark}\label{rem.penvolgrip}
In practice, we use small penalizations of the objective function $D(h,\xi)$ by the volume $\Vol(h)$ and the compliance $C(h,\xi)$ of the structure. 
The former helps removing unnecessary material islands disconnected from the main structure, while the latter enforces the connectedness of the structure,
in spite of its trend to develop thin hinges, mimicking pointwise junctions, see e.g. \cite{zhu2020design} about this common practice in the design of compliant mechanisms.
\end{remark}

We now consider a situation where the Young's modulus $E \equiv E(x,\xi)$ is unknown: its expression \cref{eq.Exxi} stems from the Karhunen-Lo\`eve expansion of the covariance kernel $\Cov(x,y)$, after truncation at $k=10$ modes. We thus solve a distributionally robust version of the optimization problem \cref{eq.target_material}, assuming a Wasserstein ambiguity set  $\calA_{\text{W}}$ based on the nominal law $\P = \delta_0$, see \cref{sec.wdro}. The problem at hand reads:
\begin{equation} \label{eq.primal_DRO_material}
\min_{h \in \Uad}\;
	\left( \sup\limits_{\Q \in \calA_{\text{W}}}\; \int_{\Xi} D(h,\xi) \:\d \Q(\xi) \right).
\end{equation}
The parameter $\e$ is set to $\e=0.01$.
Invoking \cref{prop.formulaJdir}, \cref{eq.primal_DRO_material} rewrites under the more
convenient form:
\begin{equation*}\label{eq.Dual_DRO_Material}
\min\limits_{h \in \Uad, \atop \lambda \geq 0}\calD_{\text{W}}(h,\lambda), \text{ where } \calD_{\text{W}}
	(h,\lambda) = \lambda m + \lambda\e \int_{\Xi}\log\left( \int_{\Xi}e^{\frac{D(h,\zeta)
	- c(\xi,\zeta)}{\lambda\e}}\:\d\nu_{\xi}(\zeta) \right )\:\d\P(\xi),
\end{equation*}
where the probability measure $\nu_\xi$ in \cref{eq.refcouplingW} has variance  $\sigma^2= 1e{-1}$. 

\begin{figure}[ht]
	\centering
	\begin{tabular}{ccc}
		\begin{minipage}{0.3\textwidth}\begin{overpic}[width=\linewidth]{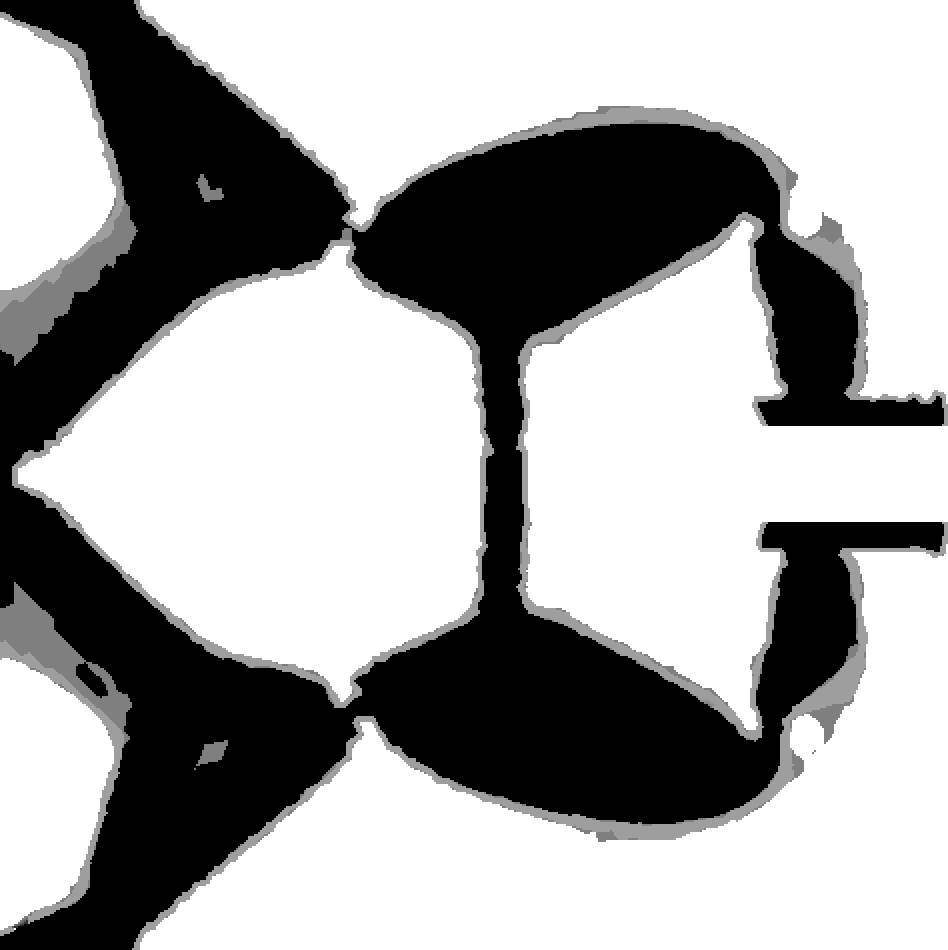} \put(2,5){\fcolorbox{black}{white}{$m = 0$}}\end{overpic}\end{minipage} & 
		\begin{minipage}{0.3\textwidth}\begin{overpic}[width=\linewidth]{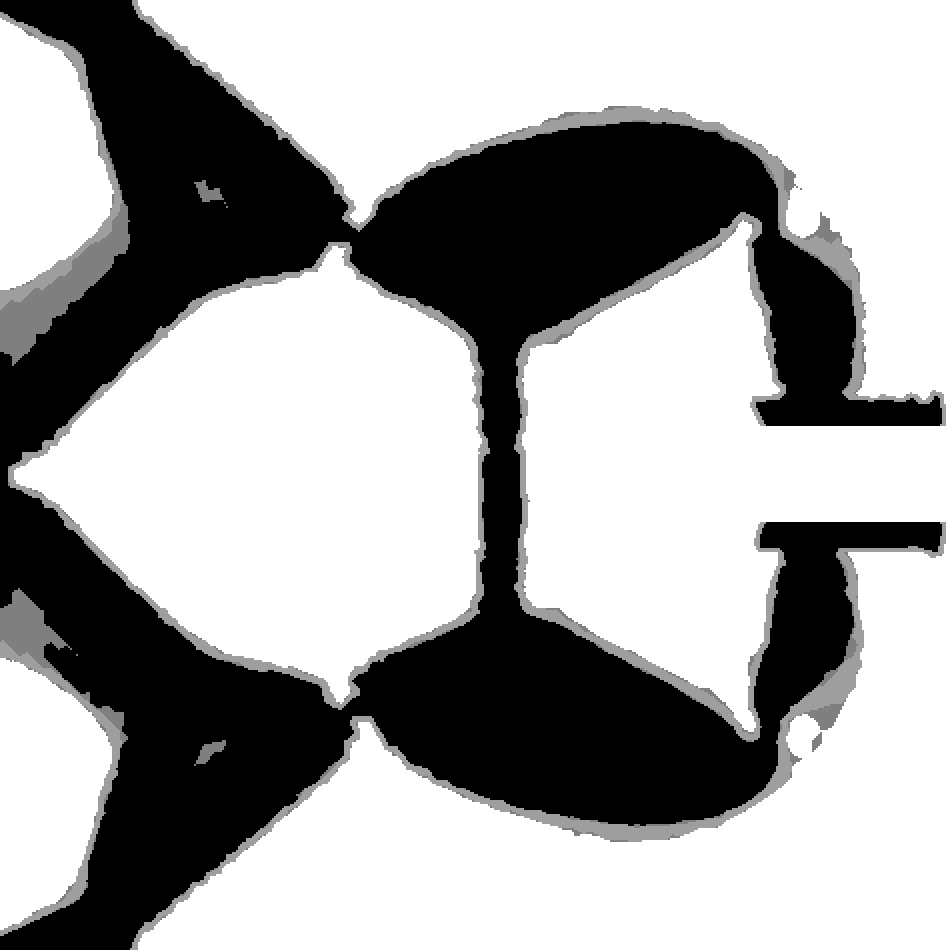} \put(2,5){\fcolorbox{black}{white}{$m = 0.4$}}\end{overpic}\end{minipage} &
		\begin{minipage}{0.3\textwidth}\begin{overpic}[width=\linewidth]{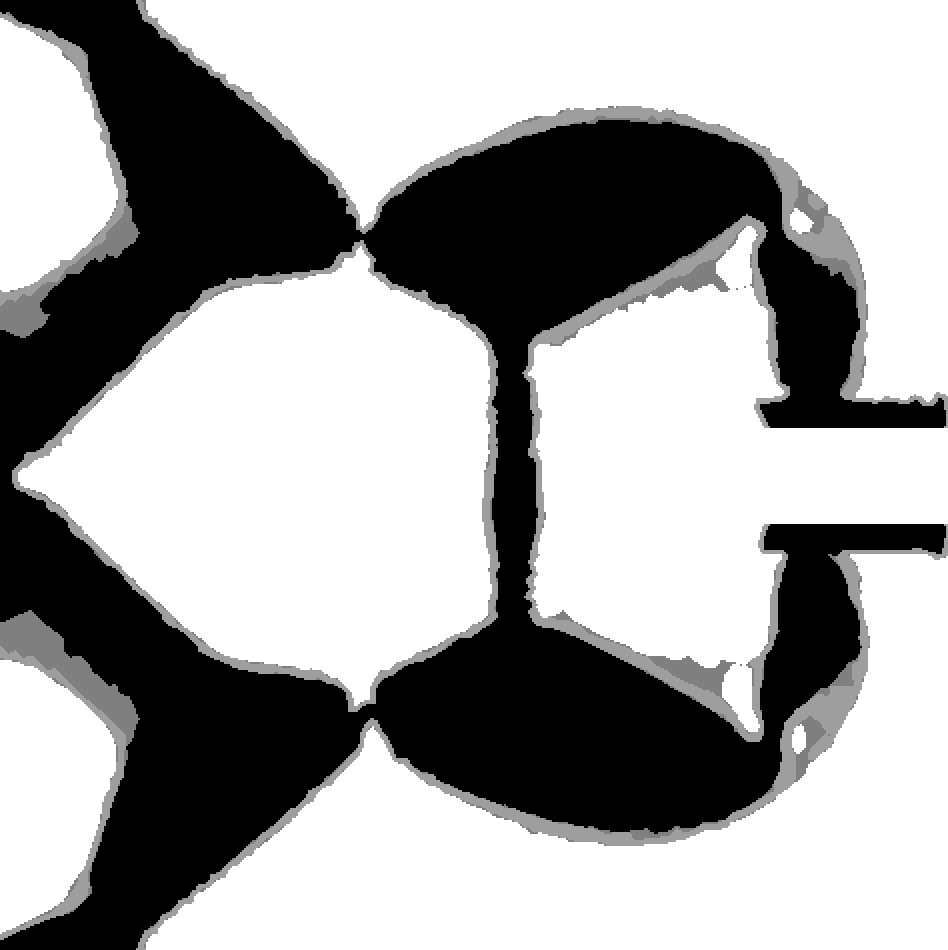} \put(2,5){\fcolorbox{black}{white}{$m = 0.5$}}\end{overpic}\end{minipage}
	\end{tabular} \par\medskip
	\begin{tabular}{ccc}
		\begin{minipage}{0.3\textwidth}\begin{overpic}[width=\linewidth]{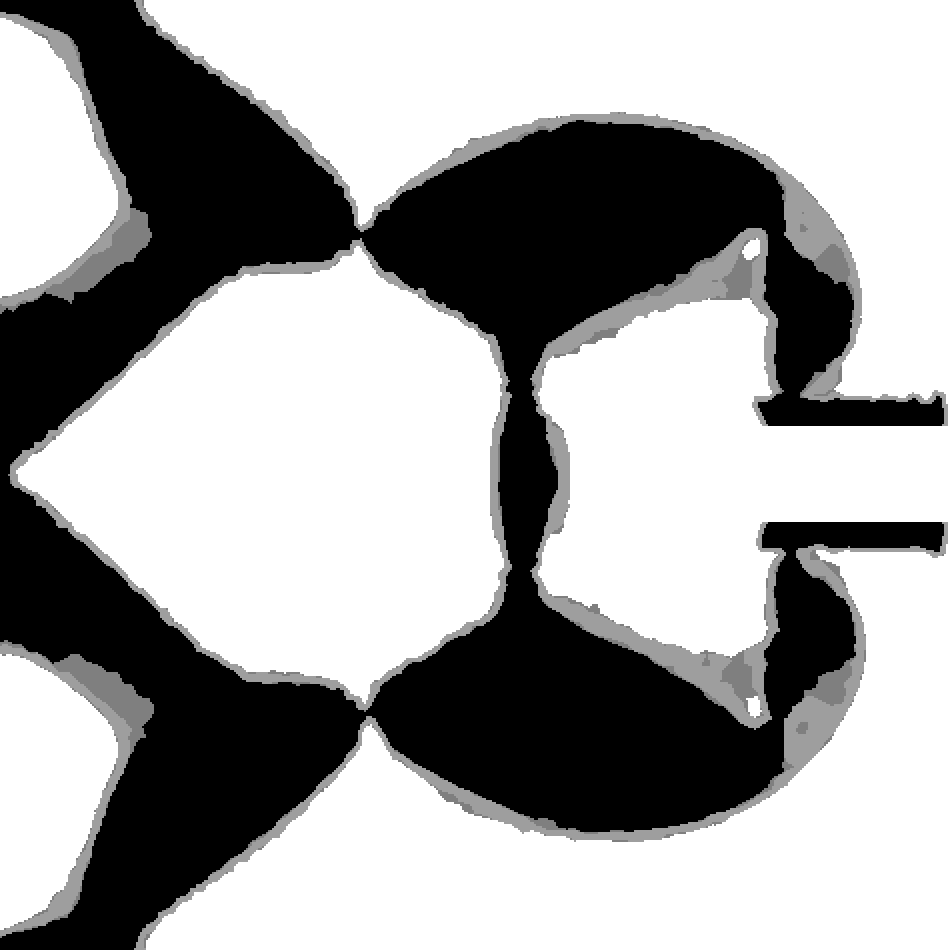} \put(2,5){\fcolorbox{black}{white}{$m = 1$}}\end{overpic}\end{minipage} &
		\begin{minipage}{0.3\textwidth}\begin{overpic}[width=\linewidth]{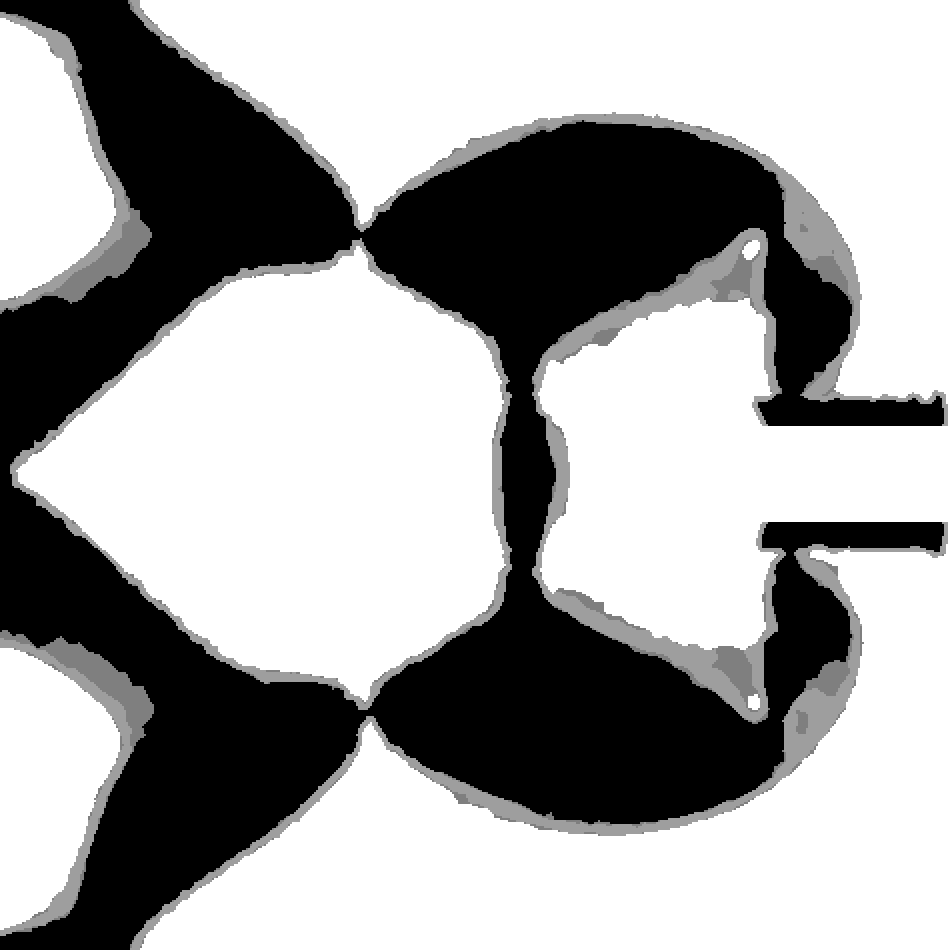} \put(2,5){\fcolorbox{black}{white}{$m = 2$}}\end{overpic}\end{minipage} &
		\begin{minipage}{0.3\textwidth}\begin{overpic}[width=\linewidth]{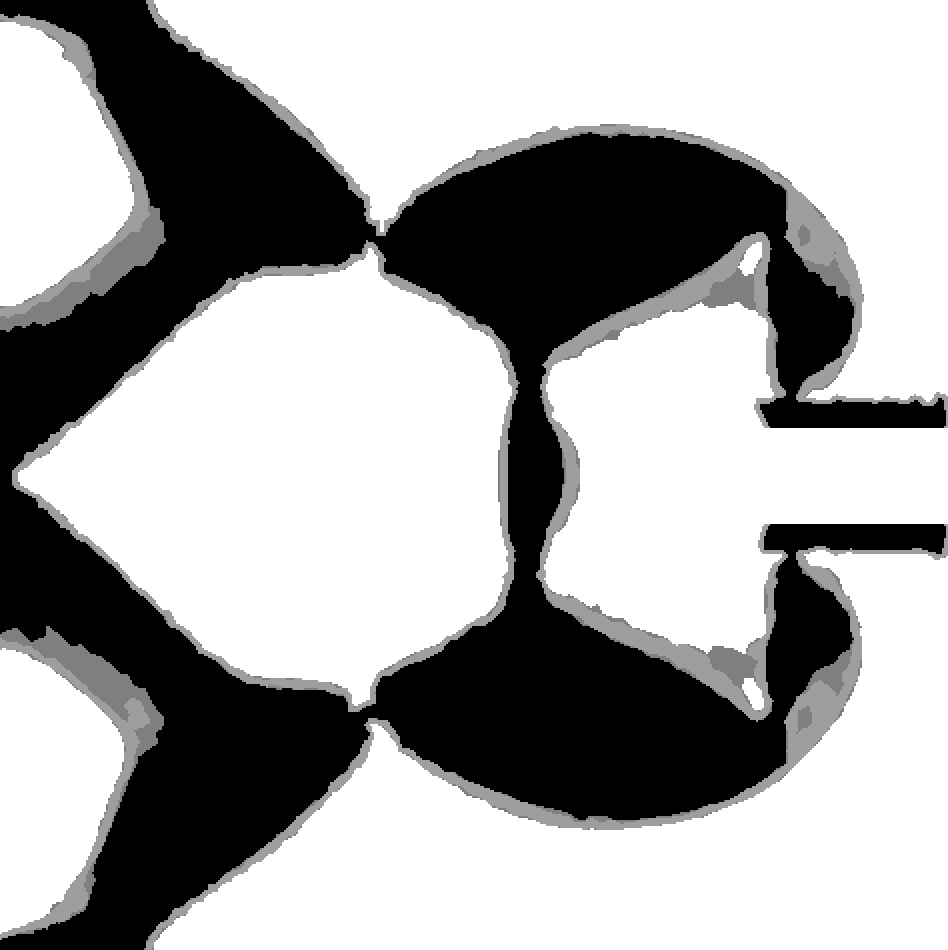} \put(2,5){\fcolorbox{black}{white}{$m = 5$}}\end{overpic}\end{minipage}
	\end{tabular}
	\caption{\it Distributionally robust optimized designs for material uncertainties of the gripping mechanism considered in \cref{subsec.matdro} associated to different values of the Wasserstein radius $m$.}
    \label{fig.optgriprob0.5_0.1}
\end{figure}

We solve this problem for several values of the radius $m$ of the Wasserstein ball; the results are reported on \cref{fig.optgriprob0.5_0.1}; each computation requires 2,000 iterations of our numerical algorithm, for a CPU time of about 13h. 
The resulting designs are quite similar from one another; note however the change in the orientation of the hinges sustaining the jaws between the cases $m=0.5$ and $m=1$.
Let us also emphasize that it is particularly difficult to achieve a ``true'' black-and-white design in this example: in spite of the penalization of the objective function by the volume hinted at in \cref{rem.penvolgrip}, the designs eliminate intermediate densities at a very slow rate.

\subsection{Reliability-based optimal design of a bridge using the conditional value at risk}\label{sec.TOcvardro}

\noindent This section illustrates the ideas of \cref{sec.relbasedOD}, regarding optimal design problems where the conditional value at risk of the cost function is used as a means to account for a failure probability. Before dealing with distributional robustness, we develop a little about this use of conditional value at risk, as it is relatively confidential in optimal design, to the best of our knowledge.

\subsubsection{Deterministic optimal design of a bridge}\label{subsec.cvarideal}

\begin{figure}[ht]
\centering
	\begin{tabular}{cc}
		\begin{minipage}{0.31\textwidth}
		\begin{overpic}[width=0.9\textwidth]{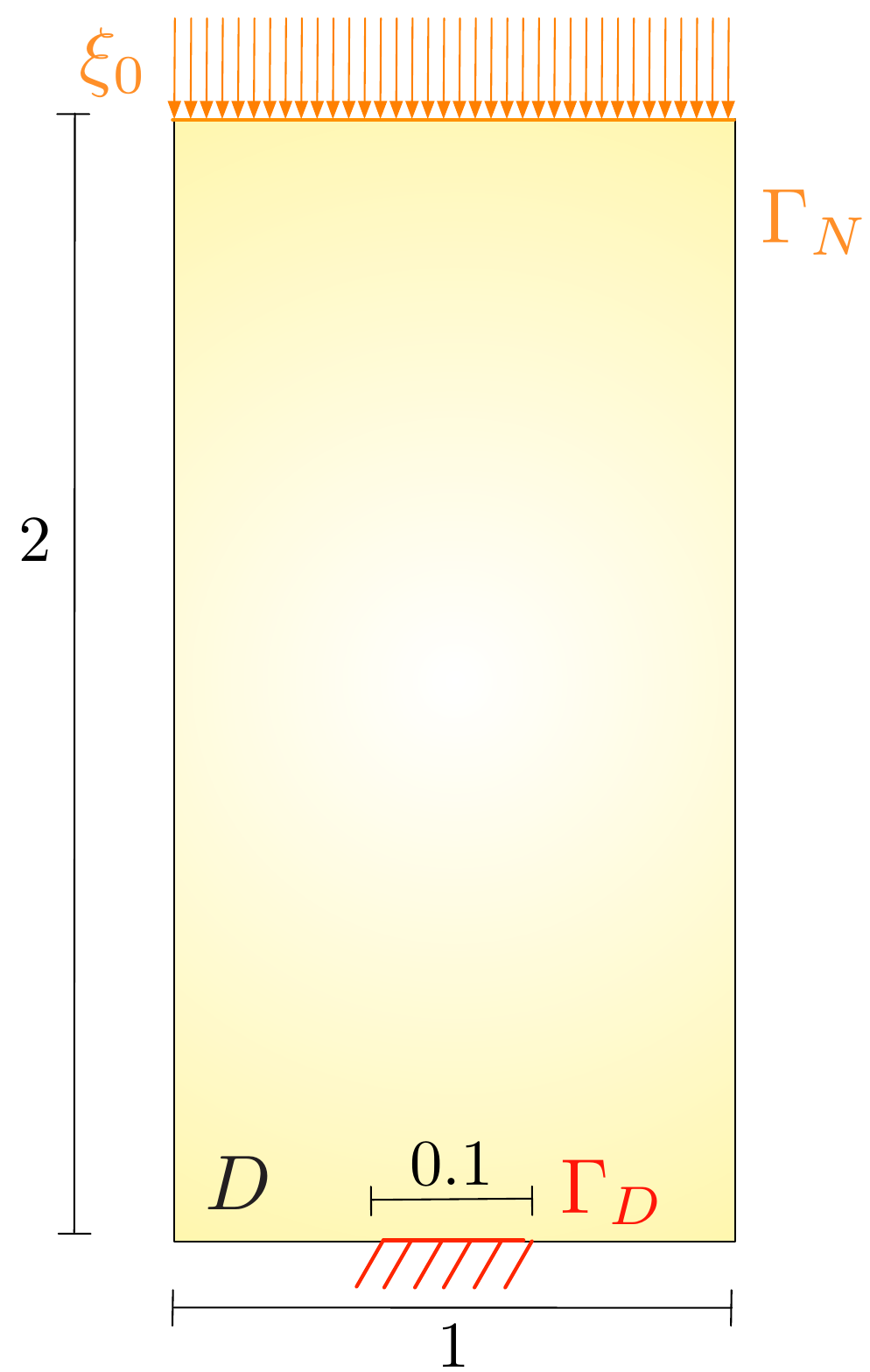} 
		\put(2,5){\fcolorbox{black}{white}{a}}
		\end{overpic}
		\end{minipage} 
		&
        \begin{minipage}{0.31\textwidth}
        \begin{overpic}[width=0.6\textwidth]{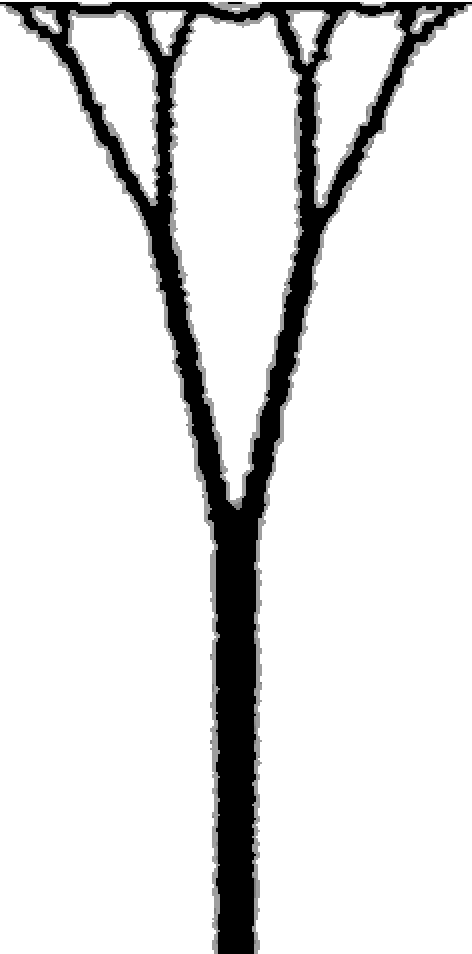}
        \put(2,5){\fcolorbox{black}{white}{b}}
        \end{overpic}            
        \end{minipage}
	\end{tabular}
	\caption{\it (a) Setting of the bridge problem considered in \cref{sec.TOcvardro}; (b) Optimized design $h^*_{\text{\rm det}}$ in ideal conditions.}
	\label{fig.complbridge}
\end{figure}

\noindent The structures under consideration are 2d bridges, contained in a rectangular domain $D$ with size $1 \times 2$, meshed with approximately $8,500$ vertices ($\approx 16,000$ triangles). They are clamped along the lower side $\Gamma_{D}$ of $\partial D$, and a load is applied on their upper boundary $\Gamma_{N}$, as illustrated in \cref{fig.complbridge} (a). In ideal conditions, this load equals $\xi^0 := (0,-1)$. 

We first optimize the compliance $C(h,\xi^0)$ of the structure, defined by \cref{eq.defcplyTO}, under a volume constraint, in the ideal situation where the applied load is exactly $\xi^0$, i.e. we solve: 
\begin{equation}\label{eq.unperturbedcvar}
    \min\limits_{h \in \Uad} \, C(h,\xi^0)\text{ s.t. } \Vol(h)=V_T, \text{ where } V_T = 0.245.
\end{equation}
The resulting design $h^*_{\text{det}}$ is depicted on \cref{fig.complbridge} (b); its compliance equals $C(h^*_{\text{det}},\xi^0)=40.3892$.

\subsubsection{Enforcement of a safety constraint via the conditional value at risk}\label{subsec.cvar}

\noindent We now assume that the load $\xi$ applied on $\Gamma_N$ is uncertain, but that its law $\P \in \calP(\Xi)$ is known perfectly: it is (the restriction to a large ball $\Xi \subset \R^2$ of) the bivariate Gaussian law centered at $\xi^0$, with variance $\sigma^2$. 
We then wish to solve the following problem: 
$$ \min\limits_{h \in \Uad} \Vol(h) \:\text{ s.t. } \:\P \Big\{\xi \in \Xi \text{ s.t. } C(h,\xi) \geq C_T \Big\} \leq 1-\beta,$$
where the threshold $C_T := 40$ is close to the optimal value $C(h^*_{\text{det}},\xi^0)$ attained in the unperturbed minimization problem \cref{eq.unperturbedcvar}. 
According to the discussion in \cref{sec.cvardro}, we replace this problem by the following conservative version:
\begin{equation}\label{eq.cvarnodro}
\min\limits_{h \in \Uad, \atop \alpha \in \mathbb{R}} \; \Vol(h) \quad \text{s.t.} \quad \alpha + \frac{1}{1-\beta} \int_{\xi \in \Xi} \left[ C(h, \xi) - \alpha \right]_+ \; \d\P(\xi) \leq C_T.
\end{equation}

\begin{figure}[ht]
	\centering
	\begin{tabular}{ccccc}
		\begin{minipage}{0.2\textwidth}
		\begin{overpic}[width=0.81\textwidth]{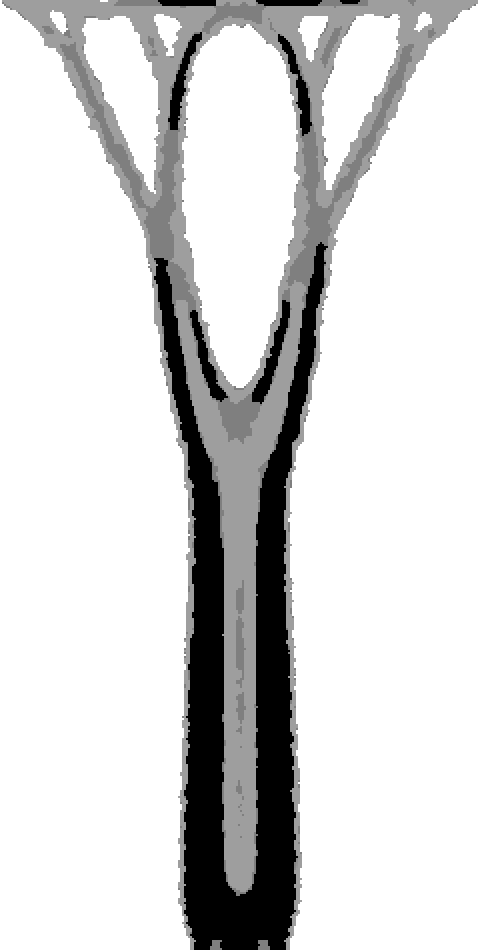 }
		\put(2,5){\fcolorbox{black}{white}{$\beta = 0.01$}}
		\end{overpic}
		\end{minipage} & 
	 \begin{minipage}{0.2\textwidth}
	 \begin{overpic}[width=0.81\textwidth]{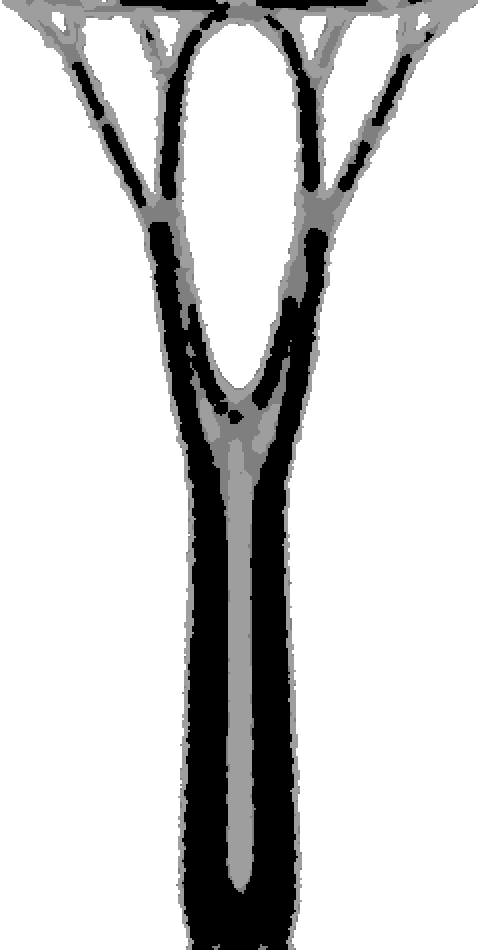 }
	 \put(2,5){\fcolorbox{black}{white}{$\beta = 0.1$}}
	 \end{overpic}
	 \end{minipage} & 
	 \begin{minipage}{0.2\textwidth}
	 \begin{overpic}[width=0.81\textwidth]{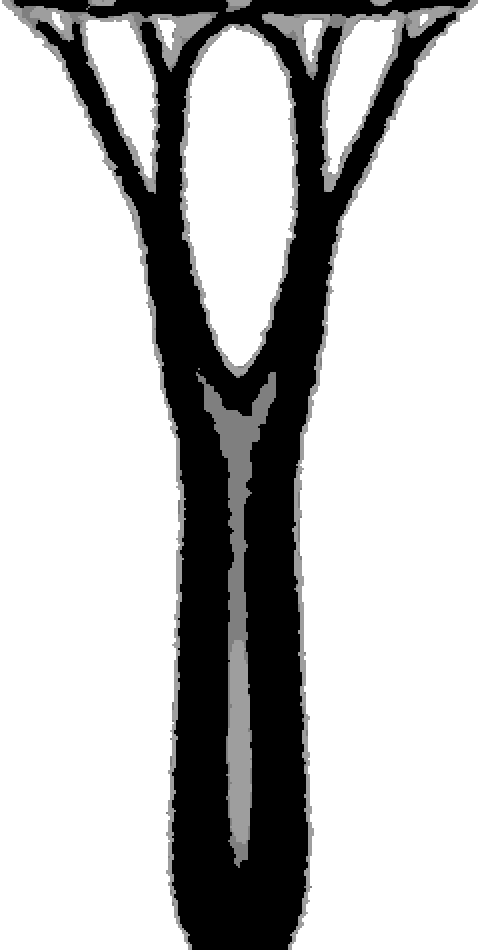} 
	 \put(2,5){\fcolorbox{black}{white}{$\beta = 0.5$}}
	 \end{overpic}
	 \end{minipage} & 
	 \begin{minipage}{0.2\textwidth}
	 \begin{overpic}[width=0.81\textwidth]{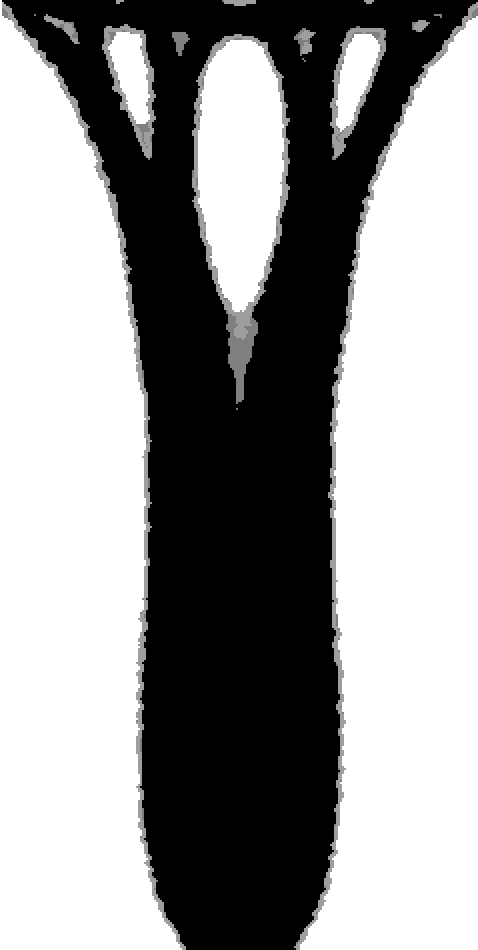}
	 \put(2,5){\fcolorbox{black}{white}{$\beta = 0.9$}}
	 \end{overpic}
	 \end{minipage} & 
	 \begin{minipage}{0.2\textwidth}
	 \begin{overpic}[width=0.81\textwidth]{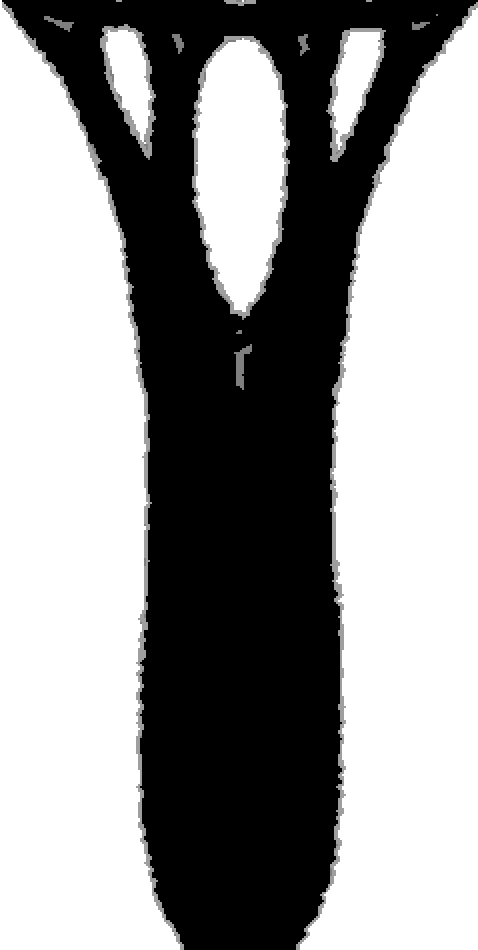} 
	 \put(2,5){\fcolorbox{black}{white}{$\beta = 0.99$}}
	 \end{overpic}
	 \end{minipage}
	\end{tabular} \par\medskip
		\begin{tabular}{ccccc}
		\begin{minipage}{0.2\textwidth}
		\begin{overpic}[width=0.81\textwidth]{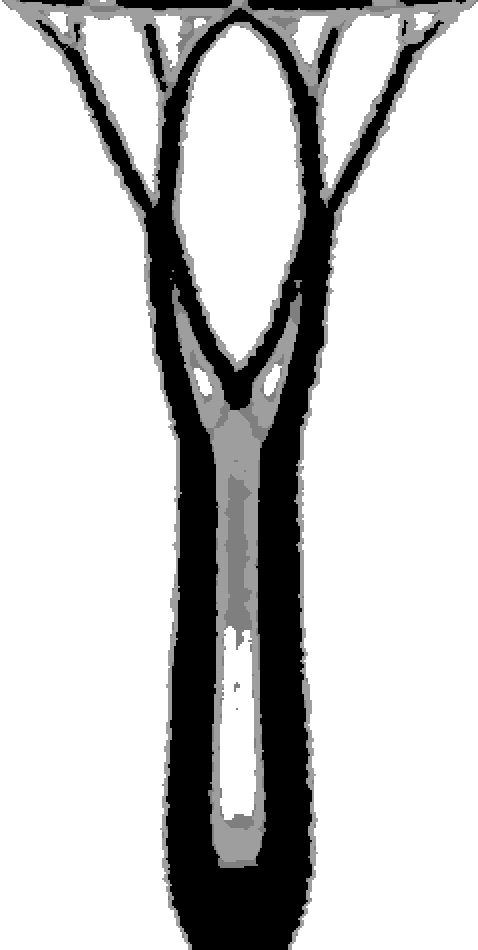} 
		\put(2,5){\fcolorbox{black}{white}{$\beta = 0.01$}}
		\end{overpic}
		\end{minipage} & 
		\begin{minipage}{0.2\textwidth}
		\begin{overpic}[width=0.81\textwidth]{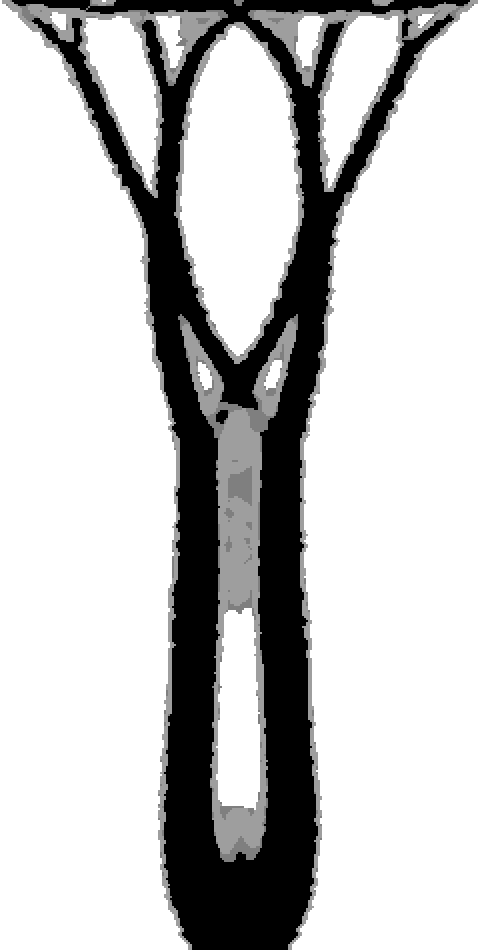} 
		\put(2,5){\fcolorbox{black}{white}{$\beta = 0.1$}}
		\end{overpic}
		\end{minipage} & 
		\begin{minipage}{0.2\textwidth}
		\begin{overpic}[width=0.81\textwidth]{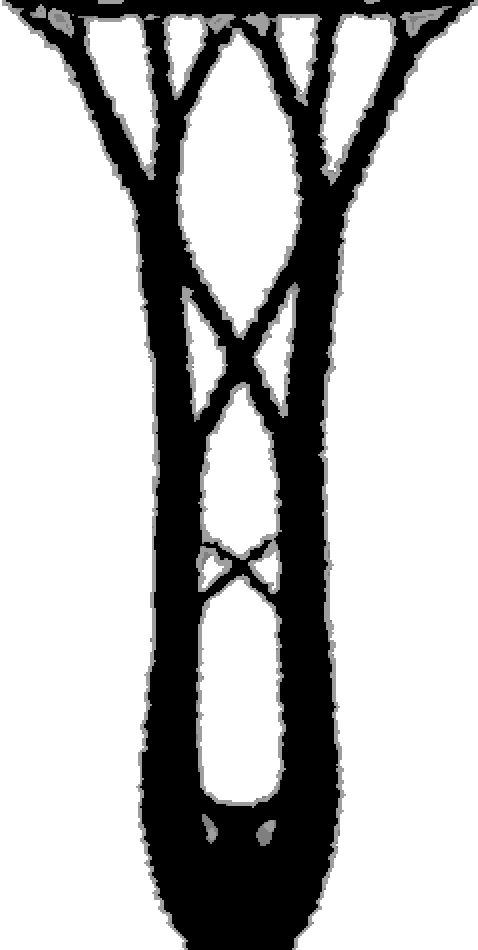} 
		\put(2,5){\fcolorbox{black}{white}{$\beta = 0.5$}}
		\end{overpic}
		\end{minipage} & 
		\begin{minipage}{0.2\textwidth}
		\begin{overpic}[width=0.81\textwidth]{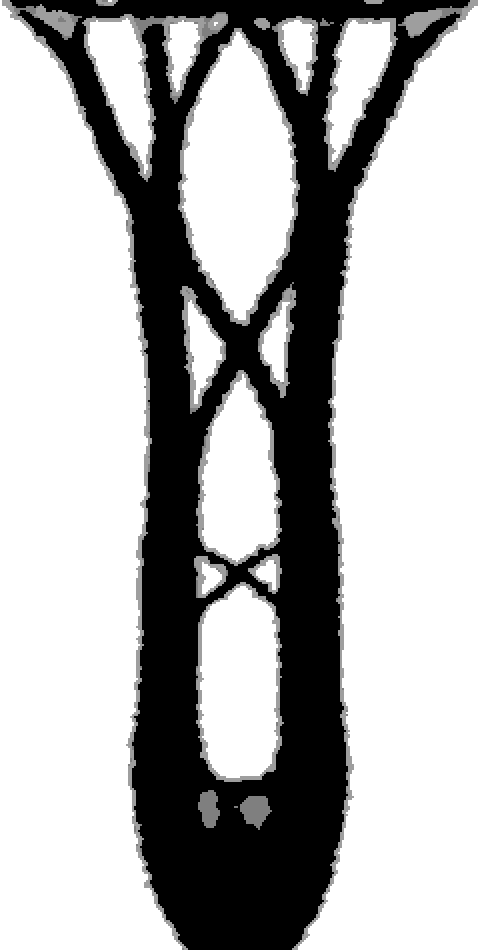 } 
		\put(2,5){\fcolorbox{black}{white}{$\beta = 0.9$}}
		\end{overpic}
		\end{minipage} & 
		\begin{minipage}{0.2\textwidth}
		\begin{overpic}[width=0.81\textwidth]{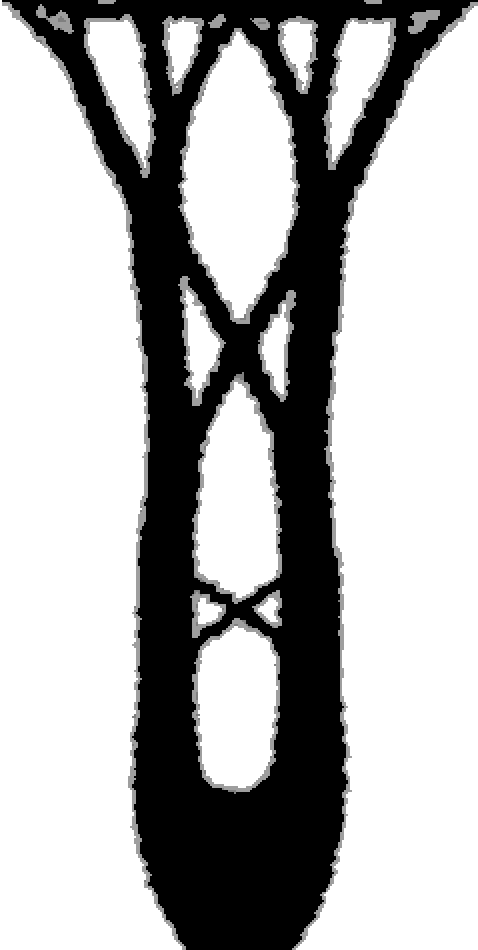} 
		\put(2,5){\fcolorbox{black}{white}{$\beta = 0.99$}}
		\end{overpic}
		\end{minipage}
	\end{tabular}
	\caption{\it Optimized designs of the bridge under a safety constraint accounted for by the conditional value at risk, as considered in \cref{subsec.cvar}; (Upper row) $\sigma^2=1e{-2.5}$; (Lower row) $\sigma^2 = 1e{-2}$.}
    \label{fig.cvarnom}
\end{figure}

We solve this problem for several values of the variance $\sigma^2$ of the law $\P$ for $\xi$ and the threshold $\beta$:
the resulting designs are presented in \cref{fig.cvarnom}, see also \cref{tab.nominalCVaR} for the associated values of the volume.
In each situation, about $600$ iterations of the optimization strategy described in \cref{sec.simpnum} are used, for an approximate CPU cost of 2h 30 mn. Again, $N=10$ samples are used for the numerical computation of each probabilistic integral featured in \cref{eq.cvarnodro}.
Understandably enough, the closer $\beta$ is to $1$, the more stringent the failure constraint, and the larger the resulting volume. Loosely speaking, 
the most significant structural modifications occur in the lower vertical part of the domain, which tends to split into two sections joined by a reinforcement crossbar.
We also note that it is more difficult for the design to converge towards clear black-and-white features (with the exact same numerical strategy), as intermediate densities seem to become more favorable.

\begin{table}[ht]
\centering
\begin{tabular}{|c|c|c|c|c|c|}
\hline
 \backslashbox[20mm]{$\sigma^2$}{$\beta$}  & $0.01$ & $0.1$ & $0.5$ & $0.9$ & $0.99$ \\
\hline
$1e{-2.5}$ & $0.365525$ & $0.398984$ & $0.523988$ & $0.835165$ & $0.840873$ \\
\hline
$1e{-2}$ & $0.438366$ & $ 0.474774$ & $0.553006$ & $0.634271$ & $0.641044$ \\
\hline
\end{tabular}
\caption{\it Volume $\Vol(h)$ of the optimized structures in the bridge example of \cref{subsec.cvar}.}
\label{tab.nominalCVaR}
\end{table}

\subsubsection{Distributionally robust optimal design considering failure probabilities}\label{subsec.cvardrobr}

\noindent We now turn to the situation where the law of the load parameter $\xi$ is itself uncertain, which raises the need for a distributionally robust formulation of the optimal design problem \cref{eq.cvarnodro}. The nominal distribution is made of a single observation, associated to the ideal load $\xi^0 = (0,-1)$: ${\P}=\delta_{\xi^0}$. 
Assuming an ambiguity set of the Wasserstein type as in \cref{sec.wdroformula}, the distributionally robust version of the problem \cref{eq.cvarnodro} reads: 
$$ \min\limits_{h\in \Uad, \atop  \alpha \in \R}\; \Vol(h) \text{ s.t. }  \sup\limits_{\Q \in
	\calA_{\text{W}}}\min\limits_{\alpha \in \R}\left(\alpha + \dfrac{1}{1-\beta}
	\int_{\xi\in \Xi}\left[\mathcal{C}(h,\xi)-\alpha\right]_{+}\;\d
	\Q(\xi)\right) \leq C_T.
$$
According to \cref{sec.cvardro}, this rewrites:
\begin{multline*}
	\min\limits_{h\in \Uad, \atop \lambda \geq 0, \alpha \in \R}\; \Vol(h) \:
	\text{ s.t. }\; \calD_{\text{C}}(h,\lambda,\alpha) \leq C_{T}, \text{
	where }\\ \calD_{\text{C}}(h,\lambda,\alpha) = \alpha+ \frac{\lambda m}{1-\beta}+ \frac{\lambda \e}{1-\beta}\log\left( \int_{\Xi}
	e^{\frac{\left[C(h,\zeta)-\alpha\right]_{+} - \lambda c(\xi^0,\zeta)}{\lambda\e}}
	\:\d\nu_{\xi^0}(\zeta) \right).
\end{multline*}

\begin{figure}[ht]
    \centering
    \begin{tabular}{ccccc}
        \begin{minipage}{0.2\textwidth}
            \begin{overpic}[width=0.81\textwidth]{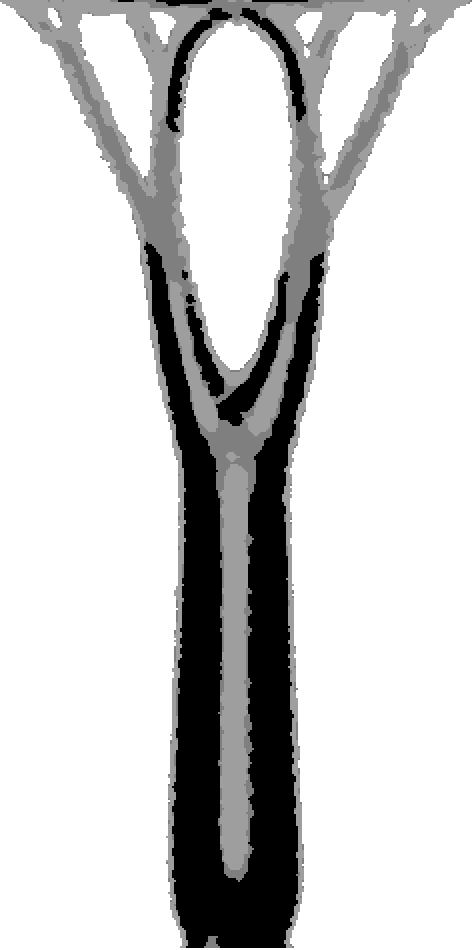}
                \put(2,5){\fcolorbox{black}{white}{$\beta = 0.01$}}
            \end{overpic}
        \end{minipage} 
        & \begin{minipage}{0.2\textwidth}
            \begin{overpic}[width=0.81\textwidth]{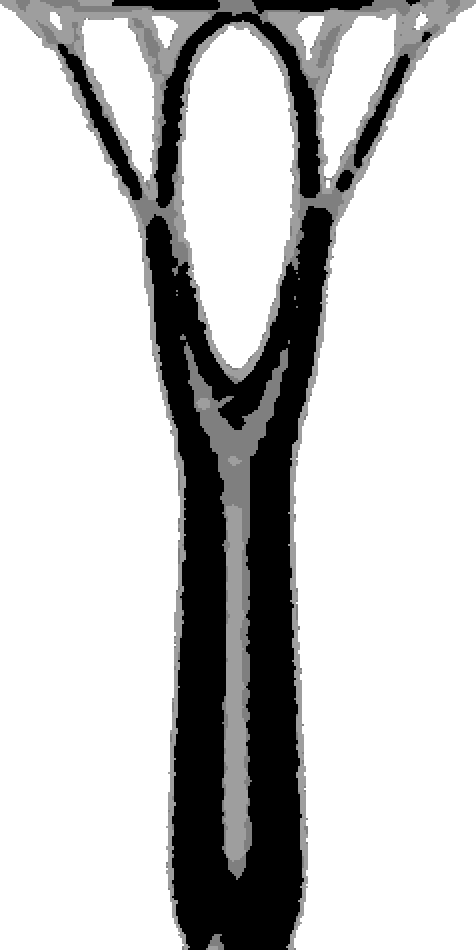}
                \put(2,5){\fcolorbox{black}{white}{$\beta = 0.1$}}
            \end{overpic}
        \end{minipage}
        & \begin{minipage}{0.2\textwidth}
            \begin{overpic}[width=0.81\textwidth]{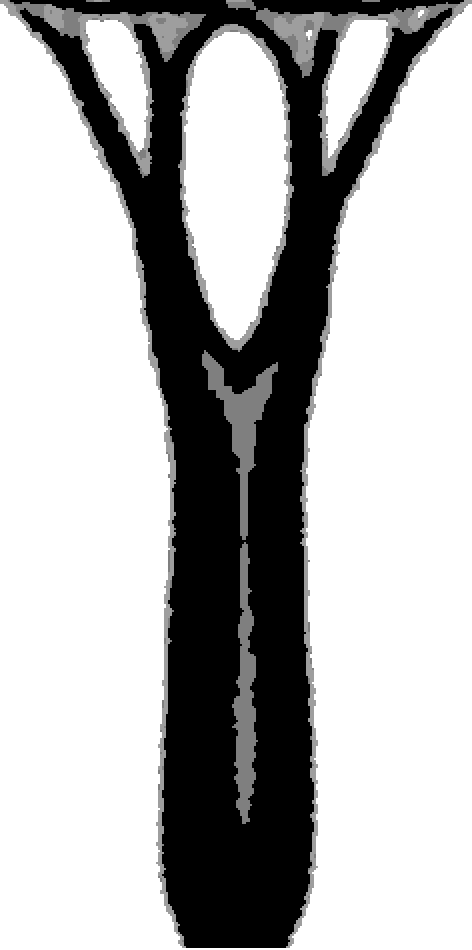}
                \put(2,5){\fcolorbox{black}{white}{$\beta = 0.5$}}
            \end{overpic}
        \end{minipage} 
        & \begin{minipage}{0.2\textwidth}
            \begin{overpic}[width=0.81\textwidth]{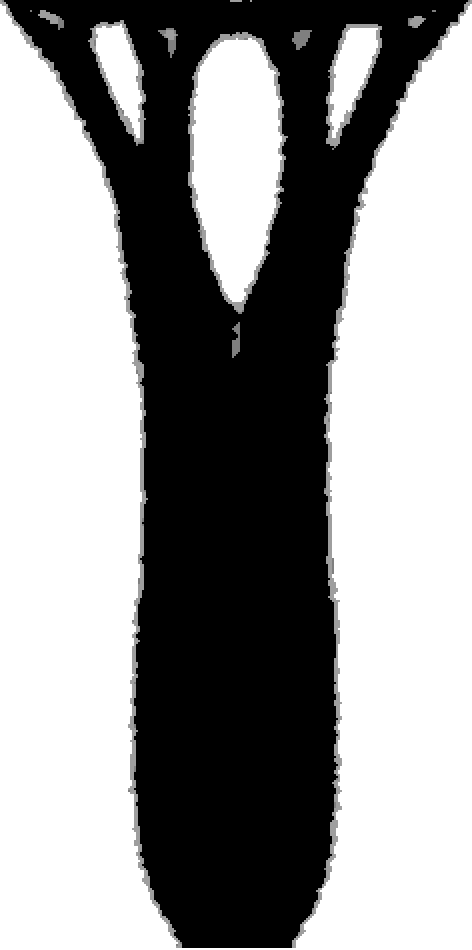}
                \put(2,5){\fcolorbox{black}{white}{$\beta = 0.9$}}
            \end{overpic}
        \end{minipage} 
        & \begin{minipage}{0.2\textwidth}
            \begin{overpic}[width=0.81\textwidth]{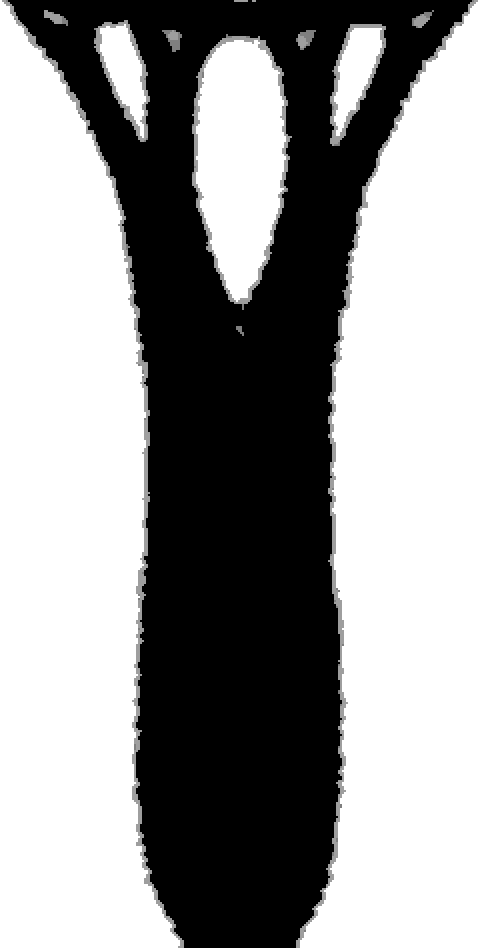}
                \put(2,5){\fcolorbox{black}{white}{$\beta = 0.99$}}
            \end{overpic}
        \end{minipage}
    \end{tabular}
    \caption{\it Distributionally robust designs of the bridge considered in \cref{subsec.cvardrobr} for a Wasserstein radius $m=0.5$, $\sigma^2=1e{-3}$ and various values of the parameter $\beta$.}
    \label{fig.cvardro3_0.5}
\end{figure}

\begin{figure}[ht]
    \centering
    \begin{tabular}{ccccc}
        \begin{minipage}{0.2\textwidth}
            \begin{overpic}[width=0.81\textwidth]{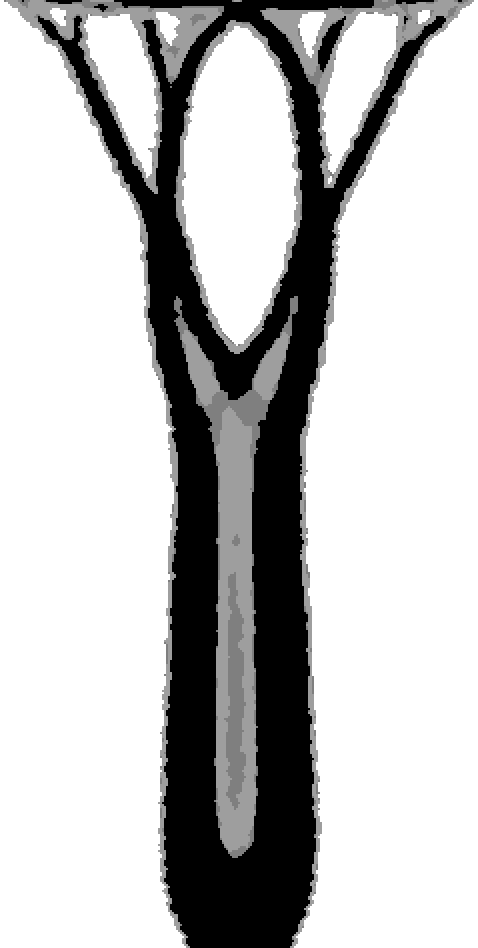}
                \put(2,5){\fcolorbox{black}{white}{$\beta = 0.01$}}
            \end{overpic}
        \end{minipage} 
        & \begin{minipage}{0.2\textwidth}
            \begin{overpic}[width=0.81\textwidth]{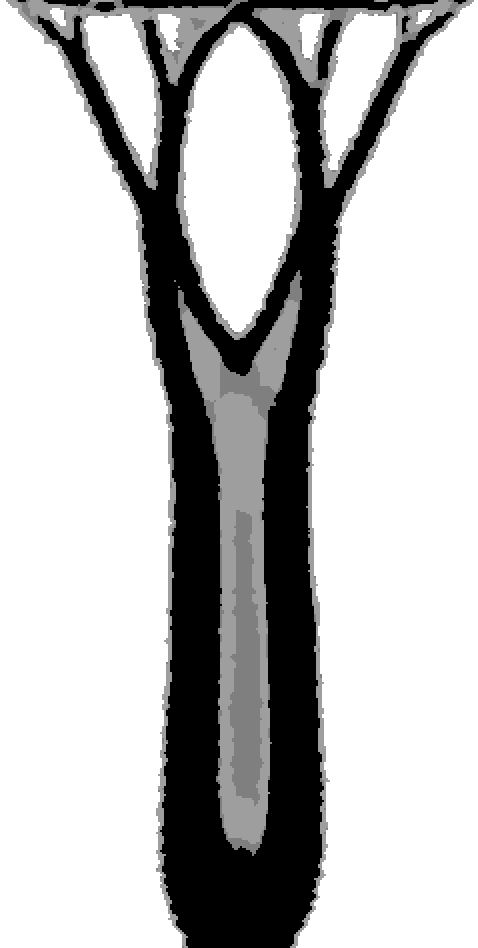}
                \put(2,5){\fcolorbox{black}{white}{$\beta = 0.1$}}
            \end{overpic}
        \end{minipage}
        & \begin{minipage}{0.2\textwidth}
            \begin{overpic}[width=0.81\textwidth]{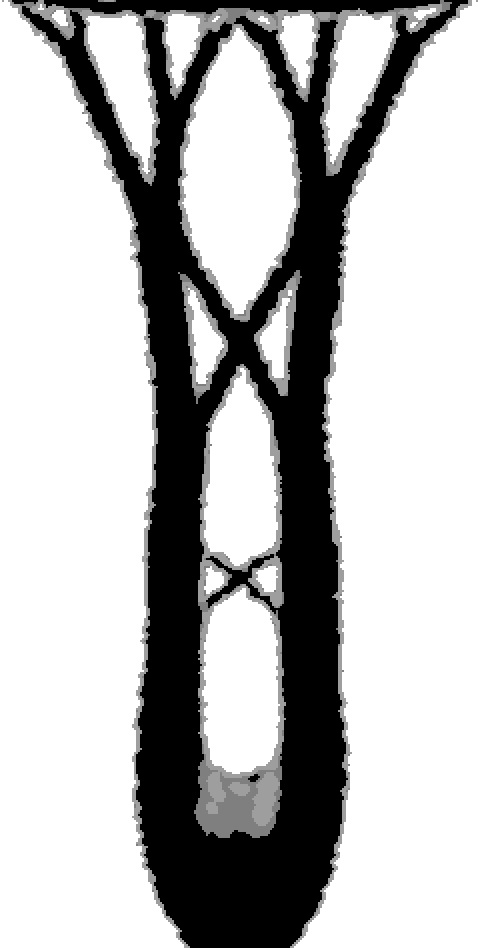}
                \put(2,5){\fcolorbox{black}{white}{$\beta = 0.5$}}
            \end{overpic}
        \end{minipage} 
        & \begin{minipage}{0.2\textwidth}
            \begin{overpic}[width=0.81\textwidth]{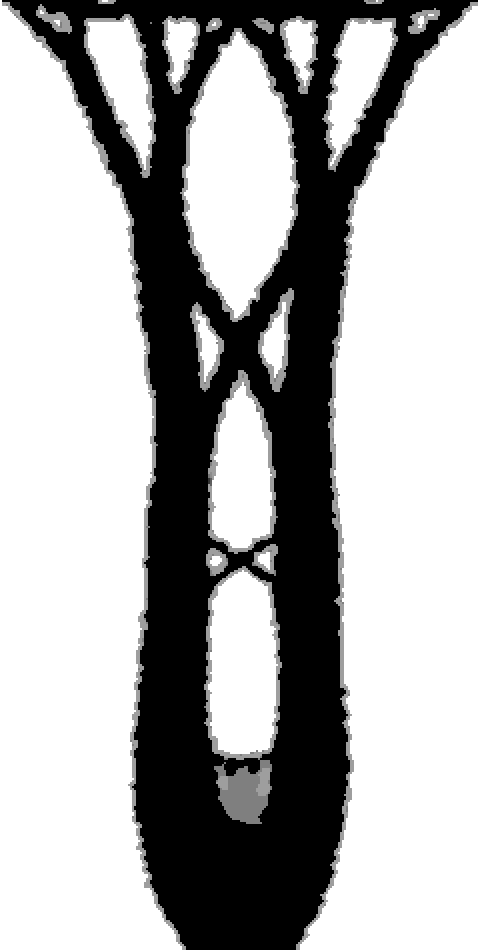}
                \put(2,5){\fcolorbox{black}{white}{$\beta = 0.9$}}
            \end{overpic}
        \end{minipage} 
        & \begin{minipage}{0.2\textwidth}
            \begin{overpic}[width=0.81\textwidth]{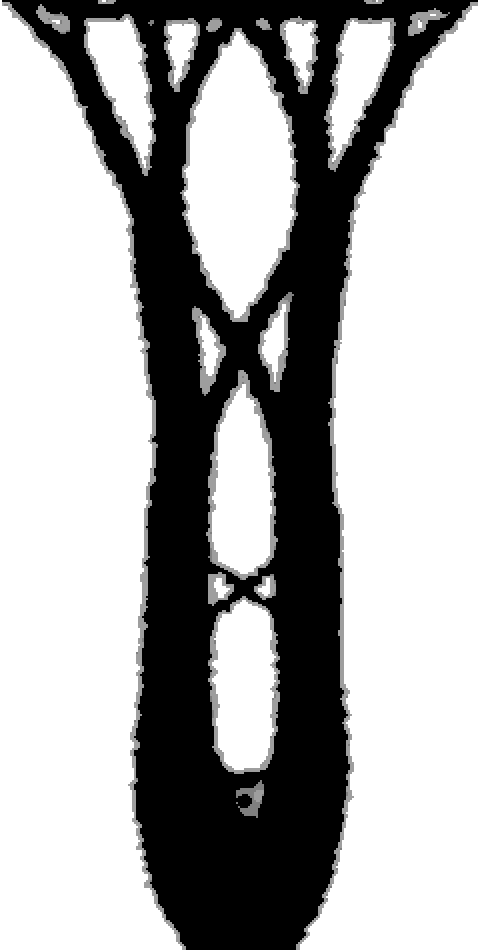}
                \put(2,5){\fcolorbox{black}{white}{$\beta = 0.99$}}
            \end{overpic}
        \end{minipage}
    \end{tabular}
    \caption{\it Distributionally robust designs of the bridge considered in \cref{subsec.cvardrobr}for a Wasserstein radius $m=0.5$, $\sigma^2=1e{-2.5}$ and various values of the parameter $\beta$.}
    \label{fig.cvardro2.5_0.5}
\end{figure}

This program is solved for several values of the threshold $\beta$, for two different values of the variance $\sigma^2$ of the reference coupling $\nu_{\xi^0}$, and for the radius $m=0.5$. The results are depicted on  \cref{fig.cvardro3_0.5,fig.cvardro2.5_0.5}. According to intuition, when perturbations become sufficiently large, the optimization algorithm gives privilege to the diagonal bars over the central vertical pillar of the structure. Also, as the variance $\sigma^2$  increases, the constraint becomes increasingly difficult to satisfy.

\section{\textbf{Distributionally Robust Geometric Shape Optimization}}
\label{sec.SO}

\noindent This section applies the generic distributionally robust formulations introduced in \cref{sec.absframework} in the context of geometric shape optimization. After presenting the mechanical and mathematical settings in \cref{sec.settingSO}, as well as the employed numerical methods for shape optimization in \cref{sec.numSO}, we conduct two numerical experiments in this context in \cref{sec.stress2dSO,sec.3dcanti}.

\subsection{Shape optimization of linear elastic structures}\label{sec.settingSO}

\noindent We first recall a few basic facts about the optimization of the shape of elastic structures.

\subsubsection{Linear elastic structures}

\noindent The designs considered henceforth are shapes, i.e. bounded Lipschitz domains $\Omega \subset \R^{d}$ ($d=2$ or $3$) representing elastic structures. The boundary $\partial \Omega$ of each shape $\Omega$ is decomposed into three disjoint regions:
$$\partial \Omega = \overline{\Gamma_D} \cup \overline{\Gamma_N} \cup \overline{\Gamma},$$
where:
\begin{itemize}
	\item The shape is clamped on the region $\Gamma_{D}$;	
	\item Surface loads $g : \Gamma_{N} \to \mathbb{R}^{d}$ are applied on $\Gamma_N$;
	\item The remaining part $\Gamma$ is traction-free and it is the only region of $\partial \Omega$ which is subject to optimization.
\end{itemize}
Omitting body forces for simplicity, the elastic displacement of $\Omega$ is the unique solution $u \in H^{1}(\Omega)^{d}$ to the linear elasticity system:
\begin{equation}\label{eq.elasshape}
\left\{
	\begin{array}{cl}
		-\dv(A e(u)) = 0 & \text{in } \Omega,   \\
		u = 0                    & \text{on } \Gamma_D, \\
		A e(u) n = g             & \text{on } \Gamma_N, \\
		A e(u) n = 0             & \text{on } \Gamma.
	\end{array}
\right.
\end{equation}
Again, we have omitted body forces for simplicity. The isotropic and homogeneous properties of the constituent material of shapes are encoded in the Hooke's tensor $A$ given by \cref{eq.Hookelaw}. 
In the sequel, we use various indices or arguments for $u$ and related quantities of interest, in order to indicate their dependence on the relevant parameters of this model, notably the loads $g$ and the shape $\Omega$ itself. 

An archetypal shape optimization problem in this context reads:
\begin{equation}\label{eq.sopb}
\min_{\Omega}\: C(\Omega) \quad \text{ s.t. } \quad \Vol(\Omega) = V_{T},
\end{equation}
where $C(\Omega)$ is the compliance of $\Omega$, i.e.:
\begin{equation}\label{eq.cplySO}
C(\Omega) = \int_{\Omega} A e(u_{\Omega}) : e(u_{\Omega}) \, \d x = \int_{\Gamma_N} g \cdot u_{\Omega} \, \d s,
\end{equation}
and $\Vol(\Omega) := \int_{\Omega} dx$ is the volume, which should match a target value $V_{T}$.

\subsubsection{Shape sensitivity by the method of Hadamard} \label{sec.shapeder}

\begin{figure}[ht]
	\centering
	\includegraphics[width=0.45\linewidth]{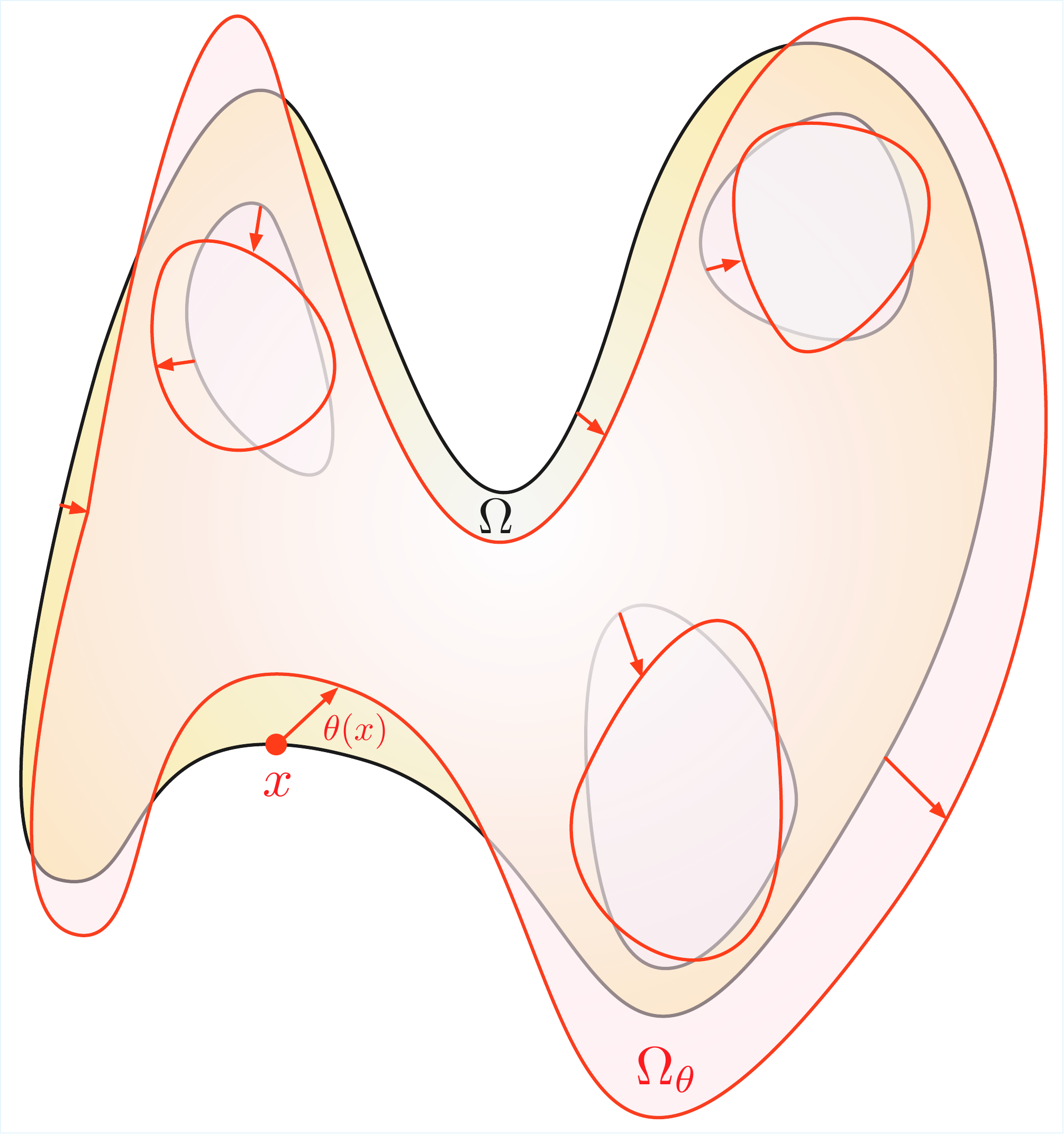}
	\caption{\it Variation $\Omega_{\theta} = (\Id + \theta)(\Omega)$ of a shape $\Omega$ by the method of Hadamard.}
	\label{fig.hadamard}
\end{figure}

\noindent When it comes to study the sensitivity of a function with respect to the domain, we rely on the boundary variation method of Hadamard, about which we refer to e.g. \cite{allaire2007conception,henrot2018shape,murat1976controle,sokolowski1992introduction}.
In brief, this method is based on variations of a reference shape $\Omega \subset \R^{d}$ of the form:
$$\Omega_{\theta} := (\Id + \theta)(\Omega), \quad \text{where } \theta \in \Winfty, \quad \| \theta \|_{\Winfty} < 1,$$
i.e. $\Omega_{\theta}$ is a version of $\Omega$ deformed by a ``small" vector field $\theta$, see \cref{fig.hadamard}.

\begin{definition}
A functional $F(\Omega)$ is called shape differentiable at a particular domain $\Omega$ if the mapping
$$\Winfty \ni \theta \mapsto F(\Omega_{\theta}) \in \R$$
is Fréchet differentiable at $0$. Its derivative $F^{\prime}(\Omega)(\theta)$ is called the shape derivative of $F$ at $\Omega$, and the following expansion holds true:
\begin{equation}\label{eq.asymhadamard}
F(\Omega_{\theta}) = F(\Omega) + F^{\prime}(\Omega)(\theta) + \o(\theta), \text{where } \frac{\o(\theta)}{\| \theta \|_{\Winfty}} \xrightarrow{\theta\to 0} 0.
\end{equation}
\end{definition}

It follows directly from this definition that a vector field $\theta$ satisfying $F^\prime(\Omega)(\theta) < 0$ is a descent direction for $F(\Omega)$, i.e. for a small enough pseudo-time step $t>0$, $F(\Omega_{t\theta}) < F(\Omega)$. More generally, in the context of a constrained optimization problem, of the form \cref{eq.cplySO}, the knowledge of the shape derivatives of the objective and constraint functions allows to infer a suitable ``descent direction'' via a constrained optimization method such as the Null-space algorithm used in this work \cite{feppon2020null}. 
The shape derivatives of the functions $F(\Omega)$ considered in this article often take the convenient form:
$$F^{\prime}(\Omega)(\theta) = \int_{\partial \Omega} v_{\Omega} \, \theta \cdot n \, \d s.$$
A descent direction for $F(\Omega)$ is immediately revealed from this structure, by taking $\theta = -v_{\Omega}n$ on $\partial \Omega$. 

\subsection{Numerical resolution of shape optimization problems by a level set based mesh evolution strategy} \label{sec.numSO}

\begin{figure}[ht]
\centering
\captionsetup{justification=centering}
\begin{tabular}{cc}
		\begin{minipage}{0.45\textwidth}\begin{overpic}[width=1.0\textwidth]{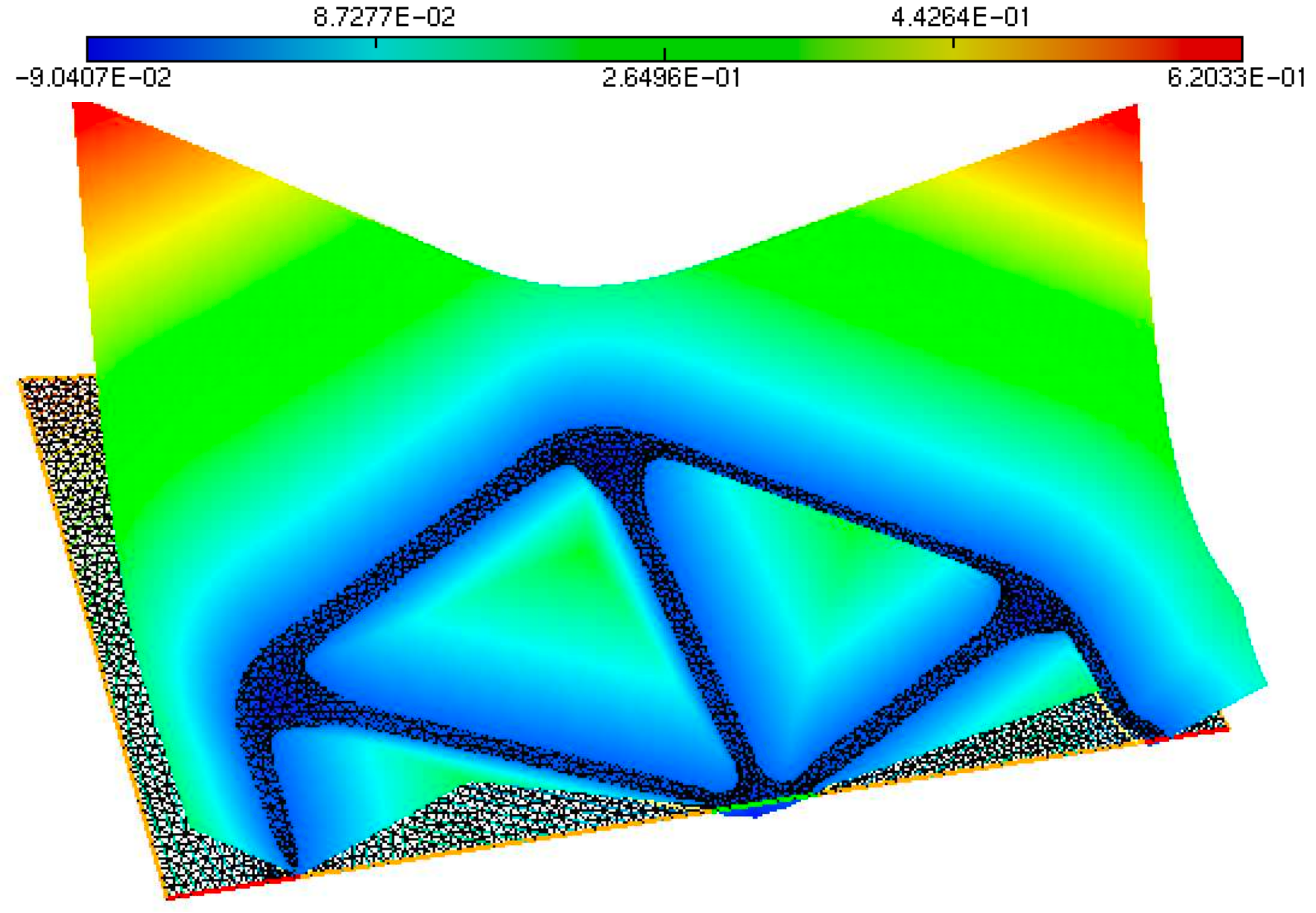} \put(2,5){\fcolorbox{black}{white}{a}}\end{overpic}\end{minipage} & \begin{minipage}{0.5\textwidth}\begin{overpic}[width=1.0\textwidth]{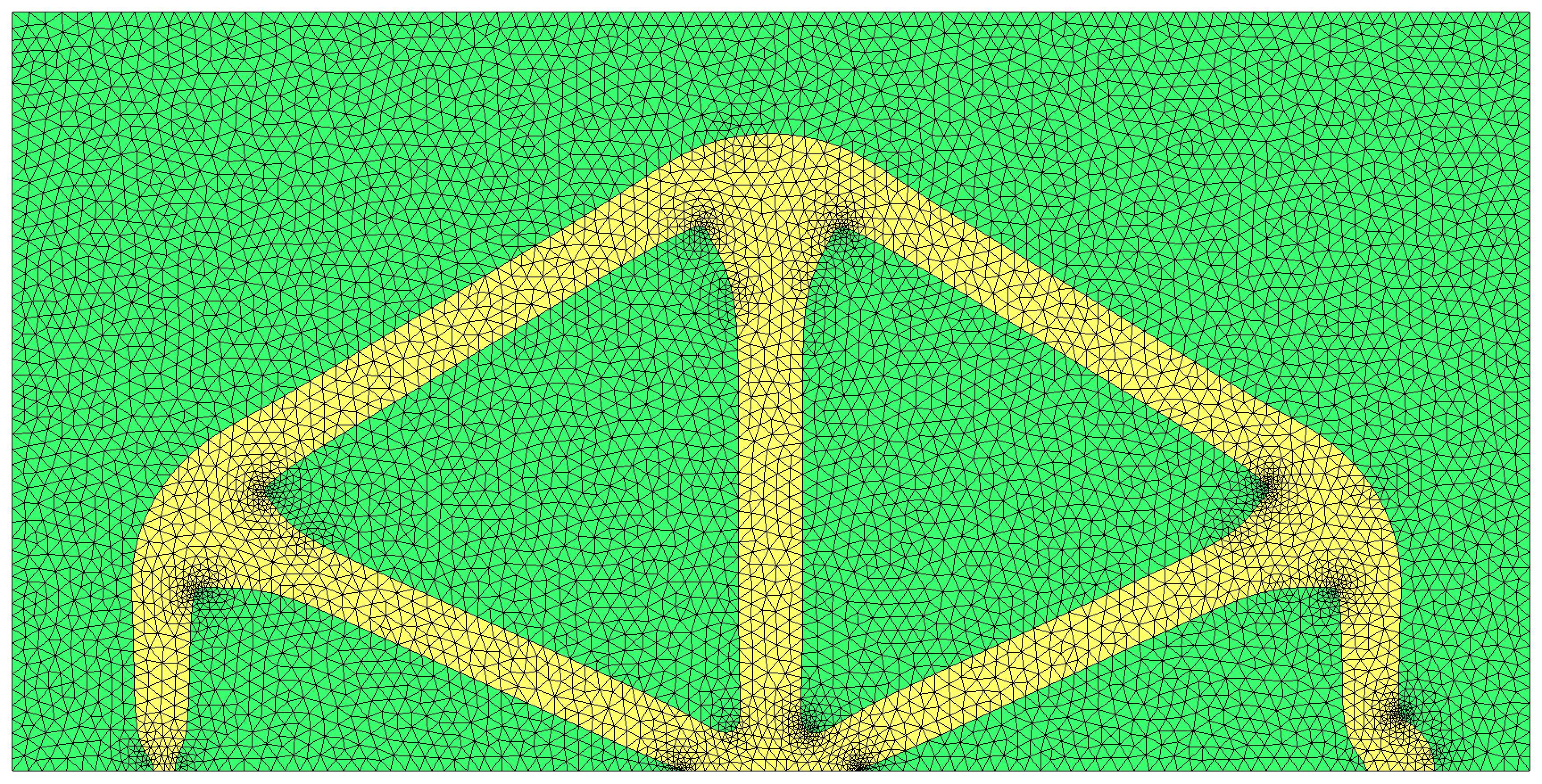} \put(2,5){\fcolorbox{black}{white}{b}}\end{overpic}\end{minipage}
	\end{tabular}
	\caption{\it Two complementary representations of a shape $\Omega \subset D$; (a) Graph of a level set function $\phi : D \to \R$; (b) Corresponding meshed representation.}
\label{fig.2reps}
\end{figure}

\noindent
Our numerical resolution of shape optimization problems of the form \cref{eq.sopb} hinges on the level set based mesh evolution method introduced in \cite{allaire2011topology,allaire2013mesh,allaire2014shape}, an open-source implementation of which is provided in \cite{dapogny2022tuto}.
Briefly, introducing a fixed computational domain $D\subset \R^d$, we rely on two complementary representations of each shape $\Omega \subset D$ arising during the optimization process:
\begin{itemize}
\item (Level set representation) As proposed in \cite{osher1988fronts} -- see also \cite{allaire2004structural,wang2003level} in the shape optimization context -- $\Omega$ is represented implicitly, as the negative subdomain of a ``level set'' function $\phi: D \to \R$, i.e.
$$\left\{
\begin{array}{cl}
\phi(x) < 0 & \text{if } x \in \Omega,          \\
\phi(x) = 0 & \text{if } x \in \partial \Omega, \\
\phi(x) > 0 & \text{if } x \in D \setminus \overline{\Omega}.
\end{array}
\right.$$
In practice, such a function $\phi$ is discretized at the vertices of a simplicial mesh $\calT$ of $D$ (i.e. made of triangles in 2d, of tetrahedra in 3d), see \cref{fig.2reps} (a).
\item (Meshed representation) The domain $D$ is equipped with a conforming, high-quality simplicial mesh $\calT$, made of two disjoint submeshes, $\calT_{\mathrm{int}}$ and $\calT_{\mathrm{ext}}$, which respectively represent $\Omega$ and its complement $D \setminus \overline{\Omega}$, see \cref{fig.2reps} (b).
\end{itemize}

Efficient numerical techniques make it possible to switch from one of these representations to the other:
\begin{itemize}
\item Given a meshed representation of a shape $\Omega \subset D$, one particular level set function for $\Omega$ is a calculated as its signed distance function, thanks to the famous Fast Marching Algorithm \cite{sethian1999fast}; this task relies on the open-source software \texttt{mshdist} from our previous work  \cite{dapogny2012computation}. 
\item Assuming a level set representation of $\Omega \subset D$, from the datum of the values of a level set function $\phi$ at the vertices of a mesh $\calT$ of $D$, a new mesh $\widetilde{\calT}$ of $D$ in which $\Omega$ is explicitly discretized is obtained thanks to a remeshing procedure, implemented in the open-source library \texttt{mmg}, see \cite{balarac2022tetrahedral,dapogny2014three}.
\end{itemize} 

Each operation of the shape optimization workflow is executed by using the most appropriate representation. In particular, ``physical'' boundary value problems posed on a shape $\Omega \subset D$ are solved by applying the Finite Element method on its meshed representation, which allows to use any numerical solver in a black-box manner; in our framework, we use the open-source platform \texttt{FreeFem} \cite{hecht2012new}. Conversely, the level set representation is particularly well-suited to track the evolution of shape $\Omega(t)$ over a time interval $(0, T)$, under the action of a prescribed velocity field $V(t, x)$. This evolution is indeed encoded in the following advection-like equation:
$$\dfrac{\partial \phi}{\partial t}(t,x) + V(t,x) \cdot \nabla \phi(t,x) = 0  \text{ for } t \in (0,T), \, x \in D,$$
which can be solved e.g. by the method of characteristics \cite{pironneau1989finite}; we rely to this end on the open-source code \texttt{advection} from \cite{bui2012accurate}.

Summarizing, this strategy features a body-fitted mesh of the shape $\Omega$ at each stage of the workflow, allowing for precise mechanical analyses based on any finite element solver, while leaving the room for arbitrarily large deformations between successive optimization iterations.

\subsection{Distributionally robust shape optimization of the stress within a 2d T-shaped beam under load uncertainties}\label{sec.stress2dSO}

\begin{figure}[ht]
\centering
	\begin{tabular}{ccc}
		\begin{minipage}{0.31\textwidth}
		\begin{overpic}[width=1.\textwidth]{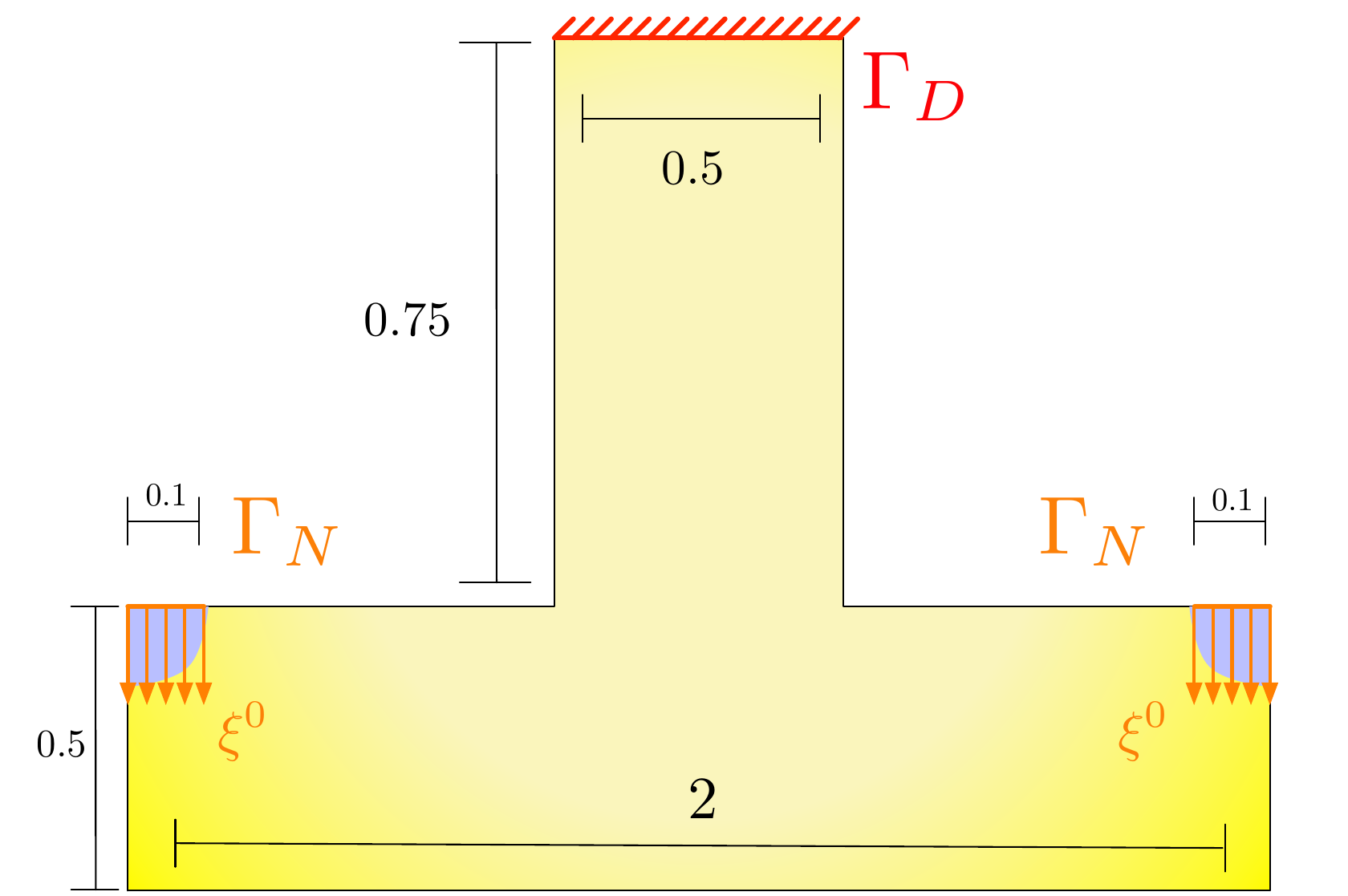} 
		\put(2,5){\fcolorbox{black}{white}{a}}
		\end{overpic}
		\end{minipage} & 
		\begin{minipage}{0.31\textwidth}
        \begin{overpic}[width=1.\textwidth]{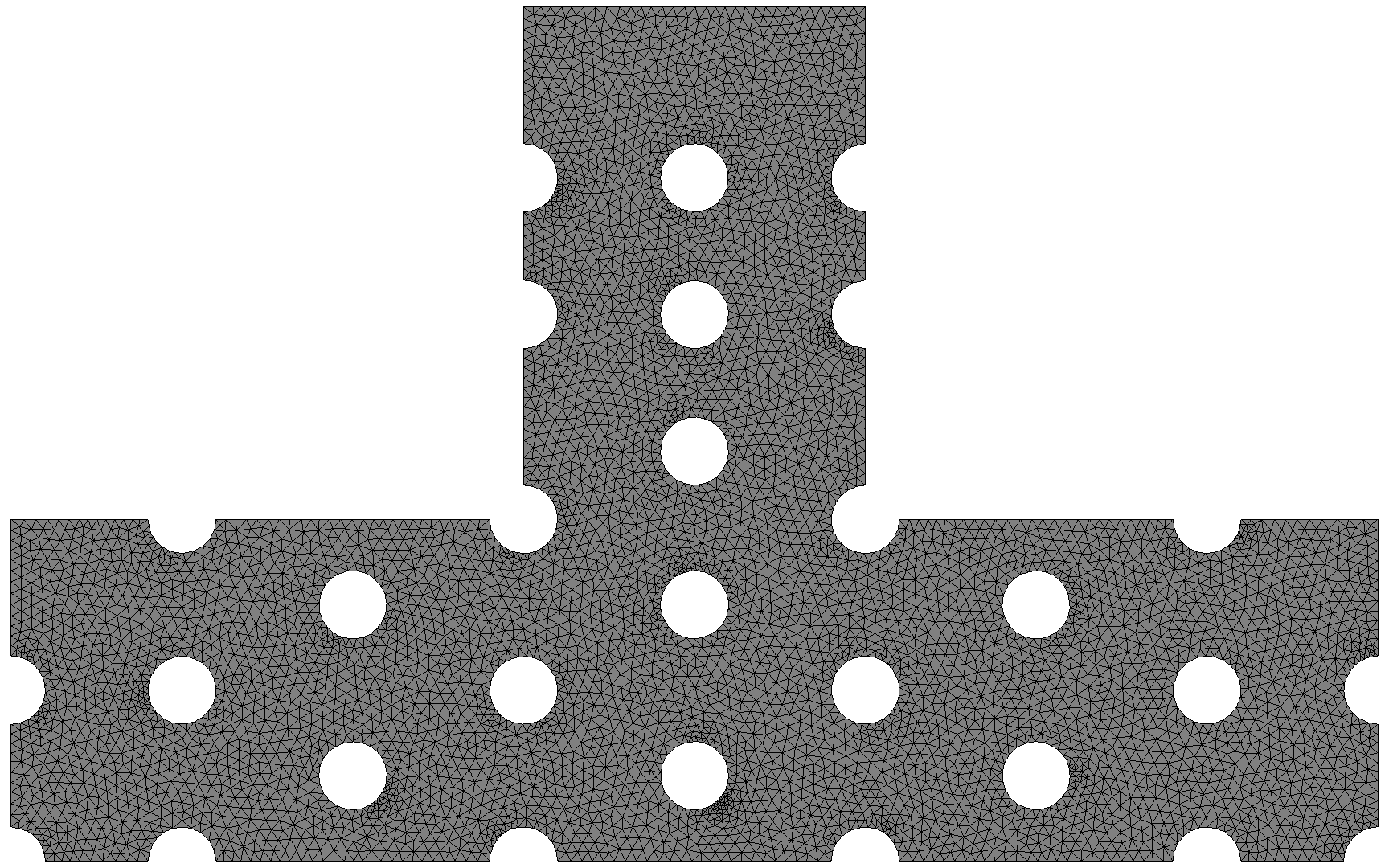} 
		\put(2,5){\fcolorbox{black}{white}{b}}
		\end{overpic}
	\end{minipage}&
    \begin{minipage}{0.31\textwidth}
    \begin{overpic}[width=1.\textwidth]{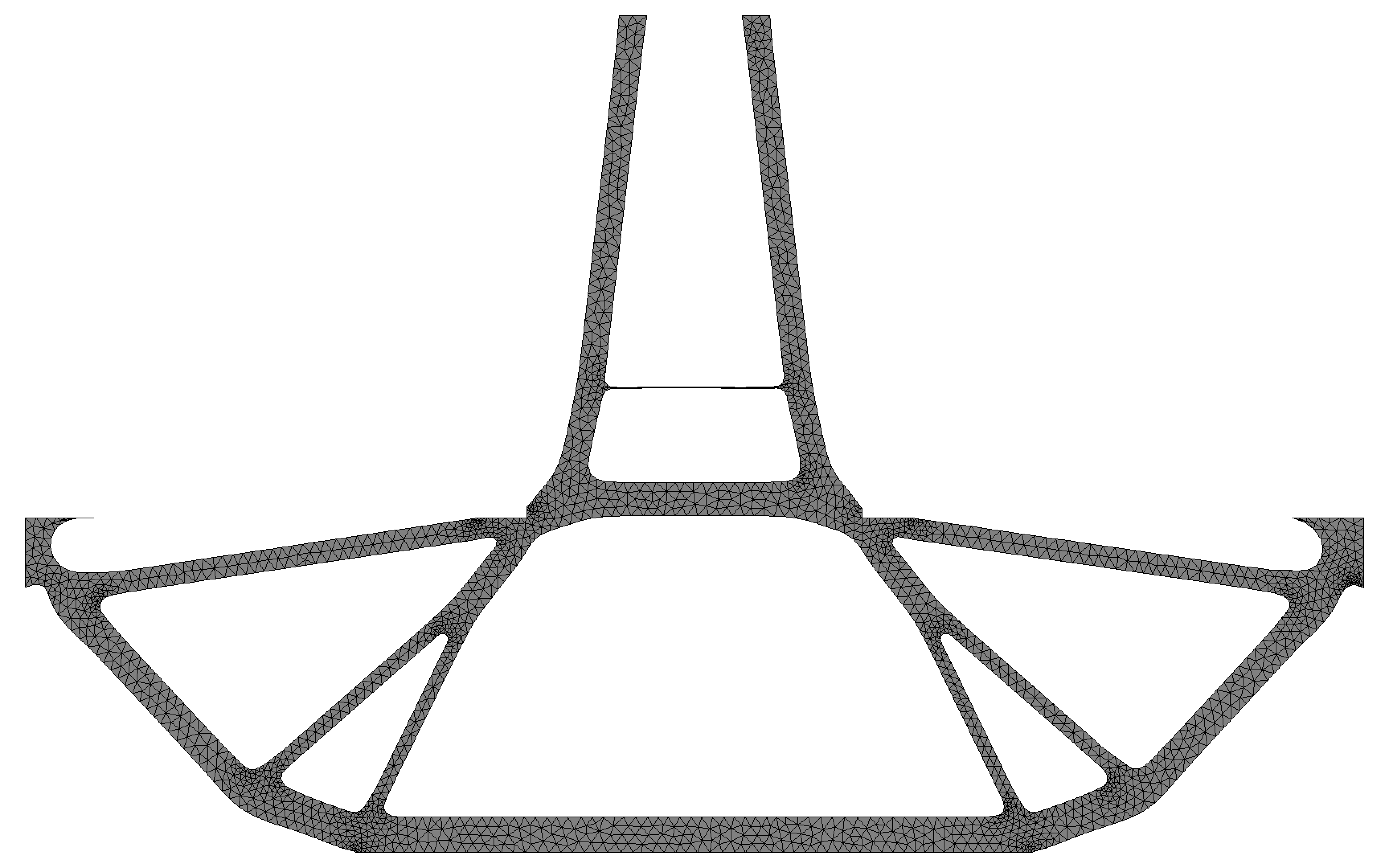}
	\put(2,5){\fcolorbox{black}{white}{c}}
	\end{overpic}        
    \end{minipage}
	\end{tabular}
	\caption{\it (a) Setting of the minimization of stress within a T-shaped beam considered in \cref{sec.stress2dSO}; (b) Mesh of the initial shape; (c) Optimized shape $\Omega^*_{\text{det}}$ for the problem \cref{eq.minpbbrso}, assuming a perfect knowledge of the loads.}
    \label{fig.TshapestressnomV2}
\end{figure}

\noindent Let the shape $\Omega$ account for a 2d T-shaped beam, enclosed in a box $D$ with size $2 \times 1$, see \cref{fig.TshapestressnomV2} (a).
 The beam is clamped on the upper region $\Gamma_D$ of $\partial D$;
under ideal conditions, it is subjected to two identical vertical loads $\xi^0 = (0,-1)$, applied on two small regions at the left and right parts of $\partial D$, whose reunion is denoted by $\Gamma_N$. In the following, we denote by $u_{\Omega,\xi}$ the displacement of the beam, solution to \cref{eq.elasshape} when the load $\xi$ is applied.

We aim to minimize the total stress within $\Omega$ under a constraint on its volume, i.e. we consider the following shape optimization problem: 
\begin{equation}\label{eq.minpbbrso}
 \min \limits_{\Omega \subset D} S(\Omega,\xi^0) \:\text{ s.t. } \Vol(\Omega) = V_T.
 \end{equation}
 Here, the stress function is defined by 
 $$ S(\Omega,\xi) = \int_{\Omega} \chi(x) ||\sigma(u_{\Omega,\xi})||^2\;\d x,$$
 where the stress $\sigma(u)$ associated to a displacement field $u$ is $\sigma(u) = Ae(u)$, and $\lvert\lvert \cdot \lvert\lvert$ denotes the usual Frobenius norm over matrices. The weight  $\chi(x)$ equals $1$ everywhere on $D$ except on two small non-optimizable neighborhoods of the region $\Gamma_N$ where it is set to $0$ (these zones are depicted in blue  in \cref{fig.TshapestressnomV2} (a)) Eventually, in \cref{eq.minpbbrso}, the target volume $V_T$ is set to $V_T = 0.3$.

The shape derivative of the stress-based functional $S(\Omega,\xi)$ is given by the following proposition, which is taken from e.g. \cite{allaire2008minimum}.

\begin{proposition}\label{prop.SOmprime}
The shape derivative of the functional $S(\Omega,\xi)$ equals: 
$$S^{\prime}(\Omega,\xi)(\theta) = \int_{\Omega}Ae(u_{\Omega,\xi}) : e(p_{\Omega,\xi}) \: \theta\cdot n \:\d s,$$ 
where the adjoint state $p_{\Omega,\xi}$ is the unique solution in $H^1(\Omega)^2$ to the following boundary-value problem:
\begin{equation} \label{eq.adjS}\left\{
	\begin{array}{cl}
		-\dv(Ae(p_{\Omega,\xi})) = -2 \dv(\chi(x) A\sigma(u_{\Omega,\xi}))& \text{in } \Omega,                                        \\
		p_{\Omega,\xi} = 0                                               & \text{on } \Gamma_D,                                 \\
		A e(p_{\Omega,\xi})n = 2 A\sigma(u_{\Omega,\xi}) n                                   & \text{on } \partial \Omega \setminus \overline{\Gamma_D}.
	\end{array}
	\right.
\end{equation}
\end{proposition}

We first solve this ideal problem with the algorithm presented in \cref{sec.numSO}. 
About $300$ iterations of our optimization algorithm are applied, corresponding to a computation of about 16 mn; the resulting shape $\Omega^*_{\text{det}}$ is depicted on \cref{fig.TshapestressnomV2} (c). The corresponding value of the stress functional is $S(\Omega^*_{\text{det}},\xi^0) = 1.99433$. The meshes of $D$ produced at each iteration of the optimization process, made of submeshes for the structure $\Omega$ and its complement, contain about $6,500$ vertices ($\approx 13,000$ triangles). 

\begin{figure}[ht]
	\centering
	\begin{tabular}{ccc}
		\begin{minipage}{0.31\textwidth}\begin{overpic}[width=1.\textwidth]{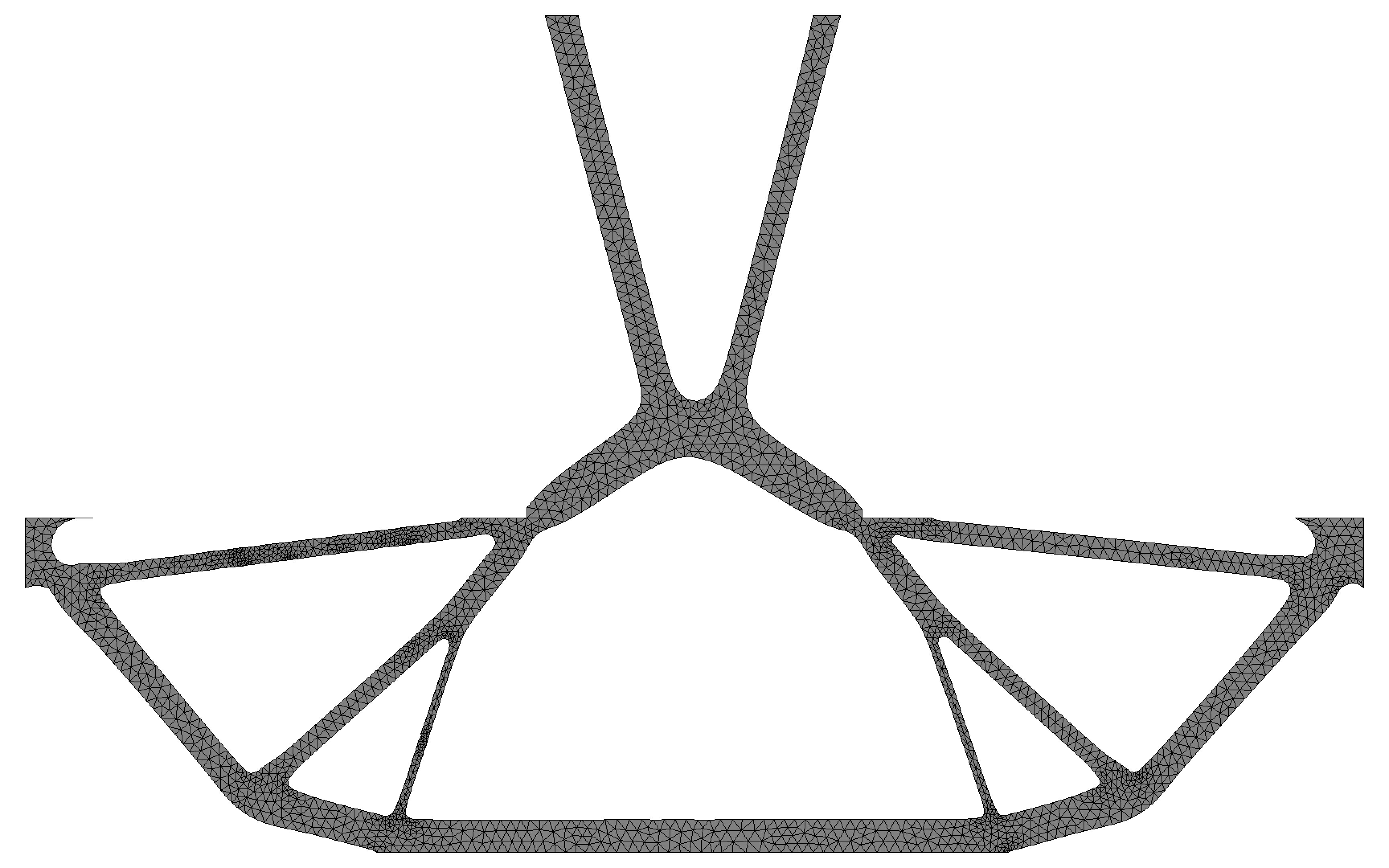} \put(2,5){\fcolorbox{black}{white}{$m=0$}}\end{overpic}\end{minipage} & \begin{minipage}{0.31\textwidth}\begin{overpic}[width=1.0\textwidth]{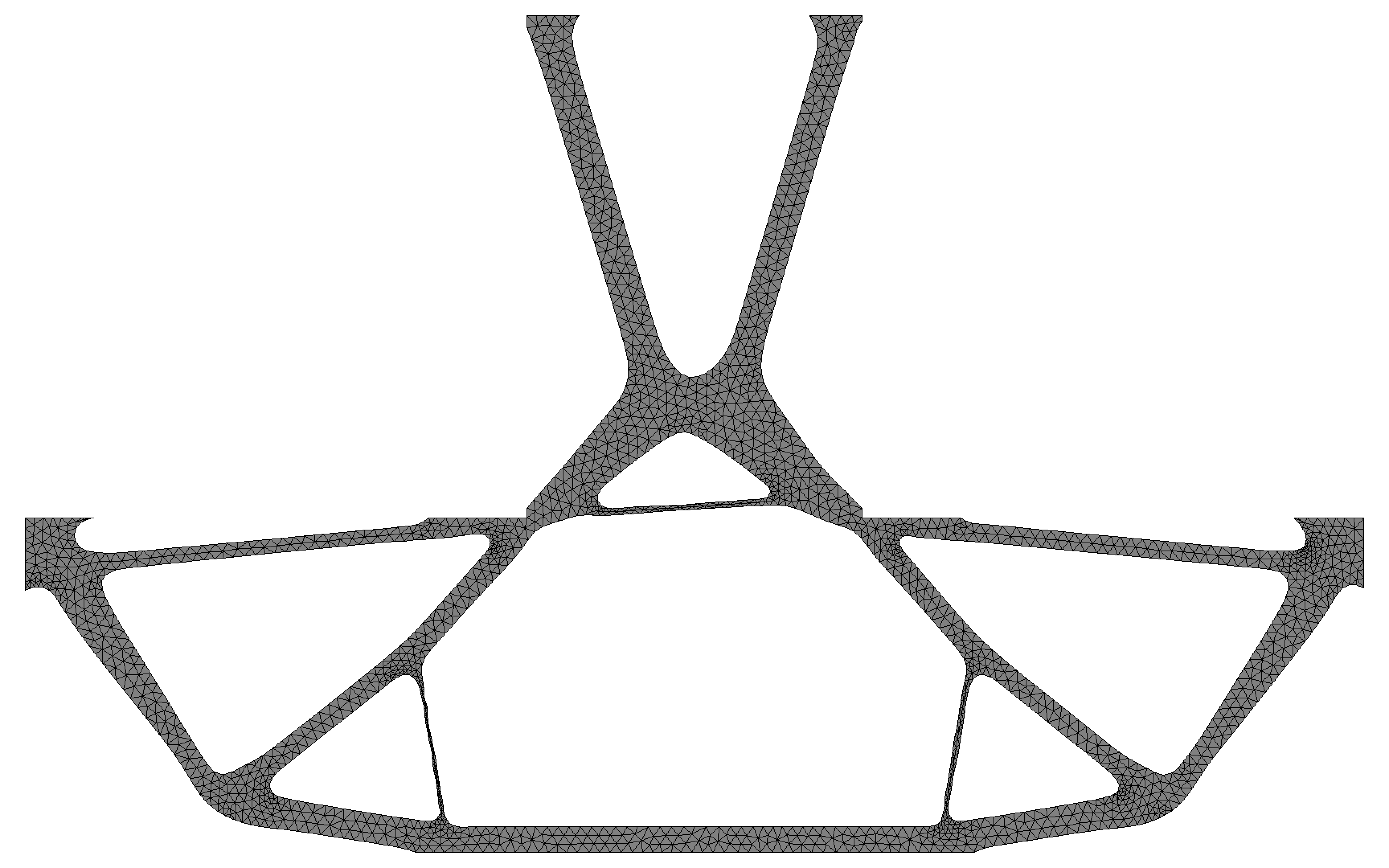} \put(2,5){\fcolorbox{black}{white}{$m=1$}}\end{overpic}\end{minipage} &
        \begin{minipage}{0.31\textwidth}\begin{overpic}[width=1.0\textwidth]{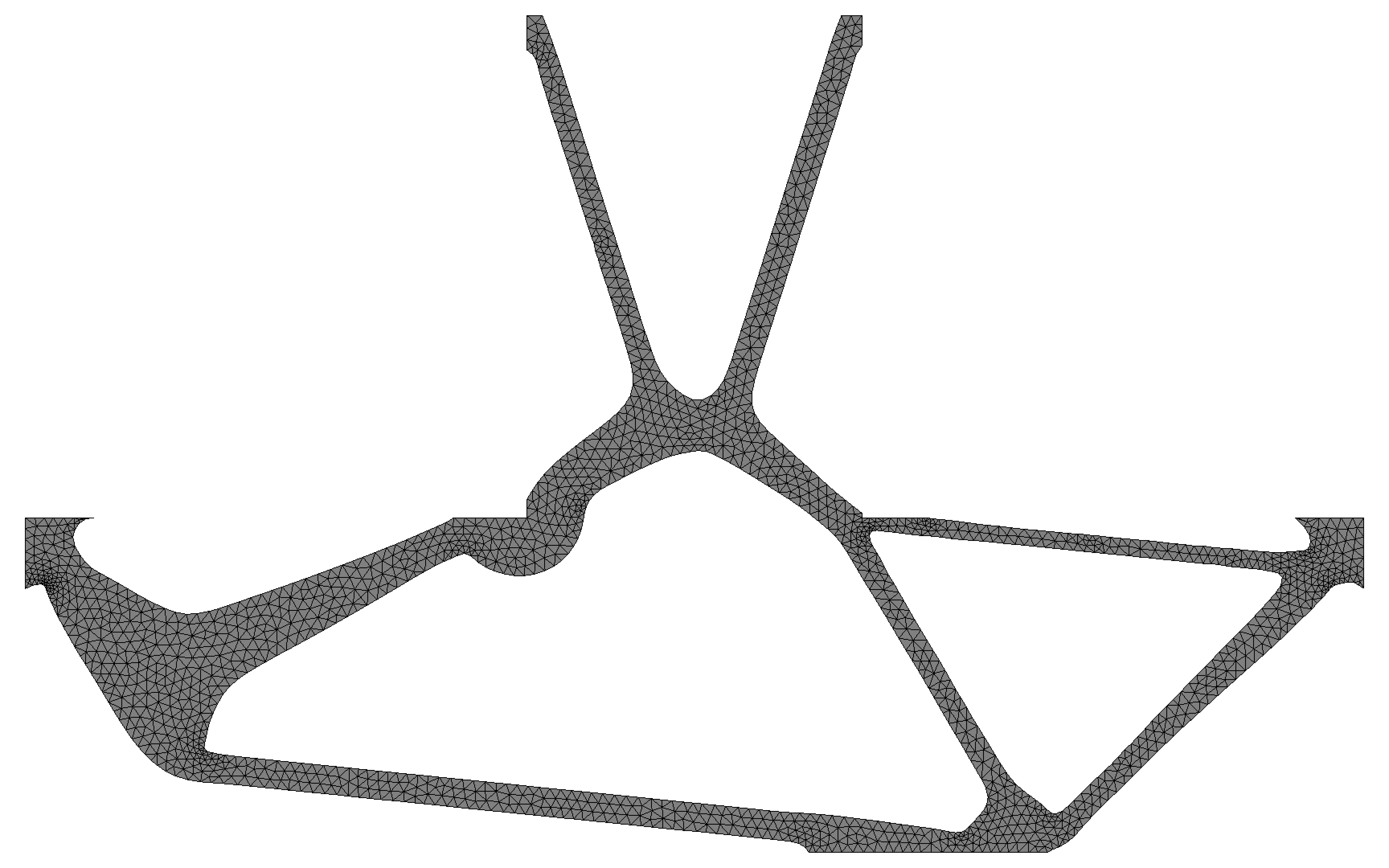} \put(2,5){\fcolorbox{black}{white}{$m=3$}}\end{overpic}\end{minipage}
	\end{tabular}
	\par
	\medskip
	\begin{tabular}{ccc}
		\begin{minipage}{0.31\textwidth}\begin{overpic}[width=1.0\textwidth]{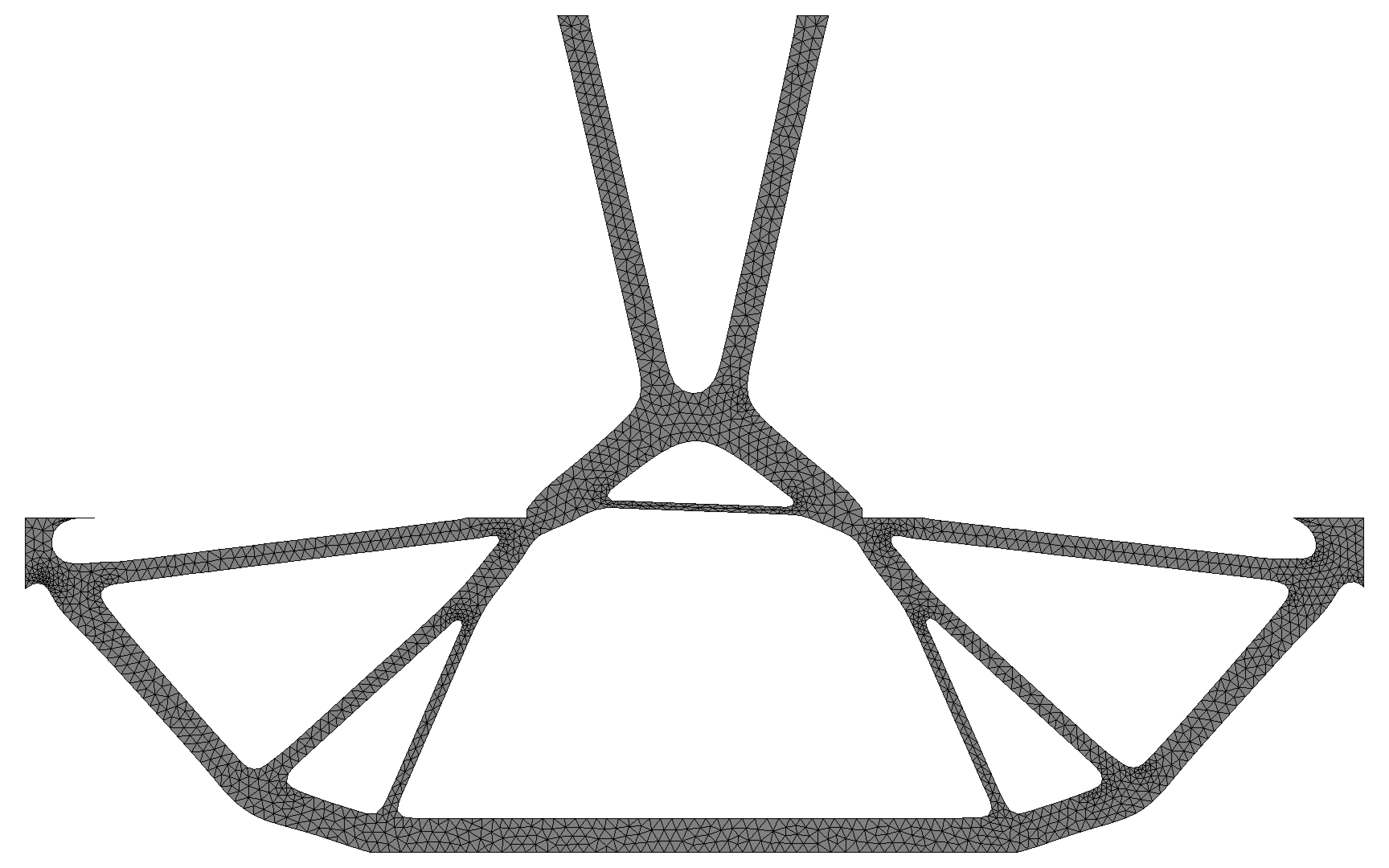} \put(2,5){\fcolorbox{black}{white}{$m=0$}}\end{overpic}\end{minipage} & \begin{minipage}{0.31\textwidth}\begin{overpic}[width=1.0\textwidth]{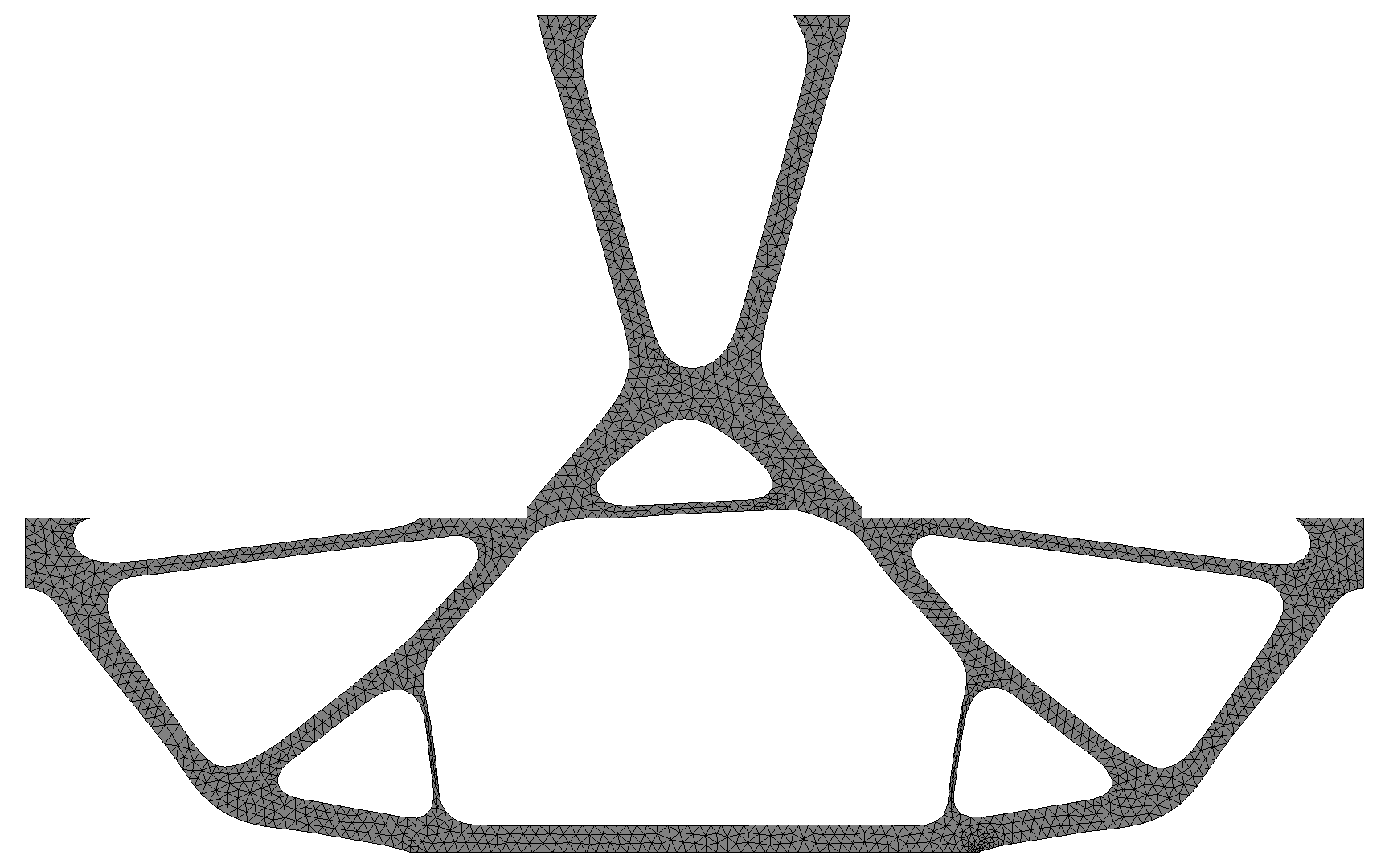} \put(2,5){\fcolorbox{black}{white}{$m=1$}}\end{overpic}\end{minipage} &
        \begin{minipage}{0.31\textwidth}\begin{overpic}[width=1.0\textwidth]{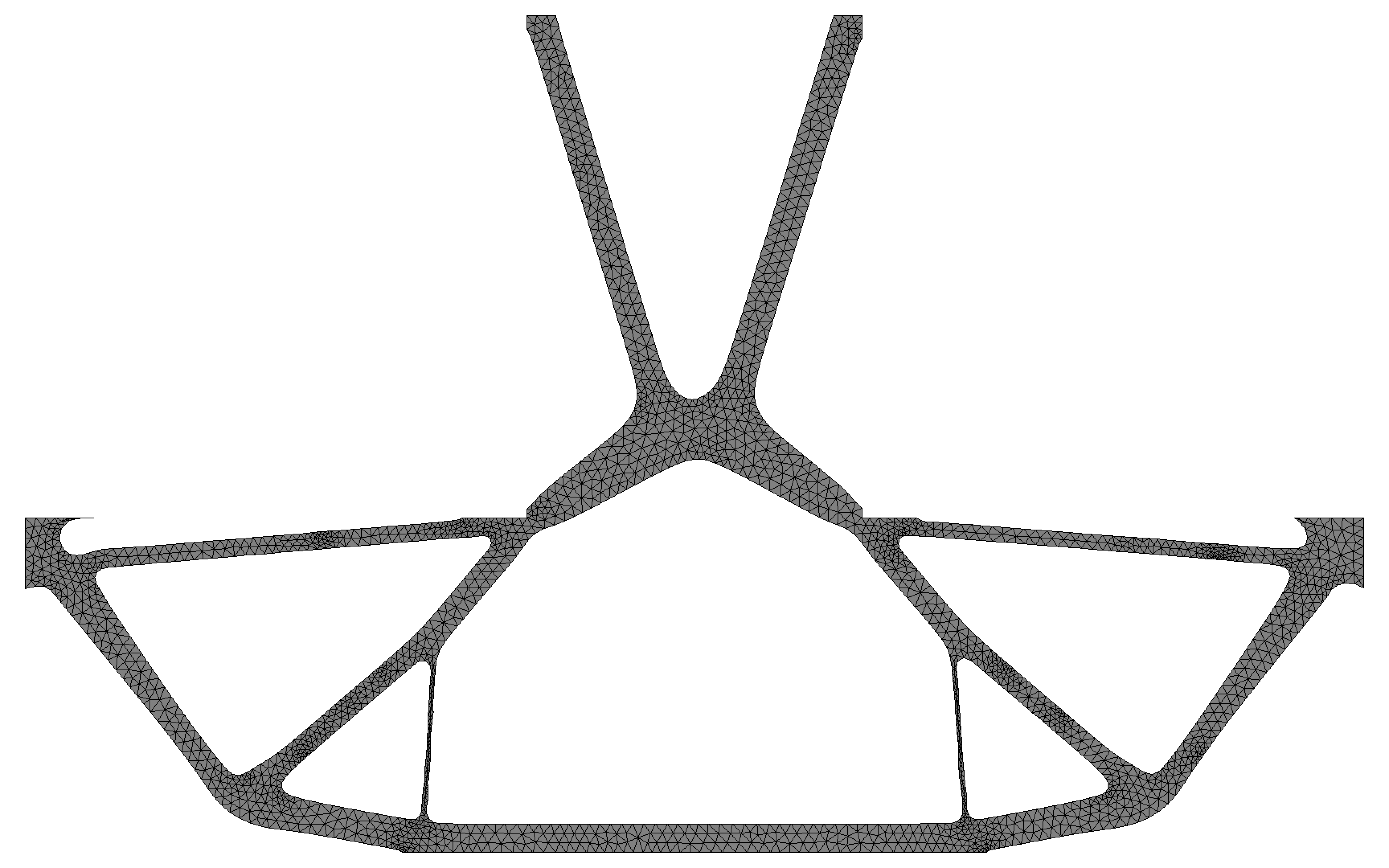} \put(2,5){\fcolorbox{black}{white}{$m=3$}}\end{overpic}\end{minipage}
	\end{tabular}
    \par\medskip
    \begin{tabular}{ccc}
		\begin{minipage}{0.31\textwidth}\begin{overpic}[width=1.0\textwidth]{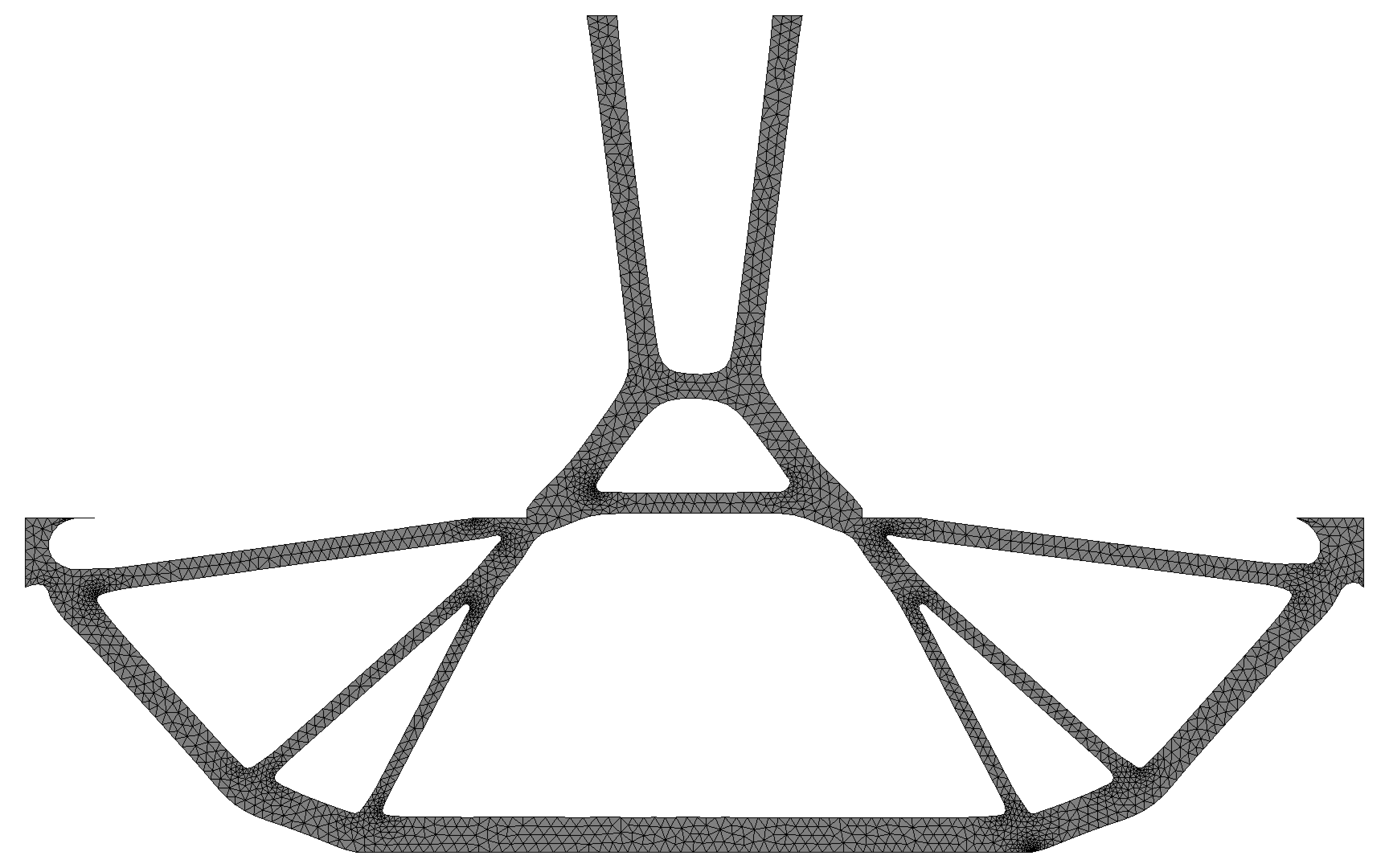} \put(2,5){\fcolorbox{black}{white}{$m=0$}}\end{overpic}\end{minipage} & \begin{minipage}{0.31\textwidth}\begin{overpic}[width=1.0\textwidth]{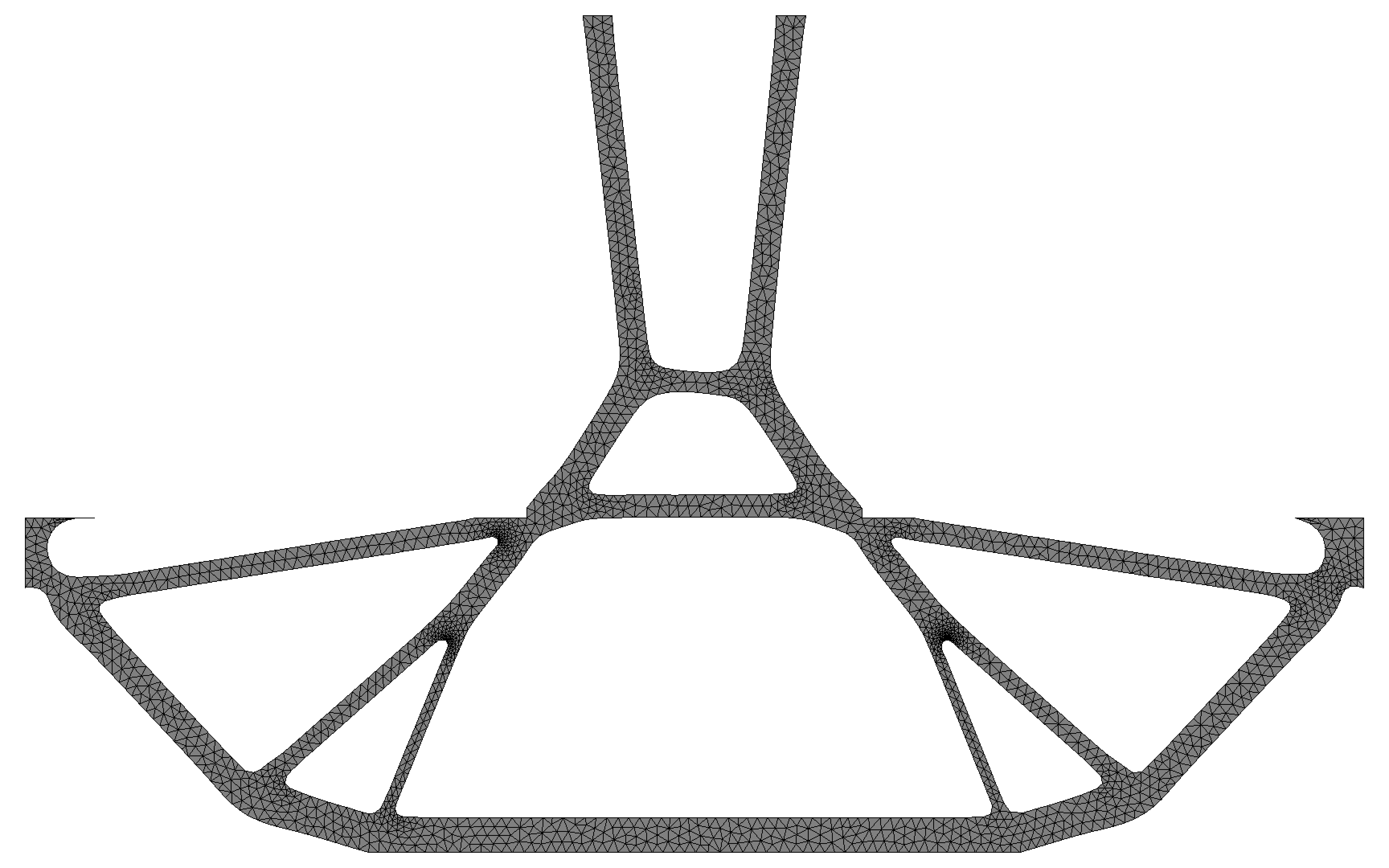} \put(2,5){\fcolorbox{black}{white}{$m=1$}}\end{overpic}\end{minipage} &
        \begin{minipage}{0.31\textwidth}\begin{overpic}[width=1.0\textwidth]{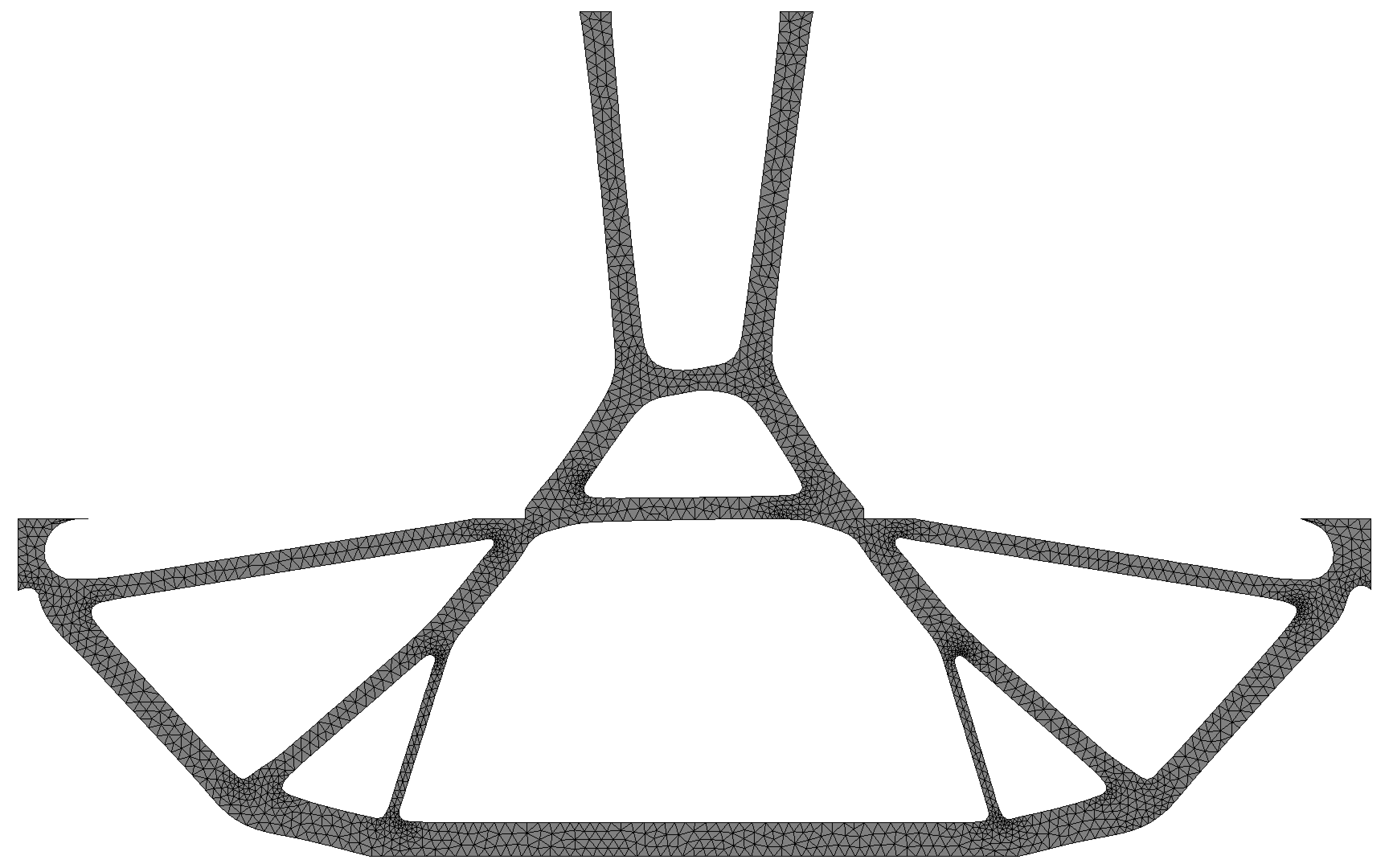} \put(2,5){\fcolorbox{black}{white}{$m=3$}}\end{overpic}\end{minipage}
	\end{tabular}
	\caption{\it Distributionally robust optimized designs obtained in the T-shaped beam example of \cref{sec.stress2dSO}, associated to various values of the parameters $m$ and $\sigma^2$;
    ($1^{st}$ row) $\sigma^2 = 1e{-0.5}$; ($2^{nd}$ row) $\sigma^2 = 1e{-1}$; ($3^{rd}$ row) $\sigma^2 = 1e{-3}$.}
    	\label{fig.TstressdroV2}
\end{figure}

We now assume that the applied load $\xi$ is uncertain; it belongs to a sufficiently large ball $\Xi \subset \R^2$, and its probability law is also uncertain -- the only available information about the latter being that it is ``close'' to the nominal distribution $\P = \delta_{\xi^0} \in \calP(\Xi)$ reconstructed from the single ideal load $\xi^0$. We thus turn to a distributionally robust version of \cref{eq.minpbbrso}: 
\begin{equation}\label{eq.drobridge}
\min\limits_{\Omega} \left(\sup\limits_{\Q \in \calA_{\text{W}}} \int_\Xi S(\Omega,\xi) \:\d\Q(\xi) \right), \text{ s.t. } \Vol(\Omega) = V_T,
\end{equation}
where $\calA_{\text{W}}$ is  the Wasserstein ambiguity set defined in \cref{eq.calAW}.
Using the material in \cref{sec.wdro}, this problem rewrites:
\begin{multline*}
	\min\limits_{\Omega, \atop \lambda \geq 0} \: \mathcal{D}_{\text{W}}(\Omega,\lambda)
	\quad \text{ s.t. } \quad \text{Vol}(\Omega) = V_T, \quad \text{where} \\
	\mathcal{D}_{\text{W}}(\Omega,\lambda) = \lambda m + \lambda \varepsilon \int_{\Xi} \log\left(
	\int_{\Xi} \exp\left(\frac{S(\Omega,\zeta) - \lambda c(\xi,\zeta)}{\lambda \varepsilon}\right)
	\,\mathrm{d}\nu_{\xi}(\zeta) \right)\, \d\P(\xi).
\end{multline*}
Here, we recall that $\nu_\xi(\zeta) \in \calP(\Xi) =  \alpha_{\xi}e^{-\frac{\lvert
	\xi-\zeta\lvert^{2}}{2\sigma^2}}\mathds{1}_{\Xi}(\zeta)$ is the probability measure featured in the definition of the reference coupling $\pi_0$, see \cref{eq.refcouplingW}.
The entropy regularization parameter $\e$ is set to $\e=0.01$.
We solve this problem for several values of the Wasserstein radius $m$ and the variance $\sigma^2$. Each computation requires $300$ iterations of our numerical strategy, for a CPU time of about 1h 20 mn. The numerical results are presented in \cref{fig.TstressdroV2}. The distributionally robust designs feature spread bars connected to the Dirichlet region $\Gamma_D$ when compared to the result of \cref{fig.TshapestressnomV2} (c). Likewise, the bars located at the center of the structure tend to adopt triangular or trapezoidal layouts, guaranteeing robustness under all perturbations.

\begin{figure}[h]
    \centering
    \begin{subfigure}{0.45\textwidth}
    \begin{overpic}[width=1.0\textwidth]
    {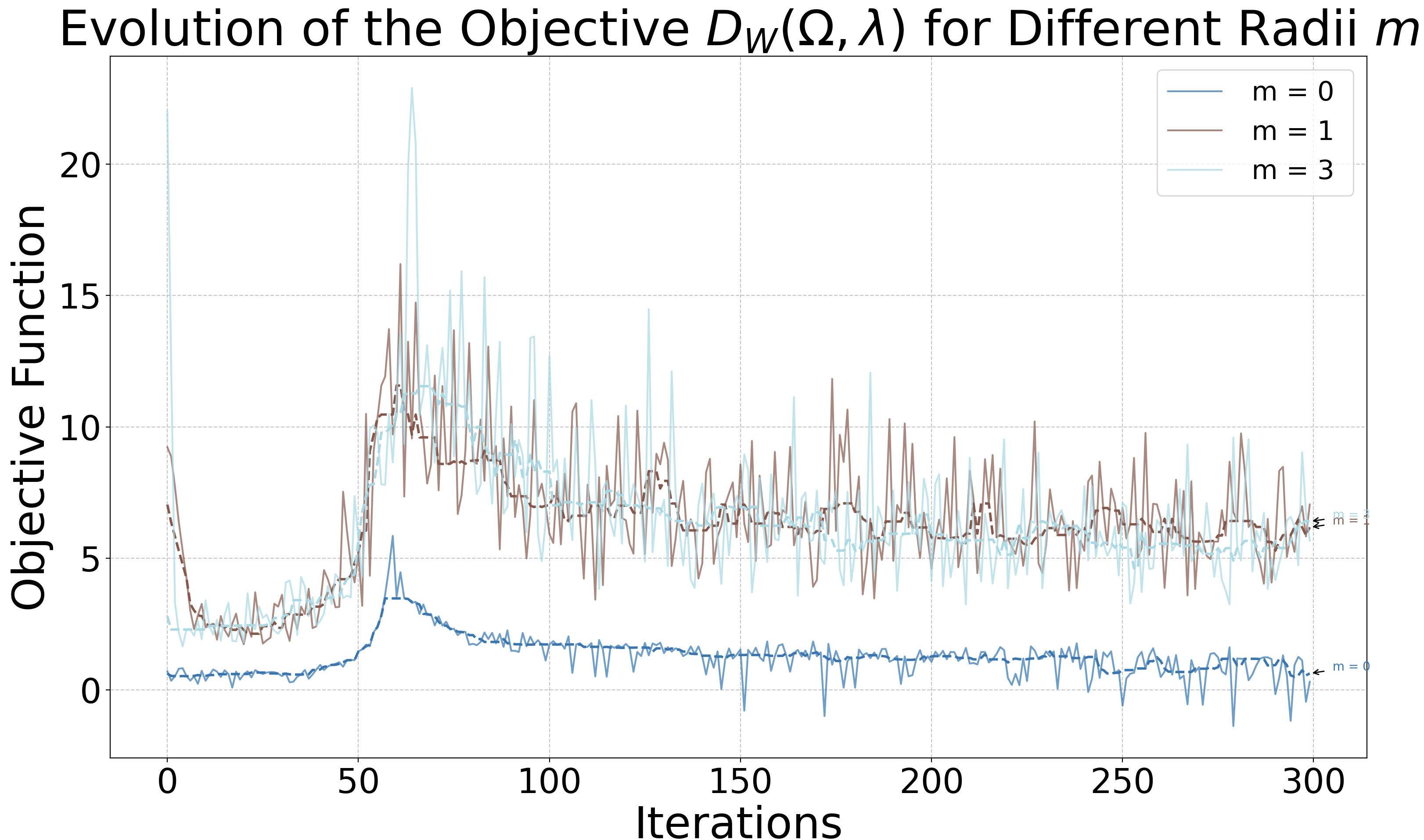}\put(2,5){\fcolorbox{black}{white}{a}}\end{overpic}
    \end{subfigure}
    \hfill
        \begin{subfigure}{0.45\textwidth}
    \begin{overpic}[width=1.0\textwidth]
    {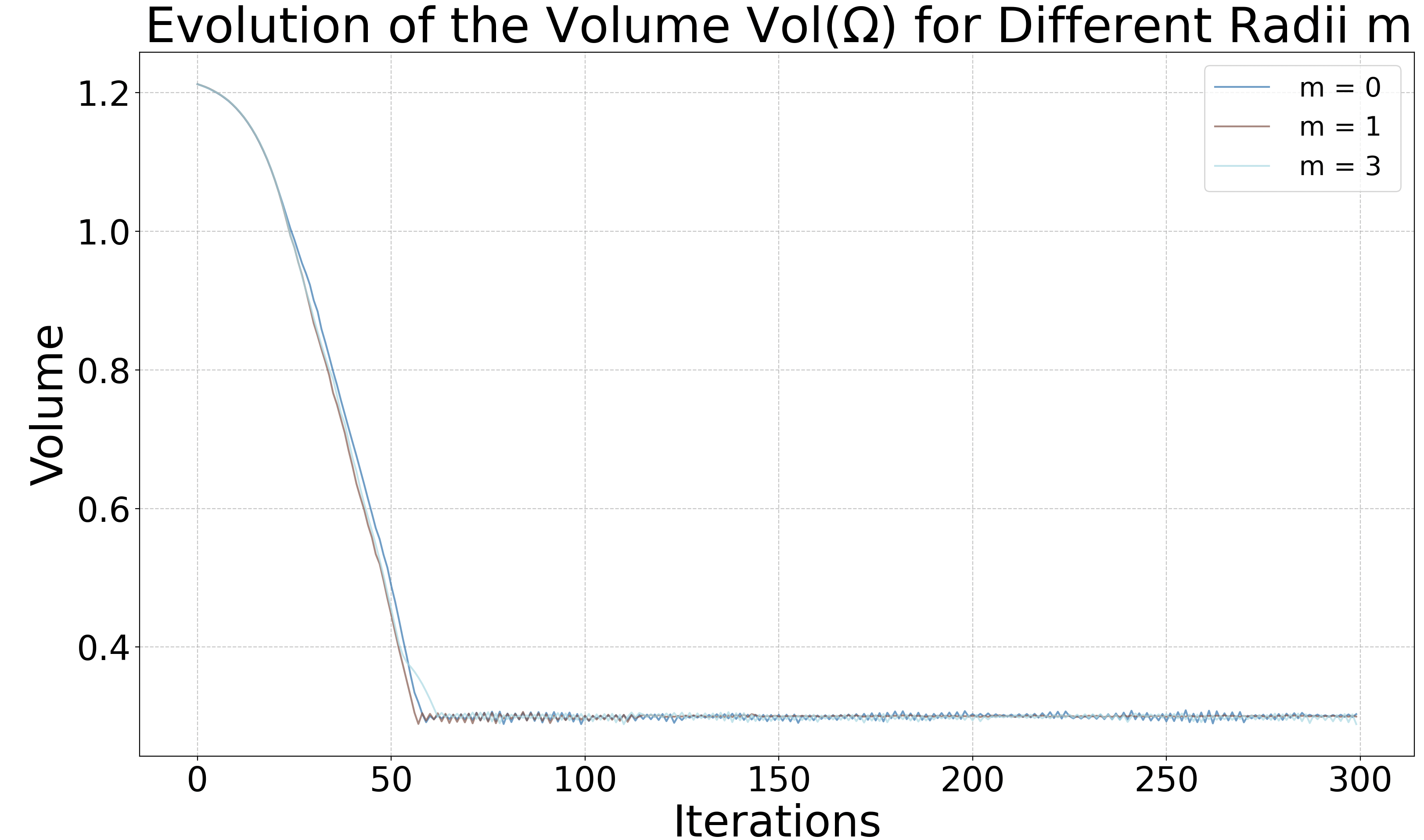}\put(2,5){\fcolorbox{black}{white}{b}}\end{overpic}
    \end{subfigure}
    \caption{\it Convergence histories of the objective function $D_W(\Omega,\lambda)$ and $\Vol(\Omega)$ of the distributionally robust optimal result of the T-shaped beam of \cref{fig.TstressdroV2} for $\sigma^2= 1e{-1}$.}
        \label{fig.cvTshape}
\end{figure}

The convergence histories of the computation are reported on \cref{fig.cvTshape}; the behavior of the objective function $\calD_{\text{W}}$ is again typical of the use of the stochastic descent algorithm: a local minimum is attained after a very noisy path, due to the use of Monte-Carlo integration for the evaluation of probabilistic integrals.

The values of the stress functional for the various optimized shapes under ideal conditions (i.e. $\xi = \xi^0$) are reported on \cref{tab.TshapestressoptiV2}. As expected, and like in the previous examples, these values deteriorate as the radius $m$ increases. Moreover, we note that for small variance $\sigma^2=1e{-3}$ and $m=0$, the stress within the distributionally robust structure is very close to that of the structure $\Omega^*_{\text{det}}$.

\begin{table}[ht]
\centering
\begin{tabular}{|c|c|c|c|}
\hline
\backslashbox{$\sigma^2$}{$m$} & $0$ & $1$ & $3$  \\
\hline
$10^{-0.5}$ & $2.12165$ & $2.32219$  & $2.81476$  \\
\hline
$10^{-1}$ & $2.02971$ & $2.27063$  & $2.27391$ \\
\hline
$10^{-3}$ & $1.99752$ & $2.01187$ & $2.02397$  \\
\hline
\end{tabular}
\caption{\it Values of the stress $S(\Omega,\xi^0)$ of the distributionally robust designs of the T-shaped beam considered in \cref{sec.stress2dSO} when the ideal load $\xi^0$ is applied.}
\label{tab.TshapestressoptiV2}
\end{table}

\subsection{Distributionally robust optimization of the shape of a 3d cantilever beam using a moment-based ambiguity set}\label{sec.3dcanti}

\begin{figure}[htbp]
\centering
	\begin{tabular}{ccc}
    \begin{minipage}{0.37\textwidth}
    \begin{overpic}[width= 1.0\textwidth]
        {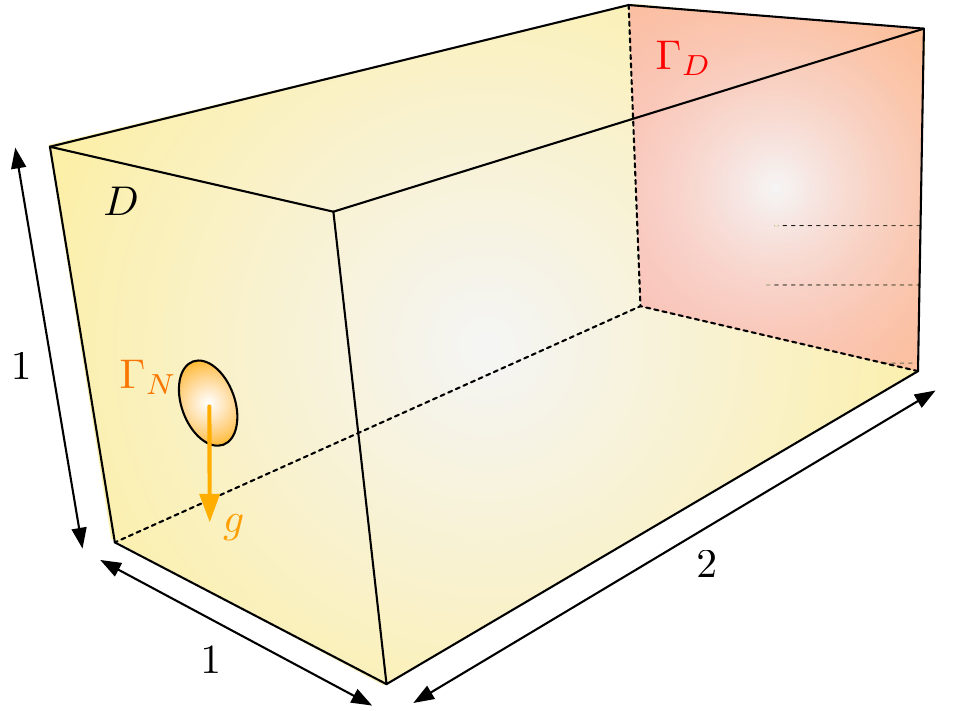}
        \put(2,5){\fcolorbox{black}{white}{a}}
    \end{overpic}
    \end{minipage}&
		\begin{minipage}{0.29\textwidth}
		\begin{overpic}[width=1.0\textwidth]{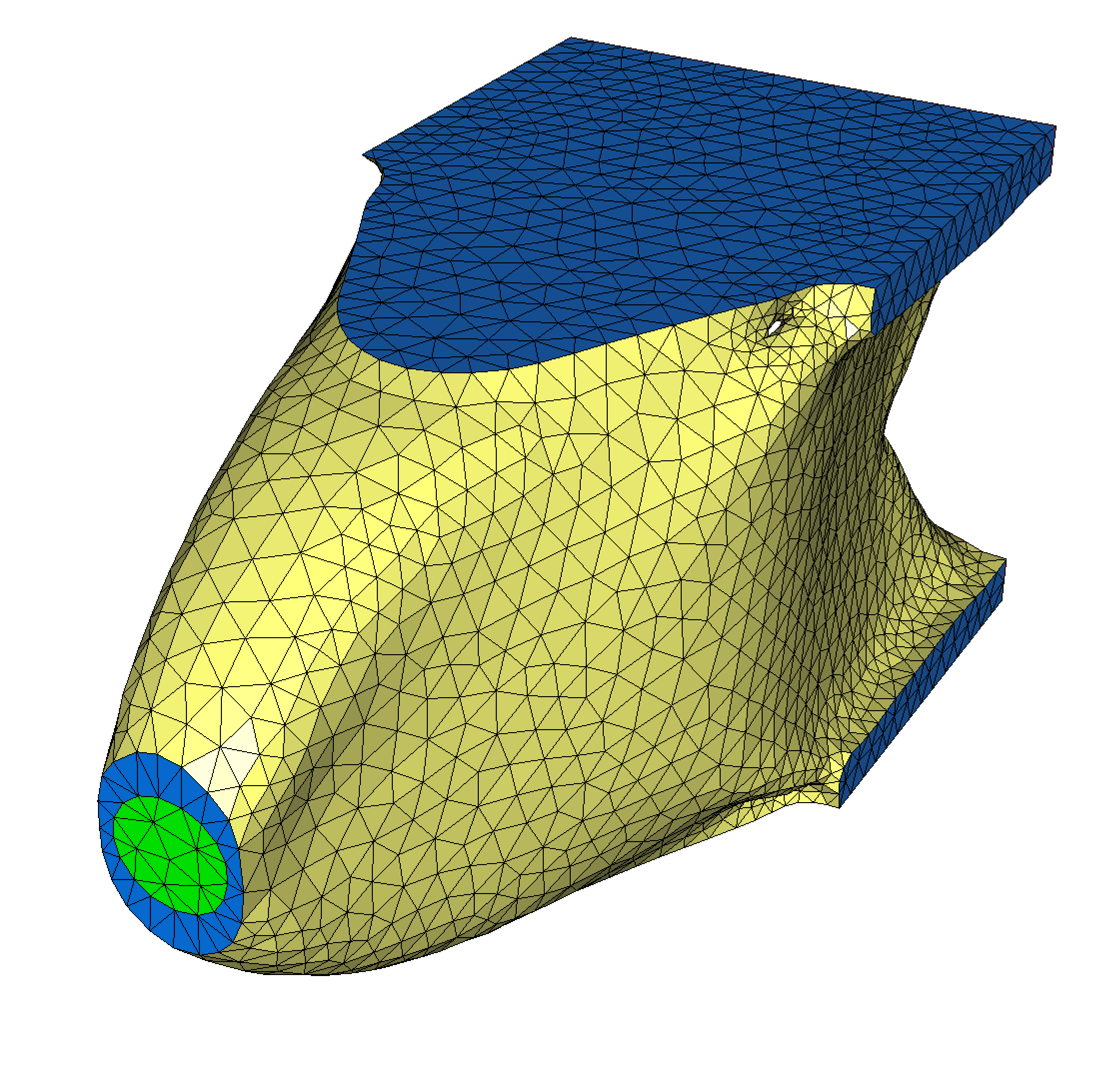}
		\put(2,5){\fcolorbox{black}{white}{b}}
		\end{overpic}
		\end{minipage} & 
		\begin{minipage}{0.32\textwidth}
        \begin{overpic}[width=1.0\textwidth]{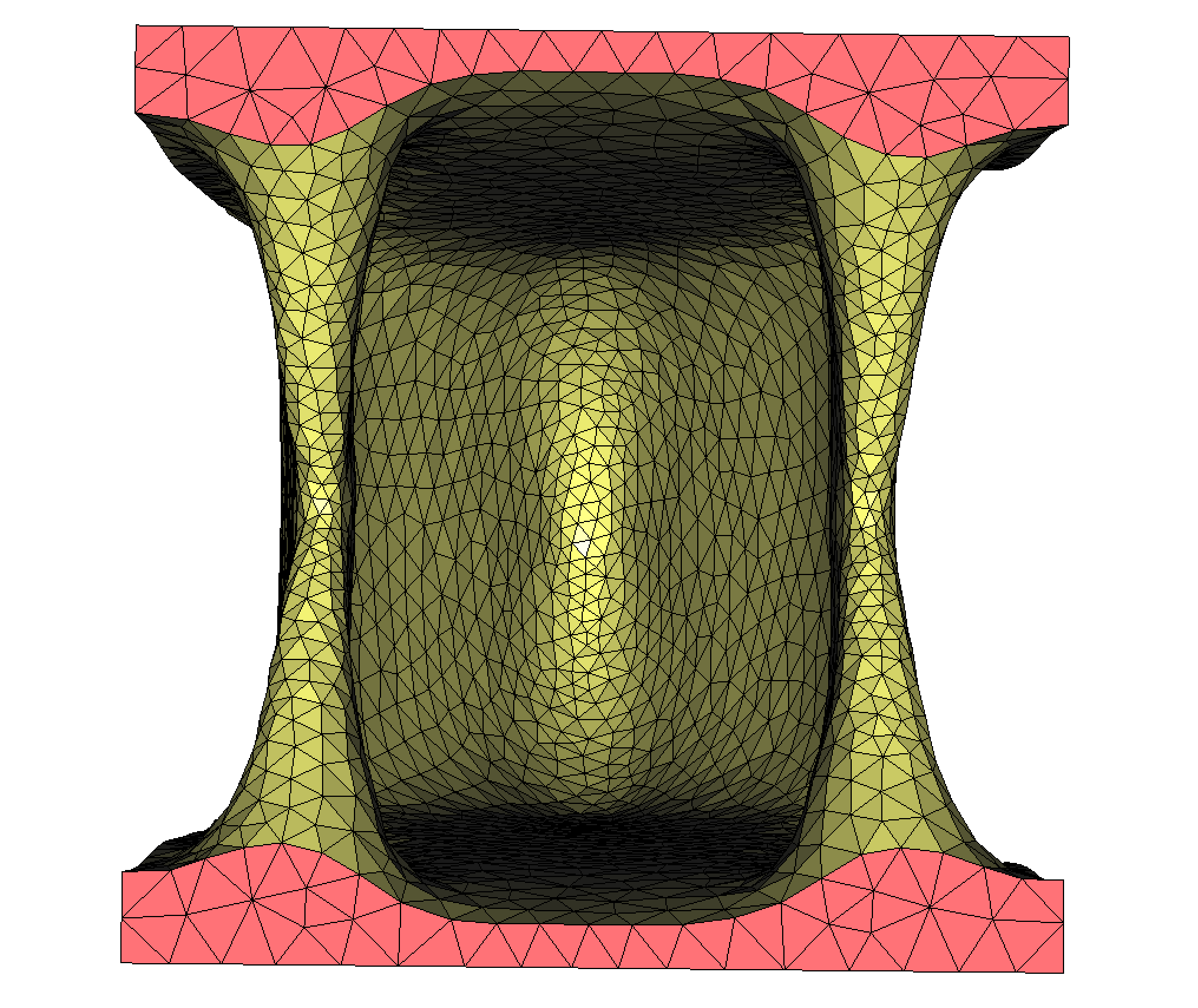}
		\put(2,5){\fcolorbox{black}{white}{c}}
		\end{overpic}
	\end{minipage}
	\end{tabular}
    \caption{\it (a) Setting of the shape optimization example of a 3d cantilever beam considered in \cref{sec.3dcanti}; (b) Front and (c) Back views of the optimized structure $\Omega^*_{\text{\rm det}}$ in the ideal situation where the applied load $\xi$ is perfectly known.}
    \label{fig.Det3DCanti}
\end{figure}

\noindent This section deals with the optimization of a 3d cantilever beam, as depicted on \cref{fig.Det3DCanti} (a). 
The structure under scrutiny is enclosed in a box $D$ with size $2 \times 1 \times 1$; it is clamped on one face $\Gamma_D$ of $\partial D$, and a load $\xi$ is applied on a small neighborhood $\Gamma_N$ of the center of the opposite face. Like in the previous \cref{sec.stress2dSO}, we denote by $u_{\Omega,\xi}$ the solution to \cref{eq.elasshape} when the load $\xi$ is applied.

In an ideal setting, $\xi$ is known perfectly, and it equals $\xi^0=(0,0,-1)$. We then consider the following problem:
\begin{equation}\label{eq.minC3dcanti}
    \min\limits_\Omega\; C(\Omega,\xi^0) \quad \text{ s.t. }\quad \Vol(\Omega)=V_T,
\end{equation}
where $C(\Omega,\xi^0)$ is the compliance \cref{eq.cplySO} of $\Omega$ and the volume target is set to $V_T = 0.45$.
For convenience, let us recall the following formula for the shape derivative of the compliance, see e.g. \cite{allaire2020shape,allaire2004structural}:
$$C^\prime(\Omega,\xi)(\theta) = -\int_\Gamma Ae(u_{\Omega,\xi}) : e(u_{\Omega,\xi}) \:\theta\cdot n \:\d s.$$
The problem \cref{eq.minC3dcanti} is solved thanks to the numerical strategy presented in \cref{sec.numSO}. Each mesh of $D$ featured during the computation is composed of about $15,000$ vertices ($\approx 80,000$ tetrahedra). The resulting shape $\Omega^*_{\text{det}}$ after 200 iterations ($\approx $ 2h 30 mn of computation) is depicted on \cref{fig.Det3DCanti} (b,c).

Now, we assume that the applied load $\xi$ is uncertain: it belongs to a sufficiently large ball $\Xi \subset \R^3$ and its probability law is also unknown. The only information at hand about the latter is about its mean value $\mu^0 \in \R^3$ and covariance matrix $\Sigma^0 \in \S^3_+$:
$$\mu^0 = \xi^0 \quad \text{ and } \Sigma^0 = \sigma^2 \I, \text{ where }\sigma^2 = 0.01.$$
As proposed in \cref{sec.momdro}, from these data, we now consider the following distributionally robust version of this problem:
$$ \min\limits_{\Omega} \left(\sup\limits_{\Q \in \calA_{\text{M}}} \int_\Xi C(\Omega,\xi) \:\d\Q(\xi) - \e H(\Q)\right) \: \text{ s.t. } \Vol(\Omega) = V_T,$$
where $\calA_{\text{M}}$ is the moment-based ambiguity-set defined in \cref{eq.calAmoments} and the entropy $H(\Q)$ of a probability law $\Q \in \calP(\Xi)$ equals \cref{eq.entropymoments}. 
The penalty parameter $\e$ is set to $0.001$. This problem has the following equivalent expression:
\begin{multline*}
	\min\limits_{\Omega, \:\lvert  \tau  \lvert\leq 1, \atop
	{\lambda \geq 0 , {S \in \mathbb{S}_+^3} } }\calD_{\text{M}}(\Omega,\lambda
	,\tau,S), \text{ where }\\[-1em] \calD_{\text{M}}(\Omega,\lambda,\tau,S):=
	\lambda m_{1}- \lambda \tau \cdot \mu_{0}+ m_{2}S : \Sigma_{0}+ \e \log
	\left( \int_{\Xi}\left(e^{\frac{C(\Omega,\xi) + \lambda \tau \cdot \xi
	- S : (\xi-\mu_{0}) \otimes (\xi-\mu_{0})}{\e}}\right) \:\d \Q_{0}(\xi
	)\right).
\end{multline*}

\begin{figure}[!ht]
\centering
	\begin{tabular}{cc}
		\begin{minipage}{0.36\textwidth}
		\begin{overpic}[width=0.9\textwidth]{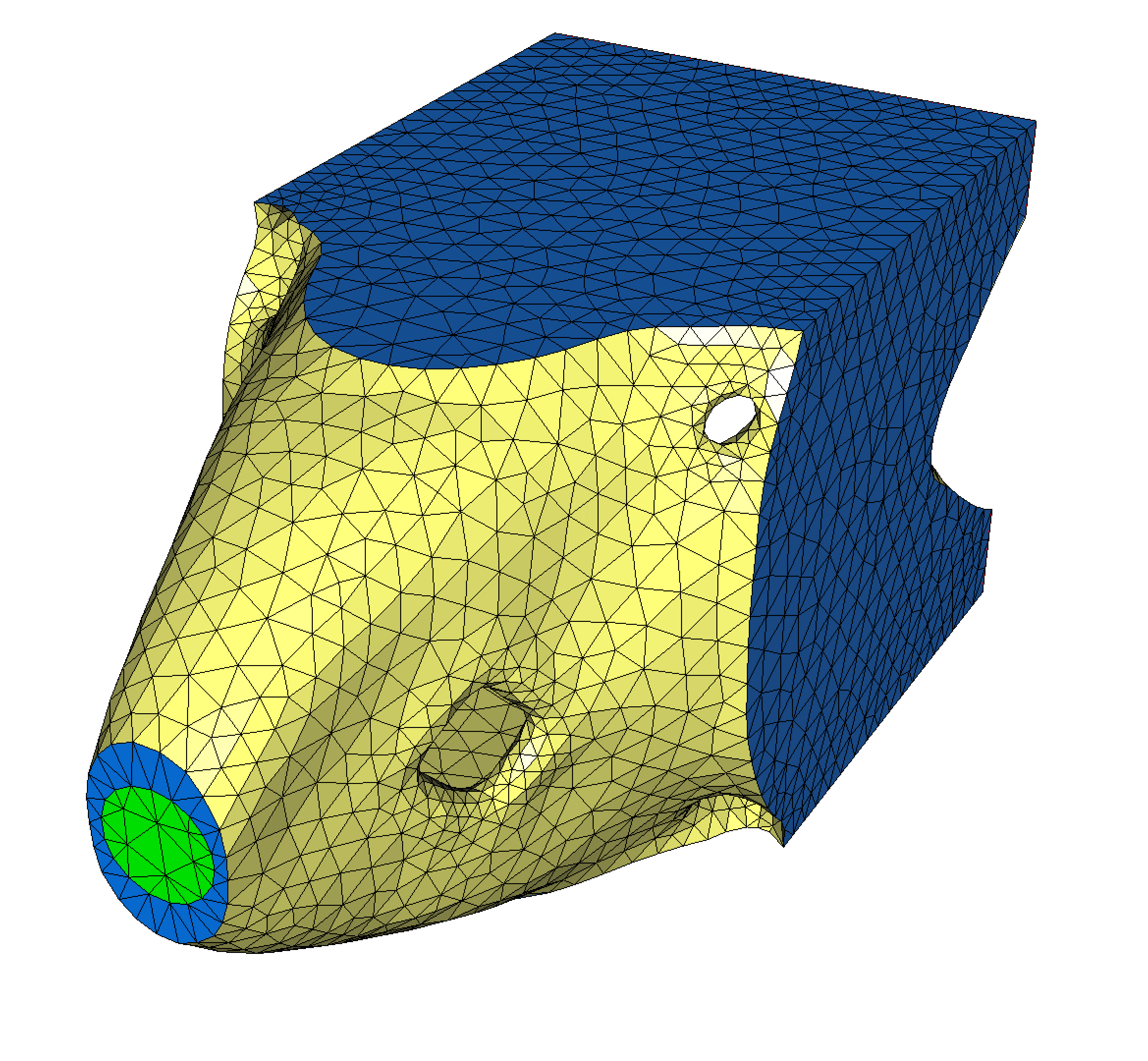}
		\put(2,5){\fcolorbox{black}{white}{a}}
		\end{overpic}
		\end{minipage} & 
		\begin{minipage}{0.4\textwidth}
        \begin{overpic}[width=0.9\textwidth]{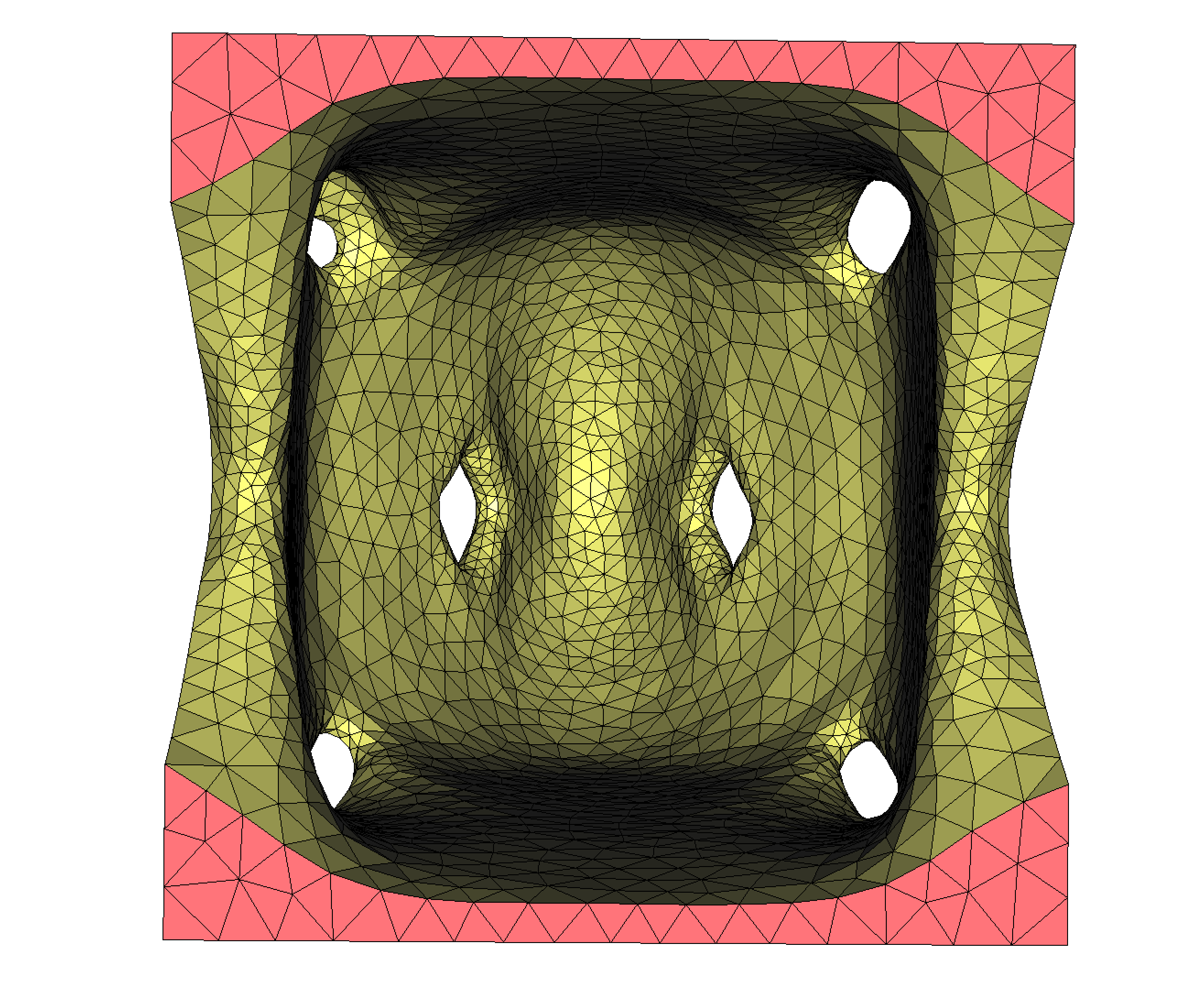}
		\put(2,5){\fcolorbox{black}{white}{a}}
		\end{overpic}
	\end{minipage}
	\end{tabular}
    \par
    \medskip
    \begin{tabular}{cc}
		\begin{minipage}{0.36\textwidth}
		\begin{overpic}[width=0.9\textwidth]{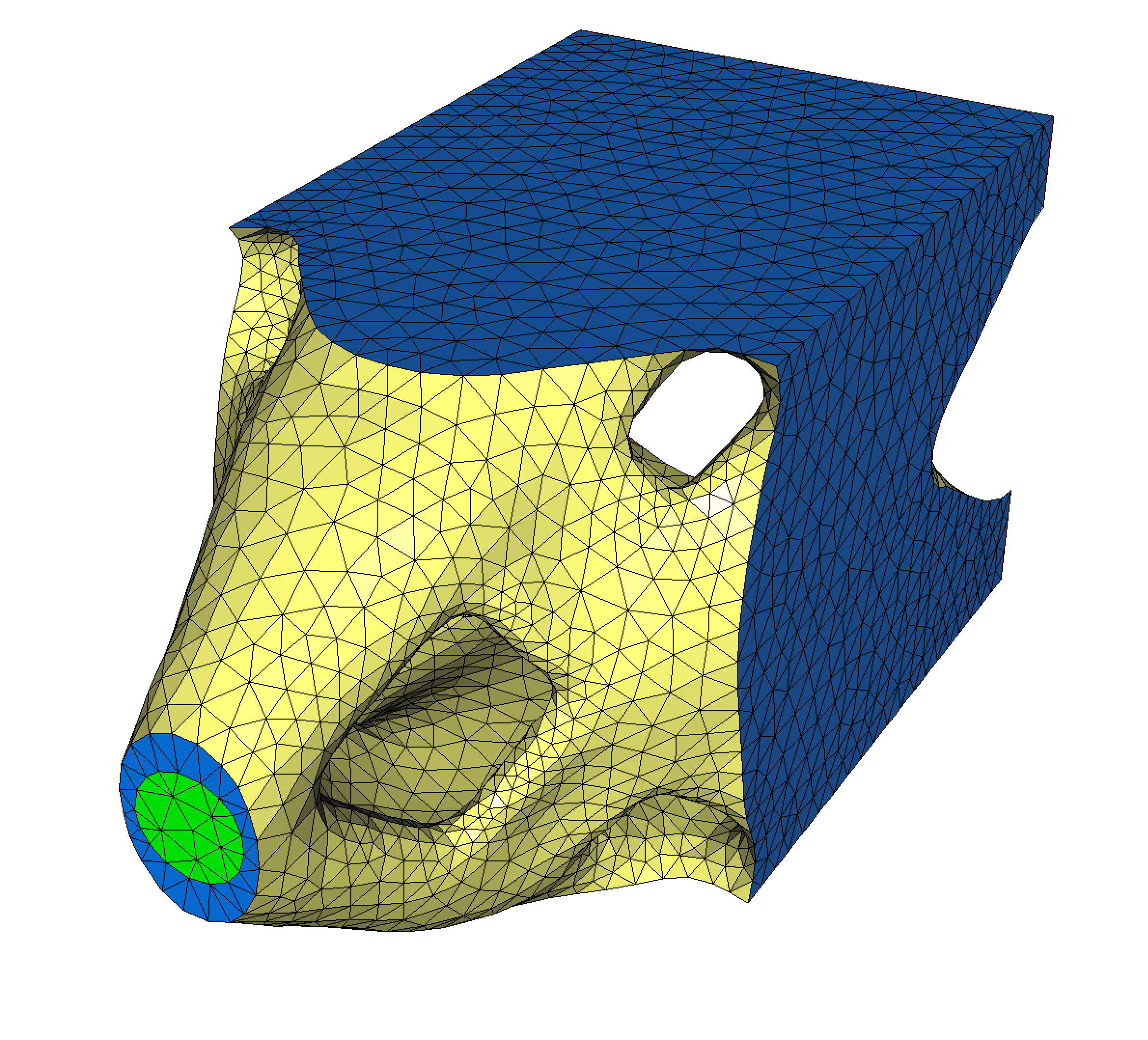}
		\put(2,5){\fcolorbox{black}{white}{b}}
		\end{overpic}
		\end{minipage} & 
		\begin{minipage}{0.4\textwidth}
        \begin{overpic}[width=0.9\textwidth]{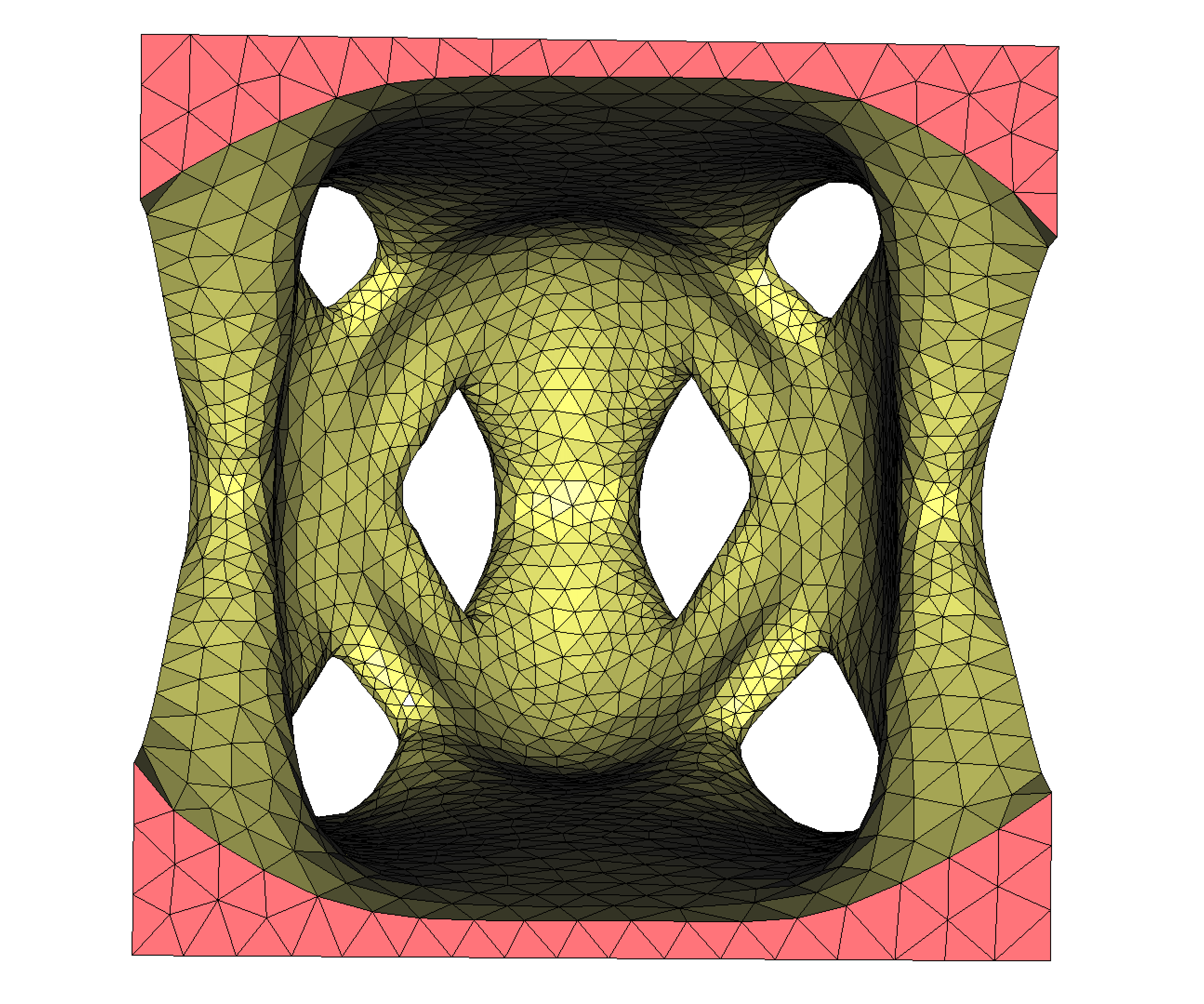}
		\put(2,5){\fcolorbox{black}{white}{b}}
		\end{overpic}
	\end{minipage}
	\end{tabular}
    \par
    \medskip
    \begin{tabular}{cc}
		\begin{minipage}{0.37\textwidth}
		\begin{overpic}[width=0.9\textwidth]{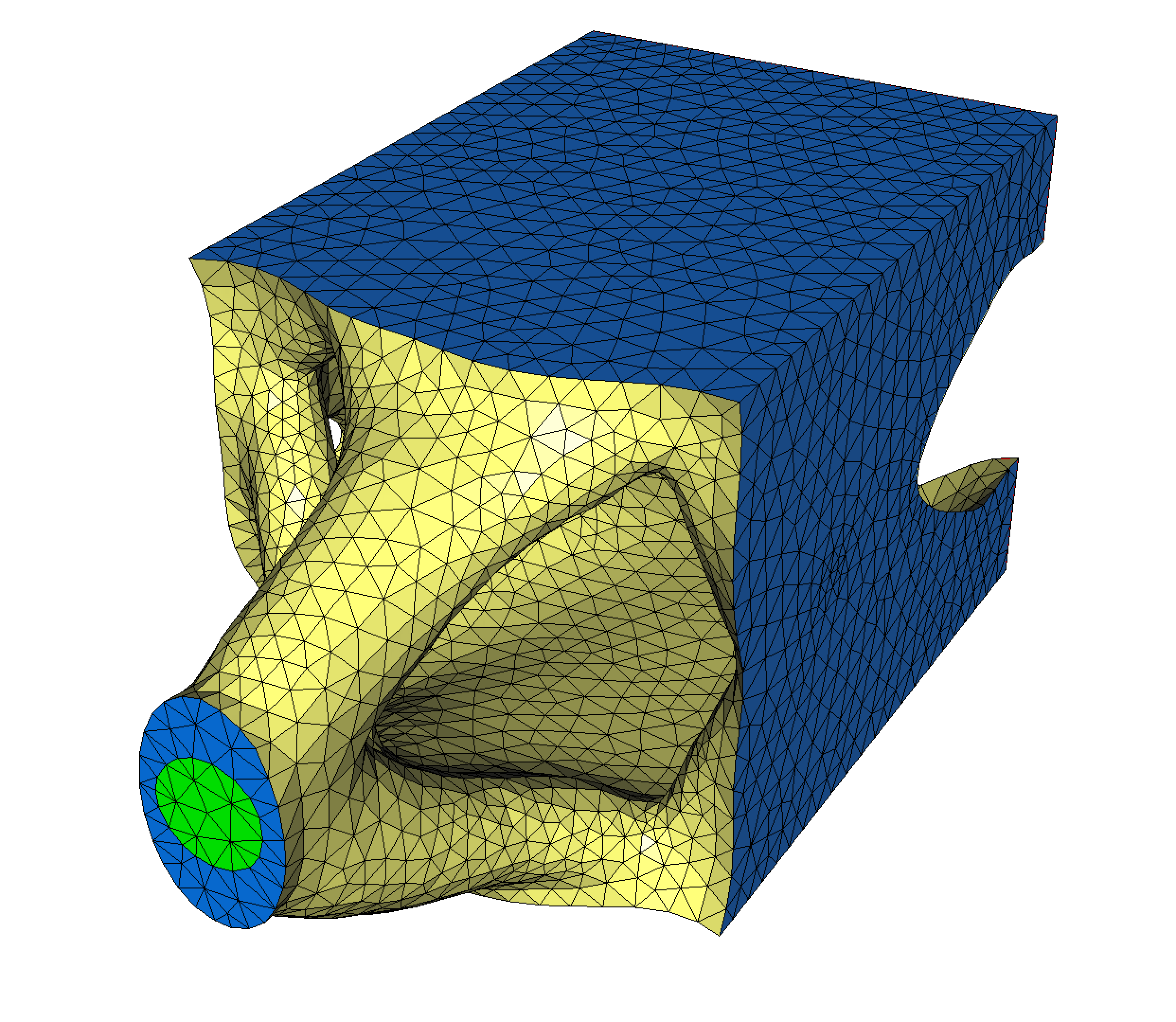}
		\put(2,5){\fcolorbox{black}{white}{c}}
		\end{overpic}
		\end{minipage} & 
		\begin{minipage}{0.4\textwidth}
        \begin{overpic}[width=0.9\textwidth]{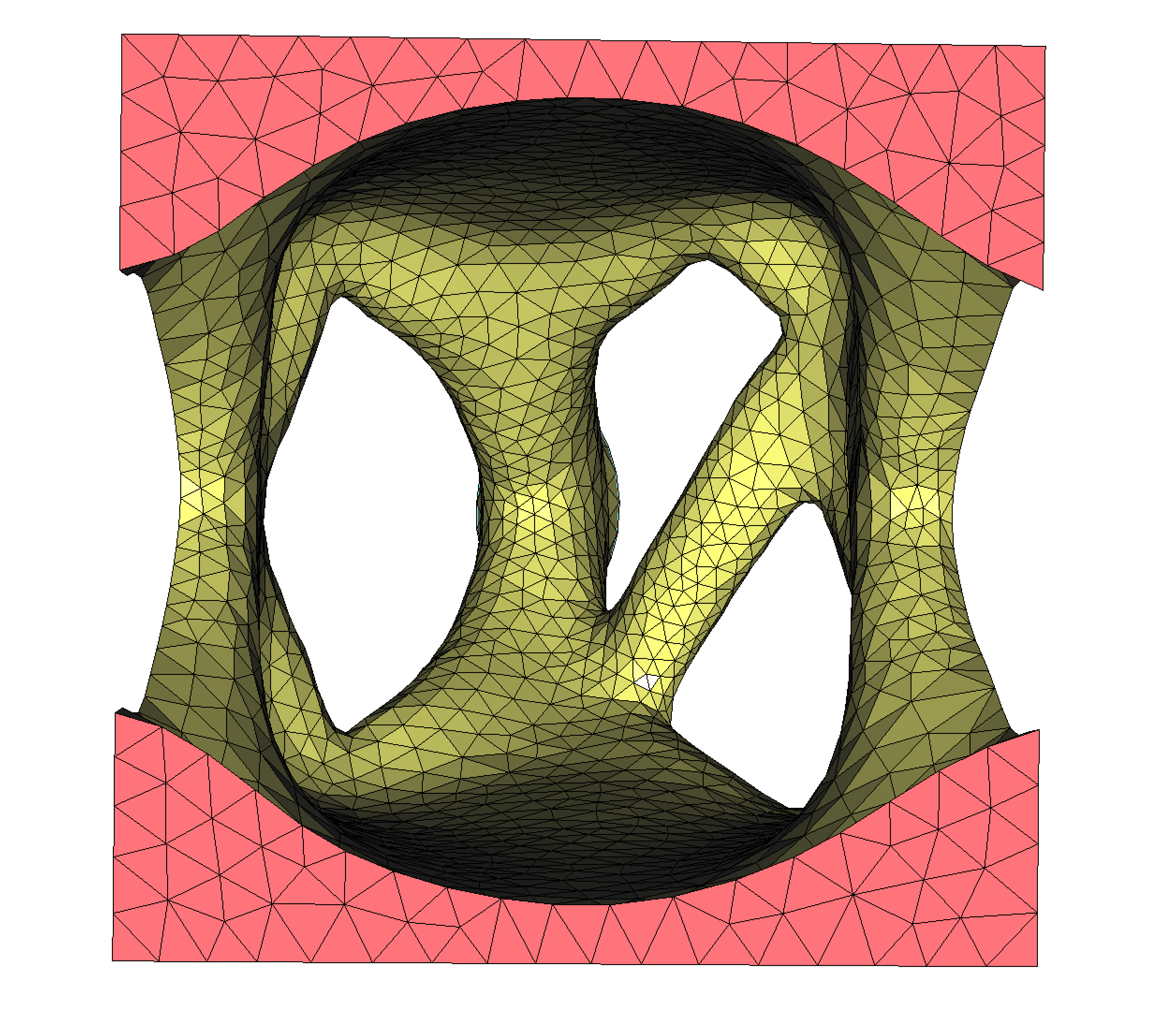}
		\put(2,5){\fcolorbox{black}{white}{c}}
		\end{overpic}
	\end{minipage}
	\end{tabular}
     \par
    \medskip
    \begin{tabular}{cc}
		\begin{minipage}{0.37\textwidth}
		\begin{overpic}[width=0.9\textwidth]{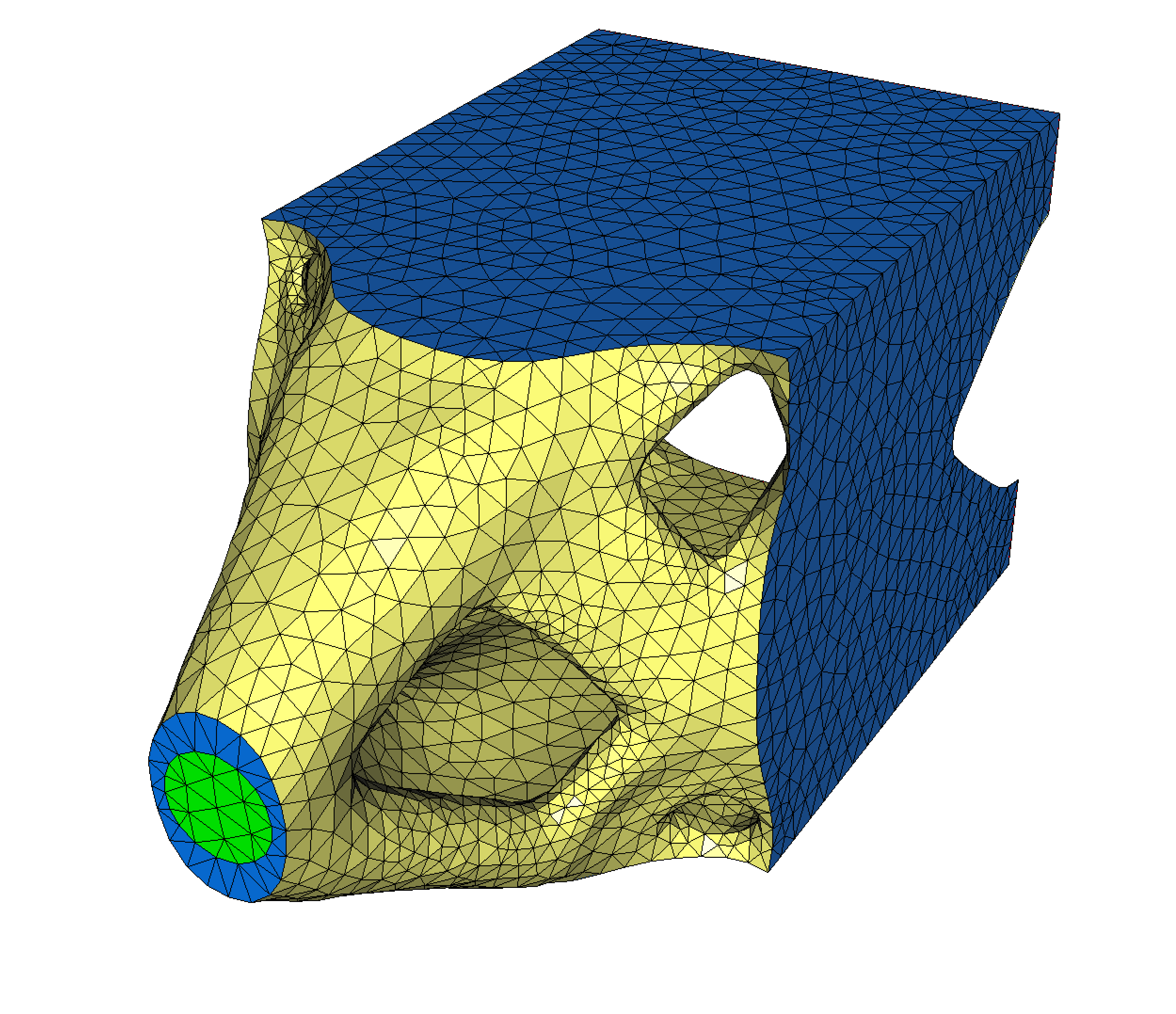}
		\put(2,5){\fcolorbox{black}{white}{d}}
		\end{overpic}
		\end{minipage} & 
		\begin{minipage}{0.4\textwidth}
        \begin{overpic}[width=0.9\textwidth]{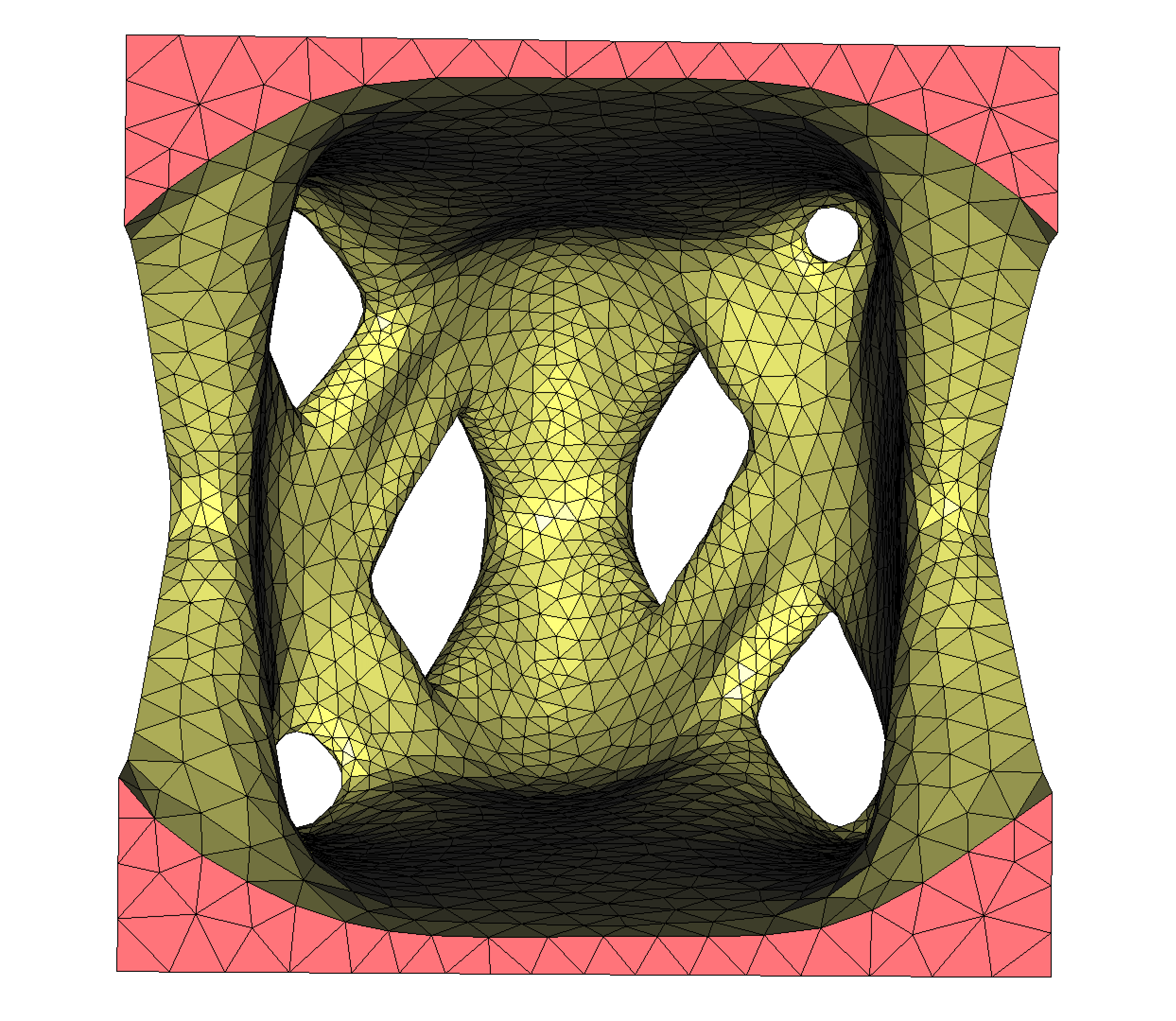}
		\put(2,5){\fcolorbox{black}{white}{d}}
		\end{overpic}
	\end{minipage}
	\end{tabular}
    \caption{\it Optimized shapes of the 3d cantilever of \cref{sec.3dcanti} under distributional uncertainties using a moment ambiguity set; front and back views for the parameters (a) $m_1=0, m_2=1$; (b) $m_1=1, m_2=1$; (c) $m_1 = 2, m_2 = 1$, and (d) $m_1 = 5, m_2 = 1$.}
    \label{fig.resmom3dcanti1}
    \end{figure}

\begin{figure}[!ht]
    \centering
    \begin{tabular}{cc}
		\begin{minipage}{0.4\textwidth}
		\begin{overpic}[width=0.9\textwidth]{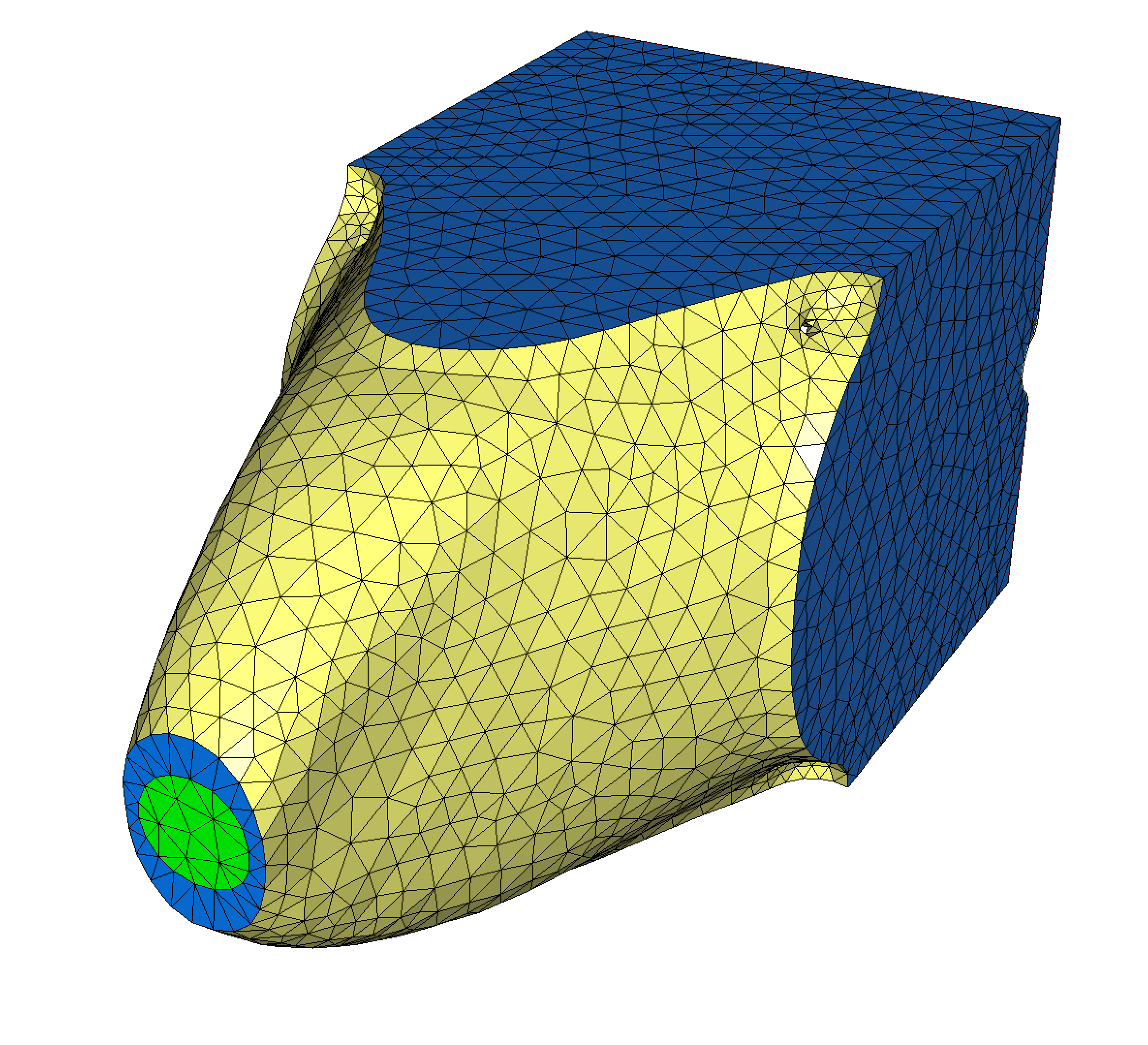}
		\put(2,5){\fcolorbox{black}{white}{e}}
		\end{overpic}
		\end{minipage} & 
		\begin{minipage}{0.45\textwidth}
        \begin{overpic}[width=0.9\textwidth]{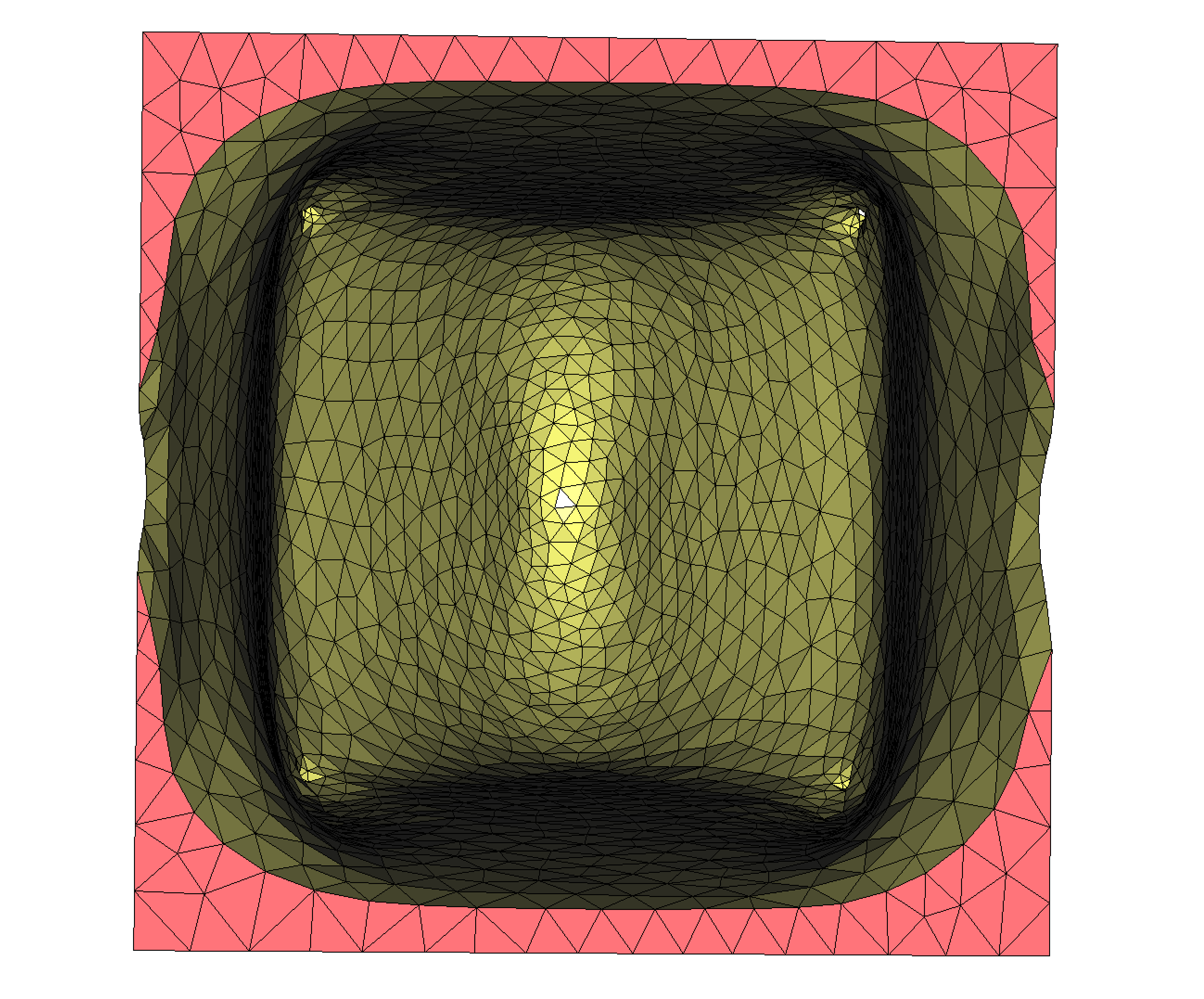}
		\put(2,5){\fcolorbox{black}{white}{e}}
		\end{overpic}
	\end{minipage}
	\end{tabular}
     \par
    \medskip
    \begin{tabular}{cc}
		\begin{minipage}{0.4\textwidth}
		\begin{overpic}[width=0.9\textwidth]{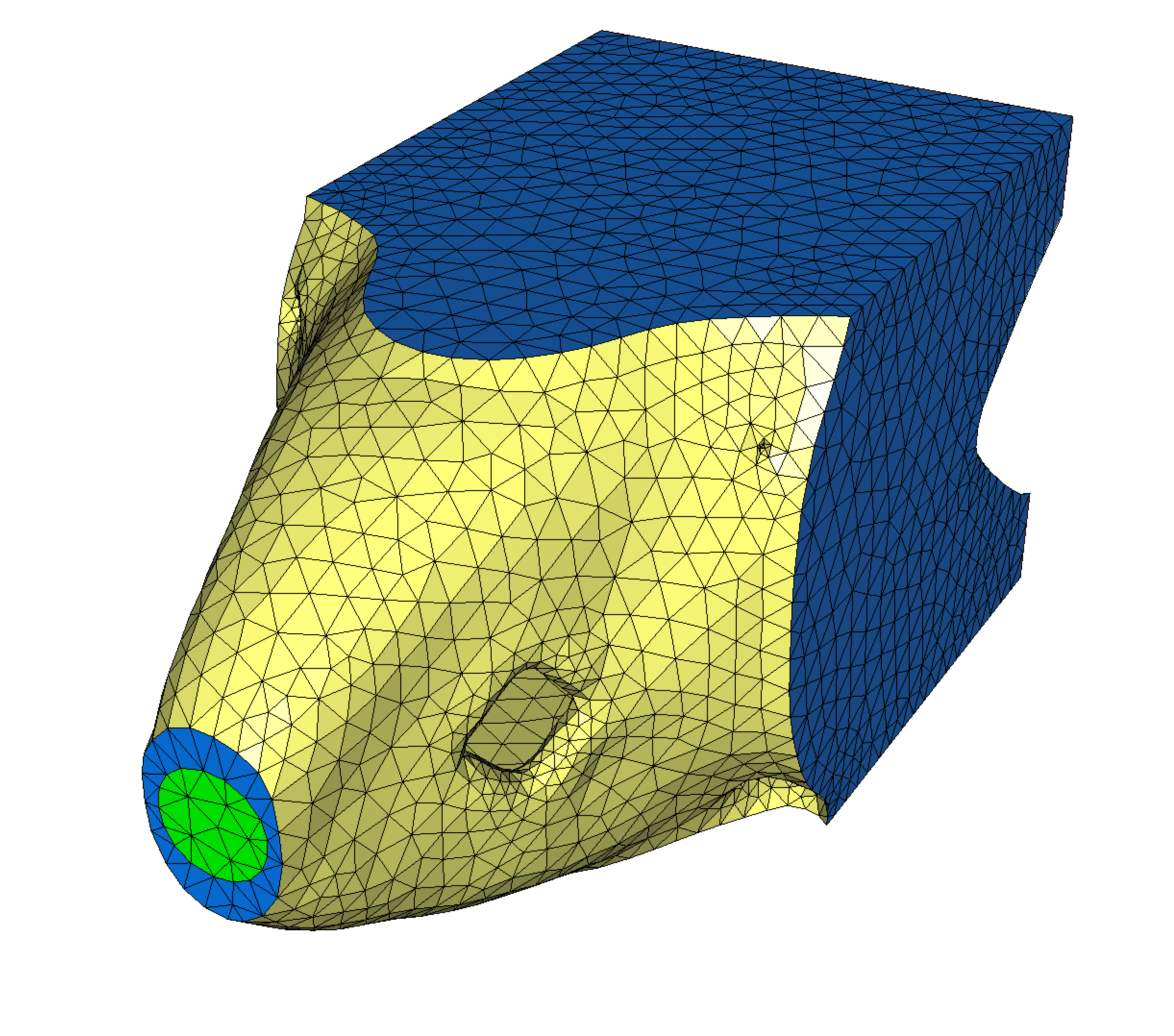}
		\put(2,5){\fcolorbox{black}{white}{f}}
		\end{overpic}
		\end{minipage} & 
		\begin{minipage}{0.45\textwidth}
        \begin{overpic}[width=0.9\textwidth]{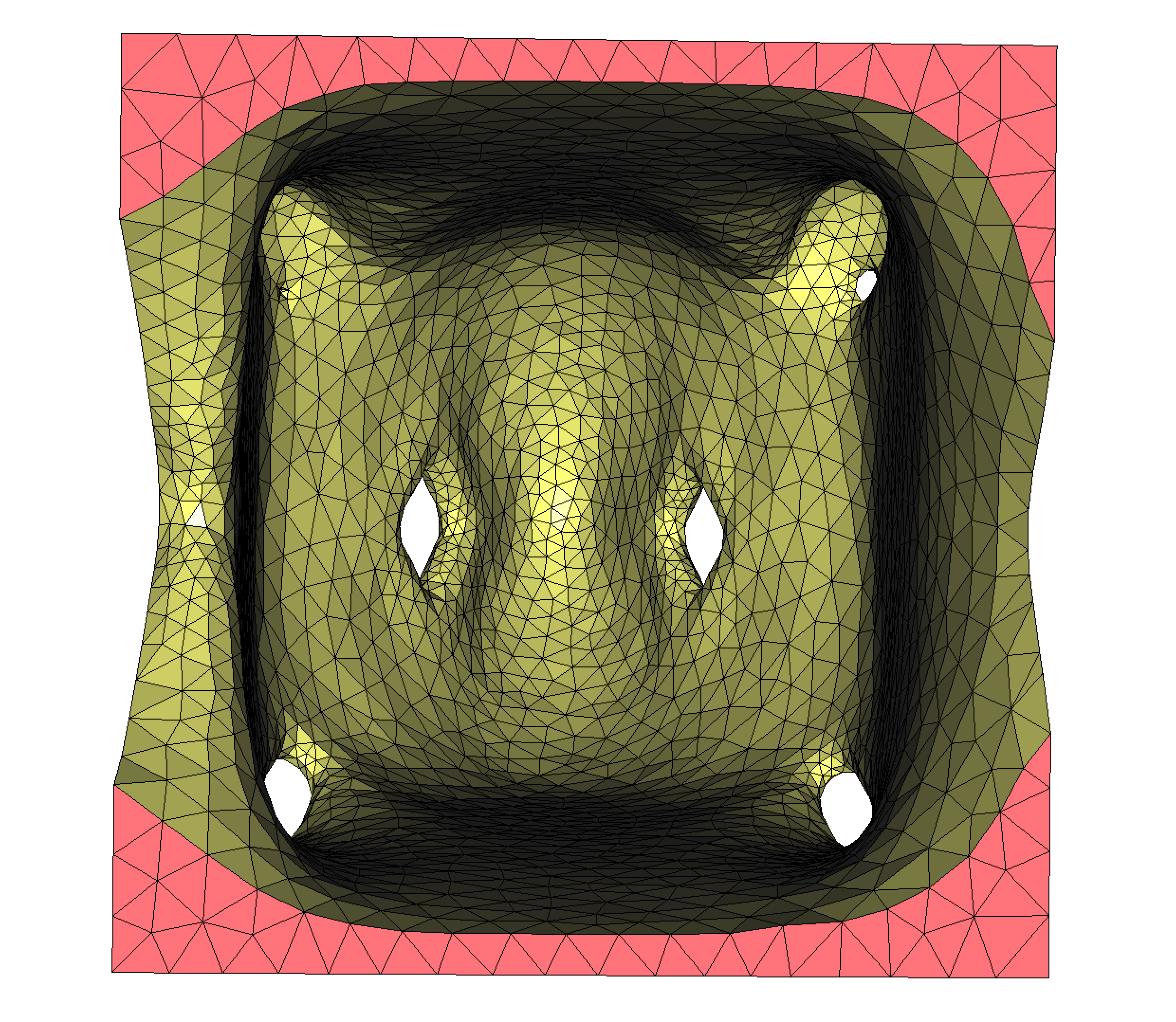}
		\put(2,5){\fcolorbox{black}{white}{f}}
		\end{overpic}
	\end{minipage}
	\end{tabular}
     \par
    \medskip
    \begin{tabular}{cc}
		\begin{minipage}{0.4\textwidth}
		\begin{overpic}[width=0.9\textwidth]{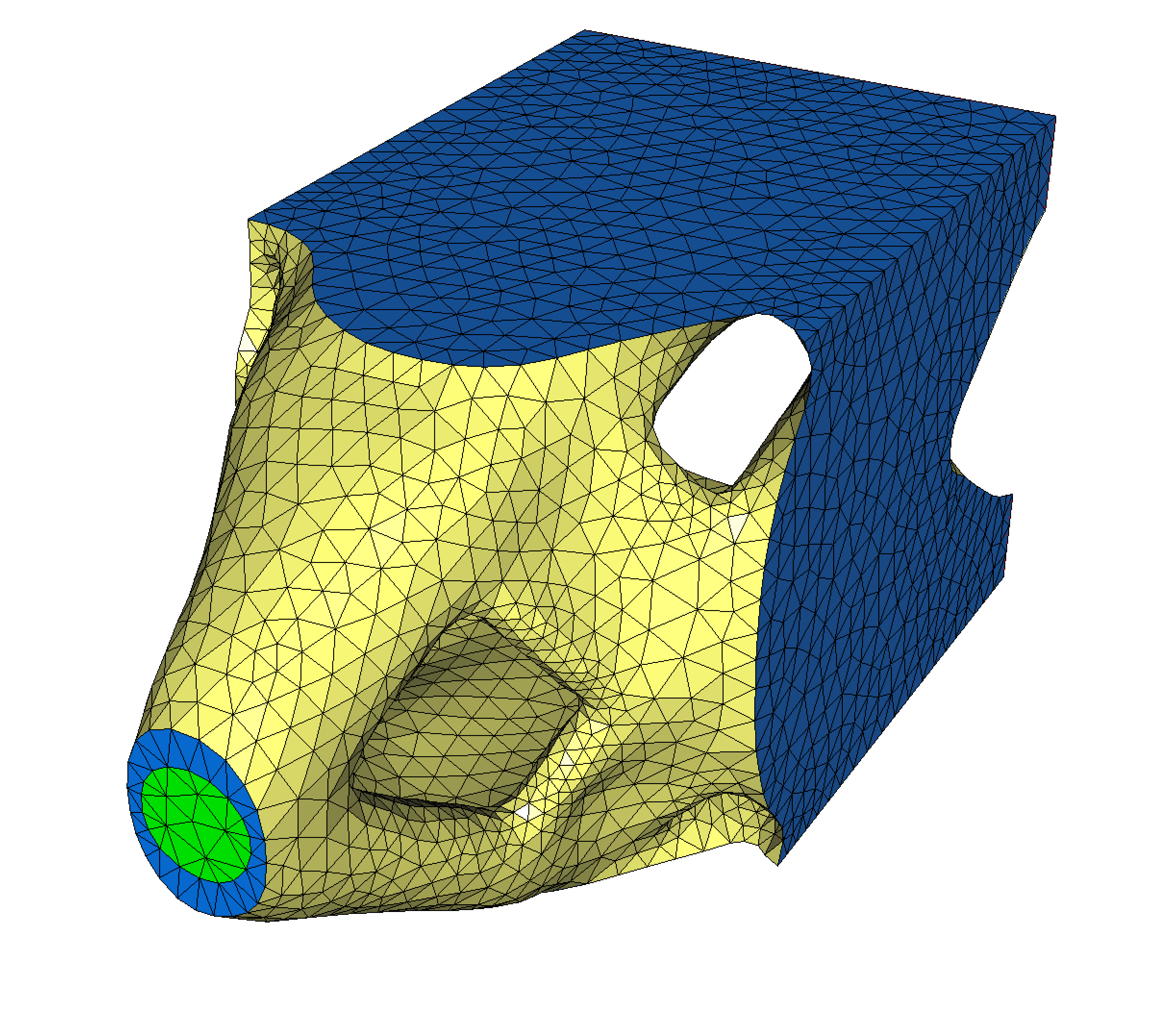}
		\put(2,5){\fcolorbox{black}{white}{g}}
		\end{overpic}
		\end{minipage} & 
		\begin{minipage}{0.45\textwidth}
        \begin{overpic}[width=0.9\textwidth]{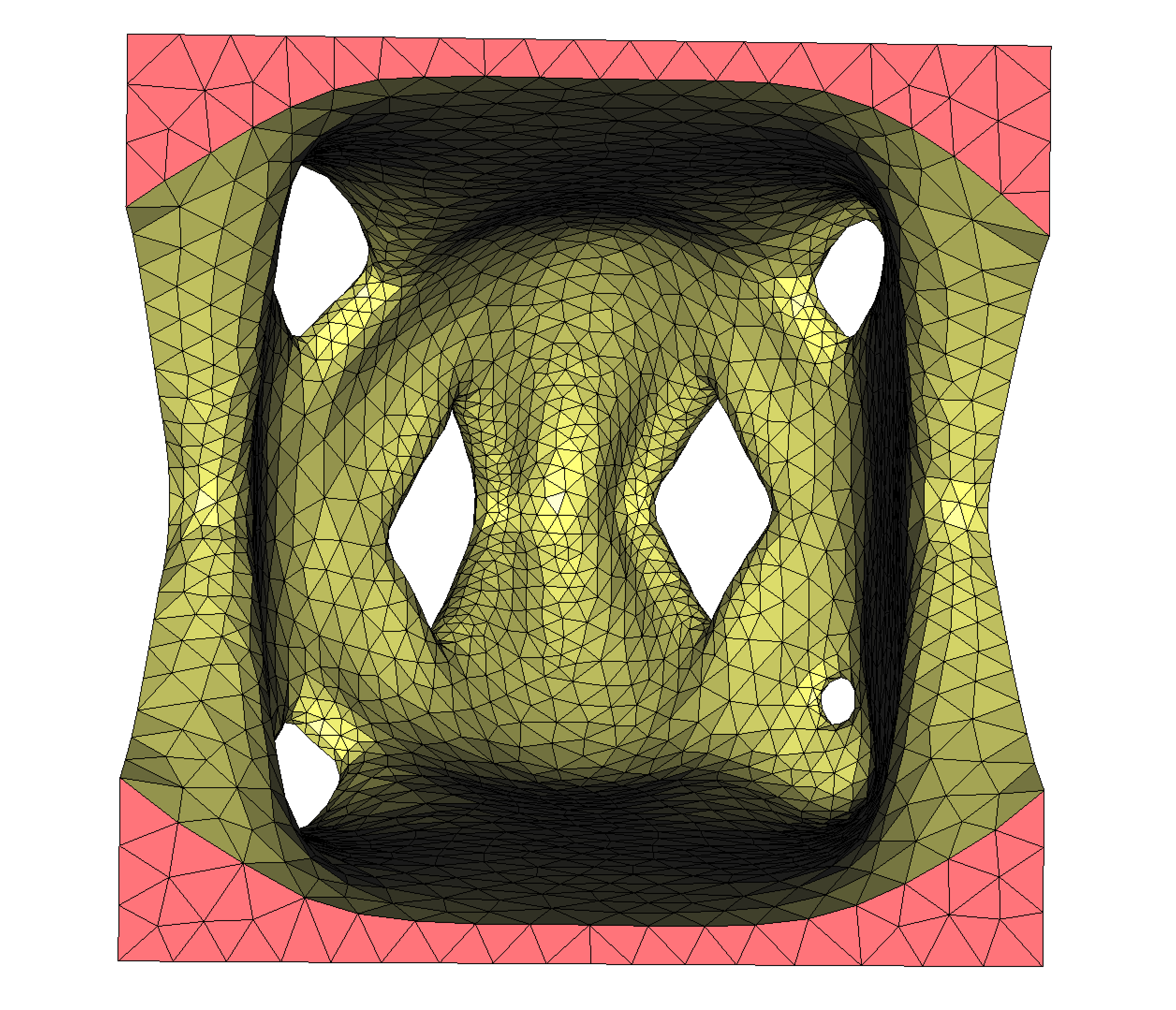}
		\put(2,5){\fcolorbox{black}{white}{g}}
		\end{overpic}
	\end{minipage}
	\end{tabular}
    \caption{\it Optimized shapes of the 3d cantilever of \cref{sec.3dcanti} under distributional uncertainties using a moment ambiguity set; front and back views for the parameters (e) $m_1=0, m_2=2$; (f) $m_1 = 2, m_2 = 1$, and (g) $m_1 = 5, m_2 = 5$.}
    \label{fig.resmom3dcanti2}
\end{figure}

We solve this problem with the numerical strategy presented in \cref{sec.numSO} for several values of the parameters $m_1$, $m_2$. The optimized shapes obtained after 300 iterations (corresponding to approximately 3h 30 mn of computation) are depicted on \cref{fig.resmom3dcanti1,fig.resmom3dcanti2}. The values of their compliance in ideal circumstances are reported on \cref{tab.3DCantiMomentopti}. 

\begin{table}[ht]
    \centering
    \resizebox{\textwidth}{!}{%
    \begin{tabular}{|c|c|c|c|c|c|c|c|c|}
        \hline
              & $\Omega^*_{\text{det}}$  & $\Omega^*_{m_{1}=0, m_{2}=1}$       & $\Omega^*_{m_{1}=1, m_{2}=1}$       & $\Omega^*_{m_{1}=2, m_{2}=1}$   & $\Omega^*_{m_{1}=5, m_{2}=1}$    & $\Omega^*_{m_{1}=0, m_{2}=2}$ &  $\Omega^*_{m_{1}=2, m_{2}=2}$     & $\Omega^*_{m_{1}=5, m_{2}=5}$         \\
        \hline
        $C(\Omega,\xi^0)$ & 0.0907753 & 0.0948306 & 0.0913598
 & 0.108606 & 0.0979698 & 0.0903963 & 0.0923416 & 0.0978554 \\
        \hline
    \end{tabular}
    }
    \caption{\it Values of the compliance $C(\Omega,\xi^0)$ of the distributionally robust designs of the 3d cantilever beam obtained in \cref{sec.3dcanti} under ideal conditions.}
    \label{tab.3DCantiMomentopti}
\end{table}

The optimized structures are not symmetric, which can be attributed to the unstructured nature of the mesh and the use of the stochastic gradient algorithm.
Also note that, contrary to intuition, the values of the compliance of the obtained designs in ideal conditions are not strictly increasing functions of the parameters $m_1$, $m_2$ defining the ambiguity set $\calA_{\text{M}}$. This phenomenon can be attributed to the fact that shape optimization problems of the form \cref{eq.minC3dcanti} generally have multiple local minima, so that the optimization path (and so the resulting optimized design) is generally very sensitive to the initial design, the mesh size, the optimization parameters, etc. In the present situation , this trend is even further exacerbated by the use of a stochastic descent algorithm. This intuition is confirmed by unreported experiments, that can be found in \cite{prando2025distributionally}.

\FloatBarrier
\section{\textbf{Conclusion and Perspectives}}
\label{sec.concl}

\noindent
This article was concerned with exploring the idea of distributional robustness in the context of shape and topology optimization: acknowledging that the physical parameters characterizing the behavior of real-life devices are inherently known with uncertainty, and that so is their probability law, we have proposed tractable robust formulations for a wide range of optimal design problems, that incorporate a degree of awareness about these sources of uncertainty. These formulations are not limited to a specific optimal design context and they can be used in both frameworks of density-based topology optimization and geometric shape optimization, in two and three space dimensions. 
As demonstrated by the numerical results, one key feature of distributionally robust optimal design formulations is their ability to anticipate scenarios ``close'', but different from the observations used for the empirical reconstruction of the probability law of the uncertain parameters.

This work paves the way to research avenues of various kinds. From the theoretical viewpoint, we would like to justify rigorously that our tractable distributionally robust formulations, based on ``small'' penalizations by suitable entropy functionals, are consistent with their exact counterparts, as the penalty parameter tends to $0$. 
From the numerical viewpoint, an obvious perspective concerns the choice of the optimization algorithm. The stochastic gradient indeed produces highly oscillatory trajectories, which can hinder convergence and degrade the quality of the final solution, particularly in high-variance settings. A different optimization algorithm, such as ADAM \cite{kingma2017adam}, could help to reduce the variance of descent direction estimates.
Furthermore, the naive sampling strategy used in this work for evaluating probabilistic integrals could be improved, for instance by using variance reduction techniques such as importance sampling.
From the modeling viewpoint, manifold physical situations suffer from uncertainty over physical parameters, and we wish to apply our techniques in e.g. thermal conduction, fluid mechanics or electromagnetism, to name a few. 
On a different note, the particular situation where the geometry of shapes itself is uncertain has been ignored in this work, in spite of its ubiquity and importance in applications. We believe that, with the help of the linearization techniques similar to those introduced in our former works \cite{allaire2014linearized,allaire2015deterministic}, it is possible to give a computationally affordable formulation of such problems. A foretaste of this important avenue for future research can be consulted in the Ph.D. manuscript \cite{prando2025distributionally}.

\par
\bigskip
\noindent
\textbf{Acknowledgements.} The work of C.D, J.P. and B.T. is partially
supported by the projects ANR-24-CE40-2216 STOIQUES and ANR-22-CE46-0006
StableProxies, financed by the French Agence Nationale de la Recherche (ANR).
This research has benefitted from multiple discussions with F. Iutzeler (Institut de Math\'ematiques de Toulouse), whose expertise and insights are gratefully acknowledged.

\appendix

\section{\textbf{Hint of Proof of \cref{prop.formulaJdir}}}\label{app.wdro}

\noindent Let $A := \sup\limits_{\Q \in \calA_{\text{W}}} \int_\Xi f(\xi) \:\d \Q(\xi)$ be the maximum value at the left-hand side of the relation to be established \cref{eq.reformspwdro}. We classically start by rewriting $A$ with the help of a Lagrange multiplier for the constraint $\Q \in \calA_{\text{W}}$; it holds:
$$A = \sup\limits_{\Q \in \calP(\Xi)}\inf\limits_{\lambda \geq 0}\left( \int_{\Xi}f(\zeta) \: \d \Q(\zeta) + \lambda \left( m - W_{\e}(\P,\Q)\right)\right).$$
In fact, the inner infimum in this formula equals $0$ when $\Q$ is such that $W_\e(\P,\Q) \leq m$ and $-\infty$ otherwise. 

Now, formally interchanging the infimum and supremum in the above formula, and invoking the definition \cref{eq.RegWass} of the regularized Wasserstein distance $W_{\e} (\cdot,\cdot)$, we obtain:
$$\begin{array}{>{\displaystyle}cc>{\displaystyle}l}
			A & = & \inf\limits_{\lambda \geq 0} \left( \lambda m + \sup\limits_{\Q \in \calP(\Xi)} \sup\limits_{\pi \in \calP(\Xi \times \Xi) \atop \pi_1 = \P, \:\pi_2 = \Q} \left(\int_\Xi f(\zeta) \:\d \Q(\zeta) - \lambda \int_{\Xi \times \Xi} c(\xi,\zeta) \:\d \pi(\xi,\zeta) - \lambda \e H(\pi)\right)\right) \\[1em]
			  & = & \inf\limits_{\lambda \geq 0} \left( \lambda m + \sup\limits_{\pi \in \calP(\Xi \times \Xi) \atop \pi_1 = \P} \int_{\Xi \times \Xi} \left( f(\zeta) - \lambda c(\xi,\zeta) \right)\:\d \pi(\xi,\zeta) - \lambda \e H(\pi)\right).
		\end{array}
$$
We now proceed to reformulate the inner supremum, that we denote by $B$. The latter ranges over all probability couplings $\pi$ whose first marginal equals $\P$. Recalling that the entropy of a coupling $\pi$ blows up when $\pi$ is not absolutely continuous with respect to $\pi_{0}$, this supremum can be recast over couplings of the form
\begin{equation} \label{eq.normalpha}
\pi = \alpha(\xi,\zeta) \:\pi_{0}, \text{ where }\alpha \in
		L^{1}(\Xi \times \Xi; \d\pi_0) \text{ satisfies }\int_{\Xi}\alpha(\xi,\zeta) \:\d\nu_{\xi}
		(\zeta) = 1 \text{ for }\P \text{-a.e. }\xi \in \Xi.
\end{equation}
Hence,
\begin{equation} \label{eq.exprBwdro}
B = \sup\limits_{\alpha \in L^1(\Xi \times \Xi; \d\pi_0) \atop \int_\Xi
		\alpha(\xi,\zeta) \:\d\nu_\xi(\zeta) = 1}\int_{\Xi \times \Xi}\Big( f(\zeta)
		- \lambda c(\xi,\zeta) - \lambda \e \log \alpha(\xi,\zeta) \Big) \alpha(\xi,\zeta
		) \:\d\pi_{0}(\xi,\zeta).
\end{equation}
It turns out that this maximization problem, involving a strictly concave objective function, can be solved explicitly. Indeed, according to the first-order optimality conditions, there exists a function $\mu(\xi)$ such that its unique solution $\alpha(\xi,\zeta)$ satisfies:
$$f(\zeta) - \lambda c(\xi,\zeta) - \lambda \e \log \alpha(\xi,\zeta) - \lambda \e + \mu(\xi) =0,$$
and so:
$$\alpha(\xi,\zeta) = e^{\frac{f(\zeta) - \lambda c(\xi,\zeta) -\lambda\e + \mu(\xi)}{\lambda\e}}.$$
Recalling the normalization \cref{eq.normalpha}, the function $\mu(\xi)$ is characterized by:
$$e^{\frac{\mu(\xi)}{\lambda\e}}= \left(\int_{\Xi}e^{\frac{f(\zeta) - \lambda c(\xi,\zeta)
		-\lambda\e}{\lambda\e}}\:\d\nu_{\xi}(\zeta) \right)^{-1},$$
whence:
$$\alpha(\xi,\zeta) = \left(\int_{\Xi}e^{\frac{f(\zeta) - \lambda c(\xi,\zeta)}{\lambda\e}}
		\:\d\nu_{\xi}(\zeta) \right)^{-1}e^{\frac{f(\zeta) - \lambda c(\xi,\zeta)}{\lambda\e}}.$$
Eventually, injecting this expression in \cref{eq.exprBwdro}, a simple calculation yields:
$$B = \lambda\e \int_{\Xi}\log\left( \int_{\Xi}e^{\frac{f(\zeta) - \lambda c(\xi,\zeta)}{\lambda\e}}
		\:\d\nu_{\xi}(\zeta) \right)\:\d\P(\xi) ,$$
and the desired formula \cref{eq.reformspwdro} follows immediately.
\begin{flushright}
$\square$
\end{flushright}

\section{\textbf{Hint of Proof of \cref{prop.formulamoment}}}\label{app.moment}

\noindent 
We rely on an argument similar to that employed in \cref{app.wdro}. Let $A := \sup\limits_{\Q \in \calA_{\text{M}}} \left(\int_\Xi f(\xi) \:\d \Q(\xi) - \e H(\Q)\right)$ be the worst-case of the (penalized) mean value of the function $f(\xi)$ when $\Q$ runs through the moment-based ambiguity set $\calA_{\text{M}}\subset \calP(\xi)$. By expressing this constraint with Lagrange multipliers, we obtain:
\begin{multline*}
		A = \sup\limits_{\Q \in \calP(\Xi)}\inf\limits_{\lambda \geq 0, \atop S \in
		\S^k_+}\left( \int_{\Xi}f(\xi) \:\d \Q(\xi) + \lambda \left(m_{1}- \left\lvert
		\mu_{0}- \int_{\Xi}\xi \:\d \Q(\xi) \right\lvert \right) \right. \\ \left.+ S
		:\left(m_{2}\Sigma_{0}- \int_{\Xi}(\xi - \mu_{0}) \otimes (\xi - \mu_{0}) \:\d
		\Q(\xi) \right) - \e H(\Q) \right).
\end{multline*}
Thanks to the simple formula
$$\forall \xi \in \R^{k}, \quad \lvert \xi \lvert = \sup\limits_{\tau \in \R^k,
		\atop \lvert \tau \lvert \leq 1}\tau \cdot \xi,$$
this rewrites:
	\begin{multline*}
		A = \sup\limits_{\Q \in \calP(\Xi)}\inf\limits_{\lambda \geq 0, \atop
		{ \lvert \tau \lvert \leq 1, \atop S \in \S^k_+}}\left( \int_{\Xi}f(\xi) \:\d \Q
		(\xi) + \lambda m_{1}- \lambda \tau \cdot \left( \mu_{0}- \int_{\Xi}\xi \:\d
		\Q(\xi) \right) \right. \\ \left.+\: S :\left(m_{2}\Sigma_{0}- \int_{\Xi}(\xi -
		\mu_{0}) \otimes (\xi - \mu_{0}) \:\d \Q(\xi) \right) - \e H(\Q) \right ).
	\end{multline*}
A formal exchange of the infimum and supremum in the above formula then yields:
\begin{multline*}
A = \inf\limits_{\lambda\geqslant 0, \atop {|\tau|\leqslant 1, \atop S \in \S^k_+}}
		\left( \lambda m_{1}- \lambda \tau \cdot \mu_{0}+ m_{2}S : \Sigma_{0}+ B \right), 
	\text{ where }\\ B := \sup\limits_{\Q \in \calP(\Xi)
		}\left( \int_{\Xi}f(\xi) \:\d \Q(\xi) + \lambda \tau \cdot \int_{\Xi}\xi \:\d
		\Q(\xi) - S : \int_{\Xi}(\xi - \mu_{0}) \otimes (\xi - \mu_{0}) \:\d \Q(\xi)
		- \e H(\Q) \right).\\ = \sup\limits_{\Q \text{ a.c. w.r.t. } \mathbb{Q}_0}
		\left( \int_{\Xi}f(\xi) \:\d \Q(\xi) + \lambda \tau \cdot \int_{\Xi}\xi \:\d
		\Q(\xi) - S : \int_{\Xi}(\xi - \mu_{0}) \otimes (\xi - \mu_{0}) \:\d \Q(\xi)
		- \e H(\Q) \right),
	\end{multline*}
where we have sued the fact that $H(\Q)$ takes infinite values when $\Q$ is not absolutely continuous with respect to $\Q_0$.
We now proceed to calculate the quantity $B$:
\begin{equation}\label{eq.defBmoms}
 B:= \sup\limits_{\alpha \in L^1(\Xi,\d\Q_0), \atop \int_\Xi
		\alpha \:\d \Q_0 = 1}\int_{\Xi}\Big(f(\xi) + \lambda \tau \cdot \xi - S : (\xi
		-\mu_{0}) \otimes (\xi-\mu_{0}) - \e \log\alpha(\xi) \Big) \alpha(\xi) \:\d \Q
		_{0}(\xi).
\end{equation}
The necessary optimality condition for this strictly concave maximization problem yields the existence of a real value $\mu$ such that the maximizer $\alpha$ satisfies:
$$f(\xi) + \lambda \tau \cdot \xi - S : (\xi-\mu_{0}) \otimes (\xi-\mu_{0}) - \e
		\log \alpha(\xi) - \e + \mu = 0.$$
By rearranging this identity and using the normalization condition for $\alpha$, we then obtain:
$$\alpha(\xi) = \frac{1}{\int_{\Xi}\left(e^{\frac{f(\xi) + \lambda \tau \cdot \xi
		- S : (\xi-\mu_0) \otimes (\xi-\mu_0)}{\e}}\right) \:\d \Q_{0}(\xi)}e^{\frac{f(\xi)
		+ \lambda \tau \cdot \xi - S : (\xi-\mu_{0}) \otimes (\xi-\mu_{0})}{\e}}.
$$
Now inserting this expression for $\alpha$ in \cref{eq.defBmoms}, we eventually end up with:
$$ B = \e \log \left( \int_{\Xi}\left(e^{\frac{f(\xi) + \lambda
		\tau \cdot \xi - S : (\xi-\mu_{0}) \otimes (\xi-\mu_{0})}{\e}}\right) \:\d \Q
		_{0}(\xi)\right),$$
which leads to the desired result \cref{eq.dualmoments}.
\begin{flushright}	
$\square$
\end{flushright}	

\section{\textbf{Hint of Proof of \cref{prop.cvarrep}}}\label{app.cvar}

\noindent 
For simplicity, we restrict ourselves to providing a sketch of the proof of this result in the case where the cumulative distribution function $\Psi(h, \cdot)$ of the law $\P$ is continuous on $\R$, i.e. there is no subset of $\Xi$ with positive measure where the cost function $\calC(h,\cdot)$ takes a given value $t \in \R$. Note however that the result holds true in the general case, see \cite{rockafellar2002conditional}.
Let us set
$$F(h, \alpha) = \alpha + \frac{1}{1-\beta}\int_{\Xi}\left[ \calC(h,\xi) - \alpha
		\right]_{+}\:\d \P(\xi).$$
A simple calculation shows that:
$$\begin{array}{>{\displaystyle}cc>{\displaystyle}l}
			F(h,\alpha) & = & \alpha + \frac{1}{1-\beta} \int_{\left\{ \xi \in \Xi, \:\: \calC(h,\xi) \geq \alpha \right\}} (\calC(h,\xi) - \alpha ) \:\d \P(\xi)\\
			& = & \frac{1}{1-\beta} \left( \int_{\left\{ \xi \in \Xi, \:\: \calC(h,\xi) \geq \alpha_\beta(h) \right\}} \alpha \:\d \P(\xi ) + \int_{\left\{ \xi \in \Xi, \:\: \calC(h,\xi) \geq \alpha \right\}} (\calC(h,\xi) - \alpha ) \:\d \P(\xi) \right),
		\end{array}$$
where we have used \cref{eq.probavar} to pass from the first line to the second one. Let us now make out two cases:

\begin{itemize}
\item If $\alpha_{\beta}(h) < \alpha$, a simple calculation yields:
$$\begin{array}{>{\displaystyle}cc>{\displaystyle}l}
			F(h,\alpha) & =    & \frac{1}{1-\beta} \left( \int_{\left\{ \xi \in \Xi, \:\: \alpha_\beta(h) \leq \calC(h,\xi) \leq \alpha \right\}} \alpha \:\d \P(\xi)+ \int_{\left\{ \xi \in \Xi, \:\: \calC(h,\xi) \geq \alpha \right\}} \calC(h,\xi) \:\d \P(\xi) \right) \\[1em]
			& =    & \frac{1}{1-\beta} \int_{\left\{ \xi \in \Xi, \:\: \calC(h,\xi) \geq \alpha_{\beta}(h) \right\}} \max\left(\calC(h,\xi),\alpha \right) \:\d \P(\xi)\\[1em]
			& \geq & \frac{1}{1-\beta} \int_{\left\{ \xi \in \Xi, \:\: \calC(h,\xi) \geq \alpha_\beta(h) \right\}} \calC(h,\xi) \:\d \P(\xi)\\
			 & =    & \phi_\beta(h).
		\end{array}
$$
\item If $\alpha < \alpha_{\beta}(h)$, it holds similarly:
$$\begin{array}{>{\displaystyle}cc>{\displaystyle}l}
			F(h,\alpha) & =    & \frac{1}{1-\beta} \Big( - \int_{\left\{ \xi \in \Xi, \:\: \alpha \leq \calC(h,\xi) \leq \alpha_\beta(h) \right\}} \alpha \:\d \P(\xi)+ \int_{\left\{ \xi \in \Xi, \:\: \calC(h,\xi) \geq \alpha_\beta(h) \right\}} \calC(h,\xi) \:\d \P(\xi) \\[1em]
			           &      & \quad\quad\quad \quad\quad\quad \quad\quad\quad \quad\quad\quad \quad\quad\quad \quad\quad\quad+ \int_{\left\{ \xi \in \Xi, \:\: \alpha \leq \calC(h,\xi) \leq \alpha_\beta(h) \right\}} \calC(h,\xi) \:\d \P(\xi) \Big)                     \\[1em]
            & =    & \phi_\beta(h) + \frac{1}{1-\beta} \int_{\left\{ \xi \in \Xi, \:\: \alpha \leq \calC(h,\xi) \leq \alpha_\beta(h) \right\}} \left( \calC(h,\xi) - \alpha\right) \:\d \P(\xi)                           \\[1em]
		    & \geq & \phi_\beta(h).
		\end{array}$$
\end{itemize}	
This proves that $F(h,\alpha) \leq \phi_{\beta}(h)$ for all $\alpha\neq \alpha_{\beta} (h)$, while the definition of $\alpha_{\beta}(h)$ immediately implies that $F(h,\alpha_{\beta}(h)) = \phi_{\beta}(h)$; this terminates the proof of the proposition. 
\begin{flushright}
$\square$
\end{flushright}

\bibliographystyle{siam}
\bibliography{genbib.bib}

\begin{thebibliography}{100}

\bibitem{acar2021modeling}
{\sc E.~Acar, G.~Bayrak, Y.~Jung, I.~Lee, P.~Ramu, and S.~S. Ravichandran}, {\em Modeling, analysis, and optimization under uncertainties: a review}, Structural and Multidisciplinary Optimization, 64 (2021), pp.~2909--2945.

\bibitem{acerbi2002coherence}
{\sc C.~Acerbi and D.~Tasche}, {\em On the coherence of expected shortfall}, Journal of banking \& finance, 26 (2002), pp.~1487--1503.

\bibitem{allaire2014linearized}
{\sc G.~Allaire and C.~Dapogny}, {\em A linearized approach to worst-case design in parametric and geometric shape optimization}, Mathematical Models and Methods in Applied Sciences, 24 (2014), pp.~2199--2257.

\bibitem{allaire2015deterministic}
\leavevmode\vrule height 2pt depth -1.6pt width 23pt, {\em A deterministic approximation method in shape optimization under random uncertainties}, SMAI Journal of computational mathematics, 1 (2015), pp.~83--143.

\bibitem{allaire2011topology}
{\sc G.~Allaire, C.~Dapogny, and P.~Frey}, {\em Topology and geometry optimization of elastic structures by exact deformation of simplicial mesh}, Comptes Rendus Mathematique, 349 (2011), pp.~999--1003.

\bibitem{allaire2013mesh}
\leavevmode\vrule height 2pt depth -1.6pt width 23pt, {\em A mesh evolution algorithm based on the level set method for geometry and topology optimization}, Structural and Multidisciplinary Optimization, 48 (2013), pp.~711--715.

\bibitem{allaire2014shape}
\leavevmode\vrule height 2pt depth -1.6pt width 23pt, {\em Shape optimization with a level set based mesh evolution method}, Computer Methods in Applied Mechanics and Engineering, 282 (2014), pp.~22--53.

\bibitem{allaire2020shape}
{\sc G.~Allaire, C.~Dapogny, and F.~Jouve}, {\em Shape and topology optimization}, in Geometric partial differential equations, part II, A. Bonito and R. Nochetto eds., Handbook of Numerical Analysis, 22 (2021), pp.~1--132.

\bibitem{allaire2008minimum}
{\sc G.~Allaire and F.~Jouve}, {\em Minimum stress optimal design with the level set method}, Engineering analysis with boundary elements, 32 (2008), pp.~909--918.

\bibitem{allaire2004structural}
{\sc G.~Allaire, F.~Jouve, and A.-M. Toader}, {\em Structural optimization using sensitivity analysis and a level-set method}, Journal of computational physics, 194 (2004), pp.~363--393.

\bibitem{allaire2007conception}
{\sc G.~Allaire and M.~Schoenauer}, {\em Conception optimale de structures}, vol.~58, Springer, 2007.

\bibitem{amstutz2016notion}
{\sc S.~Amstutz and M.~Ciligot-Travain}, {\em A notion of compliance robustness in topology optimization}, ESAIM: Control, Optimisation and Calculus of Variations, 22 (2016), pp.~64--87.

\bibitem{azizian2023exact}
{\sc W.~Azizian, F.~Iutzeler, and J.~Malick}, {\em Exact generalization guarantees for (regularized) wasserstein distributionally robust models}, Advances in Neural Information Processing Systems, 36 (2023), pp.~14584--14596.

\bibitem{azizian2022regularization}
\leavevmode\vrule height 2pt depth -1.6pt width 23pt, {\em Regularization for wasserstein distributionally robust optimization}, ESAIM: Control, Optimisation and Calculus of Variations, 29 (2023), p.~33.

\bibitem{balarac2022tetrahedral}
{\sc G.~Balarac, F.~Basile, P.~B{\'e}nard, F.~Bordeu, J.-B. Chapelier, L.~Cirrottola, G.~Caumon, C.~Dapogny, P.~Frey, A.~Froehly, et~al.}, {\em Tetrahedral remeshing in the context of large-scale numerical simulation and high performance computing}, MathematicS In Action, 11 (2022), pp.~129--164.

\bibitem{bendsoe2013topology}
{\sc M.~P. Bendsoe and O.~Sigmund}, {\em Topology optimization: theory, methods, and applications}, Springer Science \& Business Media, 2013.

\bibitem{billingsley2017probability}
{\sc P.~Billingsley}, {\em Probability and measure}, John Wiley \& Sons, 2017.

\bibitem{blanchard2021accurately}
{\sc P.~Blanchard, D.~J. Higham, and N.~J. Higham}, {\em Accurately computing the log-sum-exp and softmax functions}, IMA Journal of Numerical Analysis, 41 (2021), pp.~2311--2330.

\bibitem{blanchet2024distributionally}
{\sc J.~Blanchet, J.~Li, S.~Lin, and X.~Zhang}, {\em Distributionally robust optimization and robust statistics}, arXiv preprint arXiv:2401.14655,  (2024).

\bibitem{blanchet2021statistical}
{\sc J.~Blanchet, K.~Murthy, and V.~A. Nguyen}, {\em Statistical analysis of wasserstein distributionally robust estimators}, in Tutorials in Operations Research: Emerging optimization methods and modeling techniques with applications, INFORMS, 2021, pp.~227--254.

\bibitem{bui2012accurate}
{\sc C.~Bui, C.~Dapogny, and P.~Frey}, {\em An accurate anisotropic adaptation method for solving the level set advection equation}, International Journal for Numerical Methods in Fluids, 70 (2012), pp.~899--922.

\bibitem{cardoso2019robust}
{\sc E.~L. Cardoso, G.~Da~Silva, and A.~T. Beck}, {\em Robust topology optimization of compliant mechanisms with uncertainties in output stiffness}, International Journal for Numerical Methods in Engineering, 119 (2019), pp.~532--547.

\bibitem{cea1986conception}
{\sc J.~C{\'e}a}, {\em Conception optimale ou identification de formes, calcul rapide de la d{\'e}riv{\'e}e directionnelle de la fonction co{\^u}t}, ESAIM: Mathematical Modelling and Numerical Analysis, 20 (1986), pp.~371--402.

\bibitem{chaudhuri2020risk}
{\sc A.~Chaudhuri, M.~Norton, and B.~Kramer}, {\em Risk-based design optimization via probability of failure, conditional value-at-risk, and buffered probability of failure}, in AIAA Scitech 2020 Forum, 2020, p.~2130.

\bibitem{chen2020distributionally}
{\sc R.~Chen, I.~C. Paschalidis, et~al.}, {\em Distributionally robust learning}, Foundations and Trends{\textregistered} in Optimization, 4 (2020), pp.~1--243.

\bibitem{chen2010level}
{\sc S.~Chen, W.~Chen, and S.~Lee}, {\em Level set based robust shape and topology optimization under random field uncertainties}, Structural and Multidisciplinary Optimization, 41 (2010), pp.~507--524.

\bibitem{cherkaev1999optimal}
{\sc A.~Cherkaev and E.~Cherkaeva}, {\em Optimal design for uncertain loading condition}, in Homogenization: In Memory of Serguei Kozlov, World Scientific, 1999, pp.~193--213.

\bibitem{cuturi2013sinkhorn}
{\sc M.~Cuturi}, {\em Sinkhorn distances: Lightspeed computation of optimal transport}, Advances in neural information processing systems, 26 (2013).

\bibitem{totuto}
{\sc C.~Dapogny}, {\em \texttt{totuto} a tentative implementation of density-based topology implementation in \texttt{python}, \texttt{https://github.com/dapogny/totuto}}, 2023.

\bibitem{dapogny2014three}
{\sc C.~Dapogny, C.~Dobrzynski, and P.~Frey}, {\em Three-dimensional adaptive domain remeshing, implicit domain meshing, and applications to free and moving boundary problems}, Journal of computational physics, 262 (2014), pp.~358--378.

\bibitem{dapogny2022tuto}
{\sc C.~Dapogny and F.~Feppon}, {\em Shape optimization using a level set based mesh evolution method: an overview and tutorial}, \texttt{Hal} preprint: \texttt{https://hal.archives-ouvertes.fr/hal-03881641},  (2022).

\bibitem{dapogny2012computation}
{\sc C.~Dapogny and P.~Frey}, {\em Computation of the signed distance function to a discrete contour on adapted triangulation}, Calcolo, 49 (2012), pp.~193--219.

\bibitem{dapogny2023entropy}
{\sc C.~Dapogny, F.~Iutzeler, A.~Meda, and B.~Thibert}, {\em Entropy-regularized wasserstein distributionally robust shape and topology optimization}, Structural and Multidisciplinary Optimization, 66 (2023), p.~42.

\bibitem{de2019topologyoptimizationuncertaintyusing}
{\sc S.~De, J.~Hampton, K.~Maute, and A.~Doostan}, {\em Topology optimization under uncertainty using a stochastic gradient-based approach}, 2019.

\bibitem{de2021reliabilitybasedtopologyoptimizationusing}
{\sc S.~De, K.~Maute, and A.~Doostan}, {\em Reliability-based topology optimization using stochastic gradients}, 2021.

\bibitem{de2021topologyoptimizationmicroscaleuncertainty}
\leavevmode\vrule height 2pt depth -1.6pt width 23pt, {\em Topology optimization under microscale uncertainty using stochastic gradients}, 2021.

\bibitem{de2008shape}
{\sc F.~De~Gournay, G.~Allaire, and F.~Jouve}, {\em Shape and topology optimization of the robust compliance via the level set method}, ESAIM: Control, Optimisation and Calculus of Variations, 14 (2008), pp.~43--70.

\bibitem{dehghannasiri2018intrinsically}
{\sc R.~Dehghannasiri, X.~Qian, and E.~R. Dougherty}, {\em Intrinsically bayesian robust karhunen-lo{\`e}ve compression}, Signal Processing, 144 (2018), pp.~311--322.

\bibitem{delage2010distributionally}
{\sc E.~Delage and Y.~Ye}, {\em Distributionally robust optimization under moment uncertainty with application to data-driven problems}, Operations research, 58 (2010), pp.~595--612.

\bibitem{drusvyatskiy2023stochastic}
{\sc D.~Drusvyatskiy and L.~Xiao}, {\em Stochastic optimization with decision-dependent distributions}, Mathematics of Operations Research, 48 (2023), pp.~954--998.

\bibitem{eigel2018risk}
{\sc M.~Eigel, J.~Neumann, R.~Schneider, and S.~Wolf}, {\em Risk averse stochastic structural topology optimization}, Computer Methods in Applied Mechanics and Engineering, 334 (2018), pp.~470--482.

\bibitem{esfahani2018data}
{\sc P.~M. Esfahani and D.~Kuhn}, {\em Data-driven distributionally robust optimization using the wasserstein metric: Performance guarantees and tractable reformulations}, Mathematical Programming, 171 (2018), pp.~115--166.

\bibitem{nullspace}
{\sc F.~Feppon}, {\em \texttt{Null Space Optimizer}, \texttt{https://null-space-optimizer.readthedocs.io/en/latest/}}, 2023.

\bibitem{feppon2024tutoNS}
\leavevmode\vrule height 2pt depth -1.6pt width 23pt, {\em Density-based topology optimization with the null space optimizer: a tutorial and a comparison}, Structural and Multidisciplinary Optimization, 67 (2024), pp.~1--34.

\bibitem{feppon2020null}
{\sc F.~Feppon, G.~Allaire, and C.~Dapogny}, {\em Null space gradient flows for constrained optimization with applications to shape optimization}, ESAIM: Control, Optimisation and Calculus of Variations, 26 (2020), p.~90.

\bibitem{feydy2020analyse}
{\sc J.~Feydy}, {\em Analyse de donn{\'e}es g{\'e}om{\'e}triques, au del{\`a} des convolutions}, PhD thesis, Universit{\'e} Paris-Saclay, 2020.

\bibitem{fournier2013rateconvergencewassersteindistance}
{\sc N.~Fournier and A.~Guillin}, {\em On the rate of convergence in wasserstein distance of the empirical measure}, 2013.

\bibitem{grieshammer2024continuous}
{\sc M.~Grieshammer, L.~Pflug, M.~Stingl, and A.~Uihlein}, {\em The continuous stochastic gradient method: part i--convergence theory}, Computational Optimization and Applications, 87 (2024), pp.~935--976.

\bibitem{grieshammer2024continuous2}
\leavevmode\vrule height 2pt depth -1.6pt width 23pt, {\em The continuous stochastic gradient method: part ii--application and numerics}, Computational Optimization and Applications, 87 (2024), pp.~977--1008.

\bibitem{grigoriu1998simulation}
{\sc M.~Grigoriu}, {\em Simulation of stationary non-gaussian translation processes}, Journal of engineering mechanics, 124 (1998), pp.~121--126.

\bibitem{guest2008structural}
{\sc J.~K. Guest and T.~Igusa}, {\em Structural optimization under uncertain loads and nodal locations}, Computer Methods in Applied Mechanics and Engineering, 198 (2008), pp.~116--124.

\bibitem{guo2009confidence}
{\sc X.~Guo, W.~Bai, W.~Zhang, and X.~Gao}, {\em Confidence structural robust design and optimization under stiffness and load uncertainties}, Computer Methods in Applied Mechanics and Engineering, 198 (2009), pp.~3378--3399.

\bibitem{guo2013robust}
{\sc X.~Guo, W.~Zhang, and L.~Zhang}, {\em Robust structural topology optimization considering boundary uncertainties}, Computer Methods in Applied Mechanics and Engineering, 253 (2013), pp.~356--368.

\bibitem{hamdia2022multilevel}
{\sc K.~M. Hamdia, H.~Ghasemi, X.~Zhuang, and T.~Rabczuk}, {\em Multilevel monte carlo method for topology optimization of flexoelectric composites with uncertain material properties}, Engineering Analysis with Boundary Elements, 134 (2022), pp.~412--418.

\bibitem{hanasusanto2015distributionally}
{\sc G.~A. Hanasusanto, V.~Roitch, D.~Kuhn, and W.~Wiesemann}, {\em A distributionally robust perspective on uncertainty quantification and chance constrained programming}, Mathematical Programming, 151 (2015), pp.~35--62.

\bibitem{hasofer1974exact}
{\sc A.~M. Hasofer}, {\em An exact and invarient first order reliability format}, J. Eng. Mech. Div., Proc. ASCE, 100 (1974), pp.~111--121.

\bibitem{hecht2012new}
{\sc F.~Hecht}, {\em New development in freefem++}, Journal of numerical mathematics, 20 (2012), pp.~251--266.

\bibitem{henrot2018shape}
{\sc A.~Henrot and M.~Pierre}, {\em Shape Variation and Optimization}, EMS Tracts in Mathematics Vol. 28, 2018.

\bibitem{holmberg2015worst}
{\sc E.~Holmberg, C.-J. Thore, and A.~Klarbring}, {\em Worst-case topology optimization of self-weight loaded structures using semi-definite programming}, Structural and Multidisciplinary Optimization, 52 (2015), pp.~915--928.

\bibitem{hu2013kullback}
{\sc Z.~Hu and L.~J. Hong}, {\em Kullback-leibler divergence constrained distributionally robust optimization}, Available at Optimization Online, 1 (2013), p.~9.

\bibitem{jofre2021rapidaerodynamicshapeoptimization}
{\sc L.~Jofre and A.~Doostan}, {\em Rapid aerodynamic shape optimization under parametric and turbulence model uncertainty: A stochastic gradient approach}, 2021.

\bibitem{Kanno_2021}
{\sc Y.~Kanno}, {\em Structural reliability under uncertainty in moments: distributionally-robust reliability-based design optimization}, Japan Journal of Industrial and Applied Mathematics, 39 (2021), p.~195–226.

\bibitem{kapteyn2019distributionally}
{\sc M.~G. Kapteyn, K.~E. Willcox, and A.~B. Philpott}, {\em Distributionally robust optimization for engineering design under uncertainty}, International Journal for Numerical Methods in Engineering, 120 (2019), pp.~835--859.

\bibitem{keshavarzzadeh2017topology}
{\sc V.~Keshavarzzadeh, F.~Fernandez, and D.~A. Tortorelli}, {\em Topology optimization under uncertainty via non-intrusive polynomial chaos expansion}, Computer Methods in Applied Mechanics and Engineering, 318 (2017), pp.~120--147.

\bibitem{kingma2017adam}
{\sc D.~P. Kingma and J.~Ba}, {\em Adam: A method for stochastic optimization}, 2017.

\bibitem{kouri2016risk}
{\sc D.~P. Kouri and T.~M. Surowiec}, {\em Risk-averse pde-constrained optimization using the conditional value-at-risk}, SIAM Journal on Optimization, 26 (2016), pp.~365--396.

\bibitem{kuhn2019wasserstein}
{\sc D.~Kuhn, P.~M. Esfahani, V.~A. Nguyen, and S.~Shafieezadeh-Abadeh}, {\em Wasserstein distributionally robust optimization: Theory and applications in machine learning}, in Operations research \& management science in the age of analytics, Informs, 2019, pp.~130--166.

\bibitem{kuhn2024distributionallyrobustoptimization}
{\sc D.~Kuhn, S.~Shafiee, and W.~Wiesemann}, {\em Distributionally robust optimization}, 2024.

\bibitem{langelaar2019topology}
{\sc M.~Langelaar}, {\em Topology optimization for multi-axis machining}, Computer Methods in Applied Mechanics and Engineering, 351 (2019), pp.~226--252.

\bibitem{lazarov2012topology}
{\sc B.~S. Lazarov, M.~Schevenels, and O.~Sigmund}, {\em Topology optimization considering material and geometric uncertainties using stochastic collocation methods}, Structural and Multidisciplinary optimization, 46 (2012), pp.~597--612.

\bibitem{lin2022distributionally}
{\sc F.~Lin, X.~Fang, and Z.~Gao}, {\em Distributionally robust optimization: A review on theory and applications}, Numerical Algebra, Control \& Optimization, 12 (2022), p.~159.

\bibitem{liu2018current}
{\sc J.~Liu, A.~T. Gaynor, S.~Chen, Z.~Kang, K.~Suresh, A.~Takezawa, L.~Li, J.~Kato, J.~Tang, C.~C. Wang, et~al.}, {\em Current and future trends in topology optimization for additive manufacturing}, Structural and Multidisciplinary Optimization, pp.~1--27.

\bibitem{liu2018note}
{\sc Q.~Liu, J.~Wu, X.~Xiao, and L.~Zhang}, {\em A note on distributionally robust optimization under moment uncertainty}, Journal of Numerical Mathematics, 26 (2018), pp.~141--150.

\bibitem{loeve2017probability}
{\sc M.~Lo{\`e}ve}, {\em Probability theory}, Courier Dover Publications, 2017.

\bibitem{martinez2016large}
{\sc J.~Mart{\'i}nez-Frutos and D.~Herrero-P{\'e}rez}, {\em Large-scale robust topology optimization using multi-gpu systems}, Computer Methods in Applied Mechanics and Engineering, 311 (2016), pp.~393--414.

\bibitem{martinez2016robust}
{\sc J.~Martinez-Frutos, D.~Herrero-Perez, M.~Kessler, and F.~Periago}, {\em Robust shape optimization of continuous structures via the level set method}, Computer Methods in Applied Mechanics and Engineering, 305 (2016), pp.~271--291.

\bibitem{maute2014topology}
{\sc K.~Maute}, {\em Topology optimization under uncertainty}, in Topology optimization in structural and continuum mechanics,  (2014), pp.~457--471.

\bibitem{merigot2021optimal}
{\sc Q.~Merigot and B.~Thibert}, {\em Optimal transport: discretization and algorithms}, in Handbook of Numerical Analysis, vol.~22, Elsevier, 2021, pp.~133--212.

\bibitem{michailidis2014manufacturing}
{\sc G.~Michailidis}, {\em Manufacturing constraints and multi-phase shape and topology optimization via a level-set method}, PhD thesis, 2014.

\bibitem{mohamed2020monte}
{\sc S.~Mohamed, M.~Rosca, M.~Figurnov, and A.~Mnih}, {\em Monte carlo gradient estimation in machine learning}, Journal of Machine Learning Research, 21 (2020), pp.~1--62.

\bibitem{montesuma2024recent}
{\sc E.~F. Montesuma, F.~M.~N. Mboula, and A.~Souloumiac}, {\em Recent advances in optimal transport for machine learning}, IEEE Transactions on Pattern Analysis and Machine Intelligence,  (2024).

\bibitem{murat1976controle}
{\sc F.~Murat and J.~Simon}, {\em Sur le contr{\^o}le par un domaine g{\'e}om{\'e}trique}, Pr\'e-publication du Laboratoire d'Analyse Num\'erique,(76015),  (1976).

\bibitem{nakao2021distributionally}
{\sc H.~Nakao, R.~Jiang, and S.~Shen}, {\em Distributionally robust partially observable markov decision process with moment-based ambiguity}, SIAM Journal on Optimization, 31 (2021), pp.~461--488.

\bibitem{nguyen2022distributionally}
{\sc V.~A. Nguyen, D.~Kuhn, and P.~Mohajerin~Esfahani}, {\em Distributionally robust inverse covariance estimation: The wasserstein shrinkage estimator}, Operations research, 70 (2022), pp.~490--515.

\bibitem{nocedal2006numerical}
{\sc J.~Nocedal and S.~J. Wright}, {\em Numerical optimization 2nd}, Springer, 2006.

\bibitem{nutz2021introduction}
{\sc M.~Nutz}, {\em Introduction to entropic optimal transport}, Lecture notes, Columbia University,  (2021).

\bibitem{osher1988fronts}
{\sc S.~Osher and J.~A. Sethian}, {\em Fronts propagating with curvature-dependent speed: algorithms based on hamilton-jacobi formulations}, Journal of computational physics, 79 (1988), pp.~12--49.

\bibitem{peyre2019computational}
{\sc G.~Peyr{\'e}, M.~Cuturi, et~al.}, {\em Computational optimal transport: With applications to data science}, Foundations and Trends in Machine Learning, 11 (2019), pp.~355--607.

\bibitem{pflug2007ambiguity}
{\sc G.~Pflug and D.~Wozabal}, {\em Ambiguity in portfolio selection}, Quantitative Finance, 7 (2007), pp.~435--442.

\bibitem{pflug2024stochastic}
{\sc L.~Pflug, M.~Stingl, and A.~Uihlein}, {\em A stochastic method of moving asymptotes for topology optimization under uncertainty}, arXiv preprint arXiv:2410.19428,  (2024).

\bibitem{pironneau1989finite}
{\sc O.~Pironneau}, {\em Finite element methods for fluids}, Wiley Chichester, 1989.

\bibitem{plessix2006review}
{\sc R.-E. Plessix}, {\em A review of the adjoint-state method for computing the gradient of a functional with geophysical applications}, Geophysical Journal International, 167 (2006), pp.~495--503.

\bibitem{prando2025distributionally}
{\sc J.~Prando}, {\em Distributionally robust shape and topology optimization}, PhD thesis, Universit\'e Grenoble Alpes, 2025.

\bibitem{rahimian2019distributionally}
{\sc H.~Rahimian and S.~Mehrotra}, {\em Distributionally robust optimization: A review}, \texttt{arXiv} preprint \texttt{arXiv:1908.05659},  (2019).

\bibitem{rahimian2022frameworks}
\leavevmode\vrule height 2pt depth -1.6pt width 23pt, {\em Frameworks and results in distributionally robust optimization}, Open Journal of Mathematical Optimization, 3 (2022), pp.~1--85.

\bibitem{rockafellar2002conditional}
{\sc R.~T. Rockafellar and S.~Uryasev}, {\em Conditional value-at-risk for general loss distributions}, Journal of banking \& finance, 26 (2002), pp.~1443--1471.

\bibitem{rockafellar2000optimization}
{\sc R.~T. Rockafellar, S.~Uryasev, et~al.}, {\em Optimization of conditional value-at-risk}, Journal of risk, 2 (2000), pp.~21--42.

\bibitem{santambrogio2015optimal}
{\sc F.~Santambrogio}, {\em Optimal transport for applied mathematicians}, Birk\"auser, 2015.

\bibitem{schwab2006karhunen}
{\sc C.~Schwab and R.~A. Todor}, {\em Karhunen--lo{\`e}ve approximation of random fields by generalized fast multipole methods}, Journal of Computational Physics, 217 (2006), pp.~100--122.

\bibitem{sejourne2019sinkhorn}
{\sc T.~S{\'e}journ{\'e}, J.~Feydy, F.-X. Vialard, A.~Trouv{\'e}, and G.~Peyr{\'e}}, {\em Sinkhorn divergences for unbalanced optimal transport}, arXiv preprint arXiv:1910.12958,  (2019).

\bibitem{sethian1999fast}
{\sc J.~A. Sethian}, {\em Fast marching methods}, SIAM review, 41 (1999), pp.~199--235.

\bibitem{sigmund2013topology}
{\sc O.~Sigmund and K.~Maute}, {\em Topology optimization approaches}, Structural and Multidisciplinary Optimization, 48 (2013), pp.~1031--1055.

\bibitem{sokolowski1992introduction}
{\sc J.~Sokolowski and J.-P. Zol{\'e}sio}, {\em Introduction to shape optimization}, Springer, 1992.

\bibitem{staib2019distributionally}
{\sc M.~Staib and S.~Jegelka}, {\em Distributionally robust optimization and generalization in kernel methods}, Advances in Neural Information Processing Systems, 32 (2019).

\bibitem{tootkaboni2012topology}
{\sc M.~Tootkaboni, A.~Asadpoure, and J.~K. Guest}, {\em Topology optimization of continuum structures under uncertainty--a polynomial chaos approach}, Computer Methods in Applied Mechanics and Engineering, 201 (2012), pp.~263--275.

\bibitem{uihlein2025140}
{\sc A.~Uihlein, O.~Sigmund, and M.~Stingl}, {\em A 140 line matlab code for topology optimization problems with probabilistic parameters}, arXiv preprint arXiv:2505.10421,  (2025).

\bibitem{villani2009optimal}
{\sc C.~Villani}, {\em Optimal transport: old and new}, vol.~338, Springer, 2009.

\bibitem{vincent2024texttt}
{\sc F.~Vincent, W.~Azizian, F.~Iutzeler, and J.~Malick}, {\em \texttt{skwdro}: a library for wasserstein distributionally robust machine learning}, arXiv preprint arXiv:2410.21231,  (2024).

\bibitem{wang2020topology}
{\sc C.~Wang, B.~Xu, Q.~Meng, J.~Rong, and Y.~Zhao}, {\em Topology optimization of cast parts considering parting surface position}, Advances in Engineering Software, 149 (2020), p.~102886.

\bibitem{wang2021sinkhorn}
{\sc J.~Wang, R.~Gao, and Y.~Xie}, {\em Sinkhorn distributionally robust optimization}, arXiv preprint arXiv:2109.11926,  (2021).

\bibitem{wang2003level}
{\sc M.~Y. Wang, X.~Wang, and D.~Guo}, {\em A level set method for structural topology optimization}, Computer methods in applied mechanics and engineering, 192 (2003), pp.~227--246.

\bibitem{yue2024geometric}
{\sc M.-C. Yue, Y.~Rychener, D.~Kuhn, and V.~A. Nguyen}, {\em A geometric unification of distributionally robust covariance estimators: Shrinking the spectrum by inflating the ambiguity set}, arXiv preprint arXiv:2405.20124,  (2024).

\bibitem{zhao1999general}
{\sc Y.-G. Zhao and T.~Ono}, {\em A general procedure for first/second-order reliabilitymethod (form/sorm)}, Structural safety, 21 (1999), pp.~95--112.

\bibitem{zhou2002progress}
{\sc M.~Zhou, R.~Fleury, Y.-K. Shyy, H.~Thomas, and J.~Brennan}, {\em Progress in topology optimization with manufacturing constraints}, in 9th AIAA/ISSMO Symposium on multidisciplinary analysis and optimization, 2002, p.~5614.

\bibitem{zhu2020design}
{\sc B.~Zhu, X.~Zhang, H.~Zhang, J.~Liang, H.~Zang, H.~Li, and R.~Wang}, {\em Design of compliant mechanisms using continuum topology optimization: A review}, Mechanism and Machine Theory, 143 (2020), p.~103622.

\end{thebibliography}
\end{document}